\documentclass[11pt]{article}
\usepackage{amsmath, exscale, epsfig,amssymb, 
amsthm,
 makeidx}

\usepackage{titlesec}
\titlespacing*{\section}{0pt}{0.4\baselineskip}{0.3\baselineskip}
\titlespacing*{\subsection}{0pt}{0.4\baselineskip}{0.3\baselineskip}

\usepackage[pdftex]{hyperref} 
\hypersetup{
    unicode      = false,     
    pdftoolbar   = true,      
    pdfmenubar   = true,      
    pdffitwindow = true,      
    pdfnewwindow = true,      
    colorlinks   = true,      
    linkcolor    = blue,      
    citecolor    = red,      
    filecolor    = blue,      
    urlcolor     = green       
}

\small
\normalsize
\usepackage[utf8]{inputenc}
\usepackage[T1]{fontenc}
\usepackage[english]{babel}
\usepackage{enumerate,vmargin}

\usepackage{times}
\usepackage{amsfonts}
\usepackage{amsbsy}
\usepackage{amscd}
\usepackage{multicol}
\allowdisplaybreaks

\usepackage[all]{xy}

\usepackage{stmaryrd}
\usepackage{graphicx}
\usepackage{paralist}
\usepackage{amsfonts}
\usepackage{amssymb,mathrsfs}

\setmarginsrb{2.9cm}{2.6cm}{2.9cm}{1.7cm}{0cm}{0mm}{0cm}{9mm}

\def\build#1_#2^#3{\mathrel{\mathop{\kern 0pt#1}\limits_{#2}^{#3}}}
\def\noi{{\noindent}}

\def\cq{$\hfill \square$}
\def\un{{\bf 1}}

\newcommand{\fdelta}{{\boldsymbol{\delta}}}
\newcommand{\fLambda}{{\boldsymbol{\Lambda}}}
\newcommand{\flambda}{{\boldsymbol{\lambda}}}

\newcommand{\falpha}{{\boldsymbol{\alpha}}}

\newcommand{\fell}{{\boldsymbol{\ell}}}
\newcommand{\fmu}{{\boldsymbol{\mu}}}
\newcommand{\fnu}{{\boldsymbol{\nu}}}

\newcommand{\bbT}{\mathbb{T}}

\newcommand{\bE}{{\bf E}}

\newcommand{\bm}{\mu}
\newcommand{\bs}{{\bf s}}

\newcommand{\bM}{\mathbb{M}}
\newcommand{\bN}{\mathbb{N}}

\newcommand{\bQ}{\mathbb{Q}}
\newcommand{\bP}{{\bf P}}
\newcommand{\bR}{\mathbb{R}}

\newcommand{\bT}{\mathbb{T}}
\newcommand{\bbU}{\mathbb{U}}

\newcommand{\cB}{\mathscr{B}}
\newcommand{\cC}{\mathscr{C}}

\newcommand{\cE}{\mathcal{E}}

\newcommand{\cM}{\mathcal{M}}

\def\era{{\cal E}}

\def\bn{{\rm n}}
\def\bm{\mu}

\def\be{\begin{equation}}

\def\ee{\end{equation}}
\def\ba{\begin{eqnarray*}}
\def\ea{\end{eqnarray*}}

\def\noi{\noindent}

\def\supp{\mathtt{Supp}\, }
\newcommand{\lgeo}{[\![}
\newcommand{\rgeo}{]\!]}
\def\cqfd{ \hfill $\blacksquare$ }

\newcommand{\eqo}{\! = \! }
\newcommand{\geqo}{\! \geq \! }
\newcommand{\leqo}{\! \leq \! }
\newcommand{\ino}{\! \in \! }

\newcommand{\bxx}{\mathbf x}
\newcommand{\bhh}{\mathbf h}

\newcommand{\bbR}{\mathbb{R}}
\newcommand{\bbN}{\mathbb{N}}
\newcommand{\leko}{\! < \! }
\newcommand{\geko}{\! > \! }
\newcommand{\dHaus}{d_{\mathtt{Haus}}}
\newcommand{\dPro}{d_{\mathtt{Pro}}}
\newcommand{\deHaus}{\delta_{\mathtt{Haus}}}
\newcommand{\dePro}{\delta_{\mathtt{Pro}}}
\newcommand{\Pt}{\mathtt{P}_{\!\mathtt{t}}}
\newcommand{\Ptn}{\mathtt{P}_{\! \mathtt{t}_n}}
\newcommand{\Ptz}{\mathtt{P}_{\! \mathtt{t}_0}}

\newcommand{\ttt}{\mathtt{t}}
\newcommand{\subo}{\! \subset \!}
\newcommand{\sspaan}{\mathtt{S}\, \overline{\! \mathtt{pan}\!  }\, }
\newcommand{\spaan}{\mathtt{Span}}
\newcommand{\dGP}{\fdelta_{\mathtt{GP}}}
\newcommand{\dGH}{\fdelta_{\mathtt{GH}}}
\newcommand{\dGHP}{\fdelta_{\mathtt{GHP}}}
\newcommand{\dera}{\fdelta_{\mathtt{era}}}
\newcommand{\epp}{\varepsilon}
\newcommand{\ttP}{\mathtt{P}}

\newcommand{\esscov}{\mathtt{ess}\textrm{-}\mathtt{cov}} 
\newcommand{\esstree}{\mathtt{ess}\textrm{-}\mathtt{tree}} 

\newtheoremstyle{thmstyl}
{3.5pt} 
{2.5pt} 
{\em} 
{} 
{\bfseries} 
{.} 
{.5em} 
{} 
\theoremstyle{thmstyl}

\newtheorem{theorem}{Theorem}[section]

\newtheorem{lemma}[theorem]{Lemma}
\newtheorem{proposition}[theorem]{Proposition}
\newtheorem{corollary}[theorem]{Corollary}

\newtheoremstyle{dfstyl}
{3.5pt} 
{2.5pt} 
{} 
{} 
{\bfseries} 
{.} 
{.5em} 
{} 
\theoremstyle{dfstyl}

\newtheorem{definition}[theorem]{Definition}
\newtheorem{notation}[theorem]{Notation}

\newtheorem{remark}[theorem]{Remark}
\newtheorem{example}[theorem]{Example}

\newcommand{\eqnsection}{
\renewcommand{\theequation}{\arabic{section}.\arabic{equation}}
    \makeatletter
    \csname  @addtoreset\endcsname{equation}{section}
    \makeatother}
\eqnsection

\begin{document}

\setlength{\abovedisplayskip}{5pt}
\setlength{\belowdisplayskip}{5pt}

\title{ \bf Mass erasure on measured $\bbR$-trees,\\ applications to L\'evy forests }
\author{Thomas \textsc{Duquesne}
\thanks{Sorbonne Universit\'e, Campus Pierre et Marie Curie, 
LPSM, Case courrier 158, 4, place Jussieu, 
75252 Paris Cedex 05 
France. Email: thomas.duquesne@sorbonne-universite.fr} 
\and
Matthias {\sc Winkel}
\thanks{Department of Statistics, University of Oxford, Oxford OX1 3LB, UK.
Email: matthias.winkel@stats.ox.ac.uk} 
}
\date{\today} 

\maketitle

\begin{abstract}
Let $h\!>\!0$. For a complete and separable $\bbR$-tree $(T,d)$ equipped with a root 
$\rho$ and a finite Borel measure $\mu$, we define the $h$-mass-erased tree by removing from $T$ all fringe subtrees of mass less than $h$ and we equip it with a suitable measure such that the erasure operators $(\cE_h)_{h\ge 0}$ form a semigroup that is continuous for the Gromov-weak topology. 

Then, we say that a sequence $\fmu_n\eqo (T_n,d_n,\rho_n,\mu_n)$, $n\ino\bN$, converges in the sense of mass erasure if $(\cE_h\fmu_n)_{n\in\bN}$ converges 
Gromov-weakly for all $h\geko 0$. This notion of convergence is strictly weaker than Gromov-weak convergence and we establish criteria to relate the two notions. We define a distance function that metrizes convergence in the sense of mass erasure. By extending the notion of measured $\bbR$-trees to allow mass on the boundary (the far ends of infinite geodesics), we obtain a complete metric space. 

Next, we identify random trees of finite type (that is, discrete trees with edge lengths) satisfying the regenerative branching property as a specific class of measured (sub)critical GW forests. We then show that this class of trees is preserved by mass erasure and we compute the law of these mass-erased GW forests explicitly. 
Finally, we establish a limit theorem for these measured (sub)critical GW forests to converge to standard measured L{\' e}vy forests, i.e.~those whose total mass has the same distribution as the total population of a continuous-state branching process. This includes cases with bounded variation by crucially using the convergence in the sense of mass erasure and it 
extends the cases studied previously. 
\smallskip

{\small 
\noi

\smallskip

\noi
}

\end{abstract}

Initiated by Aldous's seminal papers \cite{Al1,Al2} in the early 1990s and with significant contributions on the underpinning analysis by Evans, Winter and co-authors \cite{EPW,EW,GPW} in the 2000s, continuum tree structures (and more general graph-like metric spaces) have become a widely used way of giving rigorous meaning to large-scale limits of discrete random structures. This includes both limits of static models of combinatorial trees and Markovian evolutions of such trees, which arise in applications in computer science, as well as in the biological and physical sciences.

Before we start making formal statements in Section \ref{sec2}, let us discuss our contributions informally in the context of the existing literature. Some of the main spaces of continuum tree structures include spaces of $\bR$-trees. These are tree-like metric spaces.

Let us recall a key notion on the Gromov--Hausdorff space of (isometry classes of) rooted compact $\mathbb{R}$-trees, $h$-erasure (or trimming \cite{EPW}), which we call here $h$-length erasure. The $h$-length-erased tree is obtained from a compact rooted $\mathbb{R}$-tree by removing all fringe subtrees of height less than $h$. In fact, this idea dates back further and was already applied in different formalisms by and Neveu \cite{Ne2} and Kesten \cite{Kes86} for (conditioned) Galton--Watson trees and by Neveu and Pitman \cite{NP89I,NP89II} and Le Gall \cite{LG89} for the tree encoded by a Brownian excursion. For compact $\mathbb{R}$-trees such as Aldous's Brownian continuum random tree \cite{Al1}, with complex features such as dense sets of branch points and leaves, $h$-length erasure for any $h\geko 0$ removes the complexity and returns a tree with a discrete branching structure. It was used to define in a projective way stable CRTs by Le Jan \cite{LJ91}. Approximation of rooted compact $\bbR$-trees via length erasure was introduced and studied systematically by Evans, Pitman and Winter \cite{EPW}, see also \cite{DuWi2} for connections with Lévy trees. The associated operators $R_h$, $h \geqo 0$, on the Gromov--Hausdorff space form a Lipschitz semigroup.  

The idea of $h$-mass erasure, which we discuss in the present paper, is to take a rooted $\mathbb{R}$-tree further equipped with a finite Borel (mass) measure and remove all fringe subtrees of mass less than $h$. A natural space for such operators is the Gromov--Prokhorov space of (measure-preserving isometry classes of) rooted Polish $\mathbb{R}$-trees equipped with a finite Borel measure, whose span of the support is the entire $\mathbb{R}$-tree. We equip the $h$-mass-erased tree with a 
measure that suitably relocates some of the mass of fringe subtrees in such a way that the resulting operators $\cE_h$ form a semigroup. We further show that this semigroup is Gromov--Prokhorov-continuous. 

Related notions appear in the literature as exit \cite{BMart,HKK,FKM,FRIME} or freezing \cite{Bertoinss,HaaMi,Ber06,Rubenthaler} mechanisms affecting small particles in fragmentation processes, by stopping small particles from fragmenting further. This may be the ultimate goal, or indeed also gives rise to a fragmentation property that allows to independently consider the further fragmentation of these fragments. See also how these ideas generalise to a context of more general recursively constructed weighted $\mathbb{R}$-trees in \cite{RW1}. Also the trimming of leaves in \cite{Geldbach} of discrete interval partition trees can be viewed as a simple instance of mass erasure. 

In the present paper we use mass erasure to introduce a new notion of convergence, in which weighted $\mathbb{R}$-trees converge to a limiting weighted $\mathbb{R}$-tree if associated $h$-mass-erased trees converge in the Gromov--Prokhorov sense for all $h>0$. We call this notion convergence in the sense of mass erasure. This turns out to be weaker than the Gromov--Prokhorov convergence of (un-mass-erased) trees. We define a distance function that metrizes this new
notion of convergence, but the resulting metric space is not complete. By extending weighted $\bR$-trees to allow mass on the boundary (at infinite distance from the root, captured by the infinite geodesics), we identify a completion of this metric space.

In previous work \cite{DuWi1,DuWi2}, building on \cite{DuLG,PW1}, we identified a large class of tree growth processes of Galton--Watson forests some of which we showed converge to limiting L\'evy forests in the Gromov--Hausdorff sense, and we established an invariance principle for Galton--Watson forests under the Grey condition that guarantees that associated continuous-state branching processes can go extinct in finite time. 
In the present paper, considering convergence in the sense of mass erasure, we  establish an invariance principle for measured Galton--Watson forests including ((sub)critical) cases where the Grey condition fails.
Furthermore, we show that the sense of convergence can be strengthened to Gromov--Prokhorov convergence if and only if the associated continuous-state branching process has unbounded variation, still including cases when the Grey condition fails, and further to Gromov--Hausdorff--Prokhorov convergence under additional assumptions.

In a companion paper, we study mass erasure for trees encoded by lower semi-continuous functions. This setting, corresponding to trees with a planar order, allows many parallel developments, where the well-established convergence of epigraphs \cite{Beer94} is a starting point that takes the role of Gromov--Prokhorov convergence here. See also L\"ohr \cite[Sec.~4]{Lohr13} 
who explored epigraph-Gromov--Prokhorov-continuity of the coding function.

\smallskip

\noi
\textbf{Structure of this article}. In Section \ref{sec2}, we state our main results. In Section \ref{sec-boundary}, we collect preliminaries on $\mathbb{R}$-trees and on the notion of tree boundary. The theory of mass erasure on a fixed $\mathbb{R}$-tree is presented in Section \ref{Merafixsec}.  In Section \ref{CVspacemerasec}, this theory is adapted and extended to the Gromov setting, and we develop relevant theory to study random measured $\mathbb{R}$-trees in the new space of measured trees. We finally demonstrate in Section \ref{secappl} how this theory can be applied to study measured Galton--Watson trees and the extended class of L\'evy forests.


\section{Statements of main definitions and results}\label{sec2}

We will use the following common abbreviations: 
i.i.d.~for \emph{independent and identically distributed}, r.v. for \emph{random variable}, w.l.o.g.~for \emph{without loss of generality}. A few other 
abbreviations will be introduced as we go along. Although, in its strict sense, the term \emph{Polish} refers only to a topology, we find convenient to use it as a synonym for a \emph{complete and separable metric} space. 
We use notation $\bbN\! :=\! \{ 0,1,2, \ldots \}$ for the set of natural numbers including $0$, and write $\bbN^*\!:=\! \bbN\backslash \{ 0\}$, $\bbR_+ \! :=\! [0, \infty)$ and $\bbR_+^*\! :=\! (0, \infty)$. Finally,  
all r.v.s are assumed to be defined on the same probability space $(\Omega, \mathscr F, \bP)$, 
which is sufficiently rich to carry as many independent r.v.s as we may need. 

To state our results precisely, we first need to recall the definition of an \textit{$\bR$-tree}. 
\begin{definition}
\label{realtrdef} A metric space $(T, d)$ is an \textit{$\bR$-tree} iff the following hold true.
\begin{compactenum} 

\smallskip

\item[$(a)$] For all $\sigma_1, \sigma_1 \! \in\!  T$, there is a unique isometry 
$f\colon [0,d(\sigma_1,\sigma_2)] \! \rightarrow \! T$ such
that $f(0)\!=\! \sigma_1$ and $f(d(\sigma_1,\sigma_2))\! =\! \sigma_2$. We shall use the notation  
$\lgeo \sigma_1,\sigma_2\rgeo \! :=\! f([0,d (\sigma_1,\sigma_2)])$. 

\smallskip

\item[$(b)$] For any continuous injective function 
$q\colon [0, 1] \! \rightarrow \! T$,  we have $q([0,1]) \! = \! \lgeo q(0), q(1)\rgeo$. 
\end{compactenum}

\smallskip

\noi
We distinguish a point $\rho\! \in \! T$ and we speak of $(T, d, \rho)$ as a \textit{rooted $\bR$-tree}.   \cq 
\end{definition}
Let us mention that we mostly focus on compact or Polish $\bbR$-trees. 
We next need to introduce basic definitions and notation attached to rooted $\bbR$-trees. 
\begin{definition}
\label{spandex} Let $(T, d, \rho)$ be a rooted $\bbR$-tree. 
Let $A \! \subseteq \! T$ be a non-empty subset.

\smallskip

\begin{compactenum} 

\item[$(a)$] The \emph{total height} of $A$ is  $\mathtt{Ht}_{d,\rho} (A) \eqo \sup_{\sigma \in A} d(\rho, \sigma) \in [0, \infty]$. 

\smallskip

\item[$(b)$] The \emph{subtree spanned by $A$} is the smallest connected subset containing $A$ and $\rho$. Namely, $\mathtt{Span} (A) \eqo \bigcup_{\sigma \in A} \lgeo \rho , \sigma \rgeo  $. We also denote by $\sspaan (A)$ the closure of $\spaan (A)$. 

\smallskip

\item[$(c)$] We say that a subtree $\ttt\! \subseteq \! T$ (i.e.~a connected subset containing $\rho$) is \emph{of finite type} if there is a finite subset $A\! \subseteq \! T$ such that $\ttt \eqo \mathtt{Span} (A)$. \cq

\end{compactenum}
\end{definition}
\noi
We observe that $(\spaan (A), d, \rho)$ and $(\sspaan (A), d, \rho)$ are rooted $\bbR$-trees. Let us mention from Lemma \ref{compacthull} that if $A$ is compact, then so is $\spaan (A)$. Finite-type $\bbR$-trees are necessarily compact: see Lemma \ref{edgelength} for a characterization. 
 
 Among metric spaces, $\bbR$-trees are characterized by the four-point condition. More precisely, we introduce the following. 
 \begin{definition}
\label{zerohypdef} A metric space $(S, d)$ is called a \textit{$0$-hyperbolic space} if for all $\sigma_1, \sigma_2, \sigma_3, \sigma_4 \ino S$,  
\begin{equation}
\label{4ptscondi}
d(\sigma_1 , \sigma_2 ) + d(\sigma_3 , \sigma_4 ) \leq \max \big( d(\sigma_1 , \sigma_3 ) + d(\sigma_2 , \sigma_4 ) 
 \, , \, d(\sigma_1 , \sigma_4 ) + d(\sigma_2 , \sigma_3 ) 
\big) \, , 
\end{equation}
which is called the \emph{four-point condition}. \cq 
\end{definition}
\begin{theorem}
\label{4ptsth} A metric space is an $\bR$-tree iff it is a connected $0$-hyperbolic space. 
The connected components of a $0$-hyperbolic space are $\bbR$-trees.
\end{theorem}
\noi
See Buneman \cite{Bun}. See also Evans \cite{EvansStF} for a self-contained proof.

Let us also recall from Proposition \ref{spanreafo} that any rooted $0$-hyperbolic space $(S,d,\rho)$ can be isometrically embedded into
an $\bbR$-tree $(\sspaan (S), d, \rho)$ that is minimal for inclusion, and if $(S,d)$ is Polish (resp.~compact), then so is $(\sspaan (S), d)$.

For any $\sigma \ino T$, we denote by 
$\bn ( \sigma, T)$ the number (possibly infinite) of connected components of $T\backslash \{ \sigma \}$. 
We define the \textit{leaves} and \textit{branch points} of $T$, as follows: 
\begin{equation}
\label{leafbranchdef}
 {\mathtt{Lf}} (T) \! :=\!   \{ \sigma \ino T \backslash \{ \rho \}  :  \bn (\sigma, T) \eqo 1 \, \}
\quad \textrm{and} \quad  {\mathtt{Br}} (T) \! :=\!  \{ \sigma \ino T \backslash \{ \rho \}  :  \bn (\sigma, T ) \geqo 3  \, 
\}\;  .
\end{equation}
For any $\sigma \ino T$ we denote by $\theta_\sigma T$ the subtree stemming from $\sigma$, i.e.~the \emph{fringe subtree} rooted at $\sigma$: 
\begin{equation}
\label{fringedef} \theta_\sigma T \eqo  \{ \sigma' \ino T  :  \sigma \ino \lgeo \rho , \sigma' \rgeo  \}\, .
\end{equation}
We observe that $\theta_\sigma T$ is closed and connected. Therefore $(\theta_\sigma T, d, \sigma)$ is a rooted $\bbR$-tree. 

For any $\sigma\ino T$ and $r\ino\bbR^*_+$, we denote open and closed balls in $T$ of radius $r$ around $\sigma$ by 
\begin{equation}
\label{ballnot}
 B^o_{T, d}(\sigma, r)\eqo  \big\{ \sigma^\prime\! \ino  T: d(\sigma,\sigma^\prime)\leko r \big\} \quad \textrm{and} \quad  B_{T,d}(\sigma, r)\eqo  \big\{\sigma^\prime\! \ino  T: d(\sigma,\sigma^\prime)\leqo r \big\} \; .
\end{equation}
When there is no ambiguity we denote these balls by $B^o_T(\sigma, r)$ and $B_T(\sigma, r)$ or even more simply by $B^o(\sigma, r)$ and $B(\sigma, r)$. 
We denote by $\cB (T)$ the sigma-field of the Borel subsets of $T$ and denote by 
$\cM_f (T)$ the set of finite positive measures $\mu\colon \cB(T)\! \to \!  \bbR_+$; when $(T,d)$ is Polish we shall always endow $\cM_f (T)$ 
with the topology of weak convergence. When equipped with $\mu\ino \cM_f(T)$, $(T,d, \rho, \mu)$ is called a \emph{rooted and (finitely) measured Polish $\bbR$-tree}, which is the main object of investigation in this article. 
 
\smallskip

We recall here the definition of \emph{length erasure} (also called \emph{trimming}) of a rooted compact $\bbR$-tree $(T,d,\rho)$ 
which is an important tool for approximating compact $\bbR$-trees by finite-type trees: for all $b\ino \bbR_+$ (we reserve $h\ino\bbR_+$ for mass erasure), the \emph{$b$-length-erased subtree} of $T$ is given by 
\begin{equation}
\label{deflengthera}
R_b (T)= \{\rho\} \cup \big\{ \sigma \ino T: \mathtt{Ht}_{d, \sigma} (\theta_\sigma T) \geqo b\big\}, 
\end{equation}
where $\mathtt{Ht}_{d, \sigma} (\theta_\sigma T)$ is the finite total height of the compact subtree $(\theta_\sigma T, d, \sigma)$. As proved by Evans, Pitman and Winter \cite[Lem.~2.6]{EPW} (and recalled in Lemma \ref{lengtheraprop} along with other important results of \cite{EPW} concerning length erasure), 
$(R_b (T), d, \rho)$ is a finite-type tree, $d_{\mathtt{Haus}} ( R_b(T), T)\leqo b$, where $d_{\mathtt{Haus}}$ stands for Hausdorff's distance on compact subsets of $T$ and length erasure enjoys the semigroup property $R_{b+b'} (T)\eqo R_b(R_{b'} (T))$, $b,b'\ino \bbR_+$.

Our goal in this article is to introduce an approximation semigroup analogous to length erasure 
but for rooted measured Polish $\bbR$-trees. To this end, we first introduce the following. 
\begin{definition}\label{deferasedtree} Let $(T, d, \rho)$ be a rooted Polish $\bR$-tree. Let $\mu \ino \cM_f (T)$ and let $h \in \bbR^*_+$. We set 
 $$T_{\! \mu, h}\! :=\!  \{ \rho \} \cup \big\{  \sigma \ino T  :   \mu (\theta_\sigma T) \geqo h \, \big\} . $$
We call $T_{\! \mu, h}$ the \emph{$h$-mass-erased tree associated with $\mu$}.  \cq
\end{definition}
As proved in Lemma \ref{erasedsubtree}, $T_{\!\mu,h}$ is a finite-type subtree of $T$. Proposition \ref{massTprop} asserts the following: 
$(1)$: $h\! \mapsto \! T_{\! \mu, h}$ is non-increasing, $d_{\mathtt{Haus}}$-left-continuous with a $d_{\mathtt{Haus}}$-right limit at all $h\ino \bbR_+$, the right limit $T_{\! \mu, h+}$ being the closure of $\bigcup_{h^\prime>h}T_{\! \mu,h^\prime}$ (here, $d_{\mathtt{Haus}}$  stands for the Hausdorff distance on the compact subsets of $T$); 
$(2)$: $\sspaan (\supp \mu) $ is the closure of $\bigcup_{h>0} T_{\! \mu,h}$ (here, $\supp \mu$ stands for the topological support of $\mu$); $(3)$: $(T_{\! \mu, h})_{h\in \bbR_+^*}$ characterizes $\mu$ (i.e.~if $T_{\! \mu, h}\eqo T_{\! \nu, h}$ for all $h\ino \bbR_+^*$ and some $\nu\ino \cM_f(T)$, then $\mu\eqo \nu$).

It is natural to equip $T_{\!\mu,h}$ with a measure. Projecting $\mu$ onto $T_{\! \mu,h}$ is a natural choice, which we explore in Section \ref{meratreesec}. However, the notion of mass erasure at the heart of our paper removes mass $h$ from every fringe subtree and arises from the following lemma.
\begin{lemma}
\label{erasedmassdefintro}
Let $(T, d, \rho)$ be a rooted Polish $\bR$-tree. Let $\bm \in \cM_f (T)$ and let $h \in \bbR^*_+$. 
Then there is a unique $\era_h \mu \ino \cM_f(T)$ such that 
\begin{equation}
\label{mainprop}
(\era_h \mu)  ( \theta_\sigma T) = \left( \mu (\theta_\sigma T) \! -\! h \right)_+ \; , \quad \sigma \ino T \; . 
\end{equation}
We call $\era_h \mu$ the \emph{$h$-erasure of $\mu$}. We also adopt the convention that $\era_0 \mu \! :=\!   \mu$. 
\end{lemma} 
\noi
We prove this in Proposition \ref{erasedmassdef}, 
starting from an explicit construction of $\era_h\mu$.

\smallskip

Let us collect some properties of mass erasure as an operation on a Polish $\mathbb{R}$-tree.
\begin{proposition}  
\label{intro:properased} Let $(T, d, \rho)$ be a rooted Polish $\bR$-tree. Then the following assertions hold true. 
\begin{compactenum}

\smallskip

\item[$(i)$] $(\era_h )_{h \in \bbR_+}$ is a semigroup on $\cM_f (T)$. Namely, $ \era_h \! \circ \! \era_{h'} \eqo  \era_{h+h'}$,  
$ h,h'\!  \ino \bbR_+ $. 

\smallskip

\item[$(ii)$] $\sup_{B \in \cB (T)} | \era_{h+h^\prime} \mu (B) \! -\! \era_{h} \mu (B) | \leqo 3 \mu(T) h'\! / h$, 
for all $h, h^\prime\!  \ino \bbR^*_+$ and all $\mu \ino \cM_f (T)$. 

\smallskip

\item[$(iii)$] $\sspaan (\mathtt{\supp}(\era_h \mu))\eqo T_{\! \mu, h+}$, $h\ino \bbR_+^*$, $\mu\ino \cM_f(T)$. 

\smallskip

\item[$(iv)$] The map $(h, \bm) \ino \bbR_+ \! \times \! \cM_f (T) \! \longmapsto \! \era_h \bm \ino \cM_f (T)$ is continuous. 

\smallskip

\item[$(v)$] Let $(h_n)_{n\in \bN}$ be an $\bbR_+^*$-valued sequence strictly decreasing to $0$ and let $(\nu_n)_{n\in \bN}$ be an $\cM_f (T)$-valued sequence such that $\era_{h_{n} \! -\! h_{n+1}} \nu_{n+1} \eqo \nu_n$ for all $n\ino \bN$. 
Then, there exists $\nu\ino\cM_f (T)$ such that $\era_{h_n}\nu\! = \! \nu_n$, $n\ino\bN$, iff  
\begin{equation}
\label{tightheight}
\sup_{n\in \bN} \nu_n \big( T\backslash B_{T,d}(\rho, r)) \underset{r\rightarrow \infty}{-\!\!\! -\!\!\!  - \!\!\! \longrightarrow} 0 \; .
\end{equation}

\end{compactenum}
\end{proposition}
\noi
The first four parts are proved in Propositions \ref{properased}--\ref{weakcvera}, the fifth part follows from Lemma \ref{descentera}.

A natural question is what happens when \eqref{tightheight} fails for a consistent family $(\nu_n)_{n\in\bN}$ of measures in $\cM_f(T)$. Clearly, this can only
happen if $T$ is unbounded. Informally, as mass erasure $\era_{h_{n} \! -\! h_{n+1}} \nu_{n+1} \! = \! \nu_n$ moves mass towards the root, passing from $\nu_n$ to $\nu_{n+1}$ moves mass away from the root 
and the supremum in \eqref{tightheight} is the mass that is eventually outside $B_{T,d}(\rho,r)$ as $n\! \to \! \infty$. If \eqref{tightheight} fails, the supremum 
decreases to a strictly positive limit as $r \! \to \! \infty$, and we can, in fact, capture such mass on the \emph{boundary of $T$}, at infinite distance from the root. 

To introduce the boundary of a rooted Polish $\bbR$-tree $(T,d,\rho)$, we need to define its \emph{rays}. Namely, 
a \emph{ray with endpoint $\rho$} is a subset $\mathbf r\!\subseteq\! T$ such that $\mathbf r\eqo c_{\mathbf r}(\bbR_+)$ for an isometry  
$c_{\mathbf r}\colon \bbR_+ \! \to \! T$ with  
$c_{\mathbf r} (0)\eqo \rho$. Then $c_{\mathbf r}$ is uniquely determined by $\mathbf r$. We note that $\mathbf r $ 
has infinite length, necessarily. 
We take the \emph{boundary $\partial T$ of $T$ as the set of rays with endpoint $\rho$}. 
If $\partial T\eqo \emptyset$, $T$ is said to be \emph{boundaryless}.

We next equip $T^*\! :=\! T \sqcup \partial T$ with a metric compatible with the \emph{cone topology}, as discussed 
e.g.~in Bridson and Haefliger \cite[Chapter II.8 and Chapter III.H]{briHae11} for the more general 
CAT($0$)- and $\delta$-hyperbolic metric spaces. Specifically, for all $\sigma\ino T$ we denote by 
$c_\sigma\colon\bbR_+ \! \to \! T$ the unique function such that $c_\sigma (0)\eqo \rho$,  
the restriction $c_\sigma\colon  [0, d(\rho, \sigma)]\! \to \! \lgeo \rho, \sigma\rgeo$ is a bijective isometry, and 
$c_\sigma (t) \eqo \sigma$ for all $t\ino [d(\rho, \sigma), \infty)$, and we set 
\begin{equation}
\label{dstarfirstdef}
\forall \, \mathbf s, \mathbf s'\! \ino T^*\! , \quad d^*(\mathbf s, \mathbf s')\! :=\! \! \int_0^\infty \! \!\! 
e^{-t} d \big(c_{\mathbf s} (t), c_{\mathbf s'}  (t) \big) \, dt . 
\end{equation}
Then Theorem \ref{bordif} asserts that $(T^ *\! ,d^*\! , \rho)$ is a boundaryless rooted Polish $\bbR$-tree, that 
$T$ is open and dense in $T^*\!$ (therefore $\partial T$ is a closed subset), that 
$\partial T \! \subseteq \! \mathtt{Lf} (T^*)$ and that $\partial T$ is totally disconnected. 
Moreover, $(T^*\! ,d^*)$ is compact iff $(T,d)$ is locally compact. Only the facts that 
$(T^*,d^*, \rho)$ is a $\bbR$-tree and that $\partial T \!\subseteq\! \mathtt{Lf} (T^*)$ may be regarded as new; 
the other statements  follow from more general results in Bridson \& Haefliger \cite[Chapter II.8 and Chapter III.H]{briHae11}
Following \cite{briHae11}, we call $(T^*\! ,d^*\! , \rho)$ the \emph{bordification} of $(T,d, \rho)$. 
Although strictly equivalent, the definition of the bordification of $\bbR$-trees provided in Section \ref{sec:boundary} 
differs from the above presentation, for various technical reasons. In particular, dealing with Gromov tree-spaces requires a more direct access to the metric $d^*$, which is rather viewed as a distortion of $d$. Theorem \ref{propdgtree} 
provides an alternative approach to the more general arguments in \cite{briHae11}, one that is specific to $\bbR$-trees and, quite understandably, somewhat simpler.
Let us also mention  Lemma \ref{embedmeas} that allows us to view $\cM_f(T)$, equipped with the topology of weak convergence on $(T,d)$, as topologically embedded in $\cM_f(T^*)$, equipped with the topology of weak convergence on $(T^*\! , d^*)$.

We observe that mass erasure makes sense on $\cM_f(T^*)$ since $(T^*\! , d^*\! , \rho)$ is a rooted Polish $\bbR$-tree and we next set 
\begin{equation}
\label{meraspace}
\cM_f^{\mathtt{era}}(T)\eqo \big\{ \mu \ino \cM_f(T^*): \forall \bs \ino \partial T, \; \mu (\{ \bs\})\eqo 0 \big\}\; , 
\end{equation}
i.e.~the set of finite Borel measures on $T^*$ that are diffuse on $\partial T$. 
Then Lemma \ref{descentera} first asserts for all $\mu \ino \cM_f(T^*)$, that $\era_h \mu \ino \cM_f(T)$ for all 
$h\ino \bbR^*_+$ iff $\mu \ino  \cM_f^{\mathtt{era}}(T)$. It also asserts that in the setting of Proposition \ref{intro:properased} (v), 
there exists a unique $\nu\ino \cM_f^{\mathtt{era}}(T)$ such that $\nu_n\eqo \era_{h_n} \nu$, $n\ino \bbN$, irrespective of \eqref{tightheight}.
\begin{example}
\label{evanescent} 
Let $(T,d, \rho)$ be 
the rooted $\bR$-tree obtained from the rooted complete infinite binary tree (i.e.~the set of finite words using a two-letter alphabet) by joining any two neighbouring vertices by a unit length geodesic. 
For all $n\ino \bbN$, we set 
 $\nu_n\! = \! 2^{-n}\sum \delta_\sigma $, where the sum is over all the 
$\sigma \ino \{ \rho \} \cup \mathtt{Br} (T)$ such that $d(\rho, \sigma ) \! \leq \! n$. We easily check that 
$\era_{2^{-n-1}} \nu_{n+1} \! = \! \nu_{n}$. 
By the above-mentioned property  (Lemma \ref{descentera} $(ii)$), there is a unique 
$\mu\ino \cM^{\mathtt{era}}_f(T)$ such that 
$\era_{2^{-n}}\mu\eqo \nu_n$ for all $n\ino \bbN$. 
This implies that $\mu(T^*)\eqo \nu_n(T)+ 2^{-n}\eqo 2^{-n} (2^{n+1} \! -\! 1)+ 2^{-n}\eqo 2$. Let $q\geko n$ and set $F_n\eqo T^* \backslash B^o_{T,d}(\rho, n)$ which is a closed subset of $T^*$ such that 
$T\cap F_n\eqo \{ \sigma \ino T: d(\rho, \sigma)\geqo n\}$. We thus see that 
$\nu_q(F_n)\eqo  2\! -\! 2^{-(q-n)}$. Since $\nu_q\! \to \! \mu$ weakly in $\cM_f(T^*)$, the Portmanteau theorem implies that $2\eqo \limsup_{q\to \infty} \nu_q(F_n)\leqo \mu(F_n) $. Thus $\mu(F_n)\eqo \mu(T^*)$ for all $n\ino \bbN$, which implies that $\mu(T)\eqo 0$, and indeed, the topological support of $\mu$ is $\partial T$. \cq 
\end{example}

It is then natural to introduce the following: let $\mu, \mu_n \ino \cM_f^{\mathtt{era}}(T)$, $n\ino \bbN$; we say that 
$(\mu_n)_{n\in \bbN}$ converges to $\mu$ \emph{in the sense of mass erasure} if $\era_h\mu_n \!\!  \to \! \era_h \mu $ weakly in $\cM_f(T)$ for all $h\ino \bbR^*_+$.  
Then, Theorem \ref{merafixpolish} $(i)$ 
asserts that the distance 
\begin{equation}
\label{dmerafix}
\delta_{\mathtt{era}} (\mu, \nu)\eqo \sum_{p\in \bbN} 2^{-p} \big( 1 \wedge \dPro^T(\era_{h_p} \mu, \era_{h_p} \nu)\big) , 
\end{equation}
metrizes the convergence in the sense of mass erasure. Here, $(h_p)_{p\in \bbN}$, is an $\bbR_+^*$-valued sequence that strictly decreases to $0$ and $\dPro^T$ stands for the Prokhorov distance on $T$. Moreover, Theorem \ref{merafixpolish} $(i)$ also shows that $( \cM_f^{\mathtt{era}}(T), \delta_{\mathtt{era}})$ is Polish and that 
$\cM_f(T)$ is $\delta_{\mathtt{era}}$-dense in $\cM_f^{\mathtt{era}}(T)$. Namely $\cM_f^{\mathtt{era}}(T)$ is the $\delta_{\mathtt{era}}$-completion of the space $(\cM_f(T), \delta_{\mathtt{era}})$. 

By Proposition \ref{intro:properased} $(iv)$, if $\mu_n\! \to \! \mu$ weakly in $\cM_f(T)$, then $\mu_n\! \to \! \mu$ in the sense of mass erasure. But the converse does not hold true, as is shown by the following example. 
\begin{example}
\label{dilutestar}
Let $(T,d, \rho)$ be a rooted Polish $\bR$-tree that is a star with a countably infinite number of unit length branches. Namely, 
$\mathtt{Lf} (T)$ is a countably infinite set listed as $(\sigma_n)_{n\in \bN}$ and we have 
$d(\rho, \sigma_n)\eqo 1$, $\lgeo \sigma_m, \rho \lgeo \, \cap\,  \lgeo \sigma_n, \rho \lgeo \,  =\! \emptyset $, for all distinct $m,n\ino \bN$ and $T\eqo \bigcup_{n\in \bN} \lgeo \rho, \sigma_n\rgeo$. For all $n\ino \bN$, we set $\mu_n \eqo \frac{1}{{n+1}} \sum_{0\leq k\leq n} \delta_{\sigma_k}$. Clearly $(\mu_n)_{n\in \bN}$ is not a weakly convergent sequence of probability measures. However, for $h\ino(\frac{1}{n+1},\infty)$, we see that 
$\era_h\mu_n\eqo(1\!-\!h)^+\delta_\rho\eqo\era_h\delta_\rho$. 
Hence $\mu_n \! \rightarrow \! \delta_\rho$ in the sense of mass erasure. \cq
\end{example}

To get a weak convergence in $\cM_f(T)$ from a convergence in the sense of mass erasure, an additional condition is needed and to formulate it we introduce for all $\mu \ino \cM_f(T)$ the \emph{height measure} $\lambda_\mu \ino \cM(\bbR_+)$ that 
captures the mass of the balls centered at $\rho$, namely it is the unique Borel measure on $\bbR_+$ such that 
\begin{equation}
\label{heightmeasdef}
\forall r\ino \bbR_+, \quad \lambda_{\mu} ([0, r])\eqo \mu \big( B_{T,d} (\rho, r)\big) 
\end{equation} 
(see Definition \ref{Lambdadef} $(b)$). Then Theorem \ref{eravsweakfix} shows for all $\mu_n, \mu \ino \cM_f(T)$, $n\ino \bbN$, that 
\begin{equation}
\label{fixedtreeeravsweak} \textrm{$\mu_n\! \to \mu$ weakly} \Longleftrightarrow \textrm{$\mu_n\! \to \mu$ in the sense of mass erasure and $ \lambda_{\mu} \!\! \to \!  \lambda_{\mu_n}$ weakly.} 
\end{equation}
To apply mass erasure to random continuum trees, in Section \ref{secappl}, we need to introduce the relavant Polish spaces of trees, namely the spaces of isometry classes of rooted Polish and compact finitely measured $\bbR$-trees equipped with their Gromov topologies. More precisely, in Section \ref{prelimGromovsec} we formally introduce the following spaces. 
\begin{compactenum} 

\smallskip

\item[$(a)$] $\bbT^0_{\! c}$ is the space isometry classes of (unmeasured) rooted compact $\bbR$-trees equipped with the \emph{Gromov--Hausdorff} distance $\dGH$.

\smallskip

\item[$(b)$] $\bbT^1_{\! c}$ is the space isometry classes of rooted finitely measured compact $\bbR$-trees 
equipped with the \emph{Gromov--Hausdorff--Prokhorov} distance $\dGHP$.

\smallskip

\item[$(c)$] $\bbT$ is the space of isometry classes of rooted finitely measured Polish $\bbR$-trees 
$(T,d,\rho, \mu)$ which are \emph{minimal}, i.e.~such that $\sspaan\,  (\supp \mu) \eqo T$ (see Definition \ref{minimal}). We endow $\bbT$ with the \emph{Gromov--Prokhorov} distance $\dGP$.

\smallskip

\end{compactenum}
\noi
Then we recall resp.~from Evans, Pitman and Winter \cite[Thm.~2]{EPW} for the first space, Abraham, Delmas and Hoscheit \cite[Cor.~3.2]{ADHoscheit} for the second one and L{ö}hr, Voisin and Winter \cite[Prop.~2.6]{LVW} for the third one, that 
\begin{equation}
\label{Treespacepolish}
\textrm{$ \big( \bbT^0_{\! c},  \dGH \big) $, $ \big( \bbT^1_{\! c}, \dGHP\big) $ and $\big( \bbT, \dGP \big)$ are Polish.} 
\end{equation}
See Definition \ref{mdef} and the comments right after for more details on $( \bbT, \dGP)$. In the whole article we write 
$(T,d,\rho, \mu) \! \equiv \! \fmu \ino \bbT$ or $\bbT^1_{\! c}$ to mean that $(T,d, \rho, \mu)$ is a representative of $\fmu$ 
(elements of $ \bbT$ or $\bbT^1_{\! c}$ are generically denoted by $\fmu$, $\fnu$, etc.); similarly 
$(T,d,\rho) \! \equiv \! \widetilde{T} \ino \bbT^0_{\! c}$ means that $(T,d, \rho)$ is a representative of $\widetilde{T}$. 
We also introduce the following subsets and functions which are useful to study the connections between $\bbT^0_{\! c}$, 
$\bbT^1_{\! c}$ and $\bbT$: 
\begin{compactenum}

\smallskip

\item[$-$] We denote by $\bbT^1_{\! c, \mathtt{min}}$ the set of $\fmu\! \equiv\! (T,d,\rho, \mu)$ in $\bbT^1_{\! c}$ that are minimal, i.e.~such that $T\eqo \sspaan (\supp \mu)$.  

\smallskip

\item[$-$] We denote by $\bbT_{\! c}$ the set of $\fmu  \! \equiv\!   (T,d,\rho, \mu) $ in $\bbT$ such that $(T,d)$ is compact. 

\smallskip

\item[$-$] For all $(T,d,\rho, \mu)  \! \equiv\!  \fmu\ino \bbT_{\! c}$, we denote by $\Phi_1 (\fmu)$ the isometry class in $\bbT^1_{\! c}$ of $(T,d, \rho, \mu)$ and we denote by $\Phi_0 (\fmu)$  the isometry class in $\bbT^0_{\! c}$ of $(T,d, \rho)$. 

\smallskip

\end{compactenum}
See Definition \ref{TcTmindef}. We mention that 
Lemma \ref{cpctsubsets} asserts that $\bbT_{\! c}$ is a Borel-measurable subset of $(\bbT, \dGP)$, 
$\bbT^1_{\! c, \mathtt{min}}$ is a Borel-measurable subset of $(\bbT^1_{\! c}, \dGHP)$, that 
$\Phi_1\colon \bbT_{\! c} \! \to \! \bbT^1_{\! c, \mathtt{min}}$ is a Borel-bi-measurable bijection and that $\Phi_0\colon\bbT_{\! c} \! \to \! \bbT^0_{\! c}$ is a Borel-measurable surjection. 

  For all $h\ino \bbR^*_+$ and for all $(T,d,\rho, \mu) \! \equiv \! \fmu \ino \bbT$ (resp.~$\bbT^1_{\! c}$), it makes sense to define $\era_h \fmu$ as the isometry class of $(T_{\! \mu, h+}, d, \rho, \era_h \mu)$ since this tree is minimal by Proposition \ref{intro:properased} $(iii)$ and since its isometry class only depends on that of $(T,d,\rho, \mu)$. Then,  
$\era_h\fmu  \ino \bbT_{\! c}$ (resp.~$\bbT^1_{\! c, \mathtt{min}}$). We set $\era_0 \fmu \! :=\! \fmu$ and observe that  
$(\era_h)_{h\in \bbR_+}$ defines a semigroup on $\bbT$ (resp.~on $\bbT^1_{\! c}$).  

Unlike what happens for a fixed tree, as mentioned above and formally stated in Proposition \ref{massTprop} $(ii)$, 
knowing  the function $h\ino \bbR_+^* \! \mapsto \! \widetilde{T}_{\! \mu, h}\! \equiv \! (T_{\! \mu, h}, d, \rho)$ 
does not generally enable us to determine $\fmu\ino \bbT$, as shown by the following example.
\begin{example}
\label{internalsym} 
We take $T\! = \! [-4, 4]$ equipped with the usual metric and the root $\rho\! = \! 0$ and we consider 
$\mu\eqo \delta_{1} +  3 \delta_{2}   +  \delta_{3}   +  3 \delta_{4} + 3 \delta_{-1} + \delta_{-2}   +  3 \delta_{-3}   +  \delta_{-4}
$ and $\nu \eqo \delta_{1} + 3 \delta_{2}   +  3 \delta_{3}   +  \delta_{4} +    3 \delta_{-1} + \delta_{-2}   +   \delta_{-3} 
  +  3 \delta_{-4} $. We denote by $\fmu$ and $\fnu$ the resulting elements of $\bT$. 
Then $\dGP (\fmu, \fnu) \! \neq \! 0$ (the average distance from the root is $\frac{1}{4}$ in $\fmu$ and 0 in $\fnu$). 
Viewed as intervals of $\bbR$, however observe that $T_{\! \mu , h}\! = \! -T_{\nu , h}$, for all $h \! \in \! (0, 3]$, and 
$T_{\! \mu , h}\! = \! T_{\nu , h}$, for all $h \! \in \! ( 3, 8]$ and trivially for all $h\ino(8,\infty)$.\cq
\end{example}
One of the main results of this article is the following. 
\begin{theorem}
\label{mainresult} If $\bbT$ is equipped with $\dGP$, $\bbT^{1}_{\! c}$ with $\dGHP$ and product spaces with product topologies, then the following holds true. 
\begin{compactenum}

\smallskip

\item[$(i)$] $(h,\fmu)\ino \bbR_+\!\! \times \! \bbT \! \mapsto \! \era_h \fmu\ino \bbT$ is continuous.

\smallskip

\item[$(ii)$] $(h,\fmu)\ino \bbR_+\!\! \times \! \bbT^1_{\! c} \! \mapsto \! \era_h \fmu\ino \bbT^1_{\! c}$ is continuous at any 
$h'\ino \bbR_+^*$ and any $(T'\! , d'\! , \rho' , \mu') \! \equiv \! \fmu'\ino   \bbT^1_{\! c}$ such that 
$T_{\! \mu'\! , h'\! +}'\eqo T_{\! \mu'\! , h'}'$. 
\smallskip
\end{compactenum}
\end{theorem}
\noi 
See Theorem \ref{eraGcont} for a more precise statement.

As in a fixed tree, this result suggests to introduce the notion of convergence in the sense of mass erasure. To this end, we introduce the space $\bbT^*$ of $\bbR$-trees equipped with their boundary.  More precisely, we introduce the following. 
\begin{compactenum}

\smallskip

\item[$-$] We say that $(T,d,\rho, \mu)$ is a $\ast$-$\bbR$-\emph{tree} if $(T,d,\rho)$ is a rooted Polish $\bbR$-tree and if $\mu\ino \cM^{\mathtt{era}}_{\! f} (T)$ such that $(T^*\! , d^*\! , \rho, \mu)$ is minimal, i.e.~$T^*$ is equal to the 
$d^*\!$-closure of the subtree of $T^*$ spanned by $\supp \mu $.

\smallskip

\item[$-$] Two $\ast$-$\bbR$-trees $(T_i,d_i,\rho_i, \mu_i)$, $i\ino \{ 1,2\}$, are 
are said to be $\ast$-equivalent if their respective bordifications $(T^*_i,d^*_i,\rho_i, \mu_i)$, $i\ino \{ 1,2\}$, are isometric.

\smallskip

\end{compactenum}
\noi
We then define $\bbT^*\! $ as the space of $\ast$-equivalence classes of $\ast$-$\bbR$-trees (see Definition \ref{Treestar}), 
and we shall view $\bbT$ as a subset of $\bbT^*\! $, as justified by Remark \ref{TinTstar}. 
For all $h\ino \bbR^*_+$, we extend $\era_h $ to a function from $\bbT^* $ to $\bbT \subo \bbT^* $ 
(see Definition \ref{masserastardef} and Remark \ref{justiferastar}). 

Let 
$\fmu_n , \fmu \ino \bbT^*$, $n\ino \bbN$. Then, we say that 
\begin{compactenum}

\smallskip

\item[$\;$]$\fmu_n \! \to \! \fmu$ \emph{in the sense of mass erasure} if $\dGP(\era_h \fmu_n, \era_h \mu) \! \to \! 0$ 
for all $h\ino \bbR_+^*$. 

\smallskip

\end{compactenum}
Mass-erasure convergence corresponds to a distance that is defined as follows: we fix 
$(h_p)_{p\in \bbN}$, an $\bbR_+^*$-valued sequence that strictly decreases to $0$ and we set 
\begin{equation}
\label{dderadef}
\forall \fmu, \fnu \! \in \! \bbT^*, \quad \dera (\fmu, \fnu)\! = \! \sum_{p\in \bN} 2^{-p} \big(
1\! \wedge\! \dGP (\era_{h_p}\fmu, \era_{h_p} \fnu )  \big)\; .
\end{equation}
Theorem \ref{cveracomplet} shows that $\dera$ is a distance which metrizes convergence in the sense of 
mass erasure, that $(\bbT^*\! , \dera)$ is Polish and that $\bbT$ is $\dera$-dense in $\bbT^*\! $. 
As (\ref{fixedtreeeravsweak}) for fixed trees, Theorem \ref{GPcvvsera} asserts for all $\fmu_n, \fmu\ino \bbT$, $n\ino \bbN$, that 
\begin{equation}
\label{eravsweak}
\Big( \dGP (\fmu_n, \fmu) \! \to \! 0 \Big) \; \Longleftrightarrow \; \Big(\,  \dera (\fmu_n, \fmu) \! \to \! 0 \; \,   
\textrm{and} \; \,   \flambda_{\fmu_n} \!\! \! \to \! \flambda_{\fmu} \, \textrm{weakly in $\cM_f(\bbR_+)$}\,  \Big),
\end{equation}
where for all $(T,d,\rho, \nu)\! \equiv\! \fnu \ino \bbT$, $\flambda_\fnu$ stands for the height measure 
$\lambda_\nu$ as defined in (\ref{heightmeasdef}), which makes sense since $\lambda_\nu$ only depends on the 
isometry class of $(T,d,\rho, \nu)$. 

 The previous deterministic results apply to the convergence in law of random continuum trees. To state our main 
 general results, let us introduce the following notation for the \emph{total mass} and \emph{total height} of $(T, d, \rho, \mu) \! \equiv \! \fmu\ino \bbT^*$: 
\begin{equation}
\label{massheight} 
\langle \fmu \rangle \eqo \mu(T^*) \quad \textrm{and} \quad \mathtt{Ht} (\fmu) = \sup_{\sigma  \in T} 
d(\rho, \sigma) \ino [0, \infty] \, , 
\end{equation}
\noi
which make sense since they only depend on the $\ast$-equivalence class of $(T, d, \rho, \mu)$. 
The following theorem collects our main general results on convergence in law in $(\bbT, \dGP)$ and $(\bbT^*\! , \dera)$. 
\begin{theorem}
\label{maincvlawera} 
Let $\mathbf m, \mathbf m',\mathbf m_n$, $n\ino \bbN$, be $\bbT^*$-valued r.v.s. 
Let $(h_p)_{p\in \bN}$ be an $\bbR^*_+$-valued sequence strictly decreasing to $0$. Then the following holds true.
\begin{compactenum}

\smallskip

\item[$(i)$] $\mathbf m$ and $\mathbf m'$ have the same law iff for all $p\ino \bbN$, $\era_{h_p} \mathbf m$ and $\era_{h_p} \mathbf m'$ have the same law. 

\smallskip

\item[$(ii)$]  Let us assume that $\era_{h_n-h_{n+1}}\mathbf m_{n+1}$ has the same law as $\mathbf m_n$. 
Then, there exists a $\bbT^*\! $-valued r.v.~$\mathbf m''$, unique in law and such that $\era_{h_n}  \mathbf m''$ has the same law as $\mathbf m_n$, for all $n\ino \bbN$.  

\item[$(iii)$] $\mathbf m_n \! \to \! \mathbf m$ in law in $(\bT^*\! , \dera)$ iff for all $p\ino \bbN$, 
$\era_{h_p} \mathbf m_n \! \to \! \era_{h_p} \mathbf m$ in law in $(\bT, \dGP)$.

\smallskip

\item[$(iv)$] The laws of $(\mathbf m_n)_{n\in \bbN}$ are tight in $(\bbT^*\! , \dera)$ iff the laws of 
$(\langle \mathbf m_n \rangle)_{n\in \bbN} $ are tight in $\bbR_+$ and for all fixed $p\ino \bbN$, 
the laws of $(\mathtt{Ht} (\era_{h_p} \mathbf m_n))_{n\in \bbN}$  are tight in $\bbR_+$. 

\smallskip

\item[$(v)$] Let us assume that $(\mathbf m_n)_{n\in \bbN}$ is 
$\bbT$-valued and that $\mathbf m_n \! \to \! \mathbf m$ in law in $(\bT^*\! , \dera)$. Then, the following assertions are equivalent.
\begin{compactenum}

\smallskip

\item[$(a)$] $\bP$-a.s.~$\mathbf m\ino \bbT$, and $\mathbf m_n \! \to \! \mathbf m$ in law in $(\bT, \dGP)$. 

\smallskip

\item[$(b)$] $\bP$-a.s.~$\flambda_{\mathbf m} \ino \cM_f(\bbR_+)$, and the random finite measures 
$(\flambda_{\mathbf m_n})_{n\in \bbN}$  on $\bbR_+$ converge to $ \flambda_{\mathbf m}$ in law in $\cM_f(\bbR_+)$. 
\smallskip
\end{compactenum}

\item[$(vi)$] For all $n\ino \bbN$, let us assume that 
a.s.~$(\mathbf T^{ n}, \mathbf d_n, \rho_n, \mathtt m_n)\! \equiv\! \mathbf m_n \ino \bbT$.  
Then the laws of $(\mathbf m_n)_{n\in \bbN}$ are tight in $(\bbT, \dGP)$ iff they are tight in $(\bbT^*\! , \dera)$ and 
\begin{equation}
\label{critightGP1bis}
\forall \epp, \eta \ino \bbR^*_+, \quad \lim_{h\to 0} \, \sup_{n\in \bbN} \bP \Big( \mathtt m_n \big( \mathbf T^{n} \backslash (\mathbf T^n_{\! \mathtt m_n , h})^{(\epp)}\big) \geko \eta \Big) = 0\; .
\end{equation}
where $(\mathbf T^n_{\! \mathtt m_n , h})^{(\epp)}\! :=\! \{ \sigma \ino \mathbf T^{n}: \mathbf d_n ( \sigma , \mathbf T^n_{\! \mathtt m_n , h})\leqo \epp\}$. 
\smallskip
\end{compactenum}
\end{theorem}
\noi
We prove this result in Lemma \ref{consistlaw}, Proposition \ref{traducvlawera}, Theorem \ref{condilawGP} and Corollary \ref{critightGP}.

Mass erasure also provides a convenient tool for studying the connections 
between convergences, in law and deterministic ones, on the spaces $(\bbT, \dGP)$, $(\bbT^0_{\! c}, \dGH)$ and $(\bbT^1_{\! c}, \dGHP)$, which is slightly more intricate than it may look at first glance, as the following example illustrates. 
\begin{example}
\label{GHGPnotGHPex} We denote by $\ell$ the Lebesgue measure on the Borel subsets of $\bbR$ and set $\mu \eqo \delta_{-1}+ 2\delta_1$. For all $n\ino \bbN$, we also set $T_n\eqo [-1,2]$ and $\mu_n\eqo \mu+ 2^{-n}\ell (\cdot \cap [1,2])$ if $n$ is even and we set 
$T_{n} \eqo [-2, 1]$ and $\mu_n\eqo \mu+2^{-n}\ell (\cdot \cap [-2,-1])$ if $n$ is odd. We equip $T_n$ with the usual distance $d$ on $\bbR$ and with the root $\rho\eqo 0$. Then $x\! \mapsto\!  -x$ is a root preserving isometry from $(T_n , d, \rho)$ onto $(T_{n+1} , d, \rho)$ and the sequence is trivially $\dGH$-converging.
We also check that $\dGP$-$\lim_{n\to \infty} (T_n, d,\rho, \mu_n)\! \equiv\!  ([-1, 1], d, \rho, \mu)$. However, $\dGHP(T_n, T_{n+1}) \geko 1/2$. Hence, $\dGP$-convergence and 
$\dGH$-convergence do not imply $\dGHP$-convergence, in general.  \cq
\end{example}
The following theorem sumarises the main general results on the connections between the convergences in law on $\bbT$, $\bbT^0_{\! c}$ and $\bbT^1_{\! c}$. 
\begin{theorem}
\label{T1T0vsTcvlaw} Let $(\mathbf m_n)_{n\in \bbN}$ and $\mathbf m$ be $\bbT$-valued r.v.s and $(h_p)_{p\in \bN}$ $\bbR^*_+\! $-valued decreasing to $0$.
\begin{compactenum}

\smallskip

\item[$(i)$] $\bP (\mathbf m\ino \bbT_{\! c})\eqo 1$ iff the laws of $(\Phi_0 (\era_{h_p} \mathbf m))_{p\in \bbN}$ are tight on $(\bbT^0_{\! c}, \dGH)$.

\smallskip

\item[$(ii)$] Let us assume that $\bP (\mathbf m_n \ino \bbT_{\! c})\eqo\bP (\mathbf m\ino \bbT_{\! c})\eqo 1$, $n\ino \bbN$. Then, 
$\Phi_1(\mathbf m_n) \! \to \!\Phi_1(\mathbf m)$ in law in $(\bbT^1_{\! c}, \dGHP)$ iff $\mathbf m_n \! \to \! \mathbf m$ in law in $(\bbT, \dGP)$ and $\Phi_0(\mathbf m_n) \! \to \!\Phi_0(\mathbf m)$ in law in $(\bbT^0_{\! c}, \dGH)$.
\end{compactenum}
\end{theorem}
\noi
We prove this in Lemma \ref{Tcrandomcriterion} and Theorem \ref{GPGH+GHP}.  

\smallskip

These results on $(\bbT^*\!, \dera)$, $(\bbT, \dGP)$ and $(\bbT^1_{\! c}, \dGHP)$ are used in Section \ref{secappl} to prove an invariance principle for measured Galton--Watson forests (GW forests for short). We consider measured GW forests characterised by a set of parameters $\xi, \varrho, c, \phi, \Lambda_n, \Lambda_{\varnothing, n}$, $n\ino \bbN$, whose meaning is as follows. 
\begin{compactenum}

\smallskip

\item[$-$] $\xi\eqo (\xi (k))_{k\in \bbN}$ is the offspring distribution, which is always assumed to be \emph{(sub)critical} and \emph{proper}, i.e.~$\sum_{k\in \bbN} k\xi (k) \leqo 1$ and $\xi (1)\eqo 0$. 

\smallskip

\item[$-$] $\varrho $ is the law of the total number of trees in the GW forest (i.e.~the degree of the root). \pagebreak[2]

\smallskip

\item[$-$] The GW forest is formally a $\bbT^1_{\! c}$-valued r.v.~$(\mathbf F, d, \rho, \mathrm m)\! \equiv \! \mathbf m$ that is 
a.s.~of \emph{finite type} (see Definition \ref{spandex} $(c)$). Conditionally given the discrete structure, branches, which are the connected components of 
$\mathbf F \backslash (\{ \rho\} \cup \mathtt{Br} (\mathbf F))$, are independent and exponentially distributed with parameter $c$. 

\smallskip

\item[$-$] Conditionally given $(\mathbf F, d, \rho)$, each branch $B\! :=\!  \, \rgeo \sigma, \sigma' \lgeo \,$ is equipped with an independent measure such that a.s. for all $\gamma \ino B$, 
$\mathrm m (\, \rgeo \sigma, \gamma \rgeo )\eqo S^{B}_{{d(\sigma, \gamma)}}$, where $(S^{B}_{s})_{s\in \bbR_+}$ is a subordinator with Laplace exponent $\phi (y)\eqo -\log \bE [\exp (-yS^{B}_{1})]$, $y\ino \bbR_+$.  

\smallskip

\item[$-$] Conditionally given $(\mathbf F, d, \rho)$ and the $\mathrm m$-measures of the branches, atoms 
$\mathrm m (\{ \sigma \})$ at $\sigma \ino \{ \rho\} \cup \mathtt{Br} (\mathbf F)$, are independent: the conditional law of $\mathrm m (\{ \sigma \})$ is $\Lambda_n (dy) $ if $\sigma \ino \mathtt{Br} (\mathbf F)$ and if $n$ is the outdegree of $\sigma$, i.e.~$n \eqo \mathtt n (\mathbf F, \sigma) \! -\! 1$; it is equal to $\Lambda_{\varnothing ,n} (dy)$ if $\sigma \eqo \rho$ and $n\eqo \mathtt n (\mathbf F, \rho)$, is the number of trees in the forest.  
\smallskip

\end{compactenum}
\noi
We refer to Definitions \ref{ellonemodel} and \ref{meaGWregedef} for more precise statements. 

The distribution of measured GW forests are precisely the laws on $\bbT^1_{\! c}$ of measured finite-type random trees satisfying the \emph{regenerative branching property}, as precisely stated and proven in Theorem \ref{branchingprop}. 
Proposition \ref{prop:eraGW} also shows that the class of measured GW forests as introduced above is stable under 
mass erasure. Namely, $\era_h \mathbf m$  is a measured GW forest whose parameters $\xi_h, \varrho_h, c_h, \phi_h, \Lambda^{\! h}_n, \Lambda^{\! h}_{\varnothing, n}$, $n\ino \bbN$, can be explicitly obtained from $h$, $\xi,\varrho, c, \phi, \Lambda_n, \Lambda_{\varnothing, n}$, $n\ino \bbN$. 

\emph{Measured Lévy forests} are the $\bbT$-valued random trees, which are not of finite type and which satisfy the regenerative branching property. In this article we do not prove such a characterization and we restrict to \emph{standard} measured Lévy forests, i.e.~measured Lévy forests whose total mass is distributed as the total population 
of a Continuous-State Branching Process (a CSBP for short). Their law is characterized by $x_0\ino \bbR_+^*$ and $\psi$ which are respectively 
the initial value and the \emph{branching mechanism} of a (sub)critical CSBP $(Z^{x_0}_s)_{s\in \bbR_+}$. Namely, $\psi$ is the Laplace exponent of a spectrally positive Lévy process $(X_s)_{s\in \bbR_+}$ that satisfies a.s.~$\liminf_{s\to \infty} X_s\eqo -\infty$, 
and it is necessarily of the following L{\'e}vy--Khintchine form, 
\begin{equation}\label{brmechintro}
\psi(y)=\mathbf{a}y+\tfrac{1}{2}\mathbf{b}y^2+\int_{\bbR^*_+}\!\! (e^{-y z}\! -\! 1+yz)\, \pi(dz)
\end{equation}
where $\mathbf{a}, \mathbf{b}\ino \bbR_+$ and where the 
Borel measure $\pi$ on $\bbR^*_+$ satisfies $\int_{\bbR_{+}^*} \! (z \! \wedge \! z^2)\, \pi(dz)\! <\! \infty$. We refer to the beginning of Section \ref{applLevy} for a more detailed account on CSBPs and spectrally positive Lévy processes. 
Following Sato \cite[Def.~11.9]{Sat99}), we introduce the following terminology.  
\begin{definition}
\label{ABCtypesdef} Let $\psi$ be a branching meachism as in (\ref{brmechintro}). We distinguish the following types. 

\begin{compactenum}

\smallskip

\item[$-$] Type A: $\;$ $\mathbf b\eqo 0$ and $\pi (\bbR^*_+) \! < \! \infty$,

\smallskip

\item[$-$] Type B: $\; $ 
$\mathbf b\eqo 0$, $\pi (\bbR^*_+) \! =\! \infty$ and $\int_{(0, 1)} \! z\, \pi(dz) \! < \! \infty$,

\item[$-$] Type C: $\; $ $\mathbf b\! >\! 0$ or $\int_{(0, 1)} \! z\, \pi(dz) \! =\! \infty$. \cq 

\smallskip

\end{compactenum}
\end{definition}
\noi
The branching mechanism $\psi$ is of Type C iff the corresponding L{\'e}vy processes or CSBPs 
have infinite-variation sample paths. 
We do not consider branching mechanisms of Type A, which correspond to compensated Poisson processes. 

In Theorems \ref{thm:growthprocconvbis} and \ref{thmIP}, we fix $x_0\ino \bbR_+^*$ and a branching mechanism $\psi$ of Types B or C, and we define a standard measured Lévy forest as a $\bbT$-valued r.v.~$ \mathbf m $, whose law is the sole  limit in distribution on $(\bbT^*\! , \dera)$ of sequences of $\bbT$-valued random trees $(\mathbf F_{\! p}, d_p, \rho_p, \mathrm m_p)\! \equiv\! \mathbf m_p$, $p\ino \bbN$, satisfying the following properties. 
\begin{compactenum}

\smallskip

\item[$(i)$] The $\mathbf m_p$ are minimal finite-type random trees satisfying the regenerative branching property (namely, they are measured GW forests as above).

\smallskip

\item[$(ii)$] The measure $\mathrm m_p$ is a deterministic fucntion of the metric space $(\mathbf F_{\! p}, d_p, \rho_p)$ and it has no mass at the root: a.s.~$\mathrm m_p (\{ \rho_p\})\eqo 0$.  

\smallskip

\item[$(iii)$] Total masses $\langle \mathbf m_p \rangle$ converge in law to the total population $\int_0^{\infty} \! Z^{x_0}_{s} ds $ of the CSBP($x_0, \psi$) $(Z^{x_0}_s)_{s\in \bbR_+}$. 

\smallskip

\end{compactenum}

More precisely, we denote by $\xi_p, \varrho_p, c_p, \phi_p, \Lambda^{{\!(p)}}_{n}, \Lambda^{{\!(p)}}_{{\varnothing, n}}$, $n\ino \bbN$, the parameters of $\mathbf m_p$ and we assume the following. 

\smallskip

\noi
$\mathbf{(1)}$ $\sum_{k\in \bbN} k \xi_p (k) \leqo 1$ and $\xi_p (1)\eqo 0$. 

\smallskip

\noi
$\mathbf{(2)}$ 
    There exists a sequence $(a_p)_{p\in\mathbb{N}}$ such that
$$ a_p\longrightarrow\infty,\quad 
\textrm{and for all $s\ino \bbR_+$} \quad 
\tfrac{1}{a_p}Z_s^{{(p)}}  \xrightarrow[p\to \infty]{\textrm{(law)}}Z^{x_0}_s ,$$
where $Z_{{ s}}^{{(p)}}\! :=\! \lim_{\epp \downarrow 0}\# \{ \sigma\ino \mathbf F_{\! p}: d_p(\rho_p, \sigma)\eqo s+\epp \}$, is the GW branching process of $\mathbf F_{\! p}$. At $s\eqo 0$, this implies that $\varrho_p (dy/a_p)\! \to \!\delta_{x_0}(dy)$ weakly on $\bbR_+$.

\smallskip

\noi
$\mathbf{(3)}$ $\Lambda^{{(p)}}_{{\varnothing, n}} (dy)\eqo \delta_{0} (dy)$ for all $n\ino \bbN$ (no mass at the root). 

\smallskip

\noi
$\mathbf{(4)}$ $\phi_p (y)\eqo \kappa_p y$, $y\ino \bbR_+$, for some (deterministic) $\kappa_p\ino\bbR_+$, 
and  $\Lambda^{{(p)}}_{n} (dy)\eqo \delta_{w_p (n-1)} (dy)$, $n\ino\bN$, for a certain weight function $w_p\colon \{-1\}\cup \bbN\! \! \to \! \bbR_+$, such that the following holds: let $(Y^{{(p)}}_{k})_{k\in \bbN^*}$ be independent r.v.s such that $\bP (Y^{{(p)}}_{k}\eqo j) \eqo \xi_p (j+1)$, for all $j\ino  \{-1\}\cup \bbN$; then 
$$a_p\kappa_p+\sum_{1\leq k\leq \lfloor a_pc_p \rfloor} w_p \big(Y^{{(p)}}_{k} \big)   \xrightarrow[p\to \infty]{\textrm{(law)}} 1.$$ 
\noi
Under $\mathbf{(1)}$--$\mathbf{(4)}$, Theorem \ref{thmIP} asserts the following. 
\begin{compactenum}

\smallskip

\item[$-$] There is a $\bbT$-valued r.v.~whose law only depends on $(x_0, \psi)$ and such that $\mathbf m_p \! \to \! \mathbf m$ in distribution in $(\bbT^*\! , \dera)$. \emph{The law of $\mathbf m$ is that of a standard measured $(x_0, \psi)$-Lévy forest}.

\smallskip

\item[$-$] $\mathbf m_p \! \to \! \mathbf m$ holds in distribution on $(\bbT , \dGP)$ iff $\psi$ is of Type C. 

\smallskip

\item[$-$]  The $\mathbf m_p$ converge in distribution on $(\bbT^1_{\! c}, \dGHP)$ iff the laws of the total height 
$(\mathtt{Ht} (\mathbf F_p))_{p\in \bbN}$ are tight on $\bbR_+$. In that case, $\bP (\mathbf m\ino \bbT_{\! c})\eqo 1$ and 
 $\mathbf m_p \! \to \! \Phi_1(\mathbf m)$ holds in distribution on $(\bbT^1_{\! c}, \dGHP)$. 

\smallskip

\end{compactenum}

Even in the compact cases, Theorem \ref{thmIP} is new in its extent. For Type C branching 
mechanisms $\psi$, standard Lévy forests correspond to the continuum random trees encoded by $\psi$-height 
processes, introduced in Le Gall and Le Jan \cite{LGLJ1} and in Duquesne and Le Gall \cite{DuLG}. We also refer to Duquesne and Rebei \cite{DuRe26prepub} for the tightness criterion in the third claim, which we state in Lemma \ref{Greydis}. For Type B branching 
mechanisms 
$(x_0, \psi)$-Lévy forests seem new. In these cases, the finite measure equipping the tree is a sum of atoms at 
branch points and the root (although the approaching GW forests $\mathbf m_p$ have no mass at their root). Moreover, in Type B cases $(\mathbf m_p)_{p\in \bbN}$ do not converge in distribution on $(\bbT, \dGP)$. Measured GW forests can converge in distribution in $(\bbT, \dGP)$ to Type B standard Lévy forests, but they must be equipped with measures which are not deterministic functions of their metric and whose extra randomness does not vanish in the limit. A direct construction of Type B Lévy forest is provided in Remark \ref{remconstructLevy} $(d)$. 
We refer to Remarks \ref{remconstructLevy} and \ref{IPcomments} for more detailed comments on Lévy forest and the invariance principle of Theorem \ref{thmIP}. 

%

\section{$\bbR$-trees and their boundary}
\label{sec-boundary}
\subsection{Basic definitions and properties of $\bbR$-trees}

We recall Definition \ref{realtrdef} of rooted $\mathbb{R}$-tree. Informally $\bR$-trees are obtained by gluing intervals of 
$\bbR$ endowed with their usual metric without creating loops (hence their name). 
Let us also mention that metric $\bR$-trees are specific instances of topological real trees that have been studied in Bowditch's article \cite{Bow99} from which we recall the following properties.
\begin{lemma} 
\label{elemBowd} Let $(T,d)$ be an $\bbR$-tree. Then the following hold true. 
\begin{compactenum}

\smallskip

\item[$(i)$] Any connected subset is path-connected, by geodesic arcs (i.e.~isometries).

\smallskip

\item[$(ii)$] The connected components of an open set are open. 

\smallskip

\item[$(iii)$] Let $\sigma_0, \sigma_1, \sigma_2 \ino T$. 
Then, there is $\sigma \ino T$ such that $\lgeo \sigma_0, \sigma \rgeo\eqo \lgeo \sigma_0, \sigma_1 \rgeo \cap \lgeo \sigma_0, \sigma_2 \rgeo$. Moreover, $\{ \sigma\} \eqo 
 \lgeo \sigma_0, \sigma_1 \rgeo \cap \lgeo \sigma_0, \sigma_2 \rgeo \cap \lgeo \sigma_1, \sigma_2 \rgeo $. 
\end{compactenum}
\end{lemma}
\noi
\textbf{Proof.} See Bowditch \cite[Lem.~1.1, 1.4 and 1.5]{Bow99}. \cqfd

\smallskip

In Lemma \ref{elemBowd} $(iii)$, we call $\sigma$ the \emph{branch point of $\sigma_0, \sigma_1, \sigma_2$} and we denote it by $\mathtt{br} (\sigma_0, \sigma_1, \sigma_2)$. 
If the tree is rooted at $\rho$, we use the notation $\sigma_1\wedge_\rho \sigma_2 \! :=\!  \mathtt{br} (\rho, \sigma_1, \sigma_2)$.
This actually depends on $\rho$, however we simply denote it by $\sigma_1\wedge \sigma_2 $ when there is no ambiguity on the choice of the root. 
If we view $T$ as a family tree whose root $\rho$ is the ancestor, then $\sigma_1\wedge \sigma_2 $ is the \emph{most recent common ancestor of $\sigma_1$ and $\sigma_2$}. 
More generally, the \emph{most recent common ancestor} of a non-empty subset of an $\bbR$-tree is defined as follows.
\begin{lemma}
\label{mrcadef} Let $(T, d, \rho)$ be a rooted $\bbR$-tree. 
Let $\emptyset  \! \neq \! A \! \subset \! T$. Then, 
there exists a unique point $\sigma_* \ino T$ such that 
$\lgeo \rho , \sigma_*  \rgeo\eqo  \bigcap_{\sigma \in A} \lgeo \rho , \sigma \rgeo $. 
Moreover, for any $\sigma \ino A$, there are $\gamma_n \ino A$, $n\ino \bN$, such that $d( \sigma_* , \gamma_n \!  \wedge \! \sigma) \! \rightarrow \! 0$ as $n\! \rightarrow \! \infty$. 
The point $\sigma_*$ is called the \emph{most recent common ancestor of $A$} and is denoted by $\mathtt{mrca} (A)$.  
\end{lemma}
\noi
\textbf{Proof.} Let $\sigma \ino A$ and set $B \eqo  \bigcap_{\gamma \in A} \lgeo \rho , \gamma \rgeo $. 
Then $\rho \ino B$ and $B$ is a compact and connected subset of $\lgeo \rho , \sigma \rgeo$ which easily 
implies the first statement of the lemma. We then prove that 
$d(\rho, \sigma_*)\eqo \inf_{\gamma \in A} d(\rho, \gamma \! \wedge \!  \sigma)$, which immediately 
entails the second statement of the lemma. 
\emph{Indeed}, we first note for all $\gamma \ino A$ that 
$\lgeo \rho, \sigma_*\rgeo \subo \lgeo \rho , \gamma \!  \wedge \!  \sigma \rgeo$ and thus   
$d(\rho, \sigma_*)\leqo \inf_{\gamma \in A} d(\rho, \gamma  \! \wedge \!  \sigma)$. If the inequality is strict, then there is 
$\sigma'\! \in \,   \rgeo \sigma_*, \sigma\rgeo $ such that $d(\rho, \sigma') \leko $ $\inf_{\gamma \in A} d(\rho, \gamma \!  \wedge \!  \sigma)$. This entails  $\lgeo \rho,  \sigma' \rgeo \subo \bigcap_{\gamma \in A} \lgeo \rho , \gamma \rgeo$, which contradicts the definition of $\sigma_*$. \cqfd 

\smallskip 
 
We recall Definition \ref{spandex} $(a)$ of the total height of a subset, from (\ref{leafbranchdef}) the notions of leaves and branch points, from (\ref{fringedef}) the definition of fringe subtrees and  from (\ref{ballnot}) the notation for open and closed balls. 
\begin{lemma} 
\label{Brcount} If an $\bbR$-tree $(T,d)$ is separable, the set $\mathtt{Br} (T)$ of branch points is at most countable. 
\end{lemma}  
\noi
\textbf{Proof.} See Bowditch \cite[Lem.~1.6]{Bow99}. \cqfd

\begin{lemma}
\label{AboTrprop} Let $(T, d, \rho)$ be a rooted $\bbR$-tree and let $\sigma , \sigma^\prime \ino T$. 
Then the following holds true. 
\begin{compactenum}

\smallskip

\item[$(i)$] The fringe subtree $\theta_\sigma T$ above $\sigma$ is closed and connected (i.e.~it is a subtree of $T$).

\smallskip

\item[$(ii)$] If $\sigma\ino \lgeo \rho, \sigma^\prime \rgeo$, then $\theta_{\sigma^\prime} T\!\subseteq \! \theta_\sigma T$. If $\sigma\wedge \sigma^\prime \! \notin \! \{ \sigma, \sigma^\prime\}$, then $\theta_\sigma T \cap \theta_{\sigma^\prime} T \eqo \emptyset$.  

\smallskip

\item[$(iii)$] For any $\sigma \ino T\backslash\{ \rho\}$, $T \backslash \theta_\sigma T$ is the connected component of $T\backslash \{ \sigma\}$ containing $\rho$.

\smallskip

\item[$(iv)$] Let $C$ be a connected component of $T\backslash \{ \sigma \}$  such that $\rho\! \notin \! C$. 
For all $\sigma'\ino C$, $\rgeo \sigma , \sigma' \rgeo \! \subseteq \!  C$, $\theta_{\sigma'} T\! \subseteq \! C$ and for all $\sigma_n \! \in \, \rgeo \sigma , \sigma' \rgeo$ such that $d(\sigma, \sigma_n) \! \to \! 0$, we also get $C\eqo \bigcup_{n\in \bN} \theta_{\sigma_n} T$. 

\smallskip

\item[$(v)$] For all $\sigma \ino T$ and $r\ino \bbR_+$, let $C$ be a connected component of $T\backslash B_{T,d}(\sigma, r)$. Then there is $\gamma \ino T$ such that $d(\rho, \gamma) \eqo r$ and $C$ is a connected component of $T\backslash \{ \gamma \} $. 

\smallskip

\item[$(vi)$] Let us suppose that $D \! \subseteq\!  T$ is countable and dense, and set $ \mathscr P \eqo 
\{ \emptyset \}\!  \cup \! \big\{ \theta_\sigma T \, ;  \sigma \ino T \backslash D \big\}$. Then, $\mathscr P$ is a pi-system generating the Borel sigma-field of $T$. 
\end{compactenum}
\end{lemma}
\noi
\textbf{Proof.} We leave the proof to the reader. \cqfd 

\smallskip

Let us recall Definition \ref{spandex} of $\mathtt{Ht}_{d,\rho} (A)$ and $\spaan (A)$, for a non-empty subset $A$ of a rooted $\bbR$-tree 
$(T,d, \rho)$. For all $\varepsilon \in \bbR^*_+$, we also 
denote by $N(A, \varepsilon)$ the minimal number of open balls centered in $A$ with radius $\varepsilon$ necessary to cover $A$. Note that $N(A, \varepsilon)$ or $\mathtt{Ht}_{d,\rho} (A)$ may be infinite. 
When $A$ is a closed subset, $\mathtt{Span} (A) $ is not necessarily closed (see Example \ref{lafuma} for instance).  
However, the cases of compact subsets have the following nice property.  
\begin{lemma} 
\label{compacthull}
Let $(T, d, \rho)$ be a rooted $\bbR$-tree, $A\! \subseteq \! T$ non-empty and $\varepsilon \in \bbR^*_+$. Then we have $N (\mathtt{Span} (A) , \varepsilon ) \leqo N(A,   \varepsilon )  ( 1+ 
\tfrac{1}{2\varepsilon}\mathtt{Ht}_{d,\rho} (A) )$, 
and  the closure of $A$ is compact iff the closure of $ \mathtt{Span} (A)$ is compact. Moreover $A$ is compact iff $A$ is closed and $ \mathtt{Span} (A)$ is compact. 
\end{lemma}
\noi
\textbf{Proof.} We leave the proof to the reader. \cqfd 


\smallskip

We recall Definition \ref{spandex} $(c)$ of $\bbR$-trees of finite type. 
\begin{lemma}
\label{edgelength} Let $(T, d, \rho)$ be a rooted $\bbR$-tree and let $\rho \ino \ttt\subo T$ be connected. 
Then $\ttt$ is of finite type iff the following three conditions are satisfied: $\ttt$ is closed, $\sup_{\sigma\in \ttt} d(\rho, \sigma) \leko \infty$ and
\begin{equation}
\label{bbound}
\!   \sup   \big\{ n \ino \bbN^*   :   \exists \,  \sigma_1, \ldots, \sigma_n \ino \ttt \; \textrm{\emph{such that}} \; \sigma_i \! \wedge \! \sigma_j \! \notin \! \{\sigma_i , \sigma_j   \} , \, 1\leqo i\! < \! j\leqo n    \;   \big\} < \infty .
 \end{equation}
\end{lemma}
\noi
\textbf{Proof.}  See Appendix \ref{Proofedgelength}. \cqfd

\subsection{$0$-hyperbolic spaces}
Let us recall Definition \ref{zerohypdef} of $0$-hyperbolic spaces and from 
Theorem \ref{4ptsth} that $\bbR$-trees are the connected $0$-hyperbolic spaces. 
We next recall mostly from Evans \cite{EvansStF} the following proposition which asserts that any $0$-hyperbolic space can be isometrically embedded into a minimal $\bR$-tree (unique up to isometry), which is called its \emph{spanning $\bR$-tree}. 
\begin{proposition}
\label{spanreafo} Let $(S, d)$ be a $0$-hyperbolic space and let $\rho \! \in \! S$.  Then the following holds.

\smallskip

\begin{compactenum}

\item[$(i)$] There is a rooted complete $\bR$-tree $(T, \delta , r)$ and an isometrical embedding $\varphi\colon (S, d, \rho) \! \hookrightarrow \! (T, \delta , r)$.   

\smallskip

\item[$(ii)$] For $i \! \in \! \{ 1, 2\}$, let $(T_i, d_i, \rho_i)$ be rooted Polish $\bR$-trees and 
let $\varphi_i\colon (S, d, \rho) \hookrightarrow (T_i, d_i, \rho_i)$ be isometrical embeddings. Let  
$\ttt_i \! =\!  \sspaan(\varphi_i(S))$ be the closure of the spanning subtree of $\varphi_i (S)$ in $T_i$. Then, there exists a bijective isometry $\phi\colon (\ttt_1, d_1, \rho_1) \rightarrow  (\ttt_2, d_2, \rho_2)$.   
\end{compactenum}

\smallskip

\noi
This allows us to define a rooted complete $\bR$-tree 
$(\sspaan (S), d, \rho)$, unique up isometry, and such that for all isometrical embedding $\varphi$ from $(S,d,\rho)$ into 
a rooted complete $\bbR$-tree $(T, \delta , r)$, the space $( \sspaan(\varphi(S)), \delta, r)$ is isometric to $(\sspaan (S), d, \rho)$.
We shall always assume that $S\subseteq \sspaan (S)$. 

We furthermore get the following. 
 
\begin{compactenum}

\smallskip

\item[$(iii)$] If $(S, d)$ is Polish, so is $(\sspaan (S), d)$. 
Moreover, the completion of $(S, d)$ is compact iff $(\sspaan (S), d)$ is compact. 
\end{compactenum}
 \end{proposition}
\noi
\textbf{Proof.} For $(i)$, see Evans \cite[Thm.~3.38]{EvansStF}. 
Assertion $(ii)$ is easy to check: we leave the details to the reader. Let us briefly prove $(iii)$: to simplify notation, we set $T\! :=\! \sspaan (S)\!\supseteq\! S$ and $T^o\!\! =\!  \bigcup_{x\in S} \lgeo \rho, x \rgeo$, whose closure is $T$. We only need to prove that $T^o$ is separable: if 
$(x_n)_{n\in \bN}$ is $S$-valued and dense, then $Q\eqo \bigcup_{n\in \bN} \lgeo \rho, x_n \rgeo$ clearly is separable, 
and we only need to prove that any point $\gamma \in T^o \backslash S$ is the limit of a $Q$-valued sequence. To this end, let  $x\ino S$ be such that $\gamma \! \in \! \lgeo \rho, x \lgeo \, $. There is an increasing $\bbN$-valued sequence $(n_p)_{p\in \bN}$ such that $d(x,x_{n_p}) \! \rightarrow \! 0$, which entails $d(x, x\! \wedge\!  x_{n_p})\!  \rightarrow \! 0$. Thus $\gamma \! \in \! \lgeo \rho, x_{n_p} \rgeo\!\subseteq\!Q$ as soon as $d(x, x\wedge x_{n_p}) \leqo  d(x, \gamma) $. The second assertion of $(iii)$ follows from Lemma \ref{compacthull}.  
\cqfd

\begin{example}
\label{lafuma} Even if a $0$-hyperbolic space $(S, d, \rho)$ is complete, $\mathtt{Span} (S)$ may not be complete. Indeed, 
let us denote by $(e_n)_{n\in \bN}$ the ``canonical'' basis of $\ell^{1} (\bN)$ and let us consider 
the complete $0$-hyperbolic space $S \! =\!  \{0 \} \cup \{ 
(1\! -\! 2^{-n})e_1+ e_n\,  ;\,  n \geqo 2\}$, when it is equipped with the $\lVert \cdot \rVert_1$-distance and rooted at the origin $0$. Then,  
$\{ e_1\} \! = \! \sspaan (S)   
\backslash \mathtt{Span} (S)$. \cq    
\end{example}

\subsection{Deformation of $\bbR$-tree metrics, $\bbR$-tree boundary}
\label{sec:boundary}

We next introduce the boundary of an $\bbR$-tree equipped with the \emph{cone topology}: for Polish $\bbR$-trees, this procedure is closely related to the \emph{bordification} of 
CAT($0$)-spaces and $\delta$-hyperbolic spaces as discussed for instance in Bridson and Haefliger \cite{briHae11}, Chapter II.8 and Chapter III.3. Here we use a direct approach, based on the following lemma which is about deformations of the metric of an $\bbR$-tree. We recall Definition \ref{spandex} $(a)$ of the total height. 
\begin{lemma} 
\label{fchangemetric} Let $a,b\ino (0, \infty]$ and let $g\colon[0, a) \! \to \! [0, b)$ be an increasing homeomorphism. We fix a rooted $\bbR$-tree $(T,d,\rho)$ such that $\mathtt{Ht}_{d, \rho} (T) \leqo a$.  We set $T^o\! \! :=\! B^o_{T,d} (\rho, a)\eqo \{ \sigma \ino T: d(\rho, \sigma) \leko a \}$ and  
\begin{equation}
\label{fmodif}
\forall \sigma, \sigma' \ino T^o, \quad d_g(\sigma, \sigma')\eqo g \big( d(\rho, \sigma) \big)+ g \big( d(\rho, \sigma') \big) \! - \! 2g \big(d(\rho, \sigma \! \wedge \! \sigma') \big).
\end{equation}
Then $(T^o\! , d_g)$ is an $\bbR$-tree and the identity function on $T^o$ is a $(d,d_g)$-homeomorphism. 
\end{lemma}
\noi
\textbf{Proof.} See Appendix \ref{Prooffchangemetric} \cqfd

\begin{notation}
\label{dgtreedef} We keep the assumptions and the notation of Lemma \ref{fchangemetric}.
\begin{compactenum}

\smallskip

\item[$(a)$] We denote by $(T_{\! g}, d_g)$ a completion of $(T^o\! , d_g)$ such that $T^o\! \subseteq\! T_{\! g}$. Consequently, $\rho\ino T_{\! g}$ and $(T_{\! g}, d_g, \rho)$ is a rooted complete $\bbR$-tree (e.g.~as a consequence of Proposition \ref{spanreafo}).     

\smallskip

\item[$(b)$] Since geodesic arcs and branch points in $(T^o\! , d)$ and $(T^o\! , d_g)$ are the same, we denote them in $(T_{\! g}, d_g)$ as in $(T,d)$, namely by $\lgeo \mathbf s,\mathbf s' \rgeo$ and $ \mathbf s\!  \wedge \! \mathbf s'$, for all $\mathbf s, \mathbf s'\ino T_{\! g}$.  
\cq

\end{compactenum}
\end{notation}

\pagebreak[2]

\begin{theorem} 
\label{propdgtree}  In the setting of Lemma \ref{fchangemetric}, let us denote by $f$ the inverse of $g$ and 
assume that $(T,d)$ is complete. 
Then the following holds true. 
\begin{compactenum}

\smallskip

\item[$(i)$] We have $T^o\!\! =\! \{\mathbf s \ino T_{\! g}\! :\!  d_g(\rho, \mathbf s) \leko b\}$, and $T^o$ is open and dense in $(T_{\! g}, d_g)$. In particular, $(T_{\! g}, d_g)$ is separable iff $(T,d)$ is separable.

\smallskip

\item[$(ii)$] Let $\mathbf s\ino T_{\! g}$. Then, $\lgeo \rho, \mathbf s \lgeo \, \subseteq \! T^o$, and if 
$\mathbf s \ino T_{\! g} \backslash T^o$, then $d_g(\rho, \mathbf s)\eqo b\leko\infty$ and $\mathbf s \ino \mathtt{Lf} (T_{\! g})$. Moreover, for all $\mathbf s, \mathbf s'\ino T_{\! g}$, if $\mathbf s\! \wedge \! \mathbf s'\ino T_{\! g} \backslash T^o$, then $\mathbf s\eqo \mathbf s'\ino T_{\! g} \backslash T^o$. In particular, $T_{\! g} \backslash T^o$ is totally disconnected.

\smallskip

\item[$(iii)$] $\mathtt{Ht}_{d_g, \rho} (T_{\! g})\leqo b$ and $\big((T_{\! g})_f, (d_g)_f, \rho \big)$ is isometric to $(T,d, \rho)$.

\smallskip

\item[$(iv)$] $(T_{\! g}, d_g)$ is locally compact iff $(T,d)$ is locally compact. In that case, if $b\leko \infty$ or $\mathtt{Ht}_{d,\rho} (T) \leko a$, then $(T_{\! g}, d_g)$ is compact. 
\end{compactenum}

\end{theorem}
\noi
\textbf{Proof.} Let us prove $(i)$. Clearly, $T^o\!\subseteq\!  \{\mathbf s \ino T_{\! g}\! :\!  d_g(\rho, \mathbf s) \leko b\}$. 
Let $\mathbf s\ino T_{\! g}$ be such that $d_g(\rho, \mathbf s) \leko b$. By definition, there is a $T^o$-valued sequence 
$(\sigma_n)_{n\in \bbN}$ such that $d_g(\sigma_n, \mathbf s) \! \to \! 0$. W.l.o.g.~we can assume that 
$\sup_{n\in \bbN} d_g(\rho, \sigma_n)$ $+$ $ \sup_{p,q\in \bbN} d_g (\sigma_p, \sigma_q) \leqo b_0$, for some 
$b_0\leko b$. We also introduce the following notation: $w_f([0, b_0], \eta)\! :=\! \max \{ |f(x)\! -\! f(y)|; x,y\ino [0, b_0] : |x\! -\! y| \leqo \eta \}$, $\eta\ino \bbR_+^*$, and $|\sigma|\! :=\! d(\rho, \sigma)$, $\sigma\ino T$. 
Since $0\leqo (g(|\sigma_p|)\! -\! g( |\sigma_p\wedge \sigma_q|)\vee (g(|\sigma_q|)\! -\! 
g(|\sigma_p\wedge \sigma_q|))  \leqo d_g (\sigma_p, \sigma_q)$, we get $0\leqo (|\sigma_p|\! -\! 
|\sigma_p\wedge \sigma_q|)\vee (|\sigma_q|\! -\! |\sigma_p\wedge \sigma_q|)  \leqo 
w_f([0, b_0], d_g (\sigma_p, \sigma_q))$ and thus $d(\sigma_p, \sigma_q) \leqo 
2w_f([0, b_0], d_g (\sigma_p, \sigma_q))$, which implies that $(\sigma_n)_{n\in \bbN}$ is a Cauchy sequence in 
$(T,d)$. Since $(T,d)$ is complete, there is $\sigma\ino  T$ such that $d(\sigma_n, \sigma)\! \to \! 0$, which implies $d(\rho, \sigma) \leqo f(b_0) \leko a$, i.e.~$\sigma\ino T^o$. Thus $d_g(\sigma_n , \sigma)\to 0$ and $\sigma\eqo \mathbf s\ino T^o$ since 
$(T^o\! , d)$ and $(T^o\! , d_g)$ are homeomorphic by Lemma \ref{fchangemetric}. This completes the proof of $T^o\!\! =\! \{\sigma \ino T_{\! g}: d_g(\rho, \sigma) \leko b\}$. Thus, $T^o$ is clearly open in $T_{\! g}$ and $d_g$-dense by definition of $T_{\! g}$. The equivalence of separability follows from Lemma \ref{fchangemetric}, which asserts that $(T^o\! , d)$ and $(T^o\! , d_g)$ are homeomorphic. 

Let us prove $(ii)$. Let $\mathbf s\ino T_{\! g}$ and $\sigma\ino \lgeo \rho, \mathbf s \lgeo \, $. Then $d_g(\rho, \sigma) \leko d_g (\rho, \mathbf s) \leqo b$ and $\sigma \ino T^o$. This also implies that $d_g(\rho, \mathbf s)\eqo b \leko \infty$ if 
$\mathbf s \ino T_{\! g} \backslash T^o$. If $\mathbf s\! \notin \! \mathtt{Lf} (T_{\! g})$, then $\mathbf s\! \in \, \rgeo \rho, \mathbf s' \lgeo \, $, where $\mathbf s'$ belongs to a connected component of $T_{\! g}\backslash \{ \mathbf s\}$ that does not contain $\rho$. Thus, $d_g (\rho, \mathbf s)\leko d_g (\rho, \mathbf s')\leqo b$ and $\mathbf s\ino T^o$. 
Let $\mathbf s, \mathbf s'\ino T_{\! g}$ be distinct. If $\mathbf s \wedge  \mathbf s'\ino \{ \mathbf s,\mathbf s'\}$, then either $\mathbf s\ino \lgeo \rho, \mathbf s' \lgeo \,\subseteq\! T^o$ or $\mathbf s'\ino \lgeo \rho, \mathbf s \lgeo \, \subseteq\! T^o$. If $\mathbf s, \mathbf s'\ino T_{\! g} \backslash T^o$, then  
$\mathbf s \! \wedge \! \mathbf s' \! \notin \!  \{ \mathbf s,\mathbf s'\}$ and we get $\mathbf s \wedge  \mathbf s' \ino  \lgeo \rho, \mathbf s \lgeo \,\subseteq\! T^o$. 
In particular, $\mathbf s$ and $\mathbf s^\prime$ cannot be in the same connected component of $T_{\! g}\backslash T^o$. 
This completes the proof of $(ii)$. 

To prove $(iii)$, we observe that the branchpoint of $\sigma_1, \sigma_2 \ino T^o$ is 
the same in $(T^o\! , d)$ and $(T^o\! , d_g)$. Thus, $(d_g)_f (\sigma_1, \sigma_2) 
\eqo f( d_g (\rho, \sigma_1)) +  f( d_g (\rho, \sigma_2))\! -\! 2 f(d_g(\rho, \sigma_1 \! \wedge \! \sigma_2))\eqo 
d(\sigma_1, \sigma_2)$, which shows that 
$(d_g)_f\eqo d$ on $T^o$ and we immediately get the desired result.

Let us prove $(iv)$. We assume that $(T,d)$ is locally compact and we prove that so is $(T_{\! g}, d_g)$, 
the converse then following from $(iii)$. First, $(T^o\! , d)$ is locally compact. If $b\eqo \infty$, then 
$T^o\eqo T_{\! g}$, and $(T^o\! , d)$ and $(T^o\! , d_g)$ being homeomorphic by Lemma \ref{fchangemetric}, 
$(T_{\! g}, d_g)$ is therefore locally compact. We now assume that $b\leko \infty$ and we prove that 
$(T_{\! g}, d_g)$ is compact. Since $(T,d)$ is complete, the
Hopf--Rinow theorem (see e.g.~Burago, Burago and Ivanov \cite[Thm 2.5.28]{BuBuIv}) shows that for all $r\ino [0, a)$, 
$B_{T,d} (\rho, r)$ is compact in $(T,d)$ and thus in $(T^o\! , d)$. Since $(T^o\! , d)$ and $(T^o\! , d_g)$ are homeomorphic 
by Lemma \ref{fchangemetric} and since $B_{T,d} (\rho, r)\eqo B_{T_g,d_g} (\rho, g(r))$, this also proves that for all 
$r\ino [0, b)$, 
$B_{T_g, d_g} (\rho, r)$ is compact in $(T^o\! , d_g)$, and thus in $(T_{\! g}, d_g)$. 
Let us prove that any $T_g$-valued sequence $(\mathbf s_n)_{n\in \bbN}$ has a $d_g$-limit point. 
If $\liminf_{n\to \infty} d_g(\rho, \mathbf s_n) \leko b$, the previous argument implies that $(\mathbf s_n)_{n\in \bbN}$ 
admits a $d_g$-converging subsequence. 
Now assume that $\lim_{n\rightarrow\infty}d_g(\rho,\mathbf s_n)=b$.
Let $(r_p)_{p\in \bbN}$ be a $[0, b]$-valued sequence that strictly increases to $b$ and such that $r_0\eqo 0$. 
Let $s(p,n)\ino \lgeo \rho, \mathbf s_n\rgeo$ be such that $d_g(\rho, s(p,n))\eqo r_p\wedge d_g(\rho,\mathbf s_n)$. 
By diagonal extraction, there are $\sigma_p\ino T_{\! g}$, $p\ino \bbN$, and a strictly increasing sequence 
of integers $(n_k)_{k\in \bbN}$ such that 
$\lim_{k\to \infty}d_g(s(p,n_k), \sigma_p)\eqo 0$ for all $p\ino\bbN$. Necessarily, 
$\sigma_p \ino \lgeo \rho, \sigma_{p+1} \rgeo$ and $b\eqo \sum_{p\in \bbN} d_g(\sigma_p, \sigma_{p+1})$. Thus $(\sigma_p)_{p\in \bbN}$ is a 
$d_g$-Cauchy sequence: its limit $\sigma\ino T_{\! g}$ must satisfy 
$d_g(\rho, \sigma)\eqo b$ and $\sigma_p\ino \lgeo \rho, \sigma \rgeo$ for all $p\in\bbN$. Then,  
$d_{g} (\mathbf s_{n_k}, \sigma ) \leqo d_g(\mathbf s_{n_k}, s(p,n_k))+ d_g(s(p,n_k), \sigma_p)+ 
d_g(\sigma_p, \sigma)
\leqo  2(b\! -\! r_p)+ d_g(s(p,n_k), \sigma_p)$. This implies $\limsup_{k\to \infty} d_{g} (\mathbf s_{n_k}, \sigma )\leqo 
2(b\! -\! r_p) \! \to \! 0$ as $p\! \to \! \infty$. Thus $\sigma$ is a limit point of $(\mathbf s_n)_{n\in \bbN}$. Hence, $(T_{\! g}, d_g)$ is compact. To complete the proof of $(iv)$, we note that $\mathtt{Ht}_{d,\rho} (T) \leko a$ 
implies that $T_{\! g}$ is bounded, and thus compact by the Hopf--Rinow theorem. \cqfd

\smallskip

Let us now define the boundary of an $\bbR$-tree. 
\begin{definition} 
\label{boundarydef}Let $(T, d, \rho)$ be a rooted complete $\bbR$-tree. 
 
 \smallskip
 
\begin{compactenum}
\item[$(a)$] A \emph{ray with endpoint $\rho$} is a subset $\mathbf r\! \subseteq\! T$ such that $\mathbf r\eqo c(\bbR_+)$ where 
$c\colon\bbR_+ \! \to \! T$ is a geodesic arc of infinite length necessarily, such that $c(0)\eqo \rho$. 

 \smallskip

\item[$(b)$] The \emph{boundary} of $(T,d , \rho)$ is the set  
$\partial T\eqo \big\{ \mathbf r \, ; \textrm{$\mathbf r$ is a ray whose endpoint is $\rho$} \big\}$.

\smallskip

\item[$(c)$] $(T, d, \rho)$ is said to be \emph{boundaryless} if $\partial T\eqo \emptyset$. 

\smallskip

\end{compactenum}
\noi
Let $\sigma \ino T$ and let $\mathbf r_1, \mathbf r_2\ino \partial T$. Let 
$c\colon [0, d(\rho, \sigma) ] \! \to \! T$, $c_i\colon \bbR_+ \!\! \to \! T$, $i\ino \{ 1, 2\}$, be geodesics such that $ c(0)\eqo \rho$, 
$c([0, d(\rho, \sigma)])\eqo \lgeo \rho, \sigma \rgeo$ and $c_i(\bbR_+)\eqo \mathbf r_i$, $i\ino \{ 1,2\}$. We next introduce the following. 
 
\medskip
 
\begin{compactenum}
\item[$(d)$] For all $i\ino \{ 1, 2\}$, we set $d(\rho, \mathbf r_i)\! :=\!  \infty$ and $\underline{\mathbf r}_i(s)\! :=\!  c_i(s)$, for all $s\ino \bbR_+$. 

\smallskip

\item[$(e)$] For all $s\ino \bbR_+$ we set $\underline{\sigma} (s) \! :=\! c(s\wedge d(\rho, \sigma))$. 

\smallskip

\item[$(f)$] We set $\mathbf r_1\! \wedge \mathbf r_2\! :=\!  c_1(s_0)$ where $s_0\eqo \max \{ s\ino \bbR_+ \! : \! c_1(s)\eqo c_2(s) \}$ if $\mathbf r_1\!\neq \! \mathbf r_2$, otherwise we set  $\mathbf r_1\! \wedge \mathbf r_2 \! :=\! \mathbf r_1$.

\smallskip

\item[$(g)$] We set  $\sigma \! \wedge \mathbf r_1\! :=\!  c(s_0)$ where $s_0\eqo \max \{ s\ino [0, d(\rho, \sigma)] : c(s)\eqo c_1(s) \}$. \cq 

\end{compactenum}
\end{definition}
\begin{theorem}
\label{bordif} Let $(T,d,\rho)$ be a rooted complete $\bbR$-tree. We set $T^*\eqo T\cup \partial T$. If $\mathbf s\ino \partial T$, we 
recall from Definition \ref{boundarydef} $(d)$ that $d(\rho, \mathbf s) \eqo \infty$ and we use the convention $e^{-\infty}\eqo 0$ in the following.
\begin{equation}
\label{conemetr}
\forall \bs_1, \bs_2\ino T^*\! , \quad d^* (\mathbf s_1, \mathbf s_2)\! :=\!  2e^{-d(\rho, \mathbf s_1 \wedge \mathbf s_2)}\!  -e^{-d(\rho, \mathbf s_1)} \! -e^{-d(\rho, \mathbf s_2)} .
\end{equation} 
Then, $(T^*\!, d^*\!, \rho)$ is isometric to $(T_{\! g}, d_g, \rho)$, where $g(x)\eqo 1\! -\! e^{-x}$, $x\ino [0, 1)$. 
We call $(T^*\!, d^*\!, \rho)$ the \emph{bordification} of $(T,d, \rho)$: it a boundaryless rooted complete $\bbR$-tree, $T$ is open and dense in $T^*\!$, $\partial T$ is closed, totally disconnected and $\partial T \! \subseteq \! \mathtt{Lf} (T^*)$. Moreover if 
 $(T,d)$ is separable (resp.~compact), then so is $(T^*\! ,d^*)$.
\end{theorem}
\noi
\textbf{Proof.} We define $\phi\colon T_{\! g} \! \to \! T^*$ as follows: we first take $\phi$ as the identity function on $T$;  
let $\mathbf s\ino T_{\! g}\backslash T$; by Theorem \ref{propdgtree} $(ii)$, 
$\lgeo \rho ,\mathbf s \lgeo \, \subseteq \!  T\eqo T^o $; thus 
$\mathbf r\! := \! \lgeo \rho , \mathbf s \lgeo\, $ has to be a ray since $\mathbf s\! \notin \! T$ and we set 
$\phi (\mathbf s)\eqo \mathbf r$. Note that $\phi$ is clearly injective. We next prove that $\phi (T_{\! g})\eqo T^*$: 
let  $\mathbf r$ be a ray in $(T,d)$; it has to be a geodesic path in $T_{\! g}$ and, since 
$\mathtt{Ht}_{d_g, \rho}(T_{\! g})\leqo 1$, its closure is of the form $\lgeo \rho, \mathbf s\rgeo$ for some $\mathbf s 
\ino T_{\! g}\backslash T$: so $\phi (\mathbf s)\eqo \mathbf r$. 
We have proved that $\phi$ is a bijection. For distinct $\mathbf s , 
\mathbf s'\ino T_{\! g}$ we set $\phi (\mathbf s)\eqo \mathbf r$ and $\phi (\mathbf s')\eqo \mathbf r'$. We then prove that 
$\mathbf s  \wedge \mathbf s'\eqo  \mathbf r  \wedge \mathbf r'$ (as in Definition \ref{boundarydef} $(f)$ and $(e)$). 
\emph{Indeed}, if $\mathbf s , \mathbf s' \ino T$, this is obvious; if 
$\mathbf s , \mathbf s' \ino  T_{\! g}\backslash T$, Theorem \ref{propdgtree} $(ii)$ entails 
$\lgeo \rho , \mathbf s  \wedge \mathbf s' \rgeo \eqo \mathbf r \cap \mathbf r'\eqo  
\lgeo \rho , \mathbf s \lgeo  \, \cap  \lgeo \rho , \mathbf s' \lgeo \, \subseteq\! T $ and thus $\mathbf s  \wedge \mathbf s' \eqo  \mathbf r  \wedge \mathbf r'$; if $\mathbf s\ino T$ and $\mathbf s' \ino T_{\! g} \backslash T$, $\mathbf s \eqo \mathbf r$, 
$\lgeo \rho , \mathbf s  \wedge \mathbf s' \rgeo \eqo \lgeo \rho , \mathbf s \rgeo  \cap \mathbf r'\!\subseteq\! T$ and 
$\mathbf s  \wedge \mathbf s' \eqo  \mathbf r  \wedge \mathbf r'$. 
We have proved for all $\mathbf s , \mathbf s'\ino T_{\! g}$, that $ \mathbf s  \wedge \mathbf s'\eqo  \phi (\mathbf s)  \wedge \phi (
\mathbf s')$, which easily implies that $d_{g} (\mathbf s, \mathbf s') \eqo d^*(\phi (\mathbf s), \phi (\mathbf s'))$. Therefore 
$(T^*\!, d^*\!, \rho)$ and $(T_{\! g}, d_g, \rho)$ are isometric and the remaining part of the proof is a direct consequence of Theorem \ref{propdgtree}. \cqfd

\begin{remark}
\label{toporem}

\noi
\textbf{(a)} To make the connection with the bordification of CAT($0$)-spaces and $\delta$-hyperbolic spaces as discussed for instance in Bridson and Haefliger \cite{briHae11}, Chapter II.8 and Chapter III.3, 
we observe that \eqref{dstarfirstdef} metrizes the topology of uniform convergence on all compact intervals, as a topology on the space $\mathbf{C}(\bbR_+,T)$ of continuous functions 
from $\bbR_+$ to $T$. The bordification is the closure of $T$, represented by the set $\{c_\sigma,\sigma\ino T\}\!\subseteq\!\mathbf{C}(\bbR_+,T)$ of possibly stopped geodesics of $T$.

\smallskip

\noi
\textbf{(b)} The definition of $(T^*\! , d^*\! , \rho)$ strongly relies on the specific choice of the function $g(x)\eqo 1\! -\! e^{-x}$, $x\ino [0, 1)$ (or on the specific distance metrizing the topology of uniform convergence of stopped geodesics). However, only the topology induced on $T^*$ by $d^*$ and its connections with the topology of $T$ actually matter.  
In separate notes (see \cite{DuWi24note}), 
we prove that the topology of the bordification $(T^*\! , d^*\! , \rho)$ introduced in Theorem \ref{bordif} 
can be characterized uniquely in terms of the topology of $T$ as the coarsest dense topological embedding (in the sense of Definition \ref{topoembedef}) of $T$ into a boundaryless topological real tree (topological real trees as defined in Bowditch \cite{Bow99}). We shall not use these results in the present paper.

\smallskip

\noi
\textbf{(c)} By Lemma \ref{fchangemetric} and Theorem \ref{bordif}, $(T,d)$ and $(T, d^*)$ are homeomorphic: the Borel sigma-field $\mathscr B(T)$ is also the Borel sigma-field generated by $d^*$ on $T$, so we can view $\mathcal M_f (T)$ as the subset of $\cM_f(T^*)$ of the finite Borel measures for which $\mu(\partial T)\eqo 0$. See Lemma \ref{embedmeas} below. \cq
\end{remark}

\begin{lemma}
\label{selfextend} Let $(T,d,\rho)$ be a rooted Polish $\bbR$-tree and let $(T^*,d^*, \rho)$ be its bordification as in Theorem \ref{bordif}. Let $\sigma_n\ino T$ with $\sigma_n \ino \lgeo \rho, \sigma_{n+1} \rgeo$, $n\ino \bbN$. Then there is $\mathbf s \ino T^*$ such that $d^*(\sigma_n , \mathbf s)\! \to \! 0$. 
\end{lemma}
\noi
\textbf{Proof.} Let us denote as in $T$ by $\lgeo \mathbf s, \mathbf s'\rgeo$ the $d^*$-geodesic arc joining $\mathbf s$ and $\mathbf s'$ in $T^*$. Then $\sigma_n \ino \lgeo \rho, \sigma_{n+1} \rgeo$ and $d^*(\rho, \sigma_{n+1})\eqo d^*(\rho, \sigma_n)+ d^*(\sigma_n, \sigma_{n+1})$. This first implies that $(d^*(\rho, \sigma_n))_{n\in \bbN}$ is a non-decreasing sequence which is  bounded by $1$, by definition of $d^*$. Therefore it converges to $\ell\ino [0, 1]$ and we then get $\ell \eqo \sum_{n\in \bbN} d^*(\sigma_{n-1}, \sigma_n)$, with the convention that $\sigma_{-1}\eqo \rho$. Thus $(\sigma_n)_{n\in \bbN}$ is a $d^*$-Cauchy sequence, which implies the desired result since $(T^*, d^*)$ is complete. \cqfd 

\smallskip

We shall use the term of \emph{topological embedding} in the following sense. 
\begin{definition}
\label{topoembedef} 
A \emph{topological embedding} $\jmath$ from a topological space $E_1$ into another topological space $E_2$ is an injective continuous function $\jmath\colon E_1 \! \to \! E_2$ such that $\jmath^{-1} \colon \jmath (E_1) \! \to \! E_1$ is continuous, when $\jmath(E_1)$ is equipped with the relative topology of $E_2$. \cq
\end{definition}
\begin{lemma}
\label{embedmeas} Let $(T,d,\rho)$ be a rooted Polish $\bbR$-tree and let $(T^*,d^*, \rho)$ be its bordification (see Theorem \ref{bordif}). We view $\cM_f(T)$ as a subset of $\cM_f(T^*)$ (see Remark \ref{toporem} $\mathbf (c)$). Then the identity function on $\cM_f(T)$ is a topological embedding from  
$\cM_f(T)$ equipped with the topology of weak convergence on $(T,d)$, into $\cM_f(T^*)$ equipped with the topology of weak convergence on $(T^*\! , d^*)$. 
\end{lemma}
\noi
\textbf{Proof.}  If $\mu_n\! \to \! \mu$ weakly in $\cM_f(T)$, then 
$\mu_n \! \to \! \mu$ weakly in $\cM_{f} (T^*)$ because the restriction to $T$ of a bounded $\bR$-valued 
continuous function on $T^*$ is also bounded and continuous on $T$. 
Conversely, suppose that $\mu_n , \mu \ino \cM_f(T^*)$ are such that $\mu_n (\partial T_n)\eqo \mu(\partial T)\eqo 0$ and $\mu_n \! \to \! \mu$ weakly in $\cM_f(T^*)$. Let $\varphi\colon T\! \to \! \bbR$ be bounded and continuous. We extend $\varphi$ to $T^*$ by setting $\varphi(\bs)\eqo 0$ for all $\bs \ino \partial T$. Then $\varphi$ is Borel-measurable and since the set of points where $\varphi$ is discontinuous is included in $\partial T$ which is $\mu$-negligible, the Portmanteau theorem implies that 
$\lim_{n\to \infty} \int_{T^*} \varphi \, d\mu_n\eqo  \int_{T^*} \varphi \, d\mu$, i.e.~$\lim_{n\to \infty} \int_{T} \varphi \, d\mu_n\eqo  \int_{T} \varphi \, d\mu$, which proves that $\mu_n \! \to \! \mu$ weakly in $\cM_f(T)$. 
\cqfd

\medskip

In the following lemma, we give properties of the bordification in connection with spanning trees and topological supports of measures. 
Let $(T,d, \rho)$ be a rooted Polish $\bbR$-tree with bordification $(T^*\! , d^*\! , \rho)$. 
For all $A\!\subseteq\! T^*$, we denote by $\, \overline{\! A}^*$ and $\sspaan^* (A)$  
the $d^*$-closures of $A$ and $\spaan (A)$ (and if $A\!\subseteq\! T$, we keep denoting by $\, \overline{\! A}$ the $d$-closure of $A$ in $T$).
Recall that arcs in $T^*$ are denoted as in $T$ by $\lgeo \bs, 
\bs'\rgeo$, $\mathbf s, \mathbf s'\ino T^*$.  
Let $\mu\ino \cM_f(T^*)$.  
We denote its topological support $\mu$ in $(T^*\! , d^*)$ by $\mathtt{\supp}_{\!\! *} \,  \mu$ (and if $\mu\ino \cM_f (T)$ we keep denoting by $\mathtt{\supp} \mu$ the topological support of $\mu$ in $T$). 
\begin{lemma}
\label{starific} Let $(T,d,\rho)$ be a rooted Polish $\bbR$-tree whose bordification is 
$(T^*\! ,d^*\! , \rho)$. Let $A\!\subseteq\! T$, and let $\mu\ino \cM_f(T)$, which is also viewed as an element of $\cM_f(T^*)$ such that $\mu (\partial T)\eqo 0$. Then the following properties hold true. 
\begin{compactenum}

\smallskip 

\item[$(i)$] $\sspaan (A)\eqo T \cap \sspaan^* (A)$ and $ \sspaan^* (A)\eqo  \sspaan^* (\, \overline{\! A}^*)$.

\smallskip 

\item[$(ii)$] $\mathtt{\supp} \mu \eqo T\cap \mathtt{\supp}_{\! \!*} \, \mu\; $ and $\; \overline{\mathtt{\supp} \mu}^{\,*} \eqo \mathtt{\supp}_{\!\! *}\,  \mu$.

\smallskip 

\item[$(iii)$] $\sspaan (\mathtt{\supp} \mu ) \eqo T\cap \sspaan^* (\mathtt{\supp}_{\!\! *} \, \mu)$. 
\end{compactenum}
\end{lemma}
\noi
\textbf{Proof.} This follows from elementary arguments, whose details are left to the reader.  \cqfd

\section{Mass erasure in fixed $\bbR$-trees}
\label{Merafixsec}
In this section we study the behaviour of $h$-erased measures on a fixed rooted $\bbR$-tree $(T,d, \rho)$, which is most of the time assumed to be Polish.
 To this end, we first need to introduce notation and recall general results on metric spaces. 

\smallskip

\noi
\textbf{Notation.} Let $(E,d)$ be a metric space. For all $x \! \in \! E$ and all $r\ino \bbR_+$, we use the notation $B^o(x, r)\eqo \big\{ y\! \in \! E: d(x,y)\leko r \big\}$ and $B(x, r)\eqo  \big\{ y\! \in \! E: d(x,y)\leqo r \big\}$ or $B_{E, d}^o(x,r)$ of $B_{E, d}(x,r)$ when there is an ambiguity .

In general, we shall say that a subset $A \! \subseteq \! E$ is \emph{precompact in $E$} if $A\! \neq \! \emptyset$ and its the closure $\overline{\! A}$ of $A$ is compact (precompact and \emph{relatively compact} are synonyms). Similarly, we say that a sequence $(x_n)_{n\in \bbN}$ is precompact in $E$ if $\{ x_n; n\ino \bbN\}$ is a precompact subset of $E$. We also use the following notation for all $\epp \ino \bbR^*_+$.
\begin{equation}
\label{entropie}
N(A, \epp)= \inf \Big\{ n\in \bN^*: \;  \exists\,  x_1, \ldots , x_n \ino A \; \textrm{such that} \;  A \! \subseteq \!  \!\! \bigcup_{1\leq i\leq n} \!\!\! B^o (x_i, \epp)  \Big\}\; ,
\end{equation}
with the convention that $\inf \emptyset \eqo \infty$.  Recall that if $(E,d)$ is complete, then $A$ is precompact iff $N(A, \epp) \leko \infty$ for all $\epp\ino \bbR^*_+$.

Let $ A\! \subseteq \! E$ be non-empty. We set $d(x, A)\eqo  \inf_{y\in A} d(x, y)$ for all $x\ino E$, and we recall that $d(\, \cdot\, , A)$ is $1$-Lipschitz, 
that $d(\, \cdot\, , A)\eqo d(\, \cdot\,  , \, \overline{\! A})$ and that $\, \overline{\! A}\eqo \{ x\ino E: d(x, A)\eqo 0\}$. 
For all $\eta \ino \bbR^*_+ $, we also set $A^{(\eta)}\eqo \{ x\ino E : d(x, A) \leqo \epp\}$ and for all non-empty $B\! \subseteq \! E$, we denote by $\dHaus(A,B)$ the \emph{Hausdorff distance} between $A$ and $B$:   
\begin{equation}
\label{Hausdis}
d_{\mathtt{Haus}} (A, B)= \inf \big\{ \eta \! \in \! \bbR^*_+ : \; B \! \subseteq \! A^{(\eta)} \;  
\textrm{and} \;  A \! \subseteq \! B^{(\eta)}\big\} \; , 
\end{equation}
with the convention that $\inf \emptyset = \infty$. We denote by $\mathcal C_E$ the space of closed subsets of $E$ and by $\mathcal K_E$ be the space of compact subsets of $E$ and we recall the following well-known results  (see e.g.~Burago et al.~\cite[Section 7.3.1]{BuBuIv}). 
\begin{compactenum}

\smallskip

\item[]$\!\!\!\!\!\!\!\!\!\textbf{Haus}$-$(i)$: \emph{$1\wedge \dHaus$ is a metric on $\mathcal C_E$ and $\dHaus$ is a metric on $\mathcal K_E$.} 

\smallskip

\item[]$\!\!\!\!\!\!\! \!\!\textbf{Haus}$-$(ii)$: \emph{If $(E, d)$ is complete, then so are $(\mathcal C_E,  1\! \wedge \! \dHaus)$ and $(\mathcal K_E ,  \dHaus)$. 
If $(E, d)$ is separable, then so is $(\mathcal K_E ,  \dHaus)$.   
If $(E, d)$ is compact, then so is $(\mathcal K_E ,  \dHaus)$. }

\smallskip

\item[]$\!\!\!\!\!\!\!\!\! \textbf{Haus}$-$(iii)$: \emph{Let us assume that $(E,d)$ is Polish. Let $\mathscr C\!\subseteq\!\mathcal K_E$ be non-empty. Then $\mathscr C$ is $\dHaus$-precompact iff for all $\varepsilon\ino \bbR_+^*$, $\sup_{K\in \mathscr C}N(K, \varepsilon)\leko \infty$.}

\smallskip

\item[]$\!\!\!\!\!\!\!\!\! \textbf{Haus}$-$(iv)$: \emph{For all $n\ino \bbN$, let $K_n \ino \mathcal K_E$ be such that $K_n\!\subseteq\! K_{n+1}$ We denote by $K$ the closure of $\, \bigcup_{n\in \bbN} K_n$. 
If $(K_n)_{n\in \bbN}$ is $\dHaus$-precompact, then $K$ is compact and $\dHaus (K_n, K)\! \to \! 0$. }

\smallskip

\item[]$\!\!\!\!\!\!\!\!\! \textbf{Haus}$-$(v)$: \emph{For all $n\ino \bbN$, let $K_n \ino \mathcal K_E$ be such that $K_{n+1}\!\subseteq\! K_n$. 
Then $K\! :=\!  \bigcap_{n\in \bbN} K_n$ is compact and $\dHaus (K_n, K)\! \to \! 0$.  }
\end{compactenum}

\medskip

  Next, we denote by $\cM_f (E)$ for the set of finite positive Borel measures on the Borel sigma-field $\cB(E)$. 
We recall here the definition of the \textit{Prokhorov distance} between $\mu$ and $\mu^\prime$ in $\cM_f (E)$: 
\begin{equation}
\label{Prokdef}
\dPro (\mu, \mu^\prime) \! = \! \inf \big\{ \eta \! \in \! \bbR_+\! : \, 
\forall C \! \in \! \mathcal C_E, \; \mu (C) \! \leq \! \mu^\prime \big( C^{(\eta)} \big)\!  + \eta \; \,  \textrm{and} \; \,  
 \mu^\prime (C)\!  \leq \! \mu \big(C^{(\eta)} \big) \! + \eta\big\} .
\end{equation}
Note that since for any subset $A$, $A^{(\eta)}= (\, \overline{\! A})^{(\eta)}$, $C$ in (\ref{Prokdef}) can range actually in $\cB(E)$. 
We next recall the following well-known results (see e.g.~Daley and Vere-Jones \cite[Prop.\ A.2.5.III]{DVJ1}). 

\smallskip

\noi
\emph{We assume that $(E,d)$ is Polish. Then, the following holds true.}

\smallskip

\begin{compactenum}

\item[]$\!\! \!\!\!\!\!\!\!\textbf{Pro}$-$(i)$: \emph{$(\cM_f (E), \dPro)$ is Polish and metrizes weak convergence. }

\smallskip

\item[]$\!\! \!\!\!\!\!\!\!\textbf{Pro}$-$(ii)$: \emph{Let $\mu, \mu^\prime\in \cM_f (E)$ be such that $\mu (E)\! = \! \mu^\prime (E)$. Then }
$$\dPro (\mu, \mu') \eqo  \inf \big\{ \eta \! \in \! \bbR_+: \; 
\forall C \! \in \! \mathcal C_E, \; \mu (C)  \leqo \mu^\prime (C^{(\eta)}) + \eta \big\}.$$ 

\smallskip

\item[]$\!\!\!\!\!\!\!\!\!\textbf{Pro}$-$(iii)$: Let $C\ino \mathcal C_E$ be non-empty. Denote by $d^{\,C}_{\mathtt{Pro}}$ the Prokhorov  distance on the space $\cM_f (C)$ of finite Borel measures on $(C, d_{| C\times C})$. 
Then $\dPro (\mu, \nu)\eqo d^{\, C}_{\mathtt{Pro}} (\mu, \nu)$ for all $\mu, \nu \ino \cM_f (C)$.

\smallskip

\end{compactenum}

\noi
We next recall Strassen's theorem. To this end we recall the following definition. 
\begin{definition}
\label{couplingdef} Let $(E,d)$ and $(E',d')$ be Polish. Let 
$\mu \ino \cM_f(E)$ and $\mu'\ino \cM_f(E')$. We equip $E\! \times \! E'$ with the product topology and 
the sigma-field $\cB(E)\! \otimes \! \cB(E')\eqo \cB(E\! \times \! E')$. 
Let $\eta \ino \bbR^*_+$. An \emph{$\eta$-coupling of $\mu$ and $\mu'$} is $\nu\ino \cM_f(E\! \times \! E')$ such that for all $B \ino \cB (E)$ and all $B'\ino \cB (E')$,  
$ 0 \leqo \mu (B) \! -\! \nu (B \! \times \! E')  \leqo \eta$ and $0 \leqo \mu^\prime (B') \! -\! \nu (E \! \times \! B')  \leqo \eta $. \cq 
\end{definition}

\noi
We refer to the following result as \emph{Strassen's theorem}. 

\smallskip

\begin{compactenum}

\item[]$\!\!\!\!\!\!\! \!\!\textbf{Pro}$-$(iv)$: \emph{Let $\mu, \mu^\prime \! \in \! \cM_f (E)$ and $\eta \ino \bbR^*_+$. Then, $ \dPro (\mu, \mu^\prime) \! \leq \!  \eta $ iff there exists an $\eta$-coupling $\nu$ of $\mu$ and $\mu'$ such that 
$ \nu \big( \big\{ (x,y)\! \in \! E^2 :  d(x,y) \geko \eta  \big\} \big)\eqo  0 $. }
\end{compactenum}

\medskip

We shall use also the following tightness result several times. 
\begin{lemma}
\label{Haustight} Let $(E,d)$ be Polish. For all $n \! \in \! \bN$, let $\mu_n \! \in \! \cM_f (E)$ be such that $\sup_{n \in \bN} \mu_n (E) \! <\!  \infty$ and let $K_n \! \subseteq \! E$ be compact subsets such that $\mu_n (E\backslash K_n)\!  \rightarrow \! 0$. If $(K_n)_{n\in \bN}$ is $\dHaus$-precompact, then the family of measures $(\mu_n)_{n\in \bN}$ is tight and therefore $\dPro$-precompact.    
\end{lemma}
\noi
\textbf{Proof.} W.l.o.g.~we can assume that $\dHaus (K_n, K)\! \rightarrow\!  0$, where $K\subo E$ is compact by $\textbf{Haus}$-$(ii)$. Let $\varepsilon \! \in \! \bbR^*_+$ and $(\eta_p)_{p\in \bN}$ be a $\bbR_+^*$-valued sequence decreasing to $0$. 
For all $p \! \in \! \bN$, let $n_p \! \in \! \bN$ be such that for all $n \! \geq \! n_p$, $K_n \subseteq K^{(\eta_p)} $ and $\mu_n (E\backslash K_n)\! \leq 2^{-p-1}\varepsilon$.  
Let $C_{p}$ be a compact set such that $\mu_n (E\backslash C_p)\! < \! 2^{-p-1} \varepsilon$ for all $n\! \in \! \{ 0, \ldots, n_p\}$. Then we set   
$Q_p \! = \! K \cup C_p$: this is a compact set such that $\mu_n ( E\backslash Q_{p}^{{(\eta_p)}}) \! < \!  2^{-p-1} \varepsilon$ for all $n\! \in \! \bN$.  
Thus, we easily see that $Q \! = \! \bigcap_{p\in \bN} Q_{p}^{{(\eta_p)}}$ is compact and such that 
$\mu_n(E\backslash Q)\! \leq \! \sum_{p\in \bN}  \mu_n ( E\backslash Q_{p}^{{(\eta_p)}}) \! < \!  \varepsilon$, which entails the desired result. \cqfd 

\medskip

We state here an elementary lemma on $\cM_f(E)$-valued r.v.s where $\cM_f(E)$ is endowed with the topology of weak convergence (which is Polish since $(E,d)$ is Polish). 
\begin{lemma}
\label{adhoclem} Let $E$ be Polish and let $\fLambda\! : \! \Omega \! \to \! \cM_f(E)$ be a random finite measure. Let $f_i \! : \! E\! \to \! \bR_+$, $i\ino \{ 1, 2\}$, be Borel-measurable. We denote by $\fLambda \circ f_i^{-1}$ the random finite measure on $\bbR_+$ that is the push-forward measure of $\fLambda$ via $f_i$. We assume the following. 
\begin{compactenum}

\smallskip

\item[$(a)$] $\bP$-a.s.~$\fLambda \big( \{ x\ino E: f_2(x) \leko f_1(x) \}\big)\eqo 0$.

\smallskip

\item[$(b)$] Under $\bP$ the two $\cM_f(\bbR_+) $-valued r.v.s $\fLambda \circ f_1^{-1}$ and $\fLambda \circ f_2^{-1}$ have the same law. 

\smallskip

\end{compactenum}
\noi
Then $\bP$-a.s.~$\fLambda \big( \{ x\ino E: f_2(x) \! \neq \!  f_1(x) \}\big)\eqo 0$. 
\end{lemma}
\noi
\textbf{Proof.} Let us first recall the following: if $X$ and $X'$ are $\bR_+$-valued r.v.s with the same law and such that $\bP$-a.s.~$X\leqo X'$, then $X\eqo X'$. (\emph{Indeed}, since $\bP (X' \leqo r \leko X)\eqo \bP (X'\leqo r)\! -\! \bP (X\leqo r) \eqo 0$, $\{ X'\leko X\}\eqo \bigcup_{r\in \bQ_+} \{ X'\leqo r \leko X\}$ is therefore $\bP$-negligible). 
This applies for all $q\ino \bQ_+$ to the r.v.s $X \eqo (\fLambda \circ f_1^{-1}) ([0, q)) $ and  
$X' \eqo (\fLambda \circ f_2^{-1}) ([0, q)) $ (\emph{indeed}, $(a)$ a.s.~entails $X'\leqo X$ and $(b)$ implies that $X$ and $X'$ have the same law). Thus there is $\Omega_0\ino \mathscr F$ such that $\bP (\Omega_0)\eqo 1$ and such that for all $\omega \ino \Omega_0$, $\fLambda_\omega \big( \{ f_2 \leko f_1 \}\big)\eqo 0$ and 
$\fLambda_\omega \big( \{ f_1 \leko q \}\big)\eqo\fLambda_\omega \big( \{ f_2 \leko q \}\big) $. 
Consequently $ \fLambda_\omega \big( \{  f_1 \leqo q \leko f_2 \}\big)\eqo \fLambda_\omega \big( \{ f_1 \leko q \}\big)\! -\! \fLambda_\omega \big( \{f_2 \leko q \}\big)\eqo 0$ and $\{ f_1 \leko f_2 \}\eqo \bigcup_{q\in \bQ_+} \{  f_1 \leqo q \leko f_2 \} $ is therefore $\fLambda_\omega$-negligible, which completes the proof of the lemma. \cqfd 

\subsection{Projections onto closed subtrees: definition and basic properties}
\label{massbase}
\begin{definition}
\label{defprojdef}
Let $(T, d, \rho)$ be a rooted \emph{separable} $\bbR$-tree. Let $\ttt$ a subtree of $T$ (i.e.~a connected subset) that contains the root $\rho$. For all $\sigma \ino T$, $\lgeo \rho, \sigma \rgeo \cap \ttt$ contains $\rho$ and  
is a compact connected subset of the arc $\lgeo \rho, \sigma \rgeo$. Therefore, there exists $\sigma^\prime\ino \ttt$ such that $\lgeo \rho, \sigma \rgeo \cap \ttt \eqo \lgeo \rho, \sigma^\prime \rgeo$.
We define the \emph{projection of $\sigma$ onto $\ttt$} as $f_\ttt (\sigma) \! :=\! \sigma'$. \cq 
\end{definition}
The following properties of projections are easy to check. 
\begin{proposition}
\label{projdef} Let $(T,d, \rho, \mu)$ be a rooted finitely measured Polish $\bbR$-tree. Let $\rho \ino \ttt\subo T $ 
be a closed subtree. Let $C_i$, $i\ino I$, be the connected components of the open set $T\backslash \ttt$ (here, the $C_i$ are non-empty open subsets and $I$ is countable since $T$ is separable). 
Then the following holds true. 
\begin{compactenum}

\smallskip

\item[$(i)$] For all $i\ino I$, there is $\sigma_i\ino \ttt$ such that $C_i$ is a connected component of $T\backslash \{ \sigma_i\}$, and $C_i\cup \{ \sigma_i\} $ is the closure of $C_i$.

\smallskip

\item[$(ii)$] For all $\sigma\ino T$, either $\sigma \ino \ttt$ and $f_\ttt (\sigma)\eqo \sigma$ or there exists $i\ino I$ such that $\sigma \ino C_i$ and $f_\ttt (\sigma) \eqo \sigma_i$. 

\smallskip

\item[$(iii)$] $f_\ttt$ is 1-Lipschitz with respect to $d$.

\smallskip

\item[$(iv)$] Let $\ttt_0$ and $\ttt_1$ be two closed subtrees of $T$ that contain $\rho$. Then, $\ttt_0\! \cap \! \ttt_1$ is a closed subtree of $T$ and  $f_{\ttt_0} \!\circ\! f_{\ttt_1}\eqo f_{\ttt_0\cap \ttt_1}$.

\smallskip

\item[$(v)$] We define the \emph{projection $\Pt \mu$ of $\mu$ onto $\ttt$} as the push-forward  
measure of $\mu$ via $f_\ttt$. Then, 
\vspace{-0.5mm}
\begin{equation}
\label{explicitproj}
\textstyle\Pt \mu = \mu (\cdot \cap \ttt) + \sum_{{i\in I}}\,  \mu (C_i)\delta_{\sigma_i} \; .
\end{equation}
Consequently we get $d_{\mathtt{var}} (\Pt\mu, \mu) := \sup_{B\in \mathscr B (T)} |\Pt \mu (B)\! -\! \mu(B)| \leqo \mu (T\backslash \ttt)$.

\smallskip

\item[$(vi)$] Let $\sigma \ino T$. If $ \sigma \ino \ttt$, then $\Pt\mu (\theta_\sigma T)\eqo \mu (\theta_\sigma T)$ and if $ \sigma \! \notin\!  \ttt$, then $\Pt\mu (\theta_\sigma T)\eqo 0$. 

\smallskip

\item[$(vii)$] Let $(T^*\! ,d^*\! , \rho)$ be the bordification of $(T,d,\rho)$. Then $\ttt$ is also a closed subtree of $T^*$. We denote by $f_{\ttt}^*\colon T^* \! \to \! \ttt$ and $\Pt^*\colon \cM_f(T^*)\! \to \! \cM_f(\ttt)$ the projections onto $\ttt$. Then $f_{\ttt}^*$ and $f_{\ttt}$ coincide on $T$ and $\Pt^*$ and $\Pt$ coincide on $\cM_f(T)$. 
\end{compactenum}
\end{proposition}
\noi
\textbf{Proof.}  We leave the details to the reader.  \cqfd 
\begin{remark}
\label{signedproj} Let us denote by $\mathcal M_\pm (T)\! :=\!  \cM_f(T)\! -\!  \cM_f(T) $, the space of \emph{signed measures}. Then 
(\ref{explicitproj}) allows to extend $\Pt$ as a \emph{linear} function from $\mathcal M_\pm (T)$ to $\mathcal M_\pm (T)$. 
Let $\mu\ino \mathcal M_\pm (T)$. We also note that $\Pt\mu$ is the push-forward  
measure of $\mu$ via $f_\ttt$. Namely, for all Borel-measurable $g\colon  T\! \to \! \bbR$, $\int_T g(f_{\mathtt t} (\sigma)) \, \mu (d\sigma) \eqo \int_{\mathtt t} \!  g(\gamma)\,  \Pt \mu (d \gamma) $. 
We see  that Proposition \ref{projdef} $(v)$ holds true and that if $\mathtt t'$ is another closed subtree containing $\rho$, 
$\mathtt{P}_{\!\ttt'} ( \mathtt{P}_{\!\ttt} \mu )\eqo \mathtt{P}_{\!\ttt\cap \ttt'} \mu$.    \cq 
\end{remark}
\begin{lemma}
\label{projprok} Let $(T, d, \rho, \mu)$ be a rooted finitely measured Polish $\bR$-tree. Let $\rho \ino \ttt_0\! \subseteq \!  \ttt \! \subseteq\! T $ be closed subtrees. Then, \vspace{-1.5mm}
\begin{equation}
\label{projprop}
\Ptz \mu \eqo\Ptz(\Pt \mu ) \quad \textrm{ and}  \quad \dPro (\Ptz \mu, \Pt\mu )  \leqo   \dPro (\Ptz \mu, \mu ) \; .
\end{equation}
\end{lemma}
\noi
\textbf{Proof.} The left-hand equality in (\ref{projprop}) is a consequence of Proposition \ref{projdef} $(iv)$. 
Let us prove the right-hand inequality. 
Let $\eta \ino \bbR^*_+$ be such that $ \dPro (\Ptz\mu,  \mu ) \leko \eta$. Then for all closed subsets $C$ of $T$, we get $\Ptz \mu (C)\leqo \eta +  \mu (C^{(\eta)})$. In particular, for all closed subsets $F$ of $\ttt$, 
we get  $\Ptz \mu (F)\leqo \eta +  \mu (F^{(\eta)})$. Let $(C_i)_{i\in I}$ be the connected components of $T\backslash \ttt$ and let $\sigma_i\ino \ttt$ be such that $C_{i} \cup \{ \sigma_i\}$ is the closure of $C_i$. Then $ \mu (F^{(\eta)})\eqo \mu(\ttt \cap F^{(\eta)})+ \sum_{i\in I} \mu (F^{(\eta)} \cap C_i)$. If $F^{(\eta)} \cap C_i \! \neq \! \emptyset$, then $\sigma_i \ino F^{(\eta)}$ since $F$ is subset of $\ttt$. Thus 
$ \mu (F^{(\eta)} \cap C_i) \leq \mu (C_i)\un_{F^{(\eta)}} (\sigma_i)$, which implies that $\Ptz \mu (F)\leqo \eta + \Pt \mu  (F^{(\eta)})$ by  Proposition \ref{projdef} $(v)$. 
Now observe that $ \Pt \mu  (F^{(\eta)})\eqo  \Pt \mu  (\ttt \cap F^{(\eta)})\eqo \Pt \mu (\{ y\ino \ttt\! : \! d(y, F) \leqo \eta \})$.
Consequently, we get $d^{\ttt}_{\mathtt{Pro}} (\Ptz \mu, \Pt \mu) \leqo \eta$, where $d^{\ttt}_{\mathtt{Pro}}$ is the Prokhorov distance on $\cM_f (\ttt)$ (here we use $\textbf{Pro}$-$(ii)$). This ends the proof of (\ref{projprop}) as $d^{\ttt}_{\mathtt{Pro}} (\Ptz \mu, \Pt \mu)\eqo d_{\mathtt{Pro}} (\Ptz \mu, \Pt \mu)$ by $\textbf{Pro}$-$(iii)$. \cqfd

\smallskip

We now discuss the convergence of measures which are consistent under projections. 
\begin{lemma}
\label{easyproj}
Let $(T, d, \rho, \mu)$ be a rooted finitely measured Polish $\bR$-tree.
For all $n\ino \bN$, let $\rho \ino \ttt_n\!\subseteq\! \ttt_{n+1}\!\subseteq\! T$ be closed subtrees. We set $\mu_n \eqo \Ptn\mu$. Then the following properties hold true. 
\begin{compactenum}

\smallskip

\item[$(i)$] For all integers $n\geqo m$, $\mathtt{P}_{\!\ttt_m}\mu_n\eqo \mu_m$. 

\smallskip

\item[$(ii)$] $\sup_{n\in \bN} \mu_n (T\backslash B(\rho, r)) \leqo \mu (T\backslash B(\rho, r)) \! \longrightarrow\!  0$ as $r\! \rightarrow \! \infty$. 

\smallskip

\item[$(iii)$] Let us denote by $\ttt $ the closure of $\, \bigcup_{n\in \bN} \ttt_n$. Then $\mu_n \! \rightarrow \! \Pt \mu$ weakly.
\end{compactenum}
\end{lemma}
\noi
\textbf{Proof.} $(i)$ follows from (\ref{projprop}). 
To prove ($ii$), we note for all $\sigma \ino T$ that $d(\rho , f_{\ttt_n} (\sigma)) \leqo d(\rho, \sigma)$, 
which yields $\mu_n (T\backslash B(\rho, r))\eqo \mu ( \{ \sigma \ino T \! :\!  d(\rho, f_{\ttt_n} (\sigma)) \geko r \}) \leqo 
\mu (T\backslash B(\rho, r))$, and thus ($ii$). To prove ($iii$), we fix $\gamma \ino \ttt$ and we first see that 
$d(\gamma, f_{\ttt_n} (\gamma))\eqo d(\gamma , \ttt_n)$ which tends to $0$. By Proposition \ref{projdef} ($iv$), we get $f_{\ttt_n} (\sigma)\eqo f_{\ttt_n} (f_{\ttt} (\sigma)) \! \rightarrow \! f_{\ttt}(\sigma)$, for all $\sigma \ino T$, which implies 
for any bounded and continuous function $\varphi \colon T\! \to \! \bR$ that  
$\int _T \varphi  \,  d\mu_n \! \eqo\!  \int_T  \varphi (f_{\ttt_n} (\sigma)) \mu(d\sigma)\! \to \! \!  \int_T  \varphi (f_{\ttt} (\sigma))\mu(d\sigma) \eqo \! \int _T \varphi \,  d \mathtt{P}_{\!\ttt}\mu $ by dominated convergence, and thus ($iii$). \cqfd 
\begin{proposition}
\label{projproj} Let $(T,d , \rho)$ be a rooted Polish $\bR$-tree with bordification $(T^*\! , d^*\! , \rho)$. For all $n\! \in \! \bN$, let $\rho\! \in \! \ttt_n \! \subseteq \! \ttt_{n+1} \!\subseteq\!T$ be closed subtrees and let $\mu_n \! \in \! \cM_{\! f} (T)$ be such that $\Ptn \mu_{n+1}\eqo  \mu_n$. Recall the definition of $\Ptn^*$ from Proposition \ref{projdef} $(vii)$. Then, the following holds. 
\begin{compactenum}

\smallskip

\item[$(i)$] There exists a unique $\mu\ino \cM_{\! f}(T^*)$ whose topological support is contained in the $d^*\! $-closure 
$\ttt^*$ of $\, \bigcup_{n\in \bbN} \ttt_n$, and such that $\Ptn^*\mu\eqo \mu_n$ for all $n\ino \bbN$ and $\mu_n \! \to \! \mu$ weakly on $\cM_{\! f} (T^*)$. 

\smallskip

\item[$(ii)$] Moreover, $\mu(\partial T)\eqo 0$ iff $ \limsup_{n\to \infty} \mu_n (\ttt_n \backslash B(\rho, r)) \! \to \! 0$ as $r\! \to \! \infty$. In this case, we view $\mu$ as an element of $\cM_f(T)$. Then its topological support in $T$ is included in $T\cap \ttt^*$ which is the closure of  $\, \bigcup_{n\in \bbN} \ttt_n$ in $T$. Moreover, $\Ptn \mu\eqo \mu_n$ for all $n\ino\bbN$ and 
 $\mu_n\! \to \! \mu$ weakly in $\cM_{\! f} (T)$. 
\end{compactenum}
\end{proposition} 
\noi
\textbf{Proof.} Observe that $c\! :=\! \mu_n(T)\eqo \mu_n (\ttt_n)$ for all $n\ino \bbN$. If $c\eqo 0$, the proposition trivially holds true. We assume that $c\geko 0$ and without loss of generality, we can also suppose that $c\eqo 1$. By Kolmogorov's extension theorem, there are r.v.s $X_n\colon \Omega \! \rightarrow \! \ttt_n$, $n\ino\bN$, 
such that  $X_n$ has law $\mu_n$ and a.s.~for all $n\! \in \! \bN$, $f_{\ttt_n} (X_{n+1})\!  = \! X_n$, which implies $X_n \ino \lgeo \rho, X_{n+1} \rgeo$. By Lemma \ref{selfextend}, there is a r.v.~$X\colon   \Omega \! \rightarrow \! T^*$ 
such that a.s~$d^*(X_n,X)\! \to \! 0$ a.s. Since a.s.~$f^*_{\ttt_n} (X_{m+n})\eqo f_{\ttt_n} (X_{m+n})\eqo X_{n}$ for all $m,n\ino \bbN$ (by Proposition \ref{projdef} ($vii$) and ($iv$)) and since $ f^*_{\ttt_n}$ is $d^*$-continuous (by Proposition \ref{projdef} $(iii)$), we get a.s.~$f^*_{\ttt_n} (X)\eqo X_n$ and thus $\Ptn^* \mu \eqo \mu_n$ for all $n\ino \bbN$, where $\mu$ is defined as the law of $X$. Since the $X_n$ a.s.~belong to $\ttt_n$, $X$ a.s.~belongs to $\ttt^*$, which implies that the topological support of $\mu$ is contained in $\ttt^*$. Thus $\mathtt{P}^*_{\ttt^*}\mu\eqo \mu$ and $\mu_n \! \to \! \mu$ weakly in $\cM_{f} (T^*)$ by Lemma \ref{easyproj} applied in $(T^*,d^*,\rho^*)$. This also entails the desired uniqueness. 

Let us prove $(ii)$. We note that $\mu_n (\ttt_n \backslash B_{T,d}(\rho, r))\leqo \mu_n( T^*\backslash B^o_{T,d}(\rho, r))$. 
Since  $T^*\backslash B^o_{T,d}(\rho, r)$ is $d^*\!$-closed. By the Portmanteau theorem $\limsup_{n\to \infty}
 \mu_n (\ttt_n \backslash B^o_{T,d}(\rho, r))\leqo \mu (T^*\backslash B^o_{T,d}(\rho, r))$, which implies the first implication in 
 $(ii)$ since 
$\lim_{r\to \infty}\mu (T^*\backslash B^o_{T,d}(\rho, r))\eqo \mu(\partial T)$. 
Conversely, since $T^*\backslash B_{T,d}(\rho, r)$ is $d^*\!$-open, the Portmanteau theorem entails 
$\mu(T^*\backslash B_{T, d}(\rho, r)) $ $\leqo$ $ \liminf_{n\to \infty}$  $\mu_n (T^* \backslash B_{T,d}(\rho, r))$ $\eqo$ 
$ \liminf_{n\to \infty}$  $\mu_n (\mathtt t_n \backslash B_{T,d}(\rho, r))$, since 
$\mathtt{Supp} \mu_n \subo \mathtt{t}_n$. Thus, 
$\mu(T^*\backslash B_{T,d}(\rho, r)) $ $\leqo$  $\limsup_{n\to \infty} $ $\mu_n (\ttt_n \backslash B_{T, d}(\rho, r))$, which easily entails the second implication in $(ii)$. 

We then suppose that $\mu(\partial T)\eqo 0$ and we view $\mu$ as an element of $\cM_f(T)$ whose topological support in $T$ is clearly included in $\ttt\! :=\! T\cap \ttt^*$, which is the $d$-closure of $\, \bigcup_{n\in \bbN} \ttt_n$. 
Then Proposition \ref{projdef} ($vii$) implies $\Ptn \mu\eqo \Ptn^*\mu\eqo \mu_n$. Lemma \ref{easyproj} implies that 
$\mu_n \! \to \!\Pt\mu\eqo \mu$ weakly in $\cM_f(T)$. \cqfd 
\begin{example}
\label{conterex} Let $(T, d, \rho)$ be the $\bR$-tree obtained from the infinite rooted binary tree by joining any two neighbouring vertices by a unit length geodesic as in Example \ref{evanescent}. 
We set  $\ttt_n\eqo  B_{T,d}(\rho, n)$ and $\mu_n \eqo \sum 2^{-n} \delta_\sigma$, where the sum extends over the $\sigma \ino \mathtt{Br} (T)$ such that $d(\rho, \sigma)\eqo n$. 
Thus $ \mu_n (T)\eqo1$ and we see that $\Ptn \mu_n \eqo \mu_{n+1}$. By Proposition \ref{projproj}, $(\mu_n)_{n\in \bbN}$ converges weakly on the bordification $T^*$ of $T$ to a measure $\mu\ino \cM_f(T^*)$. Here, $\mu (T)\eqo 0$ and $\mu (\partial T)\eqo \mu (T^*) \eqo 1$.  \cq 
\end{example}

The following lemma provides a control of projections in terms of the Hausdorff distance. It is also used to study 
Gromov--Prokhorov convergence of measure projections. 
\begin{lemma}
\label{GrC0proj} Let $(E, d)$ be Polish. Here, $\dHaus$ stands for the Hausdorff distance on $(E,d)$ and $\dPro$ for the Prokhorov distance on $\cM_f(E)$. 
Let $T$ and $T^\prime$ be closed subsets of $E$. Let $\rho \ino T$ and $\rho^\prime \ino T^\prime$. 
We assume that $(T, d, \rho)$ and that $(T', d, \rho^\prime)$ are rooted $\bR$-trees. Let $\tau$ (resp.~$\tau^\prime$) be a closed connected subset of $T$ (resp.~of $T^\prime$) which contains $\rho$ (resp.~$\rho^\prime$). Let $f_\tau$ (resp.~$f_{\tau^\prime}$) be the projection of $T$ onto $\tau$ (resp.~of $T^\prime$ onto $\tau^\prime$). 
\begin{compactenum}

\smallskip

\item[$(i)$] We assume that $\dHaus (T, T^\prime) \leko  \varepsilon $ and $\dHaus (\tau , \tau^\prime) \leko \eta $. Then \vspace{-1mm}
\begin{equation}
\label{Groproj}
\forall (\sigma , \sigma^\prime) \in T\! \times \! T^\prime , \qquad \big( \, d(\sigma, \sigma^\prime) \leq \varepsilon \;  \Longrightarrow \;  \, d( f_\tau (\sigma) , f_{\tau^\prime} (\sigma^\prime)) \leq 2\varepsilon + 3\eta \,\big). 
\end{equation}
\item[$(ii)$] Let $\mu \ino \cM_f (T)$ and $\mu^\prime \ino \cM_f (T^\prime)$, i.e.~$\mu, \mu'\ino \cM_f(E)$ and $\mu(E\backslash T)\eqo \mu'(E\backslash T')\eqo 0$. 
In the setting of $(i)$, we assume furthermore that $\dPro (\mu, \mu^\prime) \leko \varepsilon$. Then, 
\begin{equation}
\label{Groprojprok}
\dPro (\mathtt{P}_{\!\tau} \mu , \mathtt{P}_{\!\tau^\prime} \mu^\prime) \leq 2\varepsilon + 3 \eta \; .
\end{equation}
\end{compactenum}
\end{lemma}
\noi
\textbf{Proof.} Let $(\sigma, \sigma')\ino T\! \times \! T' $ be such that $d(\sigma, \sigma^\prime) \! \leq \! \varepsilon$. Let $\gamma^\prime \! \in \! \tau^\prime$ be such that $d(f_\tau (\sigma) , \gamma^\prime)\!  \leq \! \eta$. Then 
\begin{equation}
\label{Mtcor}d( \sigma, f_\tau (\sigma) )\geqo  d(\sigma^\prime , \gamma^\prime)\!  -\!  \varepsilon \! -\! \eta\geqo  \min_{x^\prime \in \tau^\prime} d(\sigma^\prime , x^\prime) \!  -\!  \varepsilon \! -\! \eta\eqo d(\sigma^\prime , f_{\tau^\prime} (\sigma^\prime)) \! -\!  \varepsilon \! -\! \eta.
\end{equation}
Similarly, $d(\sigma^\prime , f_{\tau^\prime} (\sigma^\prime)) \geqo d( \sigma, f_\tau (\sigma) ) \! -\!  \varepsilon \!-\! \eta$. 
Thus, $|d(\sigma^\prime , f_{\tau^\prime} (\sigma^\prime)) \! -\! d( \sigma, f_\tau (\sigma) )  | \leqo  \varepsilon +\eta$. 
Next we observe that 
$d(\gamma^\prime , \sigma^\prime)= d(\gamma^\prime , f_{\tau^\prime} (\sigma^\prime) ) +  d(\sigma^\prime , f_{\tau^\prime} (\sigma^\prime) )$, and by the first inequality in \eqref{Mtcor} we get 
\begin{eqnarray*}
d(f_{\tau^\prime} (\sigma^\prime), f_\tau (\sigma)) \!\! & \leq &\!\!  \eta + d(\gamma^\prime,   f_{\tau'} (\sigma')) =  
 \eta  + d(\gamma^\prime , \sigma^\prime)-d(\sigma^\prime , f_{\tau^\prime} (\sigma^\prime) )\\ 
 \!\!  &\leq & \!\! \eta+ d(\gamma', f_\tau (\sigma))+ d(f_\tau (\sigma), \sigma') -d(\sigma^\prime , f_{\tau^\prime} (\sigma^\prime) )\\ 
 &\leqo & \!\!2\eta + \varepsilon + d(f_\tau (\sigma)  , \sigma)-d(\sigma^\prime , f_{\tau^\prime} (\sigma^\prime) ) \quad\leq\quad 3 \eta + 2 \varepsilon ,
\end{eqnarray*}
which proves (\ref{Groproj}). Let us prove (\ref{Groprojprok}). Let $F$ be a closed subset of $E$. Observe that 
\begin{eqnarray*}
 \mathtt{P}_{\!\tau} \mu (F) &= &  \mathtt{P}_{\!\tau} \mu (F\cap \tau)= \mu (f^{-1}_\tau (F \cap \tau)) \\
 &\leq & \varepsilon  + \mu^\prime \big((f^{-1}_\tau (F \cap \tau))^{(\varepsilon)} \big)=\varepsilon  + \mu^\prime \big(T^\prime \cap (f^{-1}_\tau (F \cap \tau))^{(\varepsilon)} \big)  .
 \end{eqnarray*}
Let $\sigma^\prime \ino  T^\prime \cap (f^{-1}_\tau (F \cap \tau))^{(\varepsilon)}$. Then, there is $\sigma \ino T$ such that $f_\tau (\sigma ) \ino F $ and $d(\sigma, \sigma^\prime) \leqo \varepsilon$. By (\ref{Groproj}), we get $d( f_\tau (\sigma) , f_{\tau^\prime} (\sigma^\prime)) \leq 2\varepsilon + 3\eta=: \delta$. Thus, $\sigma^\prime \in f_{\tau^\prime}^{-1} (F^{(\delta)})$. This implies that 
$$  \mu^\prime \big(T^\prime \cap (f^{-1}_\tau (F \cap \tau))^{(\varepsilon)} \big) \leq \mu^\prime (f_{\tau^\prime}^{-1} (F^{(\delta)}) ) = \mathtt{P}_{\!\tau^\prime}\mu^\prime (F^{(\delta)}) \; .$$
Since $\delta \eqo 3\eta + 2\varepsilon \geq\varepsilon$, for all closed subsets $F$ of $E$, we have proved that 
$\mathtt{P}_{\!\tau} \mu (F)\leqo \delta +  \mathtt{P}_{\!\tau^\prime}\mu^\prime (F^{(\delta)})$. A similar argument entails that 
$ \mathtt{P}_{\!\tau^\prime}\mu^\prime (F)\leqo \delta +   \mathtt{P}_{\!\tau} \mu (F^{(\delta)})$, which completes the proof of (\ref{Groprojprok}). \cqfd

\subsection{Mass-erased subtrees} 
\label{meratreesec}
In this section, we study basic properties of mass-erased subtrees as introduced in Definition \ref{deferasedtree}. 
\begin{lemma}
\label{erasedsubtree} Let $(T, d, \rho, \mu)$ be a rooted finitely measured Polish $\bbR$-tree. 
Then the $h$-mass-erased tree $T_{\! \mu, h}$
is a compact subtree of $T$, contains the root and is of finite type as in Definition \ref{spandex} ($b$). 
\end{lemma}
\noi
{\bf Proof.} To prove that $T_{\! \mu, h}$ satisfies the three conditions of Lemma \ref{edgelength}, we first note that if $\sigma \ino T_{\! \mu , h}$, then $\lgeo \rho , \sigma \rgeo \!\subseteq\! T_{\! \mu, h}$,  and $T_{\! \mu, h}$ is a path-connected subset of $T$ containing $\rho$. To prove that $T_{\bm, h} $ is closed, we fix $(\sigma_n)_{n\in \bbN}$, a $T_{\! \mu, h}$-valued sequence converging to $\sigma$. As already noted, $\sigma\! \wedge\! \sigma_n \ino T_{\! \mu, h}$, which implies  
$d(\sigma\! \wedge\! \sigma_n , \sigma)\leqo d(\sigma, \sigma_n) \! \to \! 0$. Thus, $\sigma\! \wedge\! \sigma_n \! \to \!  \sigma$ and $d(\rho, \sigma_n\!\wedge\! \sigma) \! \to \! d(\rho, \sigma)$. W.l.o.g. we thus assume 
$\lgeo \rho, \sigma_n\rgeo\!\subseteq\!  \lgeo \rho, \sigma_{n+1}\rgeo\! \subseteq\!\lgeo \rho, \sigma\rgeo$ and $\sigma_n  \! \to \! \sigma $.
This implies $\theta_{\sigma_{n+1}}T\!\subseteq\! \theta_{\sigma_{n}} T$, $\theta_\sigma T\eqo \bigcap_{n\in \bbN} \theta_{\sigma_{n}} T$ and $\mu(\theta_\sigma T)\eqo \lim_{n\to \infty}\mu(\theta_{\sigma_{n}} T)\geqo h$ since $\sigma_n \ino T_{\! \mu, h}$. Thus, $\sigma \ino T_{\! \mu, h}$ and $T_{\! \mu, h}$ is closed. 

We next prove that $\mathtt{Ht}_{d,\rho}(T_{\! \mu, h}) \leko \infty$. \emph{Indeed}, let $r\ino \bbR^*_+$ be such that $\mu (T\backslash B(\rho, r))\leko h$. If $d(\rho, \sigma) \geko r$, we get $\theta_\sigma T\! \subseteq\! T\backslash B(\rho, r)$ and 
$\mu (\theta_\sigma T ) \leko h$, which shows proves that $\mathtt{Ht}_{d,\rho}(T_{\! \mu, h}) \leqo r$. 

Let $\sigma_1, \ldots, \sigma_n \ino T_{\bm, h}$ be such that $\sigma_i \wedge \sigma_j \notin \{ \sigma_i , \sigma_j \}$, for all $1\leqo i \leko j \leqo n$.  Lemma \ref{AboTrprop} $(i)$ entails that the $\theta_{\sigma_i}T$ are pairwise disjoint and thus 
$\bm (T) $  $\geqo$  $\sum_{1\leq i\leq n}$  $\bm (\theta_{\sigma_i} T ) $ $\geqo$ $nh$. Namely, $n\leqo \mu(T)/h$, which entails (\ref{bbound}). Lemma \ref{edgelength} applies and entails the desired result. \cqfd

\smallskip

The following proposition provides the main properties of mass-erased trees. 
\begin{proposition}
\label{massTprop}
Let $(T, d, \rho, \mu)$ be a rooted finitely measured Polish $\bbR$-tree. 
For all $h\ino \bbR^*_+$ we denote by $T_{\! \mu, h+}$ the closure of $T^o_{\! \mu, h+}\! := \!   \bigcup_{h^\prime >h} T_{\! \mu , h^\prime} $.
Then, the following holds true.
\begin{compactenum}

\smallskip

\item[$(i)$] Let $h\ino \bbR^*_+$, and let $\rho\ino \ttt\subo T$ be a closed subtree such that $\mu (T\backslash \ttt)\leko h$. Then, $T_{\mathtt{P}_\mathtt{t} \mu, h}\eqo T_{\! \mu, h}$. 

\smallskip

\item[$(ii)$] Let $\nu\ino \cM_f (T)$. Then $\mu \eqo \nu$ iff $T_{\! \mu, h}\eqo T_{\! \nu,h}$ for all $h\ino \bbR^*_+$.

\smallskip

\item[$(iii)$] For all $h, h' \ino \bbR^*_+$ such that $h\leko h'$, $T_{\! \mu, h'} \subo T_{\! \mu, h}$ and 
$\bigcap_{ h'' \in (0, h)} T_{\! \mu, h''} \eqo T_{\! \mu, h}$, which implies that 
$h\ino \bbR^*_+\! \mapsto \! T_{\!\mu , h}$ is $\dHaus$-left-continuous with $\dHaus$-right limits, it has at most countably many discontinuities and its $\dHaus$-right limit at $h$ is $T_{\! \mu , h+}$. 
\item[$(iv)$] Let $h\ino \bbR^*_+$. Then 
\begin{equation}
\label{discTmuh}
\bigsqcup_{\sigma \in \mathtt{Lf} (T_{\! \mu, h})}  \!\!\! \rgeo f_{T_{\! \mu, h+}}\!  (\sigma) , \sigma  \rgeo \; = T_{\! \mu, h} \backslash T_{\! \mu, h+} \subseteq 
 T_{\! \mu, h} \backslash T^o_{\! \mu, h+}\eqo \big\{ \sigma \ino T: \mu (\theta_\sigma T) \eqo h \big\}, 
 \end{equation}
where $f_{T_{\! \mu, h+}} $ stands for the projection onto $T_{\! \mu, h+}$.

\smallskip

\item[$(v)$]  If $\supp  \mu$ is connected and contains $\rho$, then $T_{\! \mu, h} \eqo T_{\! \mu, h+}$  
for all $h\ino \bbR^*_+$.

\smallskip

\item[$(vi)$] The closure of $T^o\! :=\!  \bigcup_{h\in \bbR^*_+} T_{\! \mu , h}$ is $\sspaan \, (\supp \mu)$, the closure of $\mathtt{Span} \, (\supp  \mu )$.

\end{compactenum}
\end{proposition}
\noi
{\bf Proof.} Let us prove $(i)$. Let $\sigma \ino T \backslash \ttt$. Since $\rho \! \in \! \ttt$ and since $\ttt$ is 
closed and connected, $\theta_\sigma T \! \subseteq \! T\backslash \ttt $. Thus $\mu (\theta_\sigma T) \! <\!  h$. 
Consequently, $T_{\! \mu , h} \subo \ttt$, and Proposition \ref{projdef} $(vi)$ implies the desired result. 

To prove the non-trivial implication of $(ii)$, we set $T^o \! :=\!  \bigcup_{h\in \bbR^*_+} T_{\! \mu , h}$ and  for all $\sigma\ino T^o\!$ we 
note  that $\mu (\theta_\sigma T) \! = \! \sup \{ h \! \in \! \bbR^*_+ \! : \sigma \! \in \! T_{\! \mu , h} \}$. Therefore
$\mu (\theta_{\sigma} T) \! = \! \nu (\theta_{\sigma} T)$. If $\sigma \ino T\backslash T^o$, we clearly get $\mu (\theta_{\sigma} T) \! =\!  \nu (\theta_{\sigma} T)\! =\!  0$. We complete the proof by applying 
Lemma \ref{AboTrprop} $(vi)$ and the uniqueness of the extention of finite measures on pi-systems. 

We next prove $(iii)$. The fact that $h\! \mapsto\! T_{\! \mu, h}$ is non-increasing with respect to inclusion and that 
$\bigcap_{h^\prime< h} T_{\! \mu, h^\prime} \eqo  T_{\! \mu , h}$ are immediate consequences of Definition \ref{deferasedtree}. 
By Lemma \ref{erasedsubtree}, $T_{\! \mu , h}$ is compact. 
We conclude by applying well-known results $\textbf{Haus}$-$(iv)$ and $\textbf{Haus}$-$(v)$ on the Hausdorff distance; see  
the beginning of Section \ref{Merafixsec}. 
Now a general result says that any function from $\bbR^*_+$ to a Polish space, which is left-continuous with 
right limits has at most countably many discontinuities.

Let us prove $(iv)$. On one hand, we see that $T_{\! \mu , h} \backslash T^o_{{\! \mu , h+}}\eqo \{ \sigma \ino T: 
\mu (\theta_\sigma T) \eqo h \}$. On the other hand, we note that $\, \rgeo f_{T_{\! \mu, h+}}\!  (\sigma) , \sigma  \rgeo 
\!\subseteq\! T_{\! \mu , h} \backslash T_{{\mu , h+}}$ 
for all $\sigma \ino \mathtt{Lf} (T_{\! \mu,h})$. 
We next fix $\gamma \ino T_{\! \mu , h} \backslash T_{{\! \mu , h+}}$. Then $\mu(\theta_\gamma T)\eqo h$ and 
there is $\sigma\ino \mathtt{Lf} (T_{\! \mu, h})$ such that $\sigma\ino \theta_\gamma T$. Thus, for all 
$x \ino \lgeo \gamma, \sigma\rgeo$, we get 
$\theta_{\sigma}  T \!\subseteq\! \theta_{x}  T \!\subseteq\! \theta_{\gamma}  T $ and thus 
$h\eqo \mu (\theta_{\sigma}  T) \eqo \mu ( \theta_{x}  T) \eqo \mu (\theta_{\gamma}  T )$, which first implies that  
$\lgeo \gamma , \sigma \rgeo \!\subseteq\! T_{\! \mu , h} \backslash T^o_{{\! \mu , h+}}$ and necessarily 
$\lgeo \gamma , \sigma \rgeo \!\subseteq\! T_{\! \mu , h} \backslash T_{\! {\mu , h+}}$. Therefore, we get 
$ f_{T_{\! \mu , h+}} (\sigma) \eqo  f_{T_{\! \mu , h+}} (\gamma)$. This proves that $T_{\! \mu , h} \backslash T_{{\! \mu , h+}}$ is the union of the $\, \rgeo f_{T_{\! \mu, h+}}\!  (\sigma) , \sigma  \rgeo $, $\sigma\ino \mathtt{Lf} (T_{\! \mu,h})$, which are necessarily pairwise disjoint (since at a branch point $\gamma$ of $T_{\! \mu, h}$, we get $\mu(\theta_\gamma T)\geqo 2h$). This completes the proof of $(iv)$.

Let us prove $(v)$. We suppose that $h'\ino \bbR^*_+\! \mapsto\! T_{\! \mu , h'}$ is discontinuous at $h$. By $(iv)$, there is 
$\sigma\ino \mathtt{Lf} (T_{\! \mu, h})$ such that $f_{T_{\! \mu , h+}} (\sigma)\! \neq \! \sigma$. Then for all $\gamma \! \in \, \rgeo f_{T_{\! \mu , h+}} (\sigma), \sigma \lgeo \, $, $\theta_\sigma T\!\subseteq\! \theta_\gamma T$, but $\mu(\theta_\gamma T)\eqo \mu(\theta_\sigma T)\eqo h$, which entails $\mu( \theta_\gamma T\backslash \theta_\sigma T)\eqo 0$ and $\supp \mu$ 
cannot be connected.

We now show $(vi)$. Clearly $T_{\! \mu , h} \!\subseteq\!  \mathtt{Span} \, (\supp  \mu )$ for all $h\ino \bbR_+^*$ and thus 
$T^o\!\subseteq\!  \mathtt{Span} \, (\supp \mu )$. Let 
$\sigma \ino \supp \mu$ be distinct from $\rho$ and let $\gamma \ino \lgeo \rho , \sigma \lgeo\, $. 
Then $B( \sigma , r) \!\subseteq\! \theta_\gamma T$ for all $r \leko  d(\gamma , \sigma)$. Therefore, $0\leko h\! :=\! \mu (B(\sigma, r)) \leqo \mu (\theta_\gamma T)$ and $\gamma \ino T_{\! \mu , h}$. Consequently,  
$ \lgeo \rho , \sigma \lgeo\, \subseteq \!  T^o$. We have shown that $\bigcup_{\sigma \in\supp  \mu} \lgeo \rho, \sigma \lgeo\,  \subseteq \! T^o$, which entails $(vi)$. \cqfd 

\smallskip

The three following results provide a control of mass-erased trees  in terms of the Hausdorff distance. They are used later to study Gromov--Hausdorff properties of mass-erased trees.

\begin{lemma}
\label{Gesubtr} Let $(E, d)$ be Polish. Let $\dPro$ stands for the Prokhorov distance on $\cM_f(E)$.
Let $\rho \ino T$ and $\rho'\ino T^\prime$ be closed subsets of $E$ such that $(T, d, \rho)$ and $(T', d, \rho^\prime)$ are rooted $\bR$-trees. 
\begin{compactenum}

\smallskip

\item[$(i)$] We assume that $T^\prime \!\subseteq\! T^{(\varepsilon)}$ and that $d(\rho, \rho^\prime) \leqo \varepsilon$. Then, for all $\sigma \ino T$,  
\begin{equation}
\label{frounssi}
\Big( \, T^\prime\! \cap \! \big( \theta_\sigma T\big)^{\!(\varepsilon)}\!\! \neq \! \emptyset\,  \Big) \Longrightarrow \Big( \exists \sigma^\prime \!  \in\!  T^\prime : \; d(\sigma , \sigma^\prime) \! \leq \! 4 \varepsilon  \; \,  \textrm{and} \; \,   T^\prime \! \cap \! \big( \theta_\sigma T\big)^{\!(\varepsilon)}\!\!  \subseteq \! \theta_{\sigma^\prime} T^\prime\,  \Big).  
\end{equation}
\item[$(ii)$] Let $\mu \ino \cM_f (T)$ and $\mu^\prime \ino \cM_f (T^\prime)$ (i.e.~$\mu, \mu'\ino \cM_f(E)$ and $\mu(E\backslash T)\eqo \mu'(E\backslash T')\eqo 0$).  In the setting of $(i)$, we further assume that 
$\dPro (\mu, \mu^\prime) \! <\!  \varepsilon$. Then for all $h\! \in \! (\varepsilon, \infty)$
\begin{equation}
\label{frounnsi}
 T_{\! \mu , h} \subseteq  
\big(T^\prime_{\! \mu^\prime , h-\varepsilon}\big)^{\!(4\varepsilon )}\; .  
\end{equation}
\end{compactenum}
\end{lemma}
\noi
\textbf{Proof.}  Let $\sigma \! \in \! T$ be such that $ T^\prime\! \cap \! \big( \theta_\sigma T\big)^{{\!(\varepsilon)}}_{{\!  }} \!\! \neq \! \emptyset$. 
Set $\sigma^\prime \! = \! \mathtt{mrca} \big(T^\prime\cap \big( \theta_\sigma T\big)^{{\!(\varepsilon)}}_{{\!  }}\big)$. 
 By definition of the most recent common ancestor, $ T^\prime \! \cap \! 
 \big( \theta_\sigma T\big)^{{\!(\varepsilon)}}_{{\!  }}\!  \subseteq \!  \theta_{\sigma^\prime} T^\prime$. We now prove that 
 $d(\sigma, \sigma^\prime ) \leqo 4 \varepsilon $.

We first assume that $\sigma^\prime \!  \notin \!   T^\prime \! \cap \! \big( \theta_\sigma T\big)^{{\!(\varepsilon)}}_{{\!  }}$. 
Let $a^\prime \! \in  \!  T^\prime \! \cap \!  \big( \theta_\sigma T\big)^{{\!(\varepsilon)}}_{{\!  }}$.  
We apply Lemma \ref{mrcadef} in $(T'\! , d, \rho')$ to get 
$b^\prime_n \!  \in \!  T^\prime \! \cap \! \big( \theta_\sigma T\big)^{{\!(\varepsilon)}}_{{\!  }}$ such that 
$c^\prime_n \! : = \! a^\prime \!  \wedge \! b_n^\prime\! \to \! \sigma^\prime$. We note that $c'_n\ino T'\!$ and 
w.l.o.g.~we can assume for all $n$ that 
$c^\prime_n \! \notin \!  T^\prime \! \cap \! \big( \theta_\sigma T\big)^{{\!(\varepsilon)}}_{{\!  }}$. Let $a, b_n , c_n\in T$ 
be such that $d(a, a^\prime)\vee d(b_n, b^\prime_n)\vee d(c_n, c_n^\prime) \leq \varepsilon$. We can choose 
$a, b_n \! \in\!  \theta_\sigma T$. Since $c_n^\prime \! \notin \! 
T^\prime \! \cap \! \big( \theta_\sigma T\big)^{{\!(\varepsilon)}}_{{\!  }}$ and since $c_n'\ino T'\!$,  
we get $c_n \! \notin \!  \theta_\sigma T$, which implies that 
$d(a, c_n)\! = \! d(a, \sigma) + d(\sigma , c_n)$ and $d(b_n, c_n) \! = \! d(b_n, \sigma) + d(\sigma , c_n)$. Thus, 
\begin{eqnarray*}
 d(a^\prime, b_n^\prime) = d(a^\prime , c_n^\prime) + d(c_n^\prime , b_n^\prime) & \geq & 
 d(a, c_n) + 
d(b_n, c_n) -4\varepsilon \\ 
&\geq & d(a, \sigma) + d(\sigma , b_n) + 2 d(\sigma , c_n) -4 \varepsilon .
\end{eqnarray*}
But  
$d(a, \sigma) + d(\sigma , b_n) \! \geq \! d(a,b_n) \geq d(a^\prime, b_n^\prime) \! -\! 2\varepsilon$ and so 
the previous inequalities imply that 
$d(\sigma, c_n) \! \leq \! 3 \varepsilon$ and thus $d(c^\prime_n , \sigma) \! \leq \! 4 \varepsilon$. This implies that $d(\sigma^\prime , \sigma)\!  \leq \! 4 \varepsilon$ by letting $n\! \to \! \infty$.

We next suppose that $\sigma^\prime \! \in \! T^\prime \! \cap \! \big( \theta_\sigma T\big)^{{\!(\varepsilon)}}_{{\!  }}$; 
let $\gamma \! \in \! \theta_\sigma T $ such that $d(\sigma^\prime, \gamma) \! \leq \! \varepsilon$. We first assume that 
$\sigma^\prime \! \neq \!  \rho^\prime $: let $x^\prime \! \in \! \lgeo \rho^\prime, \sigma^\prime \lgeo $ be such that 
$d(x^\prime , \sigma^\prime ) \! \leq \!  \varepsilon$ and let $x\! \in \! T$ be such that 
$d(x, x^\prime) \! \leq \! \varepsilon$; we note that $x^\prime 
\! \notin \! T^\prime \! \cap \! \big( \theta_\sigma T\big)^{{\!(\varepsilon)}}$ by definition of $\sigma^\prime \! = \! \mathtt{mrca} \big(T^\prime\cap \big( \theta_\sigma T\big)^{{\!(\varepsilon)}}\big)$. 
Since $x'\ino T'$, it entails that 
$x\! \notin \! \theta_\sigma T$ and $\sigma \! \in \! \lgeo x, \gamma \rgeo$; consequently, 
$d(\gamma , \sigma) \! \leq \! d(\gamma , x) $; but 
$d(\gamma , x)\!  \leq \! d(\sigma^\prime \! , x^\prime) + 2\varepsilon \! \leq \! 3 \varepsilon$; 
therefore, $d(\sigma, \sigma^\prime) \! \leq \! d(\gamma , \sigma) + \varepsilon \leq 4\varepsilon$. 
 
It remains to consider the case where 
$\sigma^\prime\!  \in \! T^\prime\! \cap \! \big( \theta_\sigma T\big)^{{\!(\varepsilon)}}\!$ and 
$\sigma^\prime \eqo  \rho^\prime$. Then, $d(\sigma, \sigma^\prime) \leqo d(\gamma , \sigma)+ \varepsilon$  $\leqo d(\gamma, \rho)+\varepsilon  \leqo d(\rho, \sigma^\prime) +2\varepsilon \! = \!  d(\rho, \rho^\prime) + 2\varepsilon \! \leq \! 3\varepsilon$. 
This completes the proof of (\ref{frounssi}).

Let us prove $(ii)$. Let $\sigma \! \in \! T_{\! \mu , h}$. Since $\dPro (\mu, \mu^\prime) \! <\!  \varepsilon$, 
$\mu^\prime ((\theta_\sigma T)^{\!(\varepsilon)}) \! \geq \! \mu (\theta_\sigma T) \! -\! \varepsilon \geqo 
h \! -\! \varepsilon$. This implies that $T^\prime\! \cap \! \big( \theta_\sigma T\big)^{{\!(\varepsilon)}}_{{\!  }}\!\!  \neq \! \emptyset$ and by (\ref{frounssi}) there exists $\sigma^\prime \! \in \! T^\prime$ such that $d(\sigma, \sigma^\prime) \! \leq \! 4 \varepsilon$ and $T^\prime\! \cap \! \big( \theta_\sigma T\big)^{{\!(\varepsilon)}}_{{\!  }}\!\! \subseteq \theta_{\sigma^\prime} T^\prime$ which implies that 
$\mu^\prime (\theta_{\sigma^\prime} T) \! \geq \! h-\varepsilon$. This completes the proof of (\ref{frounnsi}). \cqfd  

\smallskip

The previous lemma implies the following proposition. 
\begin{proposition}
\label{preGrera} Let $(E, d)$ be Polish. Let $\dHaus$ stand for the Hausdorff distance on $(E,d)$ and  $\dPro$ for 
the Prokhorov distance on $\cM_f(E)$. 
Let $\rho \ino T$ and $\rho'\ino T^\prime$ be closed subsets of $E$ such that $(T, d, \rho)$ and that 
$(T'\! , d, \rho^\prime)$ are rooted $\bR$-trees. Let $\mu \ino \cM_f (T)$ and $\mu^\prime \ino \cM_f (T^\prime)$ 
(i.e.~$\mu, \mu'\ino \cM_f(E)$ and $\mu(E\backslash T)\eqo \mu'(E\backslash T')\eqo 0$). 
We assume
\begin{equation}
\label{erahypo1}
\dHaus(T, T^\prime) < \varepsilon \quad \textrm{and} \quad \dPro (\mu, \mu^\prime) < \varepsilon 
\end{equation}  
Then for all $h \in (2\varepsilon , \infty)$ we get 
\begin{equation}
\label{founs}
\forall h^\prime \in [h-\varepsilon, h+ \varepsilon], \quad 
\dHaus \big( T_{\! \mu, h}, T^\prime_{\! \mu^\prime , h^\prime}\big)\leq 4\varepsilon + \dHaus (T_{\! \mu , h-2\varepsilon},T_{\! \mu , h+2\varepsilon} ) \; .
\end{equation}
\end{proposition}
\noi
\textbf{Proof.} Lemma \ref{Gesubtr} first implies 
$T_{\! \mu, h^\prime+ \varepsilon} \subo(T^\prime_{\! \mu^\prime \! , h^\prime})^{(4\varepsilon)}$ and 
$T^\prime_{\! \mu^\prime \! , h^\prime} \subo (T_{\! \mu, h^\prime- \varepsilon})^{(4\varepsilon)}$ for all 
$h^\prime\ino (2\varepsilon, \infty)$. We fix 
$u \geko \dHaus (T_{\! \mu, h^\prime- \varepsilon} , T_{\! \mu, h^\prime+ \varepsilon})$. This implies
$\dHaus (T_{\! \mu, h^\prime \pm \varepsilon} , T_{\! \mu, h}) \leko  u$ for all 
$h\ino [h^\prime \! -\! \varepsilon,  h^\prime +\varepsilon]$ and we see that
$T_{\! \mu , h} \subseteq \big(T_{\! \mu , h^\prime + \varepsilon}\big)^{(u)} \subseteq  
\big((T^\prime_{\! \mu^\prime \! , h^\prime})^{(4\varepsilon)}\big)^{(u)}
\subseteq (T^\prime_{\! \mu^\prime \! , h^\prime})^{(4\varepsilon+ u)}$. Similarly, we also find that 
$ T^\prime_{\! \mu^\prime \! , h^\prime}\! \subseteq\! (T_{\! \mu, h^\prime- \varepsilon})^{(4\varepsilon)} 
\!\subseteq\! \big(\big( T_{\! \mu , h} \big)^{(u)}  \big)^{(4\varepsilon)} \! \subseteq\! \big( T_{\! \mu , h} \big)^{(4\varepsilon+ u)}\! $. 
This yields 
$\dHaus\big( T_{\! \mu, h}, T^\prime_{\! \mu^\prime \! , h^\prime}\big) \leqo 4\varepsilon + u$ and thus (\ref{founs}) since 
$\dHaus (T_{\! \mu, h^\prime- \varepsilon} , T_{\! \mu, h^\prime+ \varepsilon}) \leqo
 \dHaus (T_{\! \mu , h-2\varepsilon},T_{\! \mu , h+2\varepsilon} )$. \cqfd 
  
\begin{proposition}
\label{fixconti} Let $(T, d, \rho)$ be a Polish $\bbR$-tree. Let $\mu, \mu_n\ino \cM_f (T)$ and 
$h, h_n \ino \bbR^*_+$, $n\ino \bN$. Suppose that $h^\prime\! \mapsto \!T_{\! \mu , h^\prime}$ is 
$\dHaus$-continuous at $h$. 
Assume that $h_n \!  \rightarrow \!  h$ and that $\mu_n \! \rightarrow \! \mu$ weakly in $\cM_f(T)$. Then \vspace{-2mm}
\begin{equation}
 \label{contsubtr}
\dHaus (T_{\! \mu_n , h_n} , T_{\! \mu , h}) \! \rightarrow \! 0 \quad \textrm{and} \quad  \mathtt{P}_{T_{\! \mu_n , h_n}} \mu_n \, \xrightarrow[n\to \infty]{\textrm{(weakly)}} \mathtt{P}_{T_{\! \mu , h}} \mu\, .
 \end{equation}
\end{proposition} 
\noi
\textbf{Proof.} This is a consequence of Proposition \ref{preGrera} and Lemma \ref{GrC0proj}. \cqfd

\subsection{Mass erasure in fixed measured $\bbR$-trees: definition and first properties} 
\label{merafixsec}
Let $(T, d, \rho)$ be a rooted Polish $\bbR$-tree and  let $\bm \ino \cM_f (T)$. We think of $\mu$ as a distribution of mass and we want to define how to erase a fixed amount $h \ino \bbR^*_+$ of mass from this measure starting from the leaves. Let us first give an informal description on a rooted finite discrete tree and to simplify, we assume here 
that at each vertex there is a positive amount of mass. 
Then, a quantity $h$ of mass is moved downwards to the root according to following local rules: 
\begin{compactenum}

\smallskip

\item[$-$] at a leaf with a mass $m$, a quantity $h\wedge m$ is handed down to the vertex right below;

\smallskip

\item[$-$] at an internal vertex with mass $m$, a quantity $m^\prime$ eventually arrives from the vertices above and a quantity $(m+m^\prime)\wedge h$ is handed down to the vertex right below.

\smallskip

\end{compactenum}

\noi
To define this process more formally in a continuous setting, we first state the following result, which in part explains what happens at leaves and branch points of an $\bbR$-tree of finite type, i.e.~how quantities that are $h$-multiples of mass are handed down.   
\begin{lemma}
\label{transpo} Let $(T, d, \rho)$ be a rooted Polish $\bbR$-tree. Let $\rho \ino \ttt\!\subseteq\! T$ be a subtree of finite type as in Definition \ref{spandex} $(b)$. 
Recall that for all $\sigma\! \in \! \ttt\backslash \{ \rho \}$, $\bn (\sigma, \ttt)$ is the finite number of connected components of $\ttt \backslash \{ \sigma \}$. We set $\mathtt{outdeg} (\sigma, \ttt)\! = \! \bn (\sigma, \ttt)\! -\! 1$ and $\mathtt{outdeg} (\rho, \ttt)\! = \! \bn (\rho, \ttt)$. 
We define the following signed measure: 
\begin{equation}
\label{signmedef}
\varpi_\ttt \! = \! \sum_{\sigma \in \ttt} \!  \big(\mathtt{outdeg} (\sigma, \ttt) \! -\! 1  \big)\delta_\sigma = 
(\bn (\rho, \ttt)\! -\! 1) \delta_{\rho}+\!\!\!  \sum_{\sigma \in \mathtt{Br} (\ttt) \backslash \{ \rho\}}\!   \big(\bn (\sigma, \ttt) \! -\! 2  \big)\delta_\sigma \, - \!  \sum_{\sigma \in \mathtt{Lf} (\ttt)} \! \delta_{\sigma}\; .
\end{equation}
\vspace{-3mm}

\noi Then, the following holds true. 
\begin{compactenum}
\smallskip
\item[$(i)$] For all $\sigma \! \in \! T$, $\varpi_{\ttt} (\theta_\sigma T) \! =\!  -\un_{\{ \ttt \, \cap \, \theta_\sigma T \, \neq \, \emptyset \}}$. 

\smallskip

\item[$(ii)$] For all $B\!\subseteq\! T$, $|\varpi_\ttt (B)| \leqo \# \mathtt{Lf} (\ttt) $. 

\smallskip
\item[$(iii)$] If $\rho \! \in \! \ttt^\prime \! \subseteq \!  \ttt$ is 
connected and closed, then $\mathtt{P}_{\!\ttt^\prime} \varpi_\ttt = \varpi_{\ttt^\prime}$ (see Remark \ref{signedproj} for the definition of projection of signed measures). 
\end{compactenum}
\end{lemma}
\noi
\textbf{Proof.} If $\sigma \! \notin \! \ttt$, then clearly $\varpi_\ttt(\theta_\sigma T)\! = \! 0$. Let $\sigma \! \in \! \ttt$. Combinatorial arguments on discrete trees imply that 
\vspace{-4mm}
 $$\sum_{ \gamma \in \{ \sigma \} \cup ( \mathtt{Br}( \ttt) \cap \, \theta_\sigma \! T) }  \!\!\!  \!\!\!  \!\!\!  \!\!\!  \!\!\! \mathtt{outdeg} 
(\gamma, \ttt)= \# \big( \mathtt{Lf} (\ttt) \! \cap \! \theta_\sigma T \big) + 
\# \big( \{ \sigma \}\!  \cup \! (\mathtt{Br} (\ttt)\!  \cap \! \theta_\sigma T) \big) -1 $$
which implies $\varpi_\ttt (\theta_\sigma T)\! = \! -1$.  
To prove $(ii)$, we note that $\bn (\rho, \ttt)\! -\! 1+ \sum_{\sigma \in \mathtt{Br} (\ttt)\backslash \{ \rho\}}\!\!   \big(\bn (\sigma, \ttt) \! -\! 2  \big)\eqo \# \mathtt{Lf} (\ttt)\! -\! 1$, which easily implies that $|\varpi_\ttt (B)| \leqo \# \mathtt{Lf} (\ttt) $ for all $B\!\subseteq\! T$. 

Let us prove $(iii)$. W.l.o.g.~we assume that $T\! = \! \ttt$. Let $\ttt_i^o$, $i\! \in \! I$, be the connected components 
of $\ttt\backslash \ttt^\prime$ and let $\sigma_i \! \in \! \ttt^\prime$ be such that $\{ \sigma_i \} \cup \ttt^o_i $ is the closure of $\ttt^o_i$. Since $\varpi_\ttt (\theta_\sigma \ttt)\! = \! -1$ for all $\sigma \! \in \! \ttt^o_i$, it easily follows from Lemma \ref{AboTrprop} $(iv)$ that $\varpi_\ttt (\ttt^o_i)\! = \! -1$ for all $i\ino I$. 
Thus, $\mathtt{P}_{\!\ttt^\prime} \varpi_\ttt\!  = \! \varpi_{\ttt} (\cdot \cap \ttt^\prime) -\! \sum_{i\in I} \delta_{\sigma_i}$ (see Remark \ref{signedproj} for the definition of projection as a linear function on the space of signed measures).  
Then for all $\sigma \! \in \! \ttt^\prime$, we see that $\# \{ i\! \in \! I : \sigma_i \! = \! \sigma\} \! = \!    \mathtt{outdeg} (\sigma, \ttt) \! -\! \mathtt{outdeg} (\sigma, \ttt^\prime)$. Thus, 
$$\mathtt{P}_{\!\ttt^\prime} \varpi_\ttt \eqo 
\varpi_{\ttt} (\cdot \cap \ttt^\prime) \! -\! \sum_{\sigma\in \ttt^\prime} \big( \mathtt{outdeg} (\sigma, \ttt) \! -\! \mathtt{outdeg} (\sigma, \ttt^\prime)\big) \delta_\sigma \eqo  \sum_{\sigma\in \ttt^\prime} (\mathtt{outdeg} (\sigma, \ttt^\prime) \! -\! 1\big) \delta_{\sigma } \eqo   \varpi_{\ttt^\prime}, $$
which completes the proof of the lemma. \cqfd     

\smallskip

The following proposition formally introduces mass erasure (see Lemma \ref{erasedmassdefintro}). We recall Definition \ref{deferasedtree} of $h$-mass-erased subtree (and also Lemma \ref{erasedsubtree}) and we recall from Proposition \ref{projdef} $(vi)$ an explicit expression of the projection of a measure onto a subtree.  
\begin{proposition}
\label{erasedmassdef}
Let $(T, d, \rho, \mu)$ be a rooted finitely measured Polish $\bbR$-tree and let $h \ino \bbR^*_+$. Then 
$\nu \! := \! \mathtt{P}_{T_{\! \mu, h}} \mu + h \varpi_{T_{\! \mu, h}}$ is the unique positive Borel measure on $T$ such that 
\begin{equation}
\label{mainprop}
 \nu ( \theta_\sigma T) = \big( \mu (\theta_\sigma T)\!  -\! h \big)_{\! +} \; , \quad \sigma \ino T \; . 
 \end{equation}
Moreover, $\nu \! = \! \mathtt{P}_{T_{\! \mu, h+}} \mu + h \varpi_{T_{\! \mu, h+}}$, where we recall that $T_{\! \mu, h+}$ is the closure of $\, \bigcup_{h^\prime >h} T_{\! \mu , h^\prime} $. 

We call $\nu$ the \emph{$h$-erasure of $\mu$} and we use notation $\cE_h \mu \eqo  \nu$. 
We also adopt the convention $\era_0 \bm \eqo \bm$. 
\end{proposition}
\noi
{\bf Proof.} To simplify notation we denote by $\pi$ the projection $\mathtt{P}_{T_{\! \mu, h}} \mu$ of 
$\mu$ onto $T_{\! \mu , h}$. We first show that $\nu \eqo \pi + h\varpi_{T_{\! \mu, h}}$ is a positive measure. 
By (\ref{signmedef}), we only need to check that $\nu (\{ \sigma \} ) \geqo 0$ for any 
$\sigma \ino \mathtt{Lf} (T_{\! \mu , h})$. \emph{Indeed}, observe that $ \pi  (\{ \sigma \})\eqo \mu (\theta_\sigma T)$ 
since $\sigma$ is a leaf of $T_{\! \mu, h}$ and by definition of the projection of $\mu$ onto $T_{\! \mu,h}$. Thus 
$\nu (\{ \sigma \})\eqo \pi  (\{ \sigma \}) \! -\! h\eqo  \mu ( \theta_\sigma T)\! -\! h\geqo 0$ since $\sigma \ino T_{\! \mu, h}$. 
Let us next prove that $\nu$ satisfies (\ref{mainprop}). Let $\sigma \ino T$. By Lemma \ref{projdef} $(vi)$  
if $\sigma \ino T_{\! \mu , h} $, then 
$\pi (\theta_\sigma T) \eqo  \mu (\theta_\sigma T) \geqo h$ and Lemma \ref{transpo} entails
$h\varpi_{T_{\! \mu , h}}(\theta_\sigma T) \! = \! -h$, which implies that 
$\nu ( \theta_\sigma T) \eqo \bm (\theta_\sigma T)\!  -\! h \geqo 0$. If $
\sigma \! \in\!  T\backslash T_{\! \mu , h}$, then $\mu(\theta_\sigma T) \leqo h$ and Lemma \ref{transpo} 
implies that $h\varpi_{T_{\! \mu , h}}(\theta_\sigma T) \! = \! 0 $. Moreover Lemma \ref{projdef} $(vi)$ entails 
$\pi (\theta_\sigma T) \eqo 0$. Therefore we get  $\nu (\theta_\sigma T) \! = \! 0 \! = \! (\mu (\theta_\sigma T) \! -\! h)_+$. 
This proves that $\nu$ is a positive measure satisfying (\ref{mainprop}). Uniqueness is a consequence of 
Lemma \ref{AboTrprop} $(vi)$ and uniqueness of the extension of a finite measure on a pi-system. 

We next set $\nu' \! = \! \mathtt{P}_{T_{\! \mu, h+}} \mu + h \varpi_{T_{\! \mu, h+}}$. We check exactly as in the previous proof that 
$\nu'$ is a positive measure such that $ \nu' ( \theta_\sigma T) \eqo  \big( \bm (\theta_\sigma T)\!  -\! h \big)_{\! +}$ 
for all $\sigma\ino T$. Therefore we get $\nu\eqo \nu'$. 
\cqfd

\medskip

Mass erasure is a nice tool to approximate measures on $\bbR$-trees. It shares many features with the trimming of $\bR$-trees (also called length erasure).  
In particular, mass erasure acts as a semigroup on $\cM_f (T)$ as shown in the following proposition which also records other basic properties. 
\begin{proposition}
\label{properased}
Let $(T, d, \rho, \mu)$ be a rooted finitely measured Polish $\bbR$-tree, let $h \ino \bbR^*_+$ and let $\rho \ino \ttt \! \subseteq \! T$ be closed and connected. Then the following holds true. 
\begin{compactenum}

\smallskip

\item[$(i)$]  If $\mu (T\backslash \ttt) \! < \!  h$, then $\era_h \mu\! = \! \era_h \Pt \mu$.  

\smallskip

\item[$(ii)$] For all connected components $C$ of $T \backslash \ttt$, 
$\era_h \mu ( C)\! =\!  (\mu (C) \! -\! h)_+$.

\smallskip

\item[$(iii)$] $\era_h \mu (T\backslash \ttt) \! \leq \! \mu (T\backslash \ttt) $.

\smallskip

\item[$(iv)$] $\mathtt{Span}\,  (\supp  \era_h \mu)$ is the compact set $ T_{\! \mu, h+}$, 
which is the closure of $\, \bigcup_{h^\prime >h} T_{\! \mu, h^\prime}$.

\smallskip

\item[$(v)$] $(\era_h )_{h \in \bbR_+}$ is a semigroup on $\cM_f (T)$, i.e.~ 
$ \era_h \!  \circ\!  \era_{h^\prime} \eqo \era_{h+h^\prime}$ for all $h,h' \ino \bbR_+ $.

\smallskip

\item[$(vi)$] For all $h, h^\prime \ino \bbR^*_+$ and all $\mu \ino \cM_f (T)$, $d_{\mathtt{var}} (\era_{h+h'}\mu, \era_h\mu ) \leqo 3 \mu(T) h^\prime / h$, where for all $\nu, \nu'\ino \cM_f(T)$ we have set $d_{\mathtt{var}} (\nu, \nu')\eqo  \sup_{B \in \cB (T)} | \nu(B)\! -\! \nu'(B) |$. 
\end{compactenum}
\end{proposition}
\noi
{\bf Proof.} Let us prove $(i)$. We suppose $\mu (T\backslash \ttt) \! < \!  h$. 
Let $\sigma \ino T \backslash \ttt$. 
Since $\ttt$ is a closed and connected subset of $T$ which contains $\rho$, we get 
$\theta_\sigma T \! \subseteq \! T\backslash \ttt $ and necessarily $\mu (\theta_\sigma T) \! <\!  h$, which implies $\era_h \mu (\theta_\sigma T)\eqo 0$. Moreover, since the support of $\Pt \mu$ is included in $\ttt$, we also get $\Pt \mu (\theta_\sigma T)\! = \! 0$. This proves that $\era_h \Pt \mu (\theta_\sigma T)\! = \! \era_h \mu (\theta_\sigma T)\! = \! 0$. 
Next suppose that $\sigma \! \in \! \ttt$. By Proposition \ref{projdef} $(vi)$, $\Pt \mu (\theta_\sigma T)\! = \! \mu (\theta_\sigma T)$. Thus, $\era_h \Pt \mu (\theta_\sigma T)\!  = \! (\mu (\theta_\sigma T)\! -\! h)_+\! = \! \era_h \mu (\theta_\sigma T)$. 
We complete the proof of $(i)$ by the uniqueness in Proposition \ref{erasedmassdef}.

 Let us prove $(ii)$ and $(iii)$. Let $C_i $, $i \ino I$, be the countably many connected components of 
$T \backslash \ttt $. Let $i \ino I$ and $\sigma_i \ino \ttt$ such that 
$\{ \sigma_i\}\!  \cup \! C_i $ is the closure of $C_i$. We note that $\rho \! \notin \! C_i$. Then, we can find 
$(\gamma_n)_{n\in \bN}$, $C_i$-valued, which converges to $\sigma_i$ and such that $\gamma_{n+1} \ino \lgeo \rho, \gamma_n\rgeo$. Therefore, $\theta_{\gamma_{n}} T \! \subseteq \! \theta_{\gamma_{n+1}} T$. By Lemma \ref{AboTrprop} 
$(iv)$, $C_i \eqo \bigcup_{n\in \bN}\theta_{\gamma_{n}} T$. Hence, $ \era_{h} (\bm) (C_i )$ $\eqo$ 
$ \lim_{n \rightarrow \infty} \! $  $ \era_{h} (\bm) ( \theta_{\gamma_n} \! T)$ $\eqo $
$ \lim_{n \rightarrow \infty} \! $  $( \bm ( \theta_{\gamma_n} \! T)\! -\! h )_{+} $ $\eqo$ $( \bm (C_i) \! -\! h )_{+}$,
which proves $(ii)$. Then, we observe that  $\era_h (\mu) (T\backslash \ttt)$ $\eqo \sum_{i\in I} \era_{h} (\bm) (C_i )\eqo \sum_{i\in I}  (\bm (C_i)\! -\!h)_{\! +} \leqo \sum_{i\in I}  \bm (C_i) \eqo \mu (T\backslash \ttt)$,
which proves $(iii)$.

To prove $(iv)$, we first show that 
$\mathtt{Span} (\supp \era_h \mu ) \! \subseteq \! T_{\! \mu , h+}$. We recall that $T_{\! \mu , h+}$ 
is the closure of $T^o_{\! \mu , h+\! } \! := \! \bigcup_{h^\prime >h} T_{\! \mu, h^\prime}$. We fix 
$\sigma \! \in \! \supp (\era_h \mu)$ distinct from $\rho$ and 
we consider $\gamma \! \in \! \lgeo \rho, \sigma \lgeo$. Then $\sigma$ belongs to a connected component $C$ 
of $T\backslash \{ \gamma\} $ that does not contain $\rho$. By Lemma \ref{AboTrprop} $(iv)$, we 
get $\sigma \ino C\! \subseteq \! \theta_\gamma T$. Since $C$ is an open neighbourhood of $\sigma$ we get 
$\era_h \mu (C) \! >\! 0$ since $\sigma \! \in \! \supp (\era_h \mu)$. Consequently, 
$(\mu (\theta_\gamma T)\! -\! h )_{+}\eqo \era_h \mu ( \theta_\gamma T) \! >\! 0$ and there exists 
$h^\prime \! >\! h$ such that $(\mu (\theta_\gamma T)\! -\!h^\prime )_{ +} \! >\! 0$, which entails 
$\gamma \ino T^o_{\! \mu, h+}$. Since $\lgeo \rho, \sigma \rgeo$ is the closure of $\lgeo \rho, \sigma \lgeo\, $, we get $\lgeo \rho, \sigma \rgeo  \! \subseteq \! T_{\! \mu , h+}$. We thus have proved that 
$\mathtt{Span} \, (\supp \era_h \mu ) \! \subseteq \! T_{\! \mu , h+}$.

To complete the proof of $(iv)$, we fix $\sigma \! \in \! T_{\! \mu , h^\prime}$ with $h^\prime \! >\! h$. 
Then, $\mu (\theta_\sigma T)\! \geq \! h^\prime$, which 
entails $\era_h \mu (\theta_\sigma T) \! \geq \! h^\prime \! -\! h \! >\! 0$. Consequently, 
$\theta_\sigma T \! \cap  \! \supp \era_h \mu \! \neq \! \emptyset $, which shows that 
$\sigma \! \in \! \mathtt{Span}\,  (\supp  \era_h \mu )$. We thus have proved $T^o_{\! \mu, h+} \! \subseteq \! 
\mathtt{Span} \, (\supp \era_h \mu )$. Since $\supp \era_h \mu$ is a closed subset of the compact subset $T_{\! \mu, h}$ (by Lemma \ref{erasedsubtree}), $\supp \era_h \mu$ is compact too, and by Lemma \ref{compacthull}, so is $\mathtt{Span} \, (\supp \era_h \mu )$. Thus, $T_{\! \mu, h+} \! \subseteq \! \mathtt{Span} \, (\supp \era_h \mu )$, which implies $(iv)$.

Let us prove $(v)$. To this end we fix $h, h^\prime \ino \bbR^*_+$. By the uniqueness in Proposition \ref{erasedmassdef}, we only need to prove for all $\sigma \ino T$ that 
$ ( \mu (\theta_\sigma T)\! -\! h\! -\! h^\prime)_{\! +} \eqo  \left(  (\mu (\theta_\sigma T)\! -\! h)_{\! +} \! -\! h^\prime\right)_{\! +} $, 
which is easy to check by considering the three possible cases: $\sigma \ino T\backslash T_{\! \mu, h}$, $\sigma \ino T_{\! \mu, h} \backslash T_{\! \mu, h+h'}$ and $\sigma \ino T_{\! \mu, h+h'}$. 

We now proceed to the proof of $(vi)$. To simplify, we set $\ttt^\prime \! :=\! T_{\! \mu , h+ h^\prime} \! \subseteq \! T_{\! \mu, h}   \! =:\! \ttt$ and $\nu\! :=\! \era_h \mu$. Let $\Gamma_i$, $i\ino I$, be the connected components of $\ttt\backslash \ttt^\prime$ in the 
compact $\bbR$-tree
$\ttt$ and let $\sigma_i \ino \ttt^\prime$ be such that the closure in $\ttt$ of $\Gamma_i$ is $\{ \sigma_i \} \cup  \Gamma_i$. 
First, note that every connected component $\Gamma_i$ contains at least one leaf of $\ttt$. Therefore,
$\# I \leqo \# \mathtt{Lf} (\ttt)$. 
Then note that $ \# \mathtt{Lf} (\ttt) \leqo \mu (T) / h$. Thus, $\# I \leqo \mu (T) / h$.

Next we recall that $\cE_{h'}\nu \eqo \cE_{h+h} \mu$, by the semigroup property of mass erasure. Then $T_{\! \nu , h^\prime}\eqo \ttt^\prime$ and by the point $(ii)$ (where $\nu$, $\ttt$ and $\ttt^\prime$ play the roles of $\mu$, $T$ and $\ttt$), we get 
$\nu  (\Gamma_i) \leqo h' $. 
Consequently, $\nu (T\backslash \ttt^\prime)\eqo \nu (\ttt\backslash \ttt^\prime)\eqo \sum_{i\in I} \nu (\Gamma_i)  \leqo h^\prime \# I \! \leq \! \mu (T) h^\prime/ h$. 
Thus, for all $B \ino \cB (T)$, we get 
$0\leqo \nu (B) \! -\! \nu ( B \cap \ttt^\prime) \leqo \nu (T\backslash \ttt^\prime) \leqo \mu (T) h^\prime/ h $. 

We next recall from Proposition \ref{projdef} $(v)$ that 
 $\mathtt{P}_{\!\ttt^\prime} \nu \eqo \nu ( \cdot \cap \ttt^\prime) + \sum_{i\in I} \nu (\Gamma_i) \delta_{\sigma_i} $. This implies that 
 $0\leqo \mathtt{P}_{\!\ttt^\prime} \nu (B)  \! -\! \nu ( B \cap \ttt^\prime) \leqo  \mu (T) h^\prime/ h  $. 
Then Lemma \ref{transpo} shows that 
$$ \big| \era_{h^\prime} \nu (B) - \mathtt{P}_{\!\ttt^\prime} \nu (B)\big| \eqo h' |\varpi_{\ttt'} (B)|\! \leq \! h'\# \mathtt{Lf} (\ttt^\prime)  \leqo h'\# \mathtt{Lf} (\ttt)\leqo   \mu (T)h^\prime/ h\; .$$ 
Thus, $|\nu (B) \! -\! \era_{h^\prime} \nu (B)| \leqo 3 \mu (T) h^\prime/h$, which proves $(vi)$. \cqfd

\begin{proposition}
\label{weakcvera} Let $(T, d, \rho)$ be a rooted Polish $\bbR$-tree. We equip $\cM_f(T)$ with the topology of weak convergence. Then $(h, \bm) \ino \bbR_+ \! \times \! \cM_f (T) \! \mapsto\!  \era_h \bm \ino \cM_f (T)$ is continuous (here $\bbR_+ \! \times \! \cM_f (T)$ is equipped with the product topology). 
\end{proposition}
\noi 
\textbf{Proof.} Let $h, h_n\ino \bbR_+$ and $\mu, \mu_n \ino \cM_f (T)$, $n\ino \bN$, be such that $h_n\! \rightarrow \! h$ and $\mu_n \! \rightarrow \! \mu$ weakly. 
Let $\eta \! \in \! \bbR^*_+$. 
By Prokhorov's tightness criterion there is $K_\eta\subo T$ compact such that $\sup_{n\in \bbN} \mu_n (T\backslash K_\eta) \leko \eta$. By Lemma \ref{compacthull}, $\ttt_\eta \! := \! \mathtt{Span} (K_\eta)$ is compact too and we also see that  
$\sup_{n\in \bN}\bm_n ( T\backslash \ttt_\eta) \! \leq \! \eta $. Then Proposition \ref{properased} $(iii)$ entails $\sup_{n\in \bN} \era_{h_n} \bm_n ( T\backslash \ttt_\eta) \! \leq \! \eta $. Since $\sup_{n\in \bbN} \era_{h_n} \mu (T)\leqo \sup_{n\in \bbN} \mu_n (T)\leko \infty $, Prokhorov's tightness criterion implies that the $(\cE_{h_n} \mu_n)_{n\in \bN}$ are tight on $T$.  
To complete the proof we only need to prove that $\era_h \mu$ is the only weak limit point of $(\cE_{h_n} \mu_n)_{n\in \bN}$.   Let $\nu \! \in \! \cM_f (T)$ be such a limit point (there is at least one by Prokhorov's tightness criterion). Namely, 
there is an increasing $\bbN$-valued sequence $(n_k)_{k\in \bN}$ such that $\era_{h_{n_k}} \mu_{n_k} \! \rightarrow \! \nu$ weakly.  
We denote by $D$ the countable set $\{ \sigma  \! \in \! T \! : \! \nu (\{ \sigma \})\! >\! 0  \; \textrm{or} \;  \mu (\{ \sigma \}) \! >\! 0  \}$. Let $\sigma \! \in \! T \backslash D$. Since 
$\theta_\sigma T$ is a closed subset whose boundary is $\{ \sigma\}$ the
Portmanteau theorem asserts that 
$$ \lim_{k \rightarrow \infty} \era_{h_{n_k}} \bm_{n_k} (\theta_\sigma T)= \nu ( \theta_\sigma T) \quad {\rm and } \quad \lim_{k \rightarrow \infty} (\bm_{n_k} (\theta_\sigma T)-h_{n_k})_+ = (\bm (\theta_\sigma T) -h)_+ . $$
Therefore, for all $\sigma \in T \backslash D$, $ \nu ( \theta_\sigma T) \! = \! 
\era_h \bm (\theta_\sigma T)$, 
which entails $\nu \! = \! \era_h \bm$ by Lemma \ref{AboTrprop} $(vi)$ and the uniqueness of the extension of finite measures on pi-systems.  \cqfd

\subsection{Convergence in the sense of mass erasure in fixed $\bbR$-trees}
\label{merafix}
We recall from Remark \ref{toporem} $\textbf{(c)}$ 
that $\cM_f(T)\!\subseteq\! \cM_f(T^*)$ and we observe that mass erasure makes sense on $\cM_f(T^*)$ since $(T^*\! , d^*\! , \rho)$ is a rooted Polish $\bbR$-tree by Theorem \ref{bordif}. In the following lemma we given a representation in $\cM_f (T^*)$ of finite measures on $T$ that are consistent under mass erasure. 
\begin{lemma}
\label{descentera}  Let $(T, d, \rho)$ be a rooted Polish $\bbR$-tree, with bordification $(T^*\! , d^*\! , \rho)$. We recall the definition of $\cM_f^{\mathtt{era}}(T)$ from (\ref{meraspace}). 
Then the following holds true. 
\begin{compactenum}

\smallskip

\item[$(i)$] $\mu \ino  \cM_f^{\mathtt{era}}(T)$ iff $\era_h \mu \ino \cM_f(T)$ for all $h\ino \bbR^*_+$. In this case, 
$T^*_{\! \mu, h}$ is a compact subset of $T$ and we simply denote it by $T_{\! \mu,h}$. 

\smallskip

\item[$(ii)$] For all $p\ino \bbN$, let $\nu_p \ino \cM_f(T)$ and let $h_p\ino \bbR^*_+$ be such that 
$h_{p+1} \leko h_p$, $\era_{h_p-h_{p+1}}\nu_{p+1}\eqo \nu_p$ and $\lim_{p\rightarrow \infty} h_p\eqo 0$. Then, there exists a unique $\mu\ino \cM_f^{\mathtt{era}}(T) $ such that $\era_{h_p} \mu \eqo \nu_p$ for all $p\ino \bbN$. Moreover $\nu_p \! \to \! \mu$ weakly in $\cM_f(T^*)$. 
\end{compactenum}
\end{lemma}
\noi
\textbf{Proof.} Let us prove $(i)$. To avoid trivial cases, we assume that $\partial T\! \neq \! \emptyset$. 
We fix $\mu\ino \cM_f(T^*)$. Let $\bs \ino \partial T$.
Since $\{ \bs \} \eqo \theta_{\bs} T^*$, we get $\era_h \mu (\{ \bs\})\eqo \big(\mu   (\{ \bs \}) \! -\! h)_{+}$ 
and if for all $h\ino \bbR^*_+$, $\era_h \mu \ino \cM_f(T)$, then $\mu(\{ \mathbf s\})\eqo 0$ and 
$\mu \ino  \cM_f^{\mathtt{era}}(T)$.  
Conversely, let us suppose that $\era_h\mu\!\notin\!\cM_f(T)$ for some $h\ino\bbR^*_+$. 
Thus, there are a ray $\bs\ino\partial T_{\! \mu,h}$ and points $\sigma_n\ino\bs$ such that $d(\rho,\sigma_n)\eqo n$, 
$n\ino\mathbb{N}$, and $\bigcap_{n\in\mathbb{N}}\theta_{\sigma_n}T^*\eqo\{\bs\}$, 
which entails $\mu(\{\bs\})\eqo\lim_{n\rightarrow\infty}\mu(\theta_{\sigma_n}T^*)\!\ge\! h$. So $\mu\!\notin\!\cM_f^{\rm era}(T)$.

We next prove $(ii)$. For all $p\ino \bbN$, we set $\ttt_p\! = \! \mathtt{Span} \, (\supp \nu_p )$, which is equal to $T_{\! \nu_{p+1}, (h_p -h_{p+1})+}$ by Proposition \ref{properased} $(iv)$. This is a subtree of $T$ of finite type, therefore compact, and contains $\rho$. 
For all $p\! \in \! \bN$, we set $\mu_p \! = \! \nu_p - h_p \varpi_{\ttt_p} $, which is a signed measure, a priori.
Let us first prove that it is actually a positive measure. 

\emph{Indeed}, the semigroup property of mass erasure (Proposition \ref{properased} $(v)$) 
implies for all $q\ino \bbN$ that $\era_{h_p -h_{p+q}} \nu_{p+q} \eqo \nu_p$. Thus $\nu_p (\theta_\sigma T) \geko 0$ 
iff $  \nu_{p+q} (\theta_\sigma T) \geko h_p-h_{p+q}$, 
which implies that $\ttt_p \eqo T_{\nu_{p+q} , (h_p-h_{p+q})+}$. Therefore Proposition \ref{erasedmassdef} implies 
that $\nu_p \eqo \mathtt{P}_{\!\ttt_p} \nu_{p+q} + (h_p\! -\! h_{p+q} ) \varpi_{\ttt_p}$. For all $B \ino \mathscr B(T)$, we get 
$\mu_p (B) + h_{p+q} \varpi_{\ttt_p} (B)  \eqo   \mathtt{P}_{\!\ttt_p} \nu_{p+q}(B) \geqo 0$ since $ \mathtt{P}_{\!\ttt_p} \nu_{p+q}$ is a 
positive measure. Thus $\mu_p(B)  \geqo - h_{p+q} \varpi_{\ttt_p} (B) \! \to \! 0 $ as $q\! \to \! \infty$, which proves that $\mu_p$ is a 
positive measure. 

We next observe the following. 
\begin{eqnarray*} 
\mathtt{P}_{\!\ttt_p} \mu_{p+1} =   \mathtt{P}_{\!\ttt_p} \nu_{p+1} -h_{p+1} \mathtt{P}_{\!\ttt_p} \varpi_{\ttt_{p+1}} &=&  
 \big(\nu_p  \! -\! (h_p \! -\! h_{p+1}) \varpi_{\ttt_p} \big) -h_{p+1} \mathtt{P}_{\!\ttt_p} \varpi_{\ttt_{p+1}} \\
\textrm{(by Proposition \ref{transpo} ($iii$))}& =&  \big(\nu_p  \! -\! (h_p \! -\! h_{p+1}) \varpi_{\ttt_p} \big) -  h_{p+1} \varpi_{\ttt_p} = \mu_p 
\end{eqnarray*}
(see Remark \ref{signedproj} for the linearity of projection of signed measures). Proposition \ref{projproj} $(i)$ shows that there is a unique $\mu\ino \cM_f(T^*)$ such that $\mathtt{P}_{\!\ttt_p}^*\mu\eqo \mu_p$ for all $p\ino \bbN$ (here, the $\mathtt{P}_{\!\ttt_p}^*\mu$ are as in Proposition \ref{projdef} ($vii$)) and such that $\mu_p \! \to \! \mu$ weakly on $\cM_{\! f} (T^*)$. 

We next check that $\era_{h_p} \mu\eqo \nu_p$. To this end, we fix $\bs \ino T^*$ and we first suppose that $\mu (\theta_\bs T^*) \eqo 0$. 
By Proposition \ref{projdef} $(vi)$, either $\bs \ino \ttt_p$ and $\mathtt{P}^*_{\ttt_p}\mu  (\theta_\bs T^*) \eqo \mu (\theta_\bs T^*) \eqo 0$ or $\bs \! \notin \! \ttt_p$ and $\mathtt{P}^*_{\ttt_p} \mu (\theta_\bs T^*) \eqo 0$. In both cases, we get  $\mu_p(\theta_\bs T^*)\eqo \mathtt{P}^*_{\ttt_p}\mu  (\theta_\bs T^*) \eqo 0$. Thus 
$\nu_p(\theta_\bs T^*) \eqo h_p \varpi_{\ttt_p} (\theta_\bs T^*)$. But Lemma \ref{transpo} $(i)$ in $T^*$ implies that 
$h_{p}\varpi_{\ttt_p} (\theta_\bs T^*) \leqo 0$, whereas $\nu_p(\theta_\bs T^*)\geqo 0$. Therefore $\nu_p(\theta_\bs T^*) \eqo h_p \varpi_{\ttt_p} (\theta_\bs T^*)\eqo  0$ and $\nu_p(\theta_\bs T^*) \eqo ( \mu(\theta_\bs T^*) \! -\! h_p)_+$.

We next assume that  $\mu(\theta_\bs T^*) \geko 0$ and that $\mu (\{ \bs \} )\eqo 0$. By the Portmanteau theorem 
$\lim_{q\to \infty}$ $ \mu_q(\theta_\bs T^*)$ $ \eqo $ $ \mu(\theta_\bs T^*) $. 
Consequently, there exists an integer $q\geko p$ such that  $\mu_q(\theta_\bs T^*) \geko 0$. Since 
$\mu_q \eqo \mathtt{P}_{\!\ttt_q} \mu$, this implies that $\bs \ino \ttt_q$ and Proposition \ref{projdef} $(vi)$
 implies that $\mu_q(\theta_\bs T^*) \eqo  \mathtt{P}_{\!\ttt_q} \mu (\theta_\bs T^*)\eqo \mu (\theta_\bs T^*)$. 
Thus, we get $0\leko \mu (\theta_\bs T^*)\eqo \nu_q(\theta_\bs T^*)\! -\! h_q \varpi_{\ttt_q} (\theta_\bs T^*)$. Since $\bs \ino \ttt_q$, Lemma \ref{transpo} $(i)$ implies $ \varpi_{\ttt_q} (\theta_\bs T^*) \eqo -1$. 
Therefore $\nu_q (\theta_\bs T^*) \eqo \mu (\theta_\bs T^*)\! -\! h_q$, which entails 
$(\mu (\theta_\bs T^*)\! -\! h_p)_+ $ $\eqo$ $ (\nu_q (\theta_\bs T^*) \! -\! (h_p\! -\! h_q))_+ $ $\eqo$ $ \era_{h_p-h_q} \nu_q (\theta_\bs T^*) $ $ \eqo$ $ \nu_p(\theta_\bs T^*) $. We thus have proved that for all $\bs \ino T^*$ 
such that $\mu(\{ \bs \} )\eqo 0$, $\nu_p(\theta_\bs T^*)\eqo (\mu(\theta_\bs T^*) \! -\! h_p)_+\eqo \era_{h_p} \mu (\theta_\bs T^*)$. Then Lemma \ref{AboTrprop} $(vi)$ combined with standard arguments entails that $\nu_p \eqo \era_{h_p} \mu$ which is the desired result. 

Finally we observe that for all $h\ino \bbR^*_+$, there is $p\ino \bbN$ such that $h_p\leko h$ and thus $\era_h \mu \eqo \era_{h-h_p} \nu_p$ by the semigroup property of mass erasure (Proposition \ref{properased} $(v)$). Thus $\era_h \mu \ino \cM_f(T)$, and Lemma \ref{descentera} (i) entails that $\mu \ino \cM_f^{\mathtt{era}}(T)$. The weak convergence of the $\nu_p$ is a direct application of Proposition \ref{weakcvera} in $T^*$. \cqfd 

\smallskip

We introduce the convergence in the sense of mass erasure in $\cM_f^{\rm era}(T)$. 
\begin{definition}
\label{meracvfixdef} Let $(T, d, \rho)$ be a rooted Polish $\bbR$-tree with bordification 
$(T^*\! , d^*\! , \rho)$. We recall $\cM_f^{\mathtt{era}}(T)$ from (\ref{meraspace}). 
Let $\mu, \mu_n \ino \cM_f^{\mathtt{era}}(T)$, $n\ino \bbN$. The sequence $(\mu_n)_{n\in \bbN}$ converges to $\mu$ \emph{in the sense of mass erasure} if $\era_h\mu_n \!\!  \to \! \era_h \mu$ weakly in $\cM_f(T)$ for all $h\ino \bbR^*_+$,  
which makes sense by Lemma \ref{descentera}. \cq 
\end{definition}
\begin{lemma}
\label{weakvsmera1} We keep the previous notation. 
Let $\mu, \mu_n\ino \cM^{\mathtt{era}}_f(T)$, $n\ino \bbN$. 
\begin{compactenum}

\smallskip

\item[$(i)$] If $\mu_n \! \to \! \mu$ weakly in $\cM_f(T^*)$, then $\mu_n \! \to \! \mu$ in the sense of mass erasure. The converse is not true in general (even if $\mu \ino \cM_f(T)$: see Example \ref{dilutestar}).

\smallskip

\item[$(ii)$] If $\mu_n \! \to \! \mu$ in the sense of mass erasure, if $h_n \! \to h \ino \bbR^*_+$ and if $h'\! \mapsto T_{\! \mu, h'}$ is $\dHaus$-continuous at $h$, then $\dHaus (T_{\! \mu_n, h_n} , T_{\! \mu, h})\! \to \! 0$ and $\mathtt{P}_{T_{\! \mu_n, h_n}} \mu_n \! \to \! \mathtt{P}_{T_{\! \mu, h}} \mu $ weakly in $\cM_f(T)$.

\end{compactenum}
\end{lemma}
\noi
\textbf{Proof.} Let us first prove $(i)$: we apply Proposition \ref{weakcvera} in $T^*$ to get $\era_h\mu_n \! \to \! \era_h \mu$ weakly in $\cM_f(T^*)$ for all $h\ino \bbR^*_+$. Now observe that $ \era_h\mu_n$ and $\era_h\mu$ belong to $\cM_f(T)$ by Lemma \ref{descentera}. Therefore $\era_h\mu_n \! \to \! \era_h \mu$ weakly in $\cM_f(T)$ by Lemma \ref{embedmeas}. This completes the proof of $(i)$. 

Let us prove $(ii)$. Let us first prove that $\dHaus (T_{\! \mu_n, h_n} , T_{\! \mu, h})\! \to \! 0$. To this end we fix 
$h'\ino (0, h)$ and we set $\mu'_n \eqo \era_{h'} \mu_n$ and $\mu'\eqo \era_{h'} \mu$. 
Then $\mu_n'\! \to \! \mu'$ weakly in $\cM_f(T)$ since $\mu_n \! \to \! \mu$ in the sense of mass erasure.  
W.l.o.g.~we assume for all 
$n\ino \bbN$ that $h_n\geko h'$ and we set $h'_n\eqo h_n\! -\! h'$.   Observe that $\mu'_n (\theta_\sigma T)\geqo h'_n$ 
iff $\mu_n(\theta_{\sigma} T^*) \geqo h'+h'_n\eqo h_n$. This shows that 
$T_{\! \mu'_n , h'_n}\eqo T_{\! \mu_n, h_n}$. Similarly, we get $T_{\! \mu'\! , h''}\eqo T_{\! \mu, h''\! +h'}$ for all $h''\ino \bbR^*_+$. 
Thus $h''\ino \bbR^*_+ \! \mapsto \! T_{\! \mu', h''}$ is $\dHaus$-continuous at $h\! -\! h'$: by Proposition \ref{fixconti},  
we get $\dHaus (T_{\! \mu_n, h_n} , T_{\! \mu, h})\! \to \! 0$ and 
$\mathtt{P}_{T_{\! \mu_n, h_n}} \mu'_n \! \to \! \mathtt{P}_{T_{\! \mu, h}} \mu'$ weakly in $\cM_f(T)$. 
  
Now recall that $\mu'_n\eqo \mathtt{P}_{T_{\! \mu_n, h'}} \mu_n + h'\varpi_{T_{\! \mu_n, h'}}$ (resp.~$\mu'\eqo \mathtt{P}_{T_{\! \mu, h'}} \mu + h'\varpi_{T_{\! \mu, h'}}$). 
By the linearity of projections of signed measures (see Remark \ref{signedproj}) and by Lemmas \ref{projprok} and \ref{transpo} 
$(iii)$, we get $\mathtt{P}_{T_{\! \mu_n, h_n}} \mu'_n\eqo \mathtt{P}_{T_{\! \mu_n, h_n}}\mu_n +h'\varpi_{T_{\! \mu_n, h_n}} $ (resp.~$\mathtt{P}_{T_{\! \mu, h}} \mu'\eqo \mathtt{P}_{T_{\! \mu, h}}\mu +h'\varpi_{T_{\! \mu, h}} $). We next observe that for all $B\ino \mathscr B(T)$, $|\varpi_{T_{\! \mu_n, h_n}}(B)|\leq \# \mathtt{Lf} (T_{\! \mu_n, h_n})\leq \mu_n(T^*)/h_n$ (resp.~$|\varpi_{T_{\! \mu, h}}(B)|\leqo \# \mathtt{Lf} (T_{\! \mu, h})\leq \mu(T^*)/h$). 
If $d_{\mathtt{var}}$ stands for the distance in variation in $\cM_f(T)$, 
then we get 
$d_{\mathtt{var}} ( \mathtt{P}_{T_{\!\mu_n, h_n}} \mu'_n, \mathtt{P}_{T_{\! \mu_n, h_n}}\mu_n) \leqo h'\mu_n(T^*)/h_n$ 
and $d_{\mathtt{var}} ( \mathtt{P}_{T_{\! \mu, h}} \mu', \mathtt{P}_{T_{\! \mu, h}}\mu) \leqo h'\mu(T^*)/h$. 
Next note that $\dPro \leqo d_{\mathtt{var}}$. Then 
\begin{eqnarray*}
\dPro \big( \mathtt{P}_{T_{\! \mu_n, h_n}}\mu_n,  \mathtt{P}_{T_{ \! \mu, h}}\mu \big) \leqo \dPro \big(  \mathtt{P}_{T_{\! \mu_n, h_n}}\mu_n  \!\!\!\! \!\! &,& \!\! \!\!\!\!  \mathtt{P}_{T_{\! \mu_n, h_n}}
 \mu'_n\big) \\
 \!\!\!\! &+ & \!\!\!\! \dPro \big( \mathtt{P}_{T_{\! \mu_n, h_n}} \mu'_n , \mathtt{P}_{T_{\! \mu, h}} \mu' \big) + \dPro
\big(  \mathtt{P}_{T_{\! \mu, h}} \mu', \mathtt{P}_{T_{\! \mu, h}}\mu \big)  .
\end{eqnarray*}
Since $\lim_{n\to \infty}\mathtt{P}_{T_{\! \mu_n, h_n}} \mu'_n \! = \! \mathtt{P}_{T_{\!\mu, h}} \mu'$ 
weakly in $\cM_f(T)$,  we get  $\limsup_{n\to \infty} 
\dPro(\mathtt{P}_{T_{\! \mu_n, h_n}}\mu_n,  \mathtt{P}_{T_{\! \mu, h}}\mu)$  $\leqo 2h'\mu(T^*)/h$ for all $h'\ino (0, h)$, which implies the desired result when we let $h'$ go to $0$. \cqfd 
\begin{theorem}
\label{merafixpolish} Let $(T, d, \rho)$ be a rooted Polish $\bbR$-tree whose bordification is $(T^*\! , d^*\! , \rho)$. 
Let $(h_p)_{p\in \bbN}$ be $\bbR_+^*$-valued and strictly decreasing to $0$. We recall $\cM_f^{\mathtt{era}}(T)$ 
from (\ref{meraspace}) and $\delta_{\mathtt{era}}$ from (\ref{dmerafix}). 
Then the following holds true. 
\begin{compactenum}

\smallskip

\item[$(i)$] $\delta_{\mathtt{era}}$ is a distance on $\cM_f^{\mathtt{era}}(T)$ which metrizes the convergence in the sense of mass erasure.

\smallskip

\item[$(ii)$] $\big( \cM_f^{\mathtt{era}}(T), \delta_{\mathtt{era}}\big)$ is Polish and $\cM_f(T)$ is $\delta_{\mathtt{era}}$-dense in $\cM_f^{\mathtt{era}}(T)$. Namely $\cM_f^{\mathtt{era}}(T)$ is the $\delta_{\mathtt{era}}$-completion of the space $(\cM_f(T), \delta_{\mathtt{era}})$. 
\end{compactenum}
\end{theorem}

\noi
\textbf{Proof.} Let us prove first that $\delta_{\mathtt{era}}$ is a distance: it is clearly nonnegative symmetric and it satisfies the triangle inequality. Let us suppose that $\delta_{\mathtt{era}} (\mu, \nu)\eqo 0$: then $\era_{h_p} \mu\eqo \era_{h_p} \nu$ for all $p\ino \bbN$. By Proposition \ref{weakcvera} in $T^*$, as $p\! \to \! \infty$, we get 
$\era_{h_p} \mu\! \to \! \mu$ and  $\era_{h_p} \nu\! \to \! \nu$ weakly in $\cM_f(T^*)$, which implies that $\mu\eqo \nu$. 

  Let $\mu, \mu_n \ino \cM_f^{\mathtt{era}}(T)$, $n\ino \bbN$. If $\mu_n\! \to \! \mu$ in the sense of mass erasure, as in Definition \ref{meracvfixdef}, then we get $\delta_{\mathtt{era}}(\mu_n, \mu)\! \to \! 0$. Conversely, let us assume that $\delta_{\mathtt{era}}(\mu_n, \mu)\! \to \! 0$. This implies that for all $p\ino \bbN$, $\era_{h_p} \mu_n \! \to \! \era_{h_p}\mu$ weakly in $\cM_f(T)$. Let $h\ino \bbR^*_+$. Then there is $p\ino \bbN$ such that $h_p\leko h$ and we get $\era_h \mu_n \eqo \era_{h-h_p} (\era_{h_p} \mu_n)\! \to \! \era_{h-h_p} (\era_{h_p} \mu)\eqo \era_h\mu$ weakly in $\cM_f(T)$ by the semigroup property in Proposition \ref{properased} ($v$) and by Proposition \ref{weakcvera} in $T$. This completes the proof of ($i$). 
  
  To prove $(ii)$, we first show that $\cM_f(T)$ is 
$\delta_{\mathtt{era}}$-dense in $\cM_f^{\mathtt{era}}(T)$. To this end, we fix $\mu \ino \cM_f^{\mathtt{era}}(T)$ and we observe that $\dPro^T(\era_{h_p} (\era_h \mu) , \era_{h_p}\mu)\eqo\dPro^T(\era_{h_p+h} \mu, \era_{h_p}\mu)\! \to \! 0$ as $h\! \to \! 0^+$ by the semigroup property of mass erasure in Proposition \ref{properased} ($v$) and by Proposition \ref{weakcvera} in $T$. Thus $\lim_{h\to 0^+}\delta_{\mathtt{era}} (\era_h \mu, \mu) \eqo 0$, which proves the desired result since $\era_h\mu\ino \cM_f(T)$ for all $h\ino \bbR^*_+$. 

Then $( \cM_f^{\mathtt{era}}(T), \delta_{\mathtt{era}})$ is separable. \emph{Indeed}, let $(\mu_n)_{n\in \bbN}$ be $\cM_f(T)$-valued and $\dPro^T$-dense; by Lemma \ref{weakvsmera1}, $(\mu_n)_{n\in \bbN}$ is 
$\delta_{\mathtt{era}}$-dense in $\cM_f(T)$ and thus in $\cM_f^{\mathtt{era}}(T)$ by the previous argument.

Let us show that $( \cM_f^{\mathtt{era}}(T), \delta_{\mathtt{era}})$ is complete. 
Let $(\mu_n)_{n\in \bbN}$ be a $\delta_{\mathtt{era}}$-Cauchy 
sequence. Then for all $p\ino \bbN$, $(\era_{h_p} \mu_n)_{n\in \bbN}$ is a $\dPro^{T}$-Cauchy 
sequence of $\cM_f(T)$. Since $(\cM_f(T), \dPro^{T})$ is complete, 
there is $\nu_p\ino \cM_f(T)$ such that $\era_{h_p} \mu_n\! \to \! \nu_p$ weakly in $\cM_f(T)$. 
By the semigroup property in Proposition \ref{properased} ($v$) and 
by Proposition \ref{weakcvera} in $T$, we get 
$\era_{h_p} \mu_n \eqo \era_{h_{p}-h_{p+1}} (\era_{h_{p+1}} \mu_n) \! \to \! \era_{h_p-h_{p+1}} \nu_{p+1}$ 
weakly in $\cM_f(T)$, which implies that $ \era_{h_p-h_{p+1}} \nu_{p+1}\eqo \nu_p$ for all $p\ino \bbN$. 
Then Lemma \ref{descentera} applies and there exists $\mu\ino  \cM_f^{\mathtt{era}}(T)$ such 
that $\era_{h_p} \mu\eqo \nu_p$. Namely,  for all $p\ino \bbN$, 
$\era_{h_p} \mu_n\! \to \! \era_{h_p} \mu$ weakly in $\cM_f(T)$, which implies that $\delta_{\mathtt{era}}(\mu_n, \mu)\! \to \! 0$.  \cqfd 

\medskip

We conclude this section by further discussing the connection between weak convergence in $\cM_f(T)$ and convergence in the sense of mass erasure. In this regard, one key feature of a given measure $\mu\ino \cM_f^{\mathtt{era}}(T)$ to look at, is the 
push-forward measure of $\mu$ via the function $\bs\ino T^*\!\!  \mapsto \! ( d(\rho, \bs), (d(\rho, f_{T_{\! \mu, h_p}}\! (\bs)))_{p\in \bbN})$, where $\bhh\eqo (h_p)_{p\in \bbN}$ is $\bbR_+^*$-valued and strictly decreasing to $0$ and where we recall that $ f_{T_{\! \mu, h_p}}$ are the projections onto the compact subtrees $T_{\! \mu, h_p}$ (we note that $d(\rho, \bs)\eqo \infty$ if $\bs\ino \partial T$). 
Before defining more precisely 
the `joint law' of theses distances from the root, we need to introduce some notation. 

In what follows, $[0,\infty]$ is equipped with the standard compact topology, and
$[0, \infty] \times \! \bbR_+^\bbN$ and $\bbR_+^\bbN$ with their respective product topologies, which makes them Polish spaces. We then set 
\begin{equation}
\label{ellspacedef}
\fell_\bullet^{\uparrow}(\bbN) = [0, \infty] \! \times \! \fell^{\uparrow} (\bbN) \quad \textrm{where} \quad  \fell^{\uparrow}(\bbN)= \big\{ (x_p)_{p\in \bbN}\ino \bbR_+^\bbN: x_p \leqo x_{p+1}, \; p\ino \bbN \big\}. 
\end{equation}
Since $ \fell^{\uparrow}(\bbN)$ is a closed subset of $\bbR_+^\bbN$, $\fell_\bullet^{\uparrow}(\bbN)$ and $\fell^{\uparrow}(\bbN)$ equipped with their respective product topologies are Polish spaces. We next denote canonical projections by $L_{p}$, $p\ino \bbN\! \cup \! \{ \infty\}$. Namely, 
\begin{equation}
\label{ellcoordef}
\forall \, \bxx\eqo \big( x_\infty, (x_p)_{p\in \bbN} \big)\ino  \fell_\bullet^{\uparrow}(\bbN) , \quad 
 L_\infty (\bxx)\! :=\!  x_\infty \quad \textrm{and} \quad  L_p (\bxx) \! :=\! x_p,  \; p\ino \bbN .
 \end{equation}
We shall also use the following conventions. Let $\Lambda \ino \cM_{\! f} (\fell_\bullet^{\uparrow}(\bbN))$. Then 
\begin{compactenum}

\smallskip

\item[$\bullet$] \emph{We denote by $\Lambda^{\! o}\ino \cM_{\! f}(\fell^{\uparrow}(\bbN))$ the push-forward measure 
of $\Lambda$ via $(L_p)_{p\in \bbN}$.}

\smallskip

\item[$\bullet$] \emph{For all $p\ino \bbN$, we denote by $\Lambda^{\! (p)}$ the measure on $[0, \infty] \! \times \! \bbR_+^{p+1}$ that is the push-forward of $\Lambda$ via $(L_\infty, (L_q)_{0\leq q \leq p})$. }

\end{compactenum}

\begin{remark}
\label{tightLambda} Let $(\Lambda_n)_{n\in \bbN}$ be an $ \cM_{\! f}(\fell^\uparrow_\bullet (\bbN))$-valued sequence. 
Since $[0, \infty]$ is compact, tightness of $(\Lambda_n^o)_{n\in \bbN}$ in $\fell^\uparrow (\bbN)$ entails tightness of $(\Lambda_n)_{n\in \bbN}$ in 
$\fell^\uparrow_\bullet (\bbN)$. In particular, if $\mathscr Q\!\subseteq\! \cM_{\! f}(\fell^\uparrow (\bbN))$ is weakly compact, then 
$\{ \Lambda \ino \cM_f(\fell^\uparrow_\bullet (\bbN)): \Lambda^{\!  o}\ino \mathscr Q\}$ is weakly compact in $\cM_{\! f}(\fell^\uparrow_\bullet (\bbN))$. \cq 
\end{remark}

\begin{definition}
\label{Lambdadef} Let $\bhh\eqo (h_p)_{p\in \bbN}$ be $\bbR^*_+$-valued and strictly decreasing to $0$. Let 
$(T, d, \rho)$ be a rooted Polish $\bbR$-tree with bordification  
$(T^*\! , d^*\! , \rho)$. Let $\mu \ino  \cM_{\! f}^{\mathtt{era}}(T)$. 

\begin{compactenum}

\smallskip

\item[$(a)$]  The \emph{$\bhh$-projected height measure} $\Lambda_{\bhh, \mu}\ino \cM_f(\fell_\bullet^{\uparrow}(\bbN))$ is defined as follows. 
\begin{compactenum}

\smallskip

\item[$-$] If $\mu(T^*)\eqo 0$, then $ \Lambda_{\bhh, \mu}$ is the null measure.

\smallskip

\item[$-$] Otherwise, let $X$ be a $T^*$-valued r.v.~with law $\mu (\cdot)/ \mu(T^*)$. We define $\Lambda_{\bhh, \mu}$ as $ \mu(T^*) \overline{\Lambda}_{\bhh, \mu}$, where $ \overline{\Lambda}_{\bhh, \mu}$ is the law of 
$$ \Big( d(\rho, X)\, , \big( d\big(\rho, f_{T_{\! \mu, h_p}}\!  (X)\big)\big)_{p\in \bbN} \Big).$$
\end{compactenum}
\noi
If $\mu(T^*)\geko 0$, we set $\Lambda_{\bhh, \mu}^{\! o} \eqo \mu(T^*) \overline{\Lambda}^o_{\bhh, \mu}$, where $\overline{\Lambda}^{o}_{\bhh, \mu}$ is the law of $\big( d(\rho, f_{T_{\! \mu, h_p}} (X))\big)_{p\in \bbN} $.

\smallskip

\item[$(b)$] The \emph{height measure} of $\mu$, which is denoted by $\lambda_\mu\ino \cM_f([0, \infty])$, 
is the push-forward measure of $\mu$ via the height function 
$\bs\ino  T^* \! \mapsto \! d(\rho, \bs)$. Namely, if $\mu(T^*)\eqo 0$, then $\lambda_\mu$ is the zero measure. Otherwise, 
$\lambda_\mu (\cdot)/ \mu(T^*)$ is the law of $d(\rho, X)$. 
Note that $\lambda_\mu$ is the push-forward measure of $L_\infty$ under $\Lambda_{\bhh, \mu}$ and that 
$\lambda_\mu\ino \cM_f(\bbR_+)$ iff $\mu\ino \cM_f(T)$. More explicitly, $\lambda_\mu ([0, r])\eqo \mu (B_{T,d} (\rho, r))$, $r\ino \bbR_+$, and $\lambda_\mu ([0, \infty])\eqo \mu(T^*)$.   \cq 
\end{compactenum}
\end{definition}
\begin{proposition}
\label{Lambdaprop} We keep the notation of Definition \ref{Lambdadef}. 
Let $\mu, \mu_n \ino  \cM_f^{\mathtt{era}}(T)$, $n\ino \bbN$. Then the following holds true. 
\begin{compactenum}

\smallskip

\item[$(i)$] $\Lambda_{\bhh, \mu}$-a.e.~$L_\infty\eqo \sup_{p\in \bbN} L_p$.

\smallskip

\item[$(ii)$] For all integers $p\geqo  q$ and all 
$r\ino \bbR^*_+$, $\dPro^{T}\! (  \ttP_{T_{\! \mu, h_p}} \mu, \ttP_{T_{\! \mu, h_q}} \mu ) \leqo r \! \vee \! 
\Lambda_{\bhh, \mu} (\{ L_p\! -\! L_q \geko r \}) $. Furthermore if $\mu \ino \cM_f(T)$, then 
$\dPro^{T}( \mu, \ttP_{T_{\! \mu, h_p}} \mu) \leqo r \! \vee\!  \Lambda_{\bhh, \mu} (\{ L_\infty\! -\! L_p \geko r \}) $.

\smallskip

\item[$(iii)$] Let $p\ino \bbN$. Let $\rho \ino  \ttt\!\subseteq\! T$ be a compact subtree such that $h_p\geko \mu (T^*\backslash \ttt)$. We set $\nu \eqo \mathtt{P}_{\!\ttt} \mu $. Then, for all Borel subsets $B$ of $[0, \infty] \! \times \!  \bbR_+^{{p+1}}\! $, we get 
$|\Lambda^{{\! (p)}}_{ \bhh, \nu} (B)\! -\! \Lambda^{{\! (p)}}_{\bhh, \mu} (B)| \leqo \mu (T^*\backslash \ttt)$.    

\smallskip

\item[$(iv)$] We suppose that $\mu_n \! \! \to \! \mu$ in the sense of mass erasure and we suppose that $h\! \mapsto \! T_{\! \mu, h}$ is $\dHaus$-continuous at each $h_p$, $p\ino \bbN$. Then the following holds.
\begin{compactenum}

\smallskip

\item[$(iv\textrm{-}a)$] $\Lambda^{\! o}_{\bhh, \mu_n} \!\! \!  \! \to \! \Lambda^{\! o}_{\bhh, \mu}$ weakly on $\cM_{\! f}(\fell^{\uparrow}(\bbN))$. 

\smallskip

\item[$(iv\textrm{-}b)$] The measures  $(\Lambda_{\bhh, \mu_n})_{n\in \bbN}$ are tight in  $\cM_{\! f}(\fell_\bullet^{\uparrow}(\bbN))$ and if $\Lambda$ is a weak limit point, then $\Lambda^{\! o}\eqo  \Lambda^{\! o}_{\bhh, \mu}$ and $\Lambda$-a.e.~$\sup_{p\in \bbN} L_p \leqo L_\infty$.   
\end{compactenum}

\smallskip

\item[$(v)$]Let $\mu_n, \mu\ino \cM_{\! f}(T)$ are such that $\mu_n \! \to \! \mu$ weakly in $\cM_f(T)$. We suppose that $h\! \mapsto \! T_{\! \mu, h}$ is $\dHaus$-continuous at each $h_p$, $p\ino \bbN$. Then, $\Lambda_{\bhh, \mu_n} \! \!\! \to \! \Lambda_{\bhh, \mu}$ weakly on $\cM_{\! f}(\fell_\bullet^{\uparrow}(\bbN))$. 
\end{compactenum}
\end{proposition}
\noi
\textbf{Proof.} If $\mu(T^*)\eqo 0$, then the whole proposition holds true trivially. So we assume that $\mu(T^*) \geko 0$ and w.l.o.g.~we also assume $\mu_n (T^*) \geko 0$ for all $n\ino \bbN$. We denote by $X$ a $T^*$-valued 
r.v.~with law $\overline{\mu} \! :=\! \mu (\cdot) / \mu(T^*)$ and by $X_n$ a $T^*$-valued r.v.~with law $\overline{\mu}_n \! :=\! \mu_n (\cdot) / \mu_n(T^*)$.

  Let us prove $(i)$. By definition of projection $d^*(X, f_{T_{\!\mu, h_p}} (X))\eqo \inf \{ d^*( X, y) ; y\ino T_{\! \mu, h_p}\}$ and by the definition of $X$, we also see that a.s.~$X\ino \mathtt{Supp}_* \mu$.  
Proposition \ref{massTprop} $(vi)$ in $T^*$ entails a.s.~that $\lim_{p\to \infty} d^*(X, f_{T_{\! \mu, h_p}} (X))\eqo 0$ and thus $d(\rho, X)\eqo \sup_{p\in \bbN}   d(\rho, f_{T_{\! \mu, h_p}} (X))$, which entails $(i)$. 

To prove $(ii)$, we set $\nu\eqo \mu(T^*) \overline{\nu}$ where $ \overline{\nu}$ is the law on $T^2$ of $(f_{T_{\! \mu, h_p}}(X), f_{T_{\! \mu, h_q}}(X))$, and we set $\nu_r\eqo \nu (\cdot \cap \Delta_r)$ where $\Delta_r\eqo \{ (\sigma, \sigma')\ino T^2\! : \! d(\sigma, \sigma') \leqo r \}$. For all $B \ino \mathscr B(T)$, we easily get 
$0\leqo \ttP_{T_{\! \mu, h_p}} \mu (B)\! -\! \nu_r (B\! \times \! T)\leqo \nu (T^2\backslash \Delta_r)$ and 
$0\leqo \ttP_{T_{\! \mu, h_q}} \mu (B)\! -\! \nu_r (T\! \times \! B)\leqo \nu (T^2\backslash \Delta_r)$. Strassen's theorem 
implies that 
$\dPro^T\! \big(  \ttP_{T_{\! \mu, h_p}} \mu, \ttP_{T_{\! \mu, h_q}} \mu \big) \leqo r \! \vee \!  \nu (T^2\backslash \Delta_r)$. 
We next observe that 
$d(f_{T_{\! \mu, h_p}}(X), f_{T_{\! \mu, h_q}}(X))$ $\eqo d(\rho, f_{T_{\! \mu, h_p}}(X))\! -\! d(\rho, f_{T_{\! \mu, h_q}}(X))$. This implies that $\nu (T^2\backslash \Delta_r)\eqo \Lambda^o_{\bhh, \mu} (\{ L_p\! -\! L_q \geko r\})$, which completes the proof of the first statement of $(ii)$. We argue similarly to get the second one. 

Let us prove $(iii)$. Set $\nu\eqo \mathtt{P}_{\!\ttt} \mu$. As 
$h_0\geko \cdots \geko  h_p \geko \mu (T^*\backslash \ttt)$, Proposition \ref{massTprop} $(i)$ (applied in $T^*$) implies that $\ttt_q\! :=\!  T_{\! \nu, h_q}\eqo T_{\! \mu, h_q}$ for all $q\ino \{ 0, \ldots, p\}$. We set $m\eqo \mu(T^*)$, which is also equal to $\nu(\ttt)$, and we set $\overline{\nu} \eqo \nu(\cdot) / m$, which is therefore 
a probability measure. We recall that $X$ has law $\overline{\mu}$ and we set $X'\! :=\! f_\ttt (X)$, which 
has law $\overline{\nu}$. We  note that $f_{\ttt_q} (X)\eqo f_{\ttt_q} (X')$. Therefore, if 
$X\ino \ttt$, and thus $X'\eqo X$, we find 
$V\! :=\! (d(\rho, X), (d(\rho, f_{\ttt_q} (X)))_{0\leq q \leq p}) $  $\eqo $ $(d(\rho, X'), (d(\rho, f_{\ttt_q} (X')))_{0\leq q \leq p})\! =:\! V' $. Since $m^{-1} \Lambda^{{\! (p)}}_{{\bhh,\mu}}$ is the law of $V$ and 
$m^{-1} \Lambda^{{\! (p)}}_{{\bhh,\nu}}$ that of $V'$, the coupling $(V, V')$ implies for any Borel subset $B$ of $[0, \infty] \! \times \! \bbR_+^{{p+1}}$ that 
$$\big| \Lambda^{\!(p)}_{\bhh, \nu} (B)\! -\! \Lambda^{\!(p)}_{\bhh, \mu} (B)\big| \leq m \big| \bP (V'\! \ino B) \! -\! \bP (V\ino B)\big|
\leq m\bP (V'\! \neq \! V) \leqo m\bP (X\! \notin \! \ttt) \eqo  \mu(T^*\backslash \ttt), $$
which is the desired result.

We next prove $(iv\textrm{-}a)$. 
To simplify notation, for all $p\ino \bbN$, we set 
$\ttt_{n,p} \! := \!T_{\! \mu_n, h_p}$, $Y_{n,p}\! := \! f_{\ttt_{n,p}} (X_n)$ and $Y_p\! := \! f_{\ttt_p} (X)$, where we recall from above that $\ttt_p\eqo T_{\! \mu, h_p}$.  
Lemma \ref{weakvsmera1} $(ii)$ implies that $Y_{n,p} \! \to \! Y_p$ in distribution in $T$ and that $\dHaus (\ttt_{n,p}, \ttt_p)\! \to \! 0$ for all $p\ino \bbN$.
Let $p,q\ino \bbN$ such that $p\geqo q$. By Proposition \ref{projdef} $(iv)$, 
$Y_q\eqo f_{\ttt_q}(Y_p)$. Since the projections are continuous (by Proposition \ref{projdef} $(iii)$) 
the convergence 
$(f_{\ttt_q} (Y_{n,p}))_{0\leq q\leq p} \! \to \! (Y_q)_{0\leq q\leq p}$ holds in distribution in $T^{p+1}$. 
By Proposition \ref{projdef} $(iv)$, we also get $Y_{n,q}\eqo f_{\ttt_{n,q}} (Y_{n,p})$. We then 
apply Lemma \ref{GrC0proj} to $T\eqo T'$, $\tau\eqo \ttt_q$, $\tau'\eqo \ttt_{n,q}$, $\sigma'\eqo \sigma\eqo Y_{n,p}$, we get $d(Y_{n,q}, f_{\ttt_q} (Y_{n,p}))\leqo 3\dHaus (\ttt_{n,q}, \ttt_q)$. This implies that 
$(Y_{n,q})_{0\leq q\leq p} \! \to \! (Y_q)_{0\leq q\leq p}$ in distribution in $T^{p+1}$. Since this holds true for all $p\ino \bbN$, we get $(Y_{n,p})_{p\in \bbN}\! \to \! (Y_{p})_{p\in \bbN}$ in distribution in $T^\bbN$ equipped with the product topology. Therefore, $((d(\rho, Y_{n,p}))_{p\in \bbN}\! \to \! d(\rho, Y_{p}))_{p\in \bbN}$ in distribution in $\fell^{\uparrow}(\bbN)$, which easily entails $(iv\textrm{-}a)$. 

  Before proving $(iv\textrm{-}b)$, we prove $(v)$: if $\mu_n\! \to \! \mu$ in 
$\cM_f(T)$, then the previous argument entails that $(X_n, (Y_{n,p})_{p\in \bbN})\! \to \! (X, (Y_{p})_{p\in \bbN})$ in distribution in $T\! \times \! T^\bbN$, equipped with the product topology. Therefore, we get $(d(\rho, X_n), (d(\rho, Y_{n,p}))_{p\in \bbN})\! \to \! (d(\rho, X), (d(\rho, Y_{p}))_{p\in \bbN})$ in distribution in $\fell_\bullet^{\uparrow}(\bbN)$ and we easily get $(v)$. 

It remains to prove $(iv\textrm{-}b)$. Tightness follows from $(iv\textrm{-}a)$ and Remark \ref{tightLambda}. 
Let $\Lambda$ be a weak limit point of the $\Lambda_{\bhh, \mu_n}$. By $(iv\textrm{-}a)$ we get 
$\Lambda^{\! o}\eqo  \Lambda^{\! o}_{\bhh, \mu}$. 
Then observe that the open set 
$U \! :=\! \{ L_\infty \leko \sup_{p\in \bbN} L_p \}$ is such that $\Lambda_{\bhh, \mu_n} (U)\eqo 0$,
for all $k\ino \bbN$. Thus, the Portmanteau theorem implies that $\Lambda (U)\eqo 0$, which completes the proof of $(iv\textrm{-}b)$. \cqfd 

\smallskip

\begin{theorem}
\label{eravsweakfix} Let $(T, d, \rho)$ be a rooted Polish $\bbR$-tree, $\mu\ino \cM_f^{\mathtt{era}}(T)$, and  
$(h_p)_{p\in \bbN}$ an $\bbR_+^*$-valued sequence which strictly decreases to $0$ and is such that 
$h\! \mapsto \! T_{\! \mu, h}$ is $\dHaus$-continuous at each $h_p$, $p\ino \bbN$ (there is always one). Keeping the notation of Definition \ref{Lambdadef},  
let $(\mu_n)_{n\in \bbN}$ be $\cM_f(T)$-valued and such that $\mu_n \! \to \mu$ in the sense of mass erasure. Then the following assertions are equivalent. 
\begin{compactenum}

\smallskip

\item[$(a)$]  $\mu\ino \cM_f(T)$ and  $ \mu_n \! \to \! \mu$ weakly in $\cM_f(T)$

\smallskip

\item[$(b)$] $\lambda_\mu\ino \cM_f(\bbR_+)$ and $\lambda_{\mu_n} \!\!\!  \to \! \lambda_{\mu} $ weakly in $\cM_f(\bbR_+)$.

\smallskip

\item[$(c)$] For all $\epp \ino \bbR^*_+$, $\lim_{p\to \infty} \sup_{n\in \bbN} \Lambda_{\bhh, \mu_n} (\{ L_\infty \! -\! L_p \geqo \epp \}) \eqo 0$. 

\end{compactenum}
\end{theorem}
\noi
\textbf{Proof.} We recall notation $L_p$, $p\ino \bbN\cup \{ \infty\}$ from (\ref{ellcoordef}) and we set 
$\mathcal L\! :=\! \sup_{p\in \bbN} L_p$, to simplify notation. Let us  first assume $(a)$ 
and prove $(a)\! \Rightarrow \! (b)$: clearly $\lambda_\mu\ino \cM_f(\bbR_+)$ and by Proposition \ref{Lambdaprop} $(v)$, $\Lambda_{\bhh, \mu_n} \! \!\! \to \! \Lambda_{\bhh, \mu}$ weakly on $\cM_{\! f}(\fell_\bullet^{\uparrow}(\bbN))$, which implies $(b)$ since 
$\lambda_{\mu_n} \eqo \Lambda_{{\bhh, \mu_n}}\circ L^{-1}_\infty$,  
$\lambda_\mu \eqo \Lambda_{{\bhh, \mu}}\circ L^{-1}_\infty$ and since $L_\infty$ is by definition continuous on $\cM_{\! f}(\fell_\bullet^{\uparrow}(\bbN))$.

We next assume $(b)$ and prove $(b)\Rightarrow (c)$. To this end, we first show that 
$\Lambda_{{\bhh, \mu_n}} \!\! \to \! \Lambda_{{\bhh, \mu}}$ weakly on $\cM_f(\fell_\bullet^{\uparrow}(\bbN))$. 
\emph{Indeed}, by Proposition \ref{Lambdaprop} $(iv\textrm{-}b)$, 
$(\Lambda_{{\bhh, \mu_n}})_{n\in \bbN}$ are tight in $\cM_f(\fell_\bullet^{\uparrow}(\bbN))$ and any weak limit point 
$\Lambda$ must satisfy $\Lambda^{\! o}\eqo  \Lambda^{\! o}_{{\bhh, \mu}}$ and 
$\Lambda$-a.e.~$\mathcal L  \leqo L_\infty$. On one hand, since 
$\lambda_{{\mu_n}} \eqo  \Lambda_{{\bhh, \mu_n}}  \circ  L_{\infty}^{{\! -1}} $, we get 
$\lambda_{\mu} \eqo \Lambda \circ L_{\infty}^{{\! -1}}$. 
On the other hand, $\mu\ino \cM_f(T)$  implies $\Lambda_{{\bhh, \mu}}$-a.e.~$L_\infty  
\eqo \mathcal L$ and 
thus $\Lambda^{\! o}_{{\bhh, \mu}}  \circ  \mathcal L^{-1}\eqo \lambda_\mu$. 
Since $\Lambda^{\! o}\eqo  \Lambda^{\! o}_{{\bhh, \mu}}$, we get $\Lambda  \circ  \mathcal L^{{-1}} \eqo \Lambda^{\! o}  \circ  \mathcal L^{{-1}}\eqo \lambda_\mu$. Then Lemma \ref{adhoclem} (for a deterministic $\Lambda$) 
implies that $\Lambda$-a.e.~$ 
L_\infty\eqo \mathcal L $, which implies $\Lambda\eqo  \Lambda_{{\bhh, \mu}}$ since 
$\Lambda^{\! o}\eqo  \Lambda^{\! o}_{{\bhh, \mu}}$.
Therefore, $\Lambda_{{\bhh, \mu}}$ is the only weak limit point of $(\Lambda_{{\bhh, \mu_n}})_{n\in \bbN}$, 
which implies the desired convergence.  

Let us complete the proof of $(c)$. 
We fix $\epp \ino \bbR_+^*$. Since $\{ L_\infty\! -\! L_p \geqo \epp \}$ is a closed subset of $\fell_\bullet^{\uparrow}(\bbN)$, the Portmanteau theorem implies $\limsup_{n\to \infty} \Lambda_{{\bhh, \mu_n}} (\{ L_\infty\! -\! L_p \geqo \epp \})\leqo 
\Lambda_{{\bhh, \mu}} (\{ L_\infty\! -\! L_p \geqo \epp  \})$, for all $p\ino \bbN$.
Since $\Lambda_{{\bhh, \mu}}$-a.e.~$\mathcal L\eqo L_\infty$, $\lim_{p\to \infty} \downarrow \Lambda_{{\bhh, \mu}}  (\{ L_\infty\! -\! L_p \geqo \epp  \})\eqo 0$, which entails $(c)$. 

Let us assume $(c)$ and prove $(c)\! \Rightarrow \! (a)$.
We first show that $\Lambda_{{\bhh, \mu_n}} \!\! \to \! \Lambda_{{\bhh, \mu}}$ weakly on $\cM_f(\fell_\bullet^{\uparrow}(\bbN))$. \emph{Indeed}, by Proposition \ref{Lambdaprop} $(iv\textrm{-}b)$, 
$(\Lambda_{{\bhh, \mu_n}})_{n\in \bbN}$ are tight in $\cM_f(\fell_\bullet^{\uparrow}(\bbN))$ and any weak limit point 
$\Lambda$ must satisfy $\Lambda^{\! o}\eqo  \Lambda^{\! o}_{{\bhh, \mu}}$ and 
$\Lambda$-a.e.~$\mathcal L  \leqo L_\infty$. Since $\{ L_\infty\! -\! L_p \geko \epp \}$ is an open 
subset of $\fell_\bullet^{\uparrow}(\bbN)$, the Portmanteau theorem implies 
$\Lambda (\{ L_\infty \! -\! L_p \geko \epp\}) \leqo \liminf_{n\to \infty} \Lambda_{{\bhh, \mu_n}} (\{ L_\infty\! -\! L_p \geko \epp \})\xrightarrow[p\to \infty]{\; } 0 $
by $(c)$. This implies that $\Lambda$-a.e.~$\mathcal L\eqo L_\infty \leko \infty$, since $\Lambda (\{ L_\infty \eqo \infty\}) \leqo \Lambda (\{ L_\infty \! -\! L_p \geko \epp\})$ for all $\epp\ino \bbR_+^*$ and all $p\ino \bbN$. Arguing as previously, we find  $\Lambda\eqo \Lambda_{{\bhh, \mu}}$. Thus, $\Lambda_{{\bhh, \mu}}$ is the only weak limit point of $(\Lambda_{{\bhh, \mu_n}})_{n\in \bbN}$, which proves the desired result. 

Let us complete the proof of $(a)$. First, we note that $\Lambda_{{\bhh, \mu}}$-a.e.~$L_\infty\leko \infty$ implies $\lambda_\mu\ino \cM_f(T)$ and thus $\mu\ino \cM_f(T)$. 
Then, by Proposition \ref{Lambdaprop} $(ii)$, for all 
$p, n\ino \bbN$ and $\epp\ino \bbR^*_+$, we get 
$$  \dPro^T( \mu_n, \mu) \leqo \epp \! \vee \! \Lambda_{\bhh, \mu_n} (\{ L_\infty\! -\! L_p \geqo \epp  \}) +  
\dPro^T( \ttP_{T_{\! \mu_n, h_p}} \mu_n , \ttP_{T_{\! \mu, h_p}} \mu) +  \epp  \! \vee \! \Lambda_{\bhh, \mu} 
(\{ L_\infty\! -\! L_p \geqo \epp  \}). $$
Lemma \ref{weakvsmera1} $(ii)$ shows that 
$\lim_{n\to \infty}\dPro^{T}( \ttP_{T_{\! \mu_n, h_p}} \mu_n , \ttP_{T_{\! \mu, h_p}} \mu)\eqo  0$. 
By the Portmanteau theorem we get 
$  \limsup_{n\to \infty}  \dPro^{T}( \mu_n, \mu) \leq 2 (  \epp   \vee  \Lambda_{{\bhh, \mu}} (\{ L_\infty\! -\! L_p \geqo \epp  \}))$, which entails $(a)$ because $\Lambda_{{\bhh, \mu}} (\{ L_\infty\! -\! L_p \geqo \epp  \}))\! \to 0$ as $p\! \to \infty$, since $\mu\ino \cM_f(T)$ and since $\epp$ can be arbitrarily small. \cqfd

\section{Convergence of $\bbR$-trees and mass erasure}
\label{CVspacemerasec}
In this section we study mass erasure for sequences of rooted finitely measured Polish $\bbR$-trees that converge in the Gromov--Prokhorov sense. 
To this end, we recall below basic definitions and known results on the Gromov-convergences of metric spaces. 
\begin{definition}
\label{GPpointdef} For $i\ino \{ 1, 2\}$, let $(E_i, d_i, \rho_i, \mu_i )$ be Polish metric spaces where 
$\rho_i \! \in \! E_i$ and $\mu_i \! \in \! \cM_{\! f} (E_i)$. We define their 
\textit{Gromov--Prokhorov pseudo-distance} by 
\begin{equation}
\label{GP}
\dGP (\mu_1, \mu_2) = \inf \big\{ \delta ( \phi_1(\rho_1), \phi_2 (\rho_2)) \vee \dePro(  
\mu_1 \! \circ \! \phi_1^{-1}\! ,  \mu_2\!  \circ \! \phi_2^{-1} )  \big\} , 
\end{equation}
and if the $E_i$ are compact, their \textit{Gromov--Hausdorff and  Gromov--Hausdorf--Prokhorov distances} by
\begin{eqnarray*}
 \dGH (E_1 , E_2) \!\!\! \!\! &=&\!\!\! \!\! \inf \! \big\{ \delta \big( \phi_1 (\rho_1), \phi_2 (\rho_2) \big) \vee \deHaus (
 \phi_1 (E_1), \phi_2(E_2) )  \big\} \quad \textrm{and} 
 \\
\dGHP\!  \big(\!(E_1, \mu_1), \! (E_2 , \mu_2)\!\big) 
   \!\!\!\! \!&=&   \!\!\!\!  \!   \inf \! \big\{\!  \delta \big( \phi_1 (\rho_1), \phi_2 (\rho_2) \big) \! \vee\!  \deHaus (
 \phi_1 (E_1), \phi_2(E_2) ) \! \vee \!  \dePro(  
\mu_1 \! \circ \! \phi_1^{-1}\! ,  \mu_2\!   \circ \! \phi_2^{-1} )\! \}
\end{eqnarray*}
where infima are taken on all Polish metric spaces $(X, \delta)$ and all isometrical embeddings $\phi_i \! : \! 
E_i \! \rightarrow \! X$, $i\ino \{ 1,2\}$; 
here, $ \deHaus$ stands for the Hausdorff distance on the space of compact subsets of $X$ and $ \dePro$ stands for the Prokhorov distance on $\cM_{\! f} (X)$. \cq 
\end{definition}

Let us recall fundamental results on $\dGH, \dGP$ and $\dGHP$. First note that these distances are invariant under bijective isometries that preserve the distinguished point and the measure, as appropriate.
\begin{compactenum}

\smallskip

\item[]$\!\!\!\!\!\!\!\!\!\textbf{Grom}$-$(i)$: \emph{$\dGH (E, E^\prime) \! =\!  0$ iff there is a bijective isometry 
$\varphi \colon E\! \rightarrow\!  E^\prime$ such that $\varphi (\rho)\! = \! \rho^\prime$. 
If we denote by $\bM^0_c$ the set of isometry classes of pointed compact metric spaces, then 
$(\bM^0_c, \dGH)$ is Polish.}\smallskip
\end{compactenum}
See Gromov, Pansu and Semmes \cite[Sec.~$3.11\frac{1}{{2}}+$]{Gro}, or Burago, Burago and Ivanov \cite[Sec. 7.3--7.4]{BuBuIv}. 

\begin{compactenum}

\smallskip

\item[]$\!\!\!\!\!\!\!\!\!\textbf{Grom}$-$(ii)$: \emph{$\dGP (\mu, \mu^\prime) \! =\!  0$ iff there is a bijective isometry $\varphi \colon
\{ \rho\}  \cup  \mathtt{\supp} \mu \rightarrow \{ \rho^\prime \}  \cup  \mathtt{\supp} \mu^\prime$ such that $\varphi (\rho)\! = \! \rho^\prime$ and $\mu \! \circ  \! \varphi^{-1}\!\! = \! \mu^\prime$. 
Let us denote by $\bM$ the set of isometry classes of pointed finitely measured Polish metric spaces $(E,d,\rho, \mu)$ such that $E\! =\! \{ \rho\}  \cup  \mathtt{\supp} \mu$. Then, $(\bM,\dGP)$ is Polish.}\smallskip
\end{compactenum}
See L{ö}hr, Voisin and Winter \cite[Prop.~2.6]{LVW}. 

\begin{compactenum}

\smallskip

\item[]$\!\!\!\!\!\!\!\!\!\textbf{Grom}$-$(iii)$: \emph{$\dGHP ((E, \mu), (E^\prime\! , \mu^\prime)) \! =\!  0$  iff there is a bijective isometry $\varphi \colon 
E\rightarrow E^\prime$ such that $\varphi (\rho)\! = \! \rho^\prime$ and $\mu \! \circ  \! \varphi^{-1}\!\! = \! \mu^\prime$ (with no assumption on the support). 
If we denote by $\bM^{1}_c$ the space of 
isometry classes of pointed finitely measured compact spaces, 
then, $(\bM^{1}_c, \dGHP)$ is Polish.}\smallskip
\end{compactenum}
\noi
See Abraham, Delmas and  Hoscheit \cite[Thm.~2.5]{ADHoscheit}. 

\smallskip

\noi
\textbf{Notation.} We shall often write $(E,d,\rho, \mu) \! \equiv\! \fmu\ino \bM$ 
(resp.~$(E,d,\rho, \mu)\! \equiv\!  \fmu \ino \bM^1_c$ or $(E,d,\rho) \! \equiv \! \widetilde{E}\ino \bM^0_c$) to mean that 
$(E,d,\rho, \mu)$ is a representative of $\fmu$ (resp.~$\fmu$ or $\widetilde{E}$) and to simplify, we shall sometimes 
confound the space and its 
isometry class. \cq 

\smallskip

It is convenient to control $\dGH$, $\dGP$ and $\dGHP$ in a more intrisic way via metric correspondences and coupling. 
To this end we recall the following definition.

\begin{definition}
Let $(E_i, d_i)$, $i\ino \{ 1,2\}$, be two metric spaces and for all $(x_1,x_2)\ino E_1\! \times \! E_2$, we set $\pi_i(x_1,x_2)\eqo x_i$, $i\ino \{ 1,2\}$. 
\begin{compactenum}

\smallskip

\item[$(a)$] Any non-empty subset $R\! \subseteq\! E_1\!  \times \! E_2$ a is (partial) \emph{correspondence}. It is furthermore a 
\emph{correspondence between} non-empty $A_1\!\subseteq\! E_1$ 
and $A_2 \!\subseteq\! E_2$ if $A_i\!\subseteq\! \pi_i (R)$, $i\ino \{ 1,2\}$.

\smallskip

\item[$(b)$] The quantity $ \mathtt{dis}_{d_1, d_2} (R) \! :=\!  \sup \big\{  \big| d_1(x_1, y_1)\! -\! d_2 (x_2, y_2) \big|\, ; \, (x_i,y_i) \ino R, i\ino \{ 1,2\}\big\}  \ino [0, \infty]$ 
is the \emph{distortion} of the partial correspondence $R$. When there is no ambiguity, we simply denote $ \mathtt{dis}_{d_1, d_2} (R)$ by $ \mathtt{dis} (R)$. \cq 
\end{compactenum}
\end{definition}
\noi
To construct distances in terms of correspondences we rely on the two following lemmas. 
\begin{lemma}
\label{extendist} Let $R\! \subseteq\! E_1\!  \times \! E_2$ be a partial correspondence on the metric spaces $(E_i,d_i)$, $i\ino \{ 1,2\}$. We assume 
that $\mathtt{dis} (R) \! <\!  \infty$. 
We set 
$X\! =\!  E _1\! \sqcup \! E_2$, we fix $r\! \in \! \bbR^*_+$ be such that $r\! \geq \! \frac{1}{2} \mathtt{dis} (R)$ and we 
define $\delta\colon X\! \times \! X \rightarrow \bbR_+$ as follows: $\delta\eqo d_i$ on $E_i\! \times \! E_i$, $i\ino \{ 1,2\}$, and 
$$ \forall (x_1,x_2) \! \in \! E_1\! \times \! E_2 , \quad    \delta (x_1,x_2) \! = \! \delta (x_2, x_1)\!  =\!   \inf \big\{ d_1(x_1,y_1) \! +\!  r\!  +\!  
d_2(x_2,y_2))\, ; \, (y_1,y_2) \! \in \! R \big\} \;  .$$
Then, $\delta$ is a metric on $X$ such that $\min_{x_1\in E_1, x_2\in E_2} \delta (x_1, x_2)\eqo r$. 
If the $E_i$ are Polish (resp.~compact), so is $X$. 
If $R$ is a correspondence between $A_1$ and $A_2$, then $\deHaus (A_1, A_2)\! = \! r$. 
\end{lemma}
\noi
\textbf{Proof.} See e.g.~Greven, Pfaffelhuber and Winter \cite[Rem.~5.5]{GPW}. \cqfd 

\begin{lemma}
\label{seqspace} Let $(X, d_X)$ and $(X_n, d_{X_n})$, $n\! \in \! \bN$, be metric spaces. For all $n\! \in \! \bN$, let $d_{X, X_n}$ be a metric on $X\! \sqcup \! X_n$ that extends $d_X$ and $d_{X_n}$ and such that $r_n\! :=\! \min_{x\in X, u\in X_n} d_{X, X_n} (x, u) \! >\! 0$.  
Then, there exists a metric $\delta$ on 
$X\sqcup \bigsqcup_{n\in \bN} X_n $ that extends the $d_{X, X_n}$.   
\end{lemma}
\noi
\textbf{Proof.} See Greven, Pfaffelhuber and Winter \cite[Lem.~A.1]{GPW}. \cqfd 

\smallskip

\noi
We next  recall $\eta$-couplings from Definition \ref{couplingdef} and we introduce the following conditions. 
\begin{definition}
\label{condiGrom} Let $\fmu_i \! \equiv \! (E_i,d_i,\rho_i, \mu_i)$ be finitely measured pointed metric spaces, $i\ino \{ 1,2\}$ and let $\eta \ino \bbR^*_+$. 
We introduce the following conditions. 
\begin{compactenum}

\smallskip

\item[$\mathtt{Cond}_{\mathtt{GH}}(\eta)$:]There is a correspondence $R$ between $E_1$ and $E_2$ such that 
\begin{equation}
\label{GHcondi}
\tfrac{1}{2} \mathtt{dis} (R) \leko \eta \quad \textrm{and} \quad (\rho_1,\rho_2)\ino R .
\end{equation}
\item[$\mathtt{Cond}_{\mathtt{GP}}(\eta)$:]There is a \emph{compact} partial correspondence $R\!\subseteq\! E_1\! \times \! E_2$ which satisfies (\ref{GHcondi}) and there is a measure $\mathtt m\ino \cM_{\! f} (E_1\! \times \! E_2)$ such that 
\begin{equation}
\label{GPcondi}
\mathtt m ((E_1\! \times \! E_2) \backslash R) \eqo 0  \quad \textrm{and $\quad \mathtt m$ is an $\eta$-coupling of $\mu_1$ and $\mu_2$.}
\end{equation}
\item[$\mathtt{Cond}_{\mathtt{GHP}}(\eta)$:]
There is a correspondence $R$ \emph{between} $E_1$ and $E_2$ satisfying (\ref{GHcondi}) and there is a measure $\mathtt m\ino \cM_{\! f} (E_1\! \times \! E_2)$ satisfying (\ref{GPcondi}) \cq
\end{compactenum}
\end{definition}
\begin{lemma}
\label{rephraseGrom} We keep the notation of Definition \ref{condiGrom}. 
\begin{compactenum}

\smallskip

\item[$(i)$] Let $E_i$, $i\ino \{ 1,2\}$, be Polish. 
Then $\fdelta_{\mathtt{GP}} (\fmu_1, \fmu_2) \leko \eta$ iff $\mathtt{Cond}_{\mathtt{GP}}(\eta)$ holds true.
\smallskip

\item[$(ii)$] Let $E_i$, $i\ino \{ 1,2\}$, be compact. Then $\fdelta_{\mathtt{GH}} (E_1, E_2) \leko \eta$
iff $\mathtt{Cond}_{\mathtt{GH}}(\eta)$ holds true. 

\smallskip

\item[$(iii)$] Let $E_i$, $i\ino \{ 1,2\}$, be compact. 
Then $\fdelta_{\mathtt{GHP}} (\fmu_1, \fmu_2) \leko \eta$ iff $\mathtt{Cond}_{\mathtt{GHP}}(\eta)$ holds true.

\end{compactenum}
\end{lemma}
\noi
\textbf{Proof.} The GH case is an immediate consequence of Lemma \ref{extendist}. 
Or see Miermont \cite[Prop.~9]{Mier09}.
The GHP case is an immediate consequence of the GH and of  GP cases. The GP case follows from Strassen's theorem and standard arguments: the details are left to the reader. See also  \cite[Prop.~6]{Mier09}, L\"ohr \cite[Thm.~3.1]{Lohr13} and Janson \cite[Prop.~3.3 and 3.5]{Jans20} for related results.
\cqfd

\smallskip

We next recall compactness criteria in  $(\bM^0_c, \dGH)$, $(\bM, \dGP)$ and $(\bM^1_c, \dGHP)$. 
To this end, we fix $\fmu\eqo(E,d, \rho, \mu)$, a pointed finitely measured Polish metric space and we recall Definition \ref{spandex} $(a)$ of the total height $\mathtt{Ht}_{d, \rho} (A)$ of a non-empty subset $A$ of $E$. To simplify, we write $\mathtt{Ht} (A)$ instead of $\mathtt{Ht}_{d, \rho} (A)$ when there is no ambiguity. Let $\epp \ino \bbR^*_+$. We also recall from (\ref{entropie}) the definition of $N(A, \epp)$ and we define the \emph{total mass} and 
the \emph{$\mu$-essential covering number} of $E$ by 
\begin{equation}
\label{esscovdef}
\!\!\!\!\! \langle \fmu\rangle \!\!  :=\! \mu(E) \; \textrm{and} \; \esscov\,  (\fmu, \epp) \! := \! \inf\!  \big\{ \! N(K, \epp)\! + \! \mathtt{Ht} (K) \, ;  K\! \subseteq \! E \, \textrm{compact s.t.~}\mu(E\backslash K) \leqo \epp \big\}\! .  \! \!\!\!
\end{equation}

\begin{compactenum}
\item[]$\!\!\!\!\!\!\!\!\!\textbf{Grom}$-$(iv)$: \emph{A sequence $(E_n , d_n , \rho_n) \ino \bM^0_c$, $n\ino \bN$, is  
$\dGH$-precompact iff }
\begin{equation}
\label{GHprecgen}
\sup_{n\in \bN} \mathtt{Ht} (E_n) \leko  \infty \quad \textit{and} \quad \forall \epp \in \bbR^*_+ , \quad \sup_{n\in \bN} N(E_n, \epp) \leko \infty \; .
\end{equation}
\item[]$\!\!\!\!\!\!\!\!\!\textbf{Grom}$-$(v)$: \emph{A sequence $\fmu_n \! \equiv\! (E_n , d_n , \rho_n, \mu_n) \ino \bM^1_c$, $n\ino \bN$, is 
$\dGHP$-precompact iff $(E_n , d_n , \rho_n)$ satisfy (\ref{GHprecgen}) and if $\sup_{n\in \bN} \langle \fmu_n \rangle \leko  \infty$.  }

\smallskip

\item[]$\!\!\!\!\!\!\!\!\!\textbf{Grom}$-$(vi)$: \emph{A sequence $\fmu_n \! \equiv\!(E_n , d_n , \rho_n, \mu_n) \ino \bM$, $n\ino \bN$, is 
$\dGP$-precompact iff} 
\begin{equation}
\label{GPprecgen}
\sup_{n \in \bN} \langle \fmu_n \rangle  < \infty \quad \textit{and} \quad   
\forall \epp \in \bbR^*_+ , \quad \sup_{n\in \bN}\, \esscov \, (\fmu_n, \epp) < \infty. 
\end{equation}
\end{compactenum}
See Gromov, Pansu and Semmes \cite[Prop.~5.5 ]{Gro} for $\textbf{Grom}$-$(iv)$; see Abraham, Delmas and Hoscheit \cite{ADHoscheit} for $\textbf{Grom}$-$(v)$; see Gromov, Pansu and Semmes \cite[Prop.~3$\frac{1}{2}$.14]{Gro} combined with L{ö}hr, Voisin and Winter \cite[Prop.~2.6]{LVW} for $\textbf{Grom}$-$(vi)$. 
\begin{remark}
\label{redondance} If $(E,d)$ is geodesic, then $\mathtt{Ht} (E) \leqo 2\varepsilon (1+ N(E, \varepsilon))$. Thus, for geodesic spaces, only the second condition in (\ref{GHprecgen}) is needed. \cq \end{remark}

\subsection{Spaces of trees: definitions and preliminary results.}
\label{prelimGromovsec}
In this section, we first introduce compact tree spaces (purely metric: $(\bbT^0_{\! c}, \dGH)$ and finitely measured: 
$(\bbT^1_{\! c}, \dGHP)$). In doing so, we review the definition and main properties of \emph{length erasure}, which is a fundamental approximation tool in $(\bbT^0_{\! c}, \dGH)$ and will be used in the applications in Section \ref{secappl}. 
We then introduce the Gromov--Prokhorov space of minimal measured $\bbR$-trees $(\bbT, \dGP)$ and briefly discuss the connections between $\bbT^0_{\! c}$, $\bbT^1_{\! c}$ and $\bbT$ (more details are given in Section \ref{connectsubsec}). 
We conclude this section with technical results used in the proof of Theorem \ref{eraGcont} on the continuity of mass erasure  and also to show the continuity of $g$-deformations of $\bbR$-tree metrics, introduced in Notation \ref{dgtreedef}, of which the bordification is a special case. This last result is used in Section \ref{Gromerasense}, where the space $\bbT^*$ of trees measured with a boundary at infinity is defined and is equipped with the mass-erasure distance $\fdelta_{\mathtt{era}}$. 

\medskip

\noi
\textbf{Spaces of compact $\bbR$-trees and length erasure.} We first recall the definition of the Gromov--Hausdorff space of pointed compact $\bbR$-trees and the Gromov--Hausdorff--Prokhorov space of such trees equipped with a finite measure.  
\begin{definition}
\label{mT01def} We denote by $\bbT^0_{\! c} \subo \bM^0_c$ 
the isometry classes of rooted compact $\bR$-trees and we denote by $\bbT^1_{\! c}\subset \bM^1_c$ the isometry classes of rooted finitely measured
compact $\bR$-trees. \cq 
\end{definition}
\noi
From resp.~Evans, Pitman and Winter \cite[Thm 2]{EPW} and Abraham, Delmas and Hoscheit \cite[Cor 3.2]{ADHoscheit}, 
\begin{equation} 
\label{Tcpolish}
\textrm{$ \big( \bbT^0_{\! c},  \dGH \big) $ and $ \big( \bbT^1_{\! c}, \dGHP\big) $ are Polish.} 
\end{equation}
We next recall the main properties of length erasure from 
Evans, Pitman and Winter \cite[Lem. 2.6]{EPW}.   
 
\begin{lemma}
\label{lengtheraprop} Let $(T,d,\rho)\! \equiv \! \widetilde{T}\ino  \bbT^0_{\! c}$ and $h\ino \bbR_+$. 
We recall from (\ref{deflengthera}) the definition of the $h$-length-erased subtree $R_h(T)$. Then, the following holds true. 
\begin{compactenum}

\smallskip

\item[$(i)$] $d_{\mathtt{Haus}} ( R_h(T), T)\leq h$, and for all $\varepsilon \ino \bbR_+^*$ we have $N(T, h+ \varepsilon) \leq  N (R_h (T), \varepsilon)$ and $N(R_h (T), h+ \varepsilon) \leq N(T, \varepsilon)$.

\smallskip

\item[$(ii)$] If $h\geko 0$, then $R_h (T)$ is a compact subtree of $T$ which is of finite type. The isometry class of 
$(R_h (T), d, \rho)$ only depends on $\widetilde{T}$ and we denote it by $R_h (\widetilde{T})$. 

\smallskip

\item[$(iii)$] $R_h \colon \bbT^0_{\! c} \! \to \! \bbT^0_{\! c} $ is $\dGH$-continuous, and for all $h'\ino \bbR_+$, $R_h \! \circ \! R_{h'}\eqo R_{h+h'}$. 

\smallskip

\item[$(iv)$] A $\bbT^0_{\! c}$-valued sequence $(\widetilde{T}_n)_{n\in \bbN}$ is $\dGH$-precompact (resp.~$\dGH$-convergent) iff for all $h\ino \bbR_+^*$, $(R_h (\widetilde{T}_n))_{n\in \bbN}$ is $\dGH$-precompact (resp.~$\dGH$-convergent). 

\smallskip

\item[$(v)$] Let $\widetilde{\mathbf T}_n$, $n\ino \bbN$, be $\bbT^0_{\! c}$-valued r.v.s. 
Their laws are tight (resp.~converge in law) in $(\bbT^0_{\! c}, \dGH)$ iff for all $h\ino \bbR_+^*$, the laws of $R_h (\widetilde{\mathbf T}_n)$, $n\ino \bbN$, are tight (resp.~converge in law) in $(\bbT^0_{\! c}, \dGH)$. 
\end{compactenum}
\end{lemma}
\noi
\textbf{Proof.} The points $(i)$ and $(ii)$ are elementary; for $(iii)$ see Evans, Pitman and Winter \cite[Lemma 2.6]{EPW}; 
$(iv)$ follows from $(i)$, $\textbf{Haus}$-$(iii)$ and $\textbf{Haus}$-$(iv)$ in Section \ref{Merafixsec}; $(v)$ is a consequence of $(iv)$ by a standard argument in Polish spaces. \cqfd 

\medskip

\noi
\textbf{The Gromov--Prokhorov space of pointed, finitely measured minimal Polish $\bbR$-trees.} We now introduce the relevant space for \emph{Gromov--Prokhorov convergence of trees} as follows. 
\begin{definition}
\label{minimal} Let $(T, d, \rho, \mu)$ be a rooted finitely measured Polish $\bR$-tree. We say that it is \emph{minimal} 
if $T \! = \! \sspaan \, (\supp  \mu)$.  \cq 
\end{definition}
Let $(T,d,\rho,\mu)$ and $(T^\prime, d^\prime, \rho^\prime, \mu^\prime)$ be  rooted finitely measured Polish $\bR$-trees. Suppose that 
$\dGP(\mu, \mu^\prime)\eqo 0$. Then the two rooted finitely measured Polish $0$-hyperbolic spaces 
$(\supp \mu , d, \rho, \mu)$ and $(\supp \mu^\prime , d^\prime, \rho^\prime, \mu^\prime)$ are isometric. By Proposition \ref{spanreafo}, this holds iff
$(\sspaan\, (\supp \mu ), d, \rho, \mu)$ and $(\sspaan\,  (\supp \, \mu^\prime) , d^\prime, \rho^\prime, \mu^\prime)$ are isometric. Therefore it makes sense to introduce the following space of $\bbR$-trees. 
\begin{definition}
\label{mdef} We denote by $\bbT$ the space of isometry classes of rooted finitely measured \emph{minimal} Polish $\bR$-trees. Let 
$\iota\colon\bbT \! \hookrightarrow \! \bM$ be the function that associates with the isometry class $\fmu$ in $\bbT$ of a rooted finitely measured minimal Polish $\bR$-tree $(T,d,\rho,\mu)$ the isometry class in $\bM$ of $(\supp\mu,d,\rho,\mu)$. 
Then, we equip $\bbT$ with the topology induced by $\iota$, i.e.~the topology associated with the metric $\dGP(\fmu,\fmu^\prime)\eqo\dGP(\iota(\fmu),\iota(\fmu^\prime))$, $\fmu,\fmu^\prime\ino\bbT$, that makes $\iota$ an isometry. \cq 
\end{definition}

We easily check that the four-point condition (\ref{4ptscondi}) which characterises $0$-hyperbolic spaces as metric spaces is a $\dGP$-closed condition. \emph{Indeed}, for all pointed finitely measured Polish metric space $(S, d, \rho, \mu)$, 
we associate the measure $q_\mu$ on $\bbR^6$ that is the push-forward of $\mu^{\otimes 4}$ via the function $(\sigma_i)_{1\leq i\leq 4}\ino T^4\! \mapsto \! (d(\sigma_i, \sigma_j))_{1\leq i<j\leq4}$. 
We then recall from Löhr, Voisin and Winter \cite[Prop.~2.6]{LVW} that 
the $\dGP$-topology on $\bM$ is the Gromov-weak topology for which in particular the function $(S,d,\rho, \mu)\! \mapsto q_\mu$ is continuous (\cite{LVW} focuses on spaces endowed with probability measures 
but the results within extend to finitely measured spaces in a straightforward way).  
Therefore we may regard $\bT$ as (isometric to) a $\dGP$-closed subspace of $\bM$ and thus 
 \begin{equation}
\label{Tpolish}
\textit{$ \big( \bbT, \dGP \big) $ is Polish.} 
\end{equation}
\noi

\smallskip

Let us briefly discuss the connections between $\bbT$, $\bbT^1_{\! c}$ and $\bbT^0_{\! c}$. We first observe the following: $(T,d,\rho, \mu) \! \equiv\!  \fmu \ino \bbT^1_{\! }$, $(T,d,\rho, \mu)$ is not necessarily minimal, i.e.~$\sspaan \, (\supp \mu )$ maybe strictly included in $T$, and therefore the isometry class of $\fmu$ may not belong to $\bbT$. This is why we introduce the following subsets and functions which are studied more precisely in Lemma \ref{cpctsubsets} at the end of Section \ref{connectsubsec} to discuss the connections between $\dGH$-, $\dGP$-, $\dGHP$-convergences and the convergence \emph{in the sense of mass erasure}, on the space $\bbT^*$ of 
trees with boundary, as introduced in Section \ref{Gromerasense}: see Theorems \ref{GPcvvsera}, \ref{GPprecpera} and \ref{GPGH+GHPdeter}. 
\begin{definition}
\label{TcTmindef} Recall the definitions of $(\bbT^0_{\! c}, \dGH)$, $(\bbT^1_{\! c}, \dGH)$ and $(\bbT, \dGP)$ from above. 
\begin{compactenum}

\smallskip

\item[$(a)$] For all $(T,d,\rho, \mu) \! \equiv\!  \fmu\ino \bbT^1_{\! c}$, we denote by $\Phi^1_0 (\fmu) \ino \bbT^0_{\! c}$ the isometry class of $(T,d, \rho)$. We note that $\Phi^1_0 \colon  \bbT^1_{\! c} \! \to \! \bbT^0_{\! c}$ is a $1$-($\dGHP, \dGH$)-Lipschitz surjective function.

\smallskip

\item[$(b)$] For all $(T,d,\rho, \mu) \! \equiv\! \fmu \ino \bT^1_{\! c}$, we denote by $\spaan \, (\fmu) $ the isometry class in 
$\bbT^1_{\! c}$ of the compact measured $\bbR$-tree 
$(\spaan \, (\supp \mu ) , d, \rho, \mu)$, which is minimal. We denote by 
$\bbT^1_{\! c, \mathtt{min}}$ the set of $ \fmu\ino \bbT^1_{\! c}$ such that $\spaan\,  (\fmu)\eqo \fmu$. 

\smallskip

\item[$(c)$] We denote by $\bbT_{\! c}$ the set of $(T,d,\rho, \mu)  \! \equiv\!  \fmu\ino \bbT$ such that $(T,d)$ is compact. 

\smallskip

\item[$(d)$] For all $(T,d,\rho, \mu)  \! \equiv\!  \fmu\ino \bbT_{\! c}$, we denote by $\Phi_1 (\fmu)$ the isometry class in $\bbT^1_{\! c}$ of $(T,d, \rho, \mu)$ and we denote by $\Phi_0 (\fmu)$  the isometry class in $\bbT^0_{\! c}$ of $(T,d, \rho)$. 

\end{compactenum}

\smallskip

\noi
In Lemma \ref{cpctsubsets}, we prove that $\bbT_{\! c}$ is a Borel-measurable subset of $(\bbT, \dGP)$, that 
$\bbT^1_{\! c, \mathtt{min}}$ is a Borel-measurable subset of $(\bbT^1_{\! c}, \dGHP)$, that $\spaan \colon  \bbT^1_{\! c} \! \to \! \bbT^1_{\! c, \mathtt{min}}$ is a Borel-measurable surjection, that 
$\Phi_1\colon  \bbT_{\! c} \! \to \! \bbT^1_{\! c, \mathtt{min}}$ is a Borel-bi-measurable bijection and that $\Phi_0\colon  \bbT_{\! c} \! \to \! \bbT^0_{\! c}$ is a Borel-measurable surjection. \cq 
\end{definition}

We next provide compactness criteria, which are derived from the general ones recalled in $\textbf{Grom}$-$(iv)$, $\textbf{Grom}$-$(v)$ and $\textbf{Grom}$-$(vi)$. 
Let $(T,d , \rho, \mu) \! \equiv \! \fmu\! \in \! \bbT$. 
We recall from (\ref{esscovdef}) the definition of $\langle \fmu \rangle$ and $\esscov (\fmu, \varepsilon)$ and we recall Definition \ref{spandex} $(a)$ of the height $\mathtt{Ht}_{d, \rho} (A)$ of a non-empty subset $A\subo T$, which is more simply denoted by $\mathtt{Ht} (A)$.  
We introduce the following quantity  
\begin{equation}
\label{spanescov}
\esstree \, ( \fmu, \varepsilon) \! = \!  
\inf \! \big\{N(\tau, \varepsilon)+ \mathtt{Ht} (\tau); \rho \! \in \! \tau \! \subseteq \! T \; \textrm{compact, connected s.t.} \; \mu (T\backslash \tau) \! \leq \! \varepsilon  
\big\} . 
\end{equation} 
Then we easily deduce from Lemma \ref{compacthull} that 
\begin{equation}
\label{spancov} 
\esscov (\fmu, \varepsilon) \leq \esstree \, ( \fmu, \varepsilon) \leq \varepsilon^{-1}\esscov (\fmu, \varepsilon)^2+\esscov (\fmu, \varepsilon). 
\end{equation}
We then get following 
$\dGP$-precompactness criterion. 
\begin{proposition}
\label{spanprecp} Let $\emptyset \! \neq \! \cC \! \subseteq \! \bbT$. Then, $\cC$ is 
$\dGP$-precompact iff 
\begin{equation}
\label{spcpexpli} 
\sup_{\fmu \in \cC} \langle \fmu \rangle < \infty \quad \textrm{and} \quad \forall \varepsilon 
\! \in \! \bbR^*_+  \quad  \sup_{\fmu \in \cC} \esstree \, ( \fmu, \varepsilon) < \infty . 
\end{equation}
\end{proposition}
\noi
\textbf{Proof.} This is an immediate consequence of  (\ref{spancov}) and (\ref{GPprecgen}) in $\textbf{Grom}$-$(vi)$.  \cqfd

\medskip

\noi
\textbf{Preliminary results on the Gromov--Prokhorov space of trees.} To prove the $\dGP$-continuity of mass erasure in the next section and Gromov-continuity results for the bordification and for the more general 
metric deformations introduced in Notation \ref{dgtreedef} (see Lemma \ref{jgtopoembed} below) 
we shall need several preliminary results which are stated below as Propositions \ref{implyspan}
and \ref{GroProproch}.  

\begin{proposition} \label{implyspan} 
Let $(S_i, d_i, \rho_i)$, $i\ino \{ 1,2\}$, be pointed compact $0$-hyperbolic spaces and let $R$ be a correspondence \emph{between} $S_1$ and $S_2$ with 
$(\rho_1,\rho_2)\ino R$ and 
finite distortion.  
Then, there exists a correspondence $R'$ \emph{between} $\mathtt{Span} (S_1) $ and $\mathtt{Span} (S_2)$ such that $R\!\subseteq\! R'$ and $\mathtt{dis} (R') \leqo 5\,  \mathtt{dis} (R)$. 
\end{proposition}
\noi
\textbf{Proof.} 
We set $T_i \! = \! \mathtt{Span} (S_i)$, with $S_i \!\subseteq\! T_i$. To simplify notation, arcs and branch points are denoted in the same way in $T_1$ and $T_2$. We set 
\begin{eqnarray}
\lefteqn{  V= \big\{ (\sigma_1 , \gamma_1 , \sigma_2, \gamma_2)  \ino S_1 \! \times \! T_1 \! \times \! S_2 \! \times \! T_2 : \; \,   \gamma_1 \! \in \!  \lgeo \rho_1 , \sigma_1 \rgeo, \; \,   \gamma_2 \! \in \!  \lgeo \rho_2 , \sigma_2 \rgeo,  } \qquad \qquad \quad 
                    \nonumber\\
 & &\!\!\!\!\!\!  (\sigma_1, \sigma_2)\ino R  \quad   \textrm{and}  \quad   d_1(\rho_1, \gamma_1) \! \wedge \! d_2 (\rho_2 , \sigma_2)
 \eqo d_2(\rho_2, \gamma_2) \! \wedge \! d_1 (\rho_1 , \sigma_1)  \big\} 
\end{eqnarray} 
 and $R' \eqo  \big\{ (\gamma_1, \gamma_2) \ino T_1 \! \times \! T_2: \exists 
(\sigma_1, \sigma_2) \ino S_1 \! \times \! S_2 \; \, \textrm{such that} \; \, (\sigma_1 , \gamma_1 , \sigma_2, \gamma_2) \ino 
V \big\} $. 

We see that $R\!\subseteq\! R'$. Let us first prove that $R'$ is a correspondence \emph{between} $T_1$ and $T_2$. 
\emph{Indeed}, let $\gamma_1 \ino T_1$. By definition of $T_1\eqo \mathtt{Span} (S_1)$, there is $\sigma_1 \ino S_1$ such that $\gamma_1 \! \in \!  \lgeo \rho_1 , \sigma_1 \rgeo$. Since $R$ is between $S_1$ and $S_2$, there is $\sigma_2 \ino S_2$ such that $(\sigma_1, \sigma_2)\ino R$. 
If $d_1(\rho_1, \gamma_1) \leqo d_2 (\rho_2, \sigma_2)$, then there is a unique point $ \gamma_2 \! \in \!  \lgeo \rho_2 , \sigma_2\rgeo$ such that $d_2 (\rho_2, \gamma_2)\eqo  d_1(\rho_1, \gamma_1)$ 
and we get $(\sigma_1, \gamma_1 , \sigma_2, \gamma_2) \ino V$. If 
$d_1(\rho_1, \gamma_1) \geko d_2 (\rho_2 , \sigma_2)$, then we set $\gamma_2 \eqo \sigma_2$ 
and we also get $(\sigma_1, \gamma_1 , \sigma_2, \gamma_2) \ino V$. In both cases, $(\gamma_1, \gamma_2)\ino R'$. Similarly, for all $\gamma_2\ino T_2$, there is $\gamma_1\ino T_1$ such that $(\gamma_1, \gamma_2)\ino R'$.  \cq 

\smallskip

To simplify notation, we set $\varepsilon \eqo \mathtt{dis} (R)$. It remains to show that $\mathtt{dis} (R')\leqo 5\varepsilon$. 
To this end, we first prove the following. 
\begin{equation}
\label{ancestor}
\forall (\sigma_1, \sigma_2),  (\sigma'_1, \sigma'_2)\ino R, \quad 
\big| d_1( \rho_1,\sigma_1 \! \wedge\! \sigma'_1) \! -\! d_2 (\rho_2 , \sigma_2 \!  \wedge \! \sigma'_2)\big| \leqo \tfrac{3}{2}\varepsilon. 
\end{equation}
\noi
\emph{Indeed}, this immediately follows from $d_i(\rho_i, \sigma_i \! \wedge \! \sigma_i')\eqo \tfrac12 \big( d_i(\rho_i, \sigma_i )+d_i(\rho_i,  \sigma_i')\! -\! d_i(\sigma_i , \sigma_i')\big)$. \cq

We next prove: 
\begin{equation}
\label{copi}
\forall (\gamma_1, \gamma_2) \ino R' , \quad |d_1(\rho_1,\gamma_1)\! -\!  d_2 ( \rho_2 , \gamma_2) | \leqo \varepsilon.  
\end{equation}
\noi
\emph{Indeed}, if $(\sigma_1, \gamma_1, \sigma_2, \gamma_2) \ino V$, then 
$d_1(\rho_1,\gamma_1)\eqo d_2(\rho_2,\gamma_2)$ unless one of the following holds
$$d_2(\rho_2,\sigma_2)\eqo d_2(\rho_2,\gamma_2)\leko d_1(\rho_1,\gamma_1)\leqo d_1(\rho_1,\sigma_1),\  
d_1(\rho_1,\sigma_1)\eqo d_1(\rho_1,\gamma_1)\leko d_2(\rho_2,\gamma_2)\leqo d_2(\rho_2,\sigma_2) $$
but since the outer terms in these inequalities differ by less than $\varepsilon$, the claim follows.\cq

\smallskip

We now suppose that $(\sigma_1 , \gamma_1 , \sigma_2, \gamma_2),(\sigma'_1 , \gamma'_1 , \sigma'_2, \gamma'_2)\ino V$ and we prove that this entails that $|d_1(\gamma_1, \gamma_1') \! -\! d_2 (\gamma_2, \gamma_2')|\leqo 5\varepsilon$. 
We consider two cases. 

\smallskip

\noi
$\bullet$ Case 1: we suppose that $\gamma_1,\gamma'_1  \ino \lgeo \rho_1,\sigma_1 \! \wedge \! \sigma'_1 \lgeo \, $.
If $\gamma_2 \! \in\,  \rgeo \sigma_2 \! \wedge \! \sigma_2^\prime,\sigma_2\rgeo$ and  
$\gamma'_2 \! \in\,  \rgeo \sigma_2 \! \wedge \! \sigma_2^\prime,\sigma'_2\rgeo$, then 
\begin{align*}
d_2(\gamma_2,\gamma'_2)
&=d_2(\gamma_2 ,\sigma_2 \! \wedge \! \sigma_2^\prime)+
d_2(\gamma_2^\prime,\sigma_2 \! \wedge \! \sigma_2^\prime)\eqo 
d_2(\rho_2,\gamma_2)
   \! -\! 2d_2(\rho_2,\sigma_2 \! \wedge \! \sigma_2^\prime)
   +d_2(\rho_2,\gamma_2^\prime)\\
&\le d_1(\rho_1,\gamma_1)\! -\! 2d_1(\rho_1,\sigma_1 \! \wedge \! \sigma'_1)+d_1(\rho_1,\gamma_1')+5\varepsilon \leqo  5 \varepsilon \leqo d_1(\gamma_1,\gamma'_1)+5\varepsilon\, , 
\end{align*}
by \eqref{ancestor}--\eqref{copi}. Otherwise,
$d_2(\gamma_2,\gamma_2')
=|d_2(\rho_2,\gamma_2^\prime)\! -\! d_2(\rho_2,\gamma_2)|
\leqo |d_1(\rho_1,\gamma_1')\! -\! d_1(\rho_1,\gamma_1)|+2\varepsilon$
$\eqo d_1(\gamma_1,\gamma_1)+2\varepsilon,
$ by \eqref{copi}.

\smallskip

\noi 
$\bullet$ Case 2: In all other configurations $d_1(\gamma_1,\gamma'_1)\eqo 
d_1(\gamma_1,\sigma_1\! \wedge \! \sigma'_1)+d_1(\gamma'_1,\sigma_1 \! \wedge \! \sigma'_1)$ 
and
\begin{align*} 
d_2(\gamma_2,\gamma'_2)&\leqo d_2(\gamma_2,\sigma_2 \! \wedge \! \sigma'_2)+d_2(\gamma'_2,\sigma_2 \! \wedge \! \sigma'_2)  \eqo 
d_2(\rho_2,\gamma_2)
   \! -\! 2d_2(\rho_2,\sigma_2 \! \wedge \! \sigma_2^\prime)
   +d_2(\rho_2,\gamma_2^\prime) \\
  & \leqo  d_1(\rho_1,\gamma_1)
   \! -\! 2d_1(\rho_1,\sigma_1 \! \wedge \! \sigma_1^\prime)
   +d_1(\rho_1,\gamma_1^\prime) +5\varepsilon 
   \eqo d_1(\gamma_1,\gamma'_1)+5\varepsilon , 
\end{align*}
by \eqref{ancestor}--\eqref{copi}.

\noi

In both cases, this proves 
$d_2(\gamma_2, \gamma'_2) \leqo d_1(\gamma_1, \gamma'_1) +5 \varepsilon $. By similar arguments we also prove that 
$d_1(\gamma_1, \gamma'_1) \leqo d_2(\gamma_2, \gamma'_2) +5 \varepsilon $, and we get the desired result, which completes the proof. \cqfd

\begin{proposition}
\label{GroProproch} 
For $i\ino \{ 1,2\}$ let $(T_i,d_i,\rho_i, \mu_i)\! \equiv \! \fmu_i$ be rooted finitely measured Polish $\bbR$-trees (not necessarily minimal) such that $\dGP (\fmu_1, \fmu_2) \leko \eta$. Then, there are compact subtrees 
$\rho_i \ino \tau_i\! \subseteq\! T_i$, $i\ino \{ 1,2\}$, a compact correspondence $R$ 
\emph{between} $\tau_1$ and $\tau_2$ such that $(\rho_1, \rho_2) \ino R$ and $\mathtt{dis} (R) \leko 10\eta$, and an 
$\eta$-coupling $\mathtt m$ of $\mu_1$ and $\mu_2$ such that $\mathtt m ((T_1\! \times \! T_2) \backslash R)\eqo 0$. 

\end{proposition}
\noi
\textbf{Proof.} By Lemma \ref{rephraseGrom} there is a partial compact correspondence $R_0\!\subseteq\! T_1\! \times \! T_2$ such that $(\rho_1, \rho_2) \ino R_0$ and $\mathtt{dis} (R_0) \leko 2\eta$, and an 
$\eta$-coupling $\mathtt m$ of $\mu_1$ and $\mu_2$ such that $\mathtt m ((T_1\! \times \! T_2) \backslash R_0)\eqo 0$. For $i\ino \{ 1,2\}$, we set $\tau_i \! :=\! \mathtt{Span} (\pi_i (R_0))$, which is a compact subtree of $T_i$, since $\pi_i(R_0)$ is compact and by Lemma \ref{compacthull}. By Proposition \ref{implyspan}, there is a correspondence \emph{between} $\tau_1$ and $\tau_2$ such that $R_0\!\subseteq\! R$ and $\mathtt{dis} (R)\leqo  10 \eta$. W.l.o.g.~$R$ can be assumed to be a closed in $\tau_1\! \times \! \tau_2$, and thus compact. \cqfd

\medskip

\noi
\textbf{Gromov--Prokhorov continuity of the bordification.} We conclude this section by discussing Gromov-continuity results for the bordification, and the more general metric deformations introduced in Notation \ref{dgtreedef}. 
To this end, we fix $a,b\ino (0, \infty]$ and $g\colon  [0, a)\! \to \! [0, b)$, an increasing homeomorphism. 
Let $a_0\ino [0, a)$, $\eta\ino \bbR^*_+$. We denote the $\eta$-modulus of continuity of $g$ on $[0, a_0]$ by  
\begin{equation}
\label{gmodulus}
\mathtt w_{g} ([0, a_0], \eta) \eqo \max \big\{ |g(x) \! -\! g(y)|\, ; \, x,y\ino [0, a_0]:  |x\! -\! y| \leqo \eta \big\} , 
\end{equation} 
The following lemma explains how correspondences are affected by the metric deformation $d\! \mapsto \! d_g$. 
\begin{lemma}
\label{distcorrdist} We keep previous notation. Let $R$ be a correspondence \emph{between} two rooted $\bbR$-trees $(T,d,\rho)$ and 
$(T'\! ,d'\! ,\rho')$ with $(\rho,\rho')\ino R$. Suppose $a_0 \eqo  \mathtt{dis}_{d,d'} (R) + \max (\mathtt{Ht}_{d, \rho} (T) ,\mathtt{Ht}_{d'\! , \rho'} (T') )\leko a $. 
Then 
$T\eqo T_{\! g}$, $T'\eqo T'_{\! g}$ and $\mathtt{dis}_{d_g,d'_g} (R) \leqo 4 \mathtt w_{g} \big([0, a_0] , \frac32 \mathtt{dis}_{d,d'} (R) \big)$. 
\end{lemma}
\noi
\textbf{Proof.} Our assumptions immediately entail the set equalities $T\eqo T_{\! g}$, $T'\eqo T'_{\! g}$. 
To simplify notation, we set $\eta \eqo \mathtt{dis}_{d,d'} (R)$ and for all $\eta'\ino \bbR^*_+$, $w(\eta') \! :=\!  \mathtt w_{g} ([0, a_0], \eta') $. Let $(\sigma_1, \sigma'_1), (\sigma_2, \sigma'_2)\ino R$. 
As in \eqref{ancestor}, 
$|d(\rho,\sigma_1\! \wedge \! \sigma_2) \! -\!  d'(\rho',\sigma'_1\! \wedge \! \sigma'_2)| 
\leqo \frac32\eta$.
Moreover, by definition of $\mathtt{dis}_{d,d'} (R)$, for all $i\ino \{ 1,2\}$, $|d(\rho, \sigma_i) \! -\! d'(\rho', \sigma'_i)| \leqo \eta$. Thus, $|g(d(\rho,\sigma_1\! \wedge \! \sigma_2)) \! -\! g(d'(\rho',\sigma'_1\! \wedge \! \sigma'_2))|\leqo  w(\frac32\eta)$ and $|g(d(\rho, \sigma_i)) \! -\! g(d'(\rho', \sigma'_i))|\leqo w(\eta)$, which easily entails the desired result. \cqfd

\smallskip

To state Gromov-continuity results for the bordification we introduce the following tree spaces.  
\begin{definition}
\label{Tadef} Let $a,b\ino (0, \infty]$ and let $g\colon  [0, a)\! \to \! [0, b)$ be an increasing homeomorphism. 
\begin{compactenum}

\smallskip

\item[$(a)$] We denote by $\bbT(a)$ the space of $(T,d,\rho, \mu)\! \equiv\!  \fmu \ino \bbT$ such that $\mathtt{Ht}_{d, \rho} (T) \leqo a$ and such that $\mu \big(\{ \sigma \ino T : d(\rho, \sigma) \eqo a \}  \big)\eqo 0$. We observe that $\bbT(\infty) \eqo \bbT$.

\smallskip

\item[$(b)$] We also denote by $\bbT^{1}_c(a)$ (resp.~by $\bbT^{0}_{c}(a)$) 
the space of $(T,d,\rho, \mu)\! \equiv\! \fmu \ino \bbT^1_{\! c}$ (resp.~of $(T,d,\rho)\! \equiv\!  \widetilde{T} \ino \bbT^0_{\! c}$) such that $\mathtt{Ht}_{d, \rho} (T) \leko a$. Note that 
$ \bbT^1_{\! c}(\infty)\eqo \bbT^1_{\! c}$ and $ \bbT^0_{\! c}(\infty)\eqo \bbT^0_{\! c}$.

\smallskip

\item[$(c)$] Let $(T,d,\rho, \mu)\! \equiv \! \fmu \ino \bbT(a)$ (resp.~$(T,d,\rho, \mu)\! \equiv\!  \fmu \ino \bbT^1_{\! c}(a)$, $(T,d,\rho)\! \equiv \! \widetilde{T} \ino \bbT^0_{\! c}(a)$). We recall $(T_{\! g}, d_g, \rho, \mu)$ from Notation \ref{dgtreedef} and Remark \ref{toporem} $\mathbf{(c)}$. By Theorem \ref{propdgtree} $(iii)$, the isometry class of $(T_{\! g}, d_g, \rho, \mu)$ belongs to $\bbT(b)$ (resp.~to $\bbT^1_{\! c}(b)$, to $\bbT^0_{\! c}(b)$) and only depends on that of $(T,d,\rho, \mu)$, so it makes sense to denote it by $\jmath_g (\fmu)$ (resp.~by $\jmath_g^1(\fmu)$,  $\jmath_g^0(\widetilde{T})$).

\smallskip

\item[$(d)$] We denote by $\jmath_*(\fmu)$ (resp.~$\jmath_*^1(\fmu)$, $\jmath^0_*(\widetilde{T})$) the isometry class 
of the bordification $(T^*\! , d^*\! , \rho, \mu)$, which corresponds to the special case where $g\colon x\ino \bbR_+ \! \mapsto \! 1\! -\! e^{-x} \ino [0, 1)$. \cq
\end{compactenum}
\end{definition}
\begin{lemma}
\label{jgtopoembed} We keep the notation of Definition \ref{Tadef}. 
Then the following holds true.
\begin{compactenum}

\smallskip

\item[$(i)$] $\jmath_g$ is a homeomorphism from $(\bbT(a), \fdelta_{\mathtt{GP}})$ onto $(\bbT(b),\fdelta_{\mathtt{GP}})$. 
In particular, $\jmath_*$ is an homeomorphism from $(\bbT,  \fdelta_{\mathtt{GP}})$  onto  $(\bbT(1), \fdelta_{\mathtt{GP}})$. 
\smallskip

\item[$(ii)$] $\jmath^1_g$ is a homeomorphism from $(\bbT^1_{\! c}(a), \fdelta_{\mathtt{GHP}})$ onto 
$(\bbT^1_{\! c}(b), \fdelta_{\mathtt{GHP}})$. 

\smallskip

\item[$(iii)$] $\jmath^0_g$ is a homeomorphism from $(\bbT^0_{\! c}(a), \fdelta_{\mathtt{GH}})$ onto 
$(\bbT^0_{\! c}(b), \fdelta_{\mathtt{GH}})$.
\end{compactenum}
\end{lemma}
\noi
\textbf{Proof.} We first prove $(i)$. By Theorem \ref{propdgtree} $(iii)$, $\jmath_g$ is a bijection from $\bbT(a)$ to 
$\bbT(b)$, whose inverse function is $\jmath_f$, where $f$ stands for the inverse function of $g$. 
So we only need to prove the  $\fdelta_{\mathtt{GP}}$-continuity of $\jmath_g$ for any $g$. 
To this end, we fix  $(T,d,\rho, \mu)\! \equiv \! \fmu \ino \bbT(a)$ and $\epp \ino \bbR^*_+$. By definition of $\bbT(a)$, there is $a_0 \ino (0, a)$, which depends on $\fmu$ and $\epp$, such that 
$\mu (T\backslash B_{a_0})\leqo \frac12 \epp$, where $B_{a_0}\eqo B_{T,d} (\rho, a_0)$. 
We next fix a positive real number $\eta$ such that 
\begin{equation}
\label{condioneta}
\eta \leko \tfrac{1}{4} \min \big( a\! -\! a_0 \, ,\,  \epp \big) \quad \textrm{and} \quad \mathtt w_{g} \big( \, [0, \tfrac{{1}}{{2}} (a+a_0)]\, , \, 15 \eta \big) \leko \tfrac14 \epp \, , 
\end{equation}
and  $(T',d',\rho', \mu')\! \equiv \! \fmu' \! \ino \bbT(a)$ such that $\fdelta_{\mathtt{GP}}(\fmu, \fmu') \leko  \eta$. 
By Lemma \ref{rephraseGrom}, there is a compact correspondence $R_0\!\subseteq\! T\! \times \! T'$ such that 
$(\rho, \rho')\ino R_0$ and $\mathtt{dis}_{d,d'\!} (R_0) \leko 2\eta$, and an $\eta$-coupling 
$\mathtt m_0 \ino \cM_{\! f} (T\! \times \! T')$ of $\mu$ and $\mu'$ such that $\mathtt m_0 \big( (T\! \times \! T') \backslash R_0 \big) \eqo 0$. We set $R_1\! :=\! R_0 \cap (B_{a_0} \! \times \! T')$, which is a compact correspondence such that $(\rho, \rho')\ino R_1$ and $\mathtt{dis}_{d,d'\!} (R_1) \leqo \mathtt{dis}_{d,d'\!} (R_0)\leko 2\eta$, and 
we set $\mathtt m\! :=\! \mathtt m_0 (\, \cdot \, \cap R_1)$. Let us check that $\mathtt m$ is an $\epp$-coupling of $\mu$ and $\mu'$. For all $B\ino \cB(T)$, we observe that $\mu(B)\! -\!  \mathtt m (B\! \times \! T') \eqo \mu (B\backslash B_{a_0})+ \mu (B\cap B_{a_0}) \! -\! \mathtt m_0 ((B_{a_0} \cap B)\! \times \! T'  ) \ino [0, \frac{1}{2} \epp + \eta )$. For all $B'\ino \cB(T')$, we next note that $ \mathtt m_0 ((T\backslash B_{a_0} )\! \times \! B'  )\leqo  \mathtt m_0 ((T\backslash B_{a_0} )\! \times \! T'  )\leqo \mu (T\backslash B_{a_0} )\leqo\frac{1}{2}  \epp$. Thus, $\mu'(B')\! -\!  \mathtt m (T\! \times \! B') \eqo 
\mu'(B')\! -\!  \mathtt m_0 (B_{a_0}\! \times \! B')\eqo  \mu'(B')\! -\!  \mathtt m_0 (T\! \times \! B')+ \mathtt m_0 ((T\backslash B_{a_0})\! \times \! B') \ino  [0,\frac{1}{2}  \epp + \eta )$. By (\ref{condioneta}), $\frac{1}{2} \epp + \eta \leko \epp$, which entails the desired result. 

We next observe that $R_1$ is a correspondence between the compact subsets 
$K\! :=\! \pi (R_1)$ and $K'\! :=\! \pi' (R_1)$, where for all $(\sigma,\sigma')\ino  T\! \times \! T'$, we have set 
$\pi (\sigma,\sigma')\eqo \sigma$ and $\pi' (\sigma,\sigma')\eqo \sigma'$. We also check that $\mathtt{Ht}_{d, \rho}( K) \leqo a_0$ and thus $\mathtt{Ht}_{d'\! , \rho'}( K')\leqo a_0+ \mathtt{dis}_{d,d'\!} (R_1) \leko a_0 +2\eta \leqo \frac12 (a+a_0)$. We next set $\tau \! :=\! \mathtt{Span} (K)$ and $\tau' \! :=\! \mathtt{Span} (K')$. These two substrees of $T$ and $T'$ 
are compact, by Lemma \ref{compacthull}. By Proposition \ref{implyspan}, there is a compact correspondence $R_2$ 
\emph{between} $\tau$ and $\tau'$ such that $R_1\!\subseteq\! R_2$ and $\mathtt{dis}_{d, d'\!} (R_2) \leqo 5\mathtt{dis}_{d, d'\!} (R_1) \leko 10 \eta$. We also see that $\mathtt{Ht}_{d, \rho} (\tau) \vee \mathtt{Ht}_{d'\! , \rho'\! } (\tau') \leko \frac12 (a+a_0)$. So Lemma \ref{distcorrdist} applies and we get $\mathtt{dis}_{d_g,d'_g} (R_2) \leqo 4 \mathtt w_g ([0, \frac{{1}}{{2} } (a+a_0)] , 15 \eta) \leko \epp$. 

We then recall from Lemma \ref{fchangemetric} that $T^o$ and $T'^o$ stand for the open balls with radius $a$ around the roots and we observe that $R_2$ is a compact subset of $T^o\! \times \! T'^o$ and thus of $T_{\! g}  \! \times \! T'_{\! g} $. 
Since $\mathtt{Supp} (\mathtt m)\!\subseteq\! R_2$, we can also view $\mathtt m$ as an element of $\cM_{\! f} (T_{\! g}  \! \times \! T'_{\! g})$. Therefore $\mathtt m$ is an $\epp$-coupling of $\mu$ and $\mu'$ such that $\mathtt m ((T_{\! g}  \! \times \! T'_{\! g})\backslash R_2 )\eqo 0$. By Lemma \ref{rephraseGrom}, we get $\fdelta_{\mathtt{GP}} (\jmath_{g} (\fmu), \jmath_{g} (\fmu')) \leko \epp$, which entails $(i)$. 
The proofs $(ii)$ and $(iii)$ are quite similar (and somewhat simpler): we leave the details to the reader. \cqfd

\subsection{Gromov--Prokhorov continuity of mass erasure}
\label{meraC0sec}

In this section, we prove a key result of the article, Theorem \ref{eraGcont}, 
which discusses continuity properties of mass erasure with respect to $ \dGP$ and $ \dGHP$. 
\begin{notation}
\label{GPmeradef} Let $(T,d,\rho, \mu) \! \equiv \! \fmu \ino \bbT$ (resp.~$(T,d,\rho, \mu) \! \equiv \! \fmu \ino\bbT^1_{\! c}$). 
\begin{compactenum}

\smallskip

\item[$(a)$] Let $h\ino \bbR^*_+$. By Proposition \ref{properased} $(iv)$, $\sspaan\,  (\mathtt{\supp} (\era_h \mu) )\eqo T_{\! \mu, h+} $ and $(T_{\! \mu, h+}, d, \rho, \era_h \mu)$ is therefore a rooted finitely measured \emph{minimal} compact  $\bbR$-tree whose isometry class only depends on that of $(T,d,\rho, \mu)$ (see Remark \ref{checkinvari} below): we denote it by $\cE_h \fmu \ino \bbT$ (resp.~$\bbT^1_{\! c}$).
We derive from Proposition \ref{properased} $(v)$ that \vspace{-0.1cm}
\begin{equation}
\label{GPmerasemi}
\forall h,h'\ino \bbR^*_+, \quad \era_h \circ \era_{h'}\eqo \era_{h+h'} \; \, \textrm{on $\bT$ or $\bbT^1_{\! c}$}\vspace{-0.1cm}
\end{equation} 
which is refered to as the \emph{semigroup property of mass erasure on} $\bbT$ (resp.~\emph{on} $\bbT^1_{\! c}$). \pagebreak

\smallskip

\item[$(b)$] Let  $\bhh\eqo (h_p)_{p\in \bbN}$ be $\bbR_+^*$-valued and strictly decreasing to $0$. 
We recall the definition of $\fell_\bullet^{\uparrow}(\bbN)$ from (\ref{ellspacedef}) and  
Definition \ref{Lambdadef} of the $\bhh$-projected height measure $\Lambda_{\bhh, \mu}$ and of the height measure $\lambda_\mu$.  
Since $\Lambda_{\bhh, \mu}$ and $\lambda_\mu$ actually depend on the equivalence class $\fmu$ of $(T, d, \rho, \mu)$ only, we use notation $\fLambda_{\bhh, \fmu}$ for $\Lambda_{\bhh, \mu}$ and $\flambda_\fmu$ for $\lambda_\mu$. \cq 
\end{compactenum}
\end{notation}
\begin{remark}
\label{checkinvari} Let $\varphi\colon (T,d,\rho, \mu) \! \to \! (T'\! ,d'\! ,\rho'\! , \mu')$ be a bijective isometry. Then 
for all $\sigma \ino T$, $\varphi(\theta_{\sigma} T)\eqo \theta_{\varphi(\sigma) } T'$ and 
$\mu'(  \theta_{\varphi(\sigma) } T') \eqo \mu (\theta_{\sigma} T)$. Therefore the restriction of $\varphi$ to $T_{\! \mu, h+}$ is a bijective isometry onto $T'_{\! \mu', h+}$ and 
$(\era_h \mu) \circ \varphi^{-1} \eqo \era_h \mu'$. Namely, $(T_{\! \mu, h+}, d, \rho, \cE_h\mu) \! \equiv \! (T'_{\! \mu', h+}, d'\! , \rho'\! , \cE_h\mu')$. \cq 
\end{remark}
\begin{theorem}
\label{eraGcont} Let $\fmu\! \equiv \! (T, d, \rho, \mu)$ and $\fmu_n \! \equiv \! (T_n, d_n , \rho_n, \mu_n)$, $n\ino \bN$, 
belong to $\bbT$ (resp.~to $\bbT^1_{\! c}$). We assume that $\lim_{n\to \infty}\dGP(\fmu_n, \fmu)\eqo 0 $ (resp.~that $\lim_{n\to \infty}\dGHP(\fmu_n, \fmu)\eqo 0 $). Let $h\ino \bbR^*_+$. Then the following holds true. 
\begin{compactenum}

\smallskip

\item[$(i)$] $\sup_{h'\in [h, \infty)} \dGP (\era_{h'} \fmu_n
, \era_{h'} \fmu) \to 0$ as $n\! \to \! \infty$. 
 
 \smallskip
 
\item[$(ii)$] Let $(h_n)_{n\in \bN}$ be $\bbR^*_+$-valued and converging to $h\ino \bbR^*_+$. Let us assume that 
$h^\prime\!  \mapsto \! T_{\! \mu, h^\prime}$ is $\dHaus$-continuous at $h$. Let us set $T^n_{\! \mu_n, h'}\eqo (T_n)_{\mu_n, h'}$ for all $h'\ino \bbR^*_+$. Then, 
\begin{align}
\label{projC0GHP}
\dGHP\big( (T^n_{\! \mu_n, h_n} , d_n, \rho_n , \era_{h_n} \mu_n)\, &,\,  (T_{\! \mu , h}, d, \rho, \era_{h} \mu) \big) \longrightarrow 0 \\
\label{eraC0GHP}
 \textrm{and}\quad \dGHP\big( (T^n_{\! \mu_n, h_n} , d_n, \rho_n ,  \mathtt{P}_{T^n_{\! \mu_n, h_n}} \mu_n)\, &, \, (T_{\! \mu , h}, d, \rho, \mathtt{P}_{T_{\! \mu , h}} \mu) \big)  \longrightarrow 0\, . 
\end{align}
\item[$(iii)$] Let $(h_n)_{n\in \bN}$ be as in $(ii)$. Let us furthermore assume that $\supp \mu \! = \! T$ (by 
Proposition \ref{massTprop} $(v)$, $h^\prime \! \mapsto \! T_{\! \mu, h^\prime}$ is therefore $\dHaus$-continuous). 
Then, 
\begin{align}
\label{eraC0uniGHP}
\sup_{h^\prime\in [h, \infty) } \dGHP\big( (T^n_{\! \mu_n, h'} , d_n, \rho_n ,  \era_{h^\prime}\mu_n)\,&,(T_{\! \mu , h^\prime},  d, \rho, \era_{h^\prime} \mu) \big) \longrightarrow 0 \nonumber \\
 \textrm{and}  \quad \sup_{h^\prime\in [h, \infty) } \dGHP\big( (T^n_{\! \mu_n, h'} , d_n, \rho_n , \mathtt{P}_{T^n_{\! \mu_n, h'}} \mu_n)\, &,\,  (T_{\! \mu , h^\prime}, d, \rho, \mathtt{P}_{T_{\! \mu , h^\prime}} \mu) \big)  \longrightarrow 0\, . 
\end{align}
\item[$(iv)$] Let $\bhh\eqo (h^*_p)_{p\in \bbN}$ be $\bbR^*_+$-valued and strictly decreasing to $0$. Let us assume 
that $h\mapsto \! T_{\! \mu, h}$ is $\dHaus$-continuous at each $h^*_p$, $p\ino \bbN$. 
Then $\fLambda_{\bhh, \fmu_n} \!\!\!   \to \! \fLambda_{\bhh, \fmu}$ weakly on $\cM_f(\fell_\bullet^{\uparrow}(\bbN))$. 

\smallskip

\item[$(v)$] $\flambda_{\fmu_n} \! \!\! \to \! \flambda_{\fmu}$ weakly on $\cM_{\! f} (\bbR_+)$. 
\end{compactenum}
\end{theorem}
\noi
\textbf{Proof.} 
We first prove (\ref{eraC0GHP}) in $(ii)$. To this end, let 
$(\varepsilon_n)_{n\in \bN}$ be a $\bbR^*_+$-valued sequence converging to 
$0$ such that $ |h_n\! -\! h|+ \frac{11}{10}\dGP (\fmu_n , \fmu)\! < \! \varepsilon_n$ 
(resp.~$ |h_n\! -\! h|+ \frac{11}{10}\dGHP (\fmu_n , \fmu)\! < \! \varepsilon_n$). 
By Proposition \ref{GroProproch}, there are compact subtrees $\rho_n \ino \tau^\prime_n \!\subseteq\! T_n$ 
and $\rho\ino \tau_n\! \subseteq\! T$, a compact correspondence $R_n$ 
\emph{between} $\tau^\prime_n$ and $\tau_n$ such that $(\rho_n, \rho) \ino R_n$ and $\mathtt{dis}_{d_n, d} (R_n) 
\leko 10 \varepsilon_n$, and a 
$(2\epp_n)$-coupling $\mathtt m_n$ of $\mu_n$ and $\mu$ such that $\mathtt m_n ((T_n\! \times \! T) \backslash R_n)\eqo 0$;
if $d_{T, T_n}$ stands for the distance on $T_n \sqcup T$ extending $d_n$ and $d$ 
obtained as in Lemma \ref{extendist} with $r\eqo 5 \epp_n$, then  
$\mu_n (T_n \backslash \tau^\prime_n) \! \vee \! \mu (T\backslash \tau_n) \leqo 2\epp_n $, and by Strassen's theorem, 
\begin{align}
\dPro^{T,T_n} \! (\mu   , \mu_n) \leqo  5\varepsilon_n  & , \quad    d_{T, T_n} (\rho, \rho_n 
 )\!  \leq\!   5\varepsilon_n \nonumber \\  
 &    \textrm{and} \quad 
0< \min \{ d_{T,T_n} (\sigma, \sigma^\prime)\, ; \,  \sigma \ino T, \sigma^\prime\!  \ino T_n \}  = 
\dHaus^{T, T_n} (\tau_n, \tau_n^\prime) \! \leq\!   5\varepsilon_n   , \label{TTen}
\end{align}
where $\dPro^{T, T_n}$  stands for the Prokhorov distance on $\cM_f (T \sqcup T_n)$ and $\dHaus^{T, T_n}$ 
stands for the Hausdorff distance on the space of compact subsets of $T \sqcup T_n$ (under the assumption of a 
$\dGHP$-convergence of the $\fmu_n$, we can take $\tau'_n \eqo T_n $ and $\tau_n\eqo T$). 
Then, Lemma \ref{seqspace} applies to $X\! = \! T$ and $X_n\! = \! T_n$ and $0\leko r_n \! \leqo  \! 5\varepsilon_n $. 
There exists a metric $\delta$ on $E^o\! = \! T \! \sqcup \bigsqcup_{n\in \bN} T_n$ that extends 
the $d_{T, T_n}$. We denote by $(E, \delta)$ the corresponding completed space that is Polish since 
$T$ and the $T_n$ are Polish. We slightly abuse notation by assuming $E^o\! \subseteq \! E$.
W.l.o.g.~we assume that $h  \! >\! 30\varepsilon_n$, which implies $h_n \! >\! 29\varepsilon_n$. 
To simplify notation, we next set $\nu^\prime_n \! = \! \mathtt{P}_{\! \tau^\prime_n}\mu_n$ and $\nu_n \! = \! \ttP_{\! \tau_n} \mu$. 
By Proposition \ref{massTprop} $(i)$,
$T^n_{\! \mu_n, h_n} \! =\! T^n_{\! \nu_n' , h_n} \! \subseteq \! \tau^\prime_n$ 
and for all $h^\prime \! \in \!  (4\varepsilon_n , \infty)$,  
$T_{\! \mu, h^\prime}\! = \! T_{\nu_n , h^\prime} \! \subseteq \! \tau_n$. 
Then, by Lemma \ref{projprok} and Proposition \ref{properased} $(i)$, 
\begin{equation}
\label{sciondub}
\ttP_{T^n_{\! \mu_n, h_n}} \mu_n \! = \! \ttP_{T^n_{\! \mu_n, h_n}} \nu^\prime_n, \quad \ttP_{T_{\! \mu, h}} \mu \! 
= \!   \ttP_{T_{\! \mu , h}} \nu_n , \quad  \era_{h_n} \mu_n = \era_{h_n} \nu^\prime_n, \quad  \era_h \mu = \era_h \nu_n. 
\end{equation}  
We prove for all $n\ino \bbN$, that 
\begin{equation}
\label{delproproj}
\dePro (\nu_n , \nu^\prime_n) \! < \! 13\varepsilon_n .
\end{equation}
\emph{Indeed}, by Proposition \ref{projdef} $(v)$, we get $\dPro^{{T_n}} (\mu_n , \nu^\prime_n) \leqo \mu_n (T_n\backslash \tau_n^\prime) \! <\!  4\varepsilon_n$. By $\textbf{Pro}$-$(iii)$, we see that $\dePro (\mu_n , \nu^\prime_n)\eqo \dPro^{{T_n}} (\mu_n , \nu^\prime_n)$. Thus $\dePro (\mu_n , \nu^\prime_n) \! < \! 4\varepsilon_n$. Similarly, we get $\dePro (\mu , \nu_n) \! < \! 4\varepsilon_n$. 
By $\textbf{Pro}$-$(iii)$ again and (\ref{TTen}), $\dePro (\mu, \mu_n)\eqo d^{{T,T_n}}_{\mathtt{Pro}} (\mu,\mu_n) \leko 5\epp_n$, which together with the previous inequalities implies (\ref{delproproj}). \cq 

\smallskip

We recall from (\ref{TTen}) that 
$\deHaus (\tau_n, \tau'_n) \eqo \dHaus^{T, T_n} (\tau_n, \tau_n^\prime) \! \leq\!   5\varepsilon_n   $. Combined with (\ref{delproproj}), this shows that  Proposition \ref{preGrera} applies with $T\eqo \tau_n$, $\mu\! = \! \nu_n$, $T^\prime\! = \! \tau_n^\prime$, $\mu^\prime\! = \! \nu_n^\prime$ and $\varepsilon \! =\! 13\varepsilon_n$, and we get 
\begin{eqnarray} 
\label{GHaussup}
\deHaus (T^n_{\! \mu_n , h_n}, T_{\! \mu, h}) &=&  \deHaus (T^n_{\! \nu_n' , h_n} , 
T_{\nu_n , h}) \nonumber \\
& \leq &  52 \varepsilon_{n} \! +\!   \deHaus \big( T_{\nu_n , h-26\varepsilon_n} ,  T_{\nu_n , h+26\varepsilon_n} \big) \nonumber  \\
& = &  52 \varepsilon_{n} \! +\!   \deHaus \big( T_{\! \mu , h-26\varepsilon_n} ,  T_{\! \mu , h+26\varepsilon_n} \big)= : \eta_n   \xrightarrow[n\to \infty]{\;} 0 . 
\end{eqnarray}
We recall that 
$T_{\! \mu, h} \! \subseteq \! \tau_n$, that $T^n_{\! \mu_n , h_n} \! \subseteq \! \tau'_n$ and that 
$\deHaus (\tau_n, \tau'_n) \leko 5\epp_n$. By (\ref{GHaussup}) and by (\ref{delproproj}), 
Lemma \ref{GrC0proj} applies to $T\eqo \tau_n$, $\tau\! = \! T_{\! \mu, h}$, $\mu\! = \! \nu_n$  
and $T^\prime\! =\!  \tau^\prime_n$, $\tau^\prime\! = \! T^n_{\! \mu_n , h_n}$, $\mu^\prime\! = \! \nu^\prime_n$, $\varepsilon\! = \! 13\varepsilon_n$ and $\eta\! = \eta_n$ as defined in (\ref{GHaussup}), and we get  
$$ \dePro \big( \ttP_{T^n_{\! \mu_n , h_n}} \mu_n , \ttP_{T_{\! \mu , h}} \mu \big) =    \dePro \big( \ttP_{T^n_{\! \mu_n , h_n}} \! \nu^\prime_n , \ttP_{T_{\! \mu , h}} \nu_n \big) \leq 26\varepsilon_n + 3\eta_n \xrightarrow[n\to \infty]{\;} 0 . $$ 
This proves (\ref{eraC0GHP}). \cq

\smallskip

It is convenient to prove $(iv)$ and $(v)$ next. 
 Recall from (\ref{ellspacedef}) the definition of $\fell_\bullet^{\uparrow}(\bbN)$ and recall from (\ref{ellcoordef}) the definition 
of the coordinate functions $L_p$, $p\ino \bbN\cup \{ \infty\}$. Then recall that if $\Lambda\ino \cM_f(\fell_\bullet^{\uparrow}(\bbN))$, then $\Lambda^{\! (p)}\! $ stands for the measure on $[0, \infty]\! \times \! \bbR_+^{p+1}\! $, which is the push-forward measure of $\Lambda$ via the function $(L_\infty, (L_q)_{0\leq q \leq p} )$.  Recall that $\flambda_\fnu \eqo \fLambda_{\bhh, \fnu} \circ L_\infty^{-1}$ for all $\fnu\ino \bbT$. Thus $(v)$ is an immediate consequence of $(iv)$. 

 To prove $(iv)$, \vspace{-0.1cm} we then note that we only need to show for all $p\ino \bbN$ that $\fLambda^{{\!(p)}}_{{\bhh,\fmu_n}} \!\!\!  \to \!  \fLambda^{{\!(p)}}_{{\bhh, \fmu}} $ weakly on $\cM_f (\bR_+ \! \times \! \bR_+^{{p+1}})$. 
To this end we fix $p\ino \bbN$ and we suppose that $h^*_p \geko 30 \epp_n$. 
We recall that $ \mu (T\backslash \tau_n) \! \vee \! \mu_n (T_n\backslash \tau_n^\prime) \! <\!  4\varepsilon_n $, 
and since $\nu^\prime_n \! = \! \mathtt{P}_{\!\tau^\prime_n}\mu_n$ and $\nu_n \! = \! \ttP_{\tau_n} \mu$, 
Proposition \ref{Lambdaprop} $(iii)$ asserts for all Borel subsets $B$ of 
$\bR_+ \! \times \! \bR_+^{{p+1}}$ that 
$|\fLambda^{{\! (p)}}_{{\bhh,\nu'_n}} (B)\! -\! \fLambda^{{\! (p)}}_{{\bhh, \mu_n}} (B)| \leqo 4\epp_n$ and 
$|\fLambda^{{\! (p)}}_{{\bhh, \nu_n}} (B)\! -\! \fLambda^{{\! (p)}}_{{\bhh, \mu}} (B)| \leqo 4\epp_n$. Thus, to prove $(iv)$, we only need to prove for all $p\ino \bbN$ that 
\begin{equation}
\label{Lambda1un}
\lim_{n\to \infty} \dPro^{(p)} \big(\Lambda^{{\!(p)}}_{\bhh, \nu'_n},   \Lambda^{{\!(p)}}_{\bhh, \nu_n} \big) =0, 
\end{equation}
where $\dPro^{(p)}$ stands for the Prokhorov distance in $\cM_f (\bR_+ \! \times \! \bR_+^{{p+1}})$, when 
$\bR_+ \! \times \! \bR_+^{{p+1}}$ is equipped with distance associated with the max-norm. 
To prove (\ref{Lambda1un}), we set 
$$ \eta_n(p)= 52 \epp_n + \max_{0\leq q\leq p} \deHaus \big( T_{\! \mu , h^*_q-26\varepsilon_n} ,  T_{\! \mu , h^*_q+26\varepsilon_n} \big) \; $$
which tends to $0$ as $n\! \to \! \infty$ since $h\mapsto T_{\! \mu, h}$ is $\deHaus$-continuous at each $h^*_q$. 
To simplify notation, we set 
$\ttt_{n, q} \eqo T^{n}_{\! \mu_n, h^*_q}$ and $\ttt_{q} \eqo T_{\! \mu, h^*_q}$. 
By Proposition \ref{properased} $(i)$, $\ttt_{n,q}\eqo T^{n}_{\! \nu'_n, h^*_q}$ and $\ttt_{q} \eqo T_{\! \nu_n, h^*_q}$.  
By (\ref{delproproj}) and (\ref{GHaussup}) applied to $h_n=h\eqo h^*_q$, for all $0\leqo q \leqo p$, we also obtain the following. 
\begin{equation}
\label{reGHaussu}
\dePro (\nu_n, \nu'_n) \leko 13 \epp_n \quad \textrm{and} \quad \max_{0\leq q\leq p} \deHaus(\ttt_{n,q}, \ttt_q) \leko \eta_n (p)  . 
\end{equation}
We extend $f_{\ttt_q} $ (resp.~$f_{\ttt_{n, q}}$) to $E$ by setting $f_{\ttt_q} (\sigma) \eqo \rho$ if $\sigma \ino E\backslash T$ (resp.~$f_{\ttt_{n,q}}(\sigma) \eqo \rho_n$ if 
$\sigma \ino E\backslash T_n$). Then for all $\sigma \ino E$,
we set 
$$\phi_{p} (\sigma) \eqo \big( \delta (\rho , \sigma), (\delta (\rho, f_{\ttt_{q}} (\sigma)))_{0\leq q \leq p}\big)\quad \textrm{and} \quad 
\phi_{n,p} (\sigma) \eqo \big( \delta (\rho , \sigma), (\delta (\rho, f_{\ttt_{n, q}} (\sigma)))_{0\leq q \leq p} \big)\; .$$ 
By definition, $\fLambda^{{\! (p)}}_{{\bhh,\nu_n}} \eqo \nu_n \! \circ\!  \phi_p^{-1}$ and that 
$\fLambda^{{\! (p)}}_{{\bhh, \nu'_n}} \eqo \nu_n' \! \circ \! \phi_{n,p}^{-1}$. Since $\dePro (\nu_n, \nu'_n) \leko 13 \epp_n$, 
and by Strassen's theorem, there is a $(13\epp_n)$-coupling $\frak{m}_n$ of $\nu_n$ and $\nu_n'$ 
(see Definition \ref{couplingdef}) that satisfies $\frak{m}_n(E^{2}\backslash (\tau_n \! \times \! \tau_n'))\eqo 0$ 
and for 
$\frak{m}_n$-almost all $(\sigma, \sigma')\ino  \tau_n \! \times \! \tau_n'$, $\delta(\sigma, \sigma')\leqo 13\epp_n$.  
We then apply Lemma \ref{GrC0proj} $(i)$ with $\epp\eqo 13\epp_n$, $\eta\eqo \eta_n (p)$, $T\eqo \tau_n$, 
$T'\eqo \tau_n'$, and successively $\tau\eqo \ttt_q$ and $\tau'\eqo \ttt_{n,q}$, and we get 
$\delta(f_{\ttt_q} (\sigma), f_{\ttt_{n,q} }(\sigma') ) \leqo 26 \epp_n + 3 \eta_n (p)$, for all $q\ino \{ 0, \ldots, p\}$ and all $(\sigma, \sigma')\ino  \tau_n \! \times \! \tau_n'$ such that $\delta(\sigma, \sigma')\leqo 13\epp_n$.  
We also recall that $\delta (\rho, \rho_n) \leko 5\epp_n$. This proves that for $\frak{m}_n$-almost all $(\sigma, \sigma')\ino  \tau_n \! \times \! \tau_n'$,
\begin{equation}\label{eq:strassenbound} \big| \delta (\rho, \sigma) \! -\! \delta (\rho_n, \sigma')  \big| \vee \max_{0\leq q \leq p}  \big| \delta (\rho, f_{\ttt_q}(\sigma)) \! -\! \delta (\rho_n, f_{\ttt_{n,q}}(\sigma'))  \big| \leqo 31 \epp_n +   3 \eta_n (p) \; .
\end{equation}
For all $(\sigma, \sigma')\ino E^2$, we then set $\psi_{n,p}(\sigma,\sigma') \eqo (\phi_p(\sigma), \phi_{n,p} (\sigma'))$ and we denote by $\pi_n$ the measure on $(\bR_+\!\!  \times \! \bR^{{p+1}}_{+})^2$ that is the push-forward measure \vspace{-0.1cm} of $\frak{m}_n$ via $\psi_{n,p}$: namely $\pi_n\eqo \frak{m}_n \circ \psi_{n,p}^{-1}$. Then 
$\pi_n$ is a 
$(31\epp_n + 3\eta_n (p))$-coupling of $\fLambda^{{\! (p)}}_{{\bhh ,\nu_n} } $ and $\fLambda^{{\! (p)}}_{{\bhh, \nu'_n}} $, and by \eqref{eq:strassenbound}, Strassen's theorem applies and yields (\ref{Lambda1un}), which therefore completes the proof of $(iv)$. \cq

\smallskip

Let us prove (\ref{projC0GHP}). We recall that Proposition \ref{properased} $(iv)$ asserts that 
$\supp \era_{h_n} \mu_n \! \subseteq \! T^n_{\! \mu_n , h_n}$ and 
we recall (\ref{GHaussup}). Since $\sup_{n\in \bN} \mu_n (T_n) \! < \! \infty$, 
Lemma \ref{Haustight} applies and asserts that the measures $\era_{h_n} \mu_n $ are tight on $(E, d)$ so we only need to prove that $\era_h \mu$ is their only weak limit. Namely, let $\pi \! \in \! \cM_f (T)$ and 
$(n_k)_{k\in \bN}$ be a strictly increasing sequence of integers such that $\era_{h_{n_k}} \mu_{n_k} \!\!  \rightarrow \! \pi$ weakly in $\cM_f(E)$. We want to prove that $\pi \! = \! \era_h \mu$, necessarily. To this end, we set for any $k \! \in \! \bN$, 
$$ s_k \! :=\! 13 \varepsilon_{n_k}+ |h\! -\! h_{n_k}|+ \dePro ( \era_{h_{n_k}} \mu_{n_k}, \pi) + 
\deHaus (T^{n_k}_{\! \mu_{n_k} , h_{n_k}} , T_{\! \mu, h}) \xrightarrow[k\to \infty]{\;} 0 \, .$$  
We also introduce the following notation: for all $\sigma \ino T$ and all $r \ino \bbR_+$, we denote by $\sigma (r)$ the unique $\gamma \! \in \! \lgeo \rho, \sigma \rgeo_T$ such that $d(\sigma, \gamma) \! = \! r \! \wedge \! d(\rho, \sigma)$. Here the notation $\theta_\sigma T$ still means $\{ \sigma'\ino T\! : \! \sigma \ino \lgeo \rho, \sigma'\rgeo_T\}$. We first claim that 
\begin{equation}
\label{subsupp}
\supp \pi \! \subseteq \! T_{\! \mu, h} \quad \textrm{and} \quad \forall \sigma \in   T\backslash T_{\! \mu, h}, \quad  \pi (\theta_\sigma T)\! = \! \era_h \mu (\theta_\sigma T)= 0\, .
\end{equation}
\textit{Proof of (\ref{subsupp}).} Observe that for all $n\! \in \! \bN$ and all $\varepsilon \! \in \! \bbR^*_+$ such that $\deHaus (T^n_{\! \mu_n , h_n}, T_{\! \mu, h})\! < \! \varepsilon$, we get $\era_{h_n} \mu_n (E \backslash T_{{\mu, h}}^{{(\varepsilon)}}) \! = \! 0$. The Portmanteau theorem implies that $\pi (  E \backslash T_{{\mu, h}}^{{(\varepsilon)}})\! = \! 0$, for all $\varepsilon\! \in \! \bbR^*_+$, which entails $\supp \pi \! \subseteq \! T_{\! \mu, h}$. Since $\theta_\sigma T \! \subseteq \!  T\backslash T_{\! \mu, h}$ for all $\sigma \ino T\backslash T_{\! \mu, h}$, the proof of (\ref{subsupp}) is complete. \cq 

\smallskip

For all $\sigma\! \in \! T_{\! \mu, h} $ such that 
$\era_h \mu (\theta_\sigma T)\! >\! 0$, we next prove that 
\begin{equation}
\label{eratopi}
\forall r \in (0, \era_h \mu (\theta_\sigma T)), \quad  \era_h \mu (\theta_\sigma T)\leq  \pi \big( \theta_{\sigma (r)} T\big) + r. 
\end{equation}  
\textit{Proof of (\ref{eratopi}).} Let $k$ be such that $8 s_k \! \leq \! r < \! \era_h \mu (\theta_\sigma T)$. Then 
\begin{eqnarray}
8s_k \leko \era_h \mu (\theta_\sigma T) \!\!\!\! &= &\!\!\!\!  \era_h \nu_{n_k} (\theta_\sigma T) \quad \textrm{(by Prop.~\ref{properased} $(i)$, since $\nu_{n_k}\! = \! \ttP_{\tau_{n_k}} \mu$ and $h\!> \! \mu (T\backslash \tau_{n_k})$)} \nonumber \\
\!\!\!\! & = &\!\!\!\!  \nu_{n_k} (\theta_\sigma T) \! -\! h =  \nu_{n_k} (\theta_\sigma \tau_{n_k}) \! -\! h \nonumber  \\
\!\!\!\! & \leq &\!\!\!\!  \label{glurnsi} \nu^\prime_{n_k}\!  \big( (\theta_\sigma \tau_{n_k})^{(s_k)}\big) + s_k \! -\! h \!= \!   \nu^\prime_{n_k} \! \big( \tau^\prime_{n_k} \!  \cap \!  (\theta_\sigma \tau_{n_k})^{(s_k)}\big) + s_k \! -\! h,  
\end{eqnarray}
since $s_k \! >\! 13\varepsilon_{n_k} \! >\! \dePro (\nu_{n_k} , \nu^\prime_{n_k})$ by (\ref{delproproj}).  
We next observe that $h\!  -\!  s_k+ \era_h \mu (\theta_\sigma T)  $ $ >$ $ h+ 7s_k $ $ >$ $0$. By (\ref{glurnsi})
we get $\tau^\prime_{n_k}\!  \cap (\theta_\sigma \tau_{n_k})^{(s_k)}\! \neq \! \emptyset $. We then recall from (\ref{TTen}) that 
$\deHaus(\tau_{n_k}, \tau'_{n_k})\leko 5\epp_{n_k} \leko s_k$. Therefore (\ref{frounssi}) in Lemma \ref{Gesubtr} applies with $T\! = \! \tau_{n_k}$, $T^\prime\! = \! \tau^\prime_{n_k}$  
and $\varepsilon\! = \! s_k$: there exists $\sigma^\prime\! \in \! \tau^\prime_{n_k}$ such that $\delta (\sigma, \sigma^\prime) \! \leq \! 4 s_k$ and $ \tau^\prime_{n_k}\!  \cap (\theta_\sigma \tau_{n_k})^{(s_k)}\! \subseteq  \! \theta_{\sigma^\prime} \tau^\prime_{n_k} $. Consequently, 
\begin{eqnarray*}
 8s_k  \leko \era_h \mu (\theta_\sigma T) \!\!\!\! &\leq  &\!\!\!
\nu^\prime_{n_k} \! \big(\theta_{\sigma^\prime} \tau^\prime_{n_k}  \big) + s_k \! -\! h  =  \nu^\prime_{n_k} \! \big(\theta_{\sigma^\prime} T_{n_k}  \big) + s_k \! -\! h \\
 & =  &  \!\!\!\mu_{n_k} \! \big(\theta_{\sigma^\prime} T_{n_k}  \big) + s_k \! -\! h 
  \quad \textrm{(by Prop.~\ref{projdef} since $\ttP_{\tau_{n_k}^\prime} \mu_{n_k}\! = \! \nu^\prime_{n_k}$ and since  $\sigma^\prime \! \in\! \tau^\prime_{n_k}$)}   \\
& \leq &  \!\!\!  \mu_{n_k} \! \big(\theta_{\sigma^\prime} T_{n_k}  \big) \! -\! h_{n_k}+ 2s_k  \quad \textrm{(since $|h_{n_k}\! -\! h| \! \leq \! s_k$)} \\
& \leq & \!\!\! \era_{h_{n_k}} \mu_{n_k}  (\theta_{\sigma^\prime} T_{n_k} ) + 2s_k \! = \!  
 \era_{h_{n_k}} \mu_{n_k}  (\theta_{\sigma^\prime} T^{n_k}_{\! \mu_{n_k}, h_{n_k} } ) + 2s_k  \\
& \leq & \!\!\!\pi \big( (\theta_{\sigma^\prime} T^{n_k}_{\! \mu_{n_k}, h_{n_k}} )^{(s_k)} ) + 3s_k  \quad 
\textrm{(since $\dePro ( \era_{h_{n_k}} \mu_{n_k}, \pi) \! \leq \! s_k$)}  \\  
& = & \!\!\! \pi \big( T_{\! \mu, h} \! \cap \! (\theta_{\sigma^\prime} T^{n_k}_{\! \mu_{n_k}, h_{n_k}} )^{(s_k)} ) + 3s_k. \end{eqnarray*} 
Consequently, $T_{\! \mu, h} \cap (\theta_{\sigma^\prime} T^{n_k}_{\! \mu_{n_k}, h_{n_k}} )^{(s_k)} \! \neq \! \emptyset $. Since 
$\deHaus (T^{n_k}_{\! \mu_{n_k} , h_{n_k}} , T_{\! \mu, h})\leko s_k$, by definition of $s_k$, 
(\ref{frounssi}) in Lemma \ref{Gesubtr} applies with 
$T\! = \! T^{n_k}_{\! \mu_{n_k}, h_{n_k}}$, $T^\prime\! = \! T_{\! \mu, h}$ and $\varepsilon\! = \! s_k$: there is 
$\sigma^{\prime \prime} \! \! \in \! T_{\! \mu, h}$ such that $\delta (\sigma^\prime, \sigma^{\prime \prime}) \! \leq \! 4 s_k $ and $T_{\! \mu, h}  \cap (\theta_{\sigma^\prime} T^{n_k}_{\! \mu_{n_k}, h_{n_k}} )^{(s_k)} \! \subseteq \! \theta_{\sigma^{\prime \prime}} T_{\! \mu, h} $. Observe that $\delta (\sigma, \sigma^{\prime \prime}) \! \leq \! 8 s_k$, therefore $  \theta_{\sigma^{\prime \prime}} T_{\! \mu, h} \! \subseteq \! \theta_{\sigma (8s_k)} T_{\! \mu, h} \! \subseteq \! \theta_{\sigma (r)} T_{\! \mu, h} \! \subseteq \! \theta_{\sigma (r)} T $. Consequently,  
$$  \era_h \mu (\theta_\sigma T) \leq  \pi \big( T_{\! \mu, h} \! \cap \! (\theta_{\sigma^\prime} T_{\! \mu_{n_k}, h_{n_k}} )^{(s_k)} ) + 3s_k  \leq \pi \big( \theta_{\sigma (r)} T\big) + 3s_k  $$
which completes the proof of (\ref{eratopi}).   \cq 
 
\smallskip 

We next prove the following: for all $\sigma\! \in \! T_{\! \mu, h} $ such that 
$\pi(\theta_\sigma T)\! >\! 0$, we have 
\begin{equation}
\label{pitoera}
\forall r \in (0, \pi(\theta_\sigma T)), \quad \pi (\theta_\sigma T)\leq   \era_h \mu \big( \theta_{\sigma (r)} T\big) + r. 
\end{equation}  
\textit{Proof of (\ref{pitoera}).} Let $k$ be such that $8 s_k \! \leq \! r \! < \! \pi (\theta_\sigma T)$. 
By (\ref{subsupp}), this implies that $\sigma \! \in \! T_{\! \mu, h}$. Thus, $T_{\! \mu_, h} \cap \theta_\sigma T = \theta_\sigma T_{\! \mu_, h}$ and 
\begin{eqnarray*}
8s_k < \pi (\theta_\sigma T)  \!\!\! &= & \!\!\! \pi (\theta_\sigma T_{\! \mu_, h}) \\
 \!\!\! &\leq  & \!\!\! \era_{h_{n_k}}\mu_{n_k} \!  \big( (\theta_\sigma T_{\! \mu, h})^{(s_k)} \big) + s_k  \quad 
\textrm{(since $s_k \! > \!  \dePro ( \era_{h_{n_k}} \mu_{n_k}, \pi) $ )}  \\
 \!\!\! & \leq & \!\!\!  \era_{h_{n_k}}\mu_{n_k} \!  \big( T^{n_k}_{\! \mu_{n_k}, h_{n_k} } 
 \!\!  \cap \! (\theta_\sigma T_{\! \mu, h})^{(s_k)} \big) + s_k .
\end{eqnarray*}
This implies $ T^{{n_k}}_{{\!\mu_{n_k}, h_{n_k}} } \! \cap  (\theta_\sigma T_{\! \mu, h})^{(s_k)} \! \neq \! \emptyset$ and we can apply  (\ref{frounssi}) in Lemma \ref{Gesubtr} with $T\! = \! T_{\! \mu, h}$, $T^\prime\! = \! T^{{n_k}}_{{\!\mu_{n_k}, h_{n_k}} }$ and $\varepsilon\! = \! s_k$: there is $\sigma^{\prime } \! \! \in \! T^{{n_k}}_{{\!\mu_{n_k}, h_{n_k}} }$ such that $\delta (\sigma, \sigma^{\prime }) \! \leq \! 4 s_k $ and 
$ T^{{n_k}}_{{\!\mu_{n_k}, h_{n_k}} } \! \cap  (\theta_\sigma T_{\! \mu, h})^{(s_k)} \! \subseteq \! 
\theta_{\sigma^\prime} T^{{n_k}}_{{\!\mu_{n_k}, h_{n_k}} }$. Consequently, 
$$  8s_k <\pi(\theta_\sigma T) \! \leq  \!
 \era_{h_{n_k}}\mu_{n_k} \!  \big(  \theta_{\sigma^\prime} T^{{n_k}}_{{\!\mu_{n_k}, h_{n_k}} } \big) + s_k =    \era_{h_{n_k}}\nu^\prime_{n_k} \!  \big(  \theta_{\sigma^\prime}T^{{n_k}}_{{\!\mu_{n_k}, h_{n_k}} }\big) + s_k , $$
by Prop.~\ref{properased} $(i)$, since $\nu^\prime_{n_k}\! = \! \ttP_{\tau^\prime_{n_k}} \mu_{n_k}$ and $h_{n_k}\!> \! 29 \epp_{n_k} \geko 4\epp_{n_k} \geko \mu_{n_k} (T_{n_k}\backslash \tau^\prime_{n_k})$. Thus, we get:  
\begin{eqnarray*}
 8s_k < \pi (\theta_\sigma T) \!\!\!\! &\leq  &\!\!\! 
\nu^\prime_{n_k} \!  \big(  \theta_{\sigma^\prime} \tau^\prime_{n_k} \big) -h_{n_k} + s_k \\
& \leq & \nu_{n_k} \!  \big(  (\theta_{\sigma^\prime} \tau^\prime_{n_k})^{(s_k)} \big)-h_{n_k} + 2s_k
\quad \textrm{(since $s_k \! >\! 13\varepsilon_{n_k} \! >\! \dePro (\nu_{n_k} , \nu^\prime_{n_k})$)} \\
& \leq & \nu_{n_k} \!  \big(  \tau_{n_k} \! \cap \! (\theta_{\sigma^\prime} \tau^\prime_{n_k})^{(s_k)} \big)-h_{n_k} + 2s_k .
 \end{eqnarray*} 
This shows $ \tau_{n_k} \! \cap \! (\theta_{\sigma^\prime} \tau^\prime_{n_k})^{(s_k)} \! \neq \! \emptyset$ and we can apply   (\ref{frounssi}) in Lemma \ref{Gesubtr} with $T\! = \! \tau_{n_k}^\prime$, $T^\prime\! = \! \tau_{n_k}$ and $\varepsilon\! = \! s_k$: there exists $\sigma^{\prime \prime} \! \in \! \tau_{n_k}$ such that $\delta(\sigma^\prime, \sigma^{\prime \prime})\! \leq \! 4 s_k$ and 
$ \tau_{n_k} \! \cap \! (\theta_{\sigma^\prime} \tau^\prime_{n_k})^{(s_k)} \! \subseteq  \!  \theta_{\sigma^{\prime \prime}} \tau_{n_k}$. Thus,
 \begin{eqnarray*}
 8s_k < \pi  (\theta_\sigma T) \!\!\!\! &\leq  &\!\!\! \nu_{{n_k}} (\theta_{\sigma^{\prime \prime}} \tau_{n_k}) -h_{n_k} + 2s_k = \nu_{{n_k}} (\theta_{\sigma^{\prime \prime}}T)  -h_{n_k} + 2s_k   \\
 & \leq &  \mu (\theta_{\sigma^{\prime \prime}}T)  -h_{n_k} + 2s_k 
\; \,  \textrm{(by Prop.~\ref{projdef} $(vi)$ since $\ttP_{\tau_{n_k}} \mu\! = \! \nu_{n_k}$ and since 
$\sigma^{\prime \prime} \! \in\! \tau_{n_k}$)} \\
& \leq &  \mu (\theta_{\sigma^{\prime \prime}}T) -h + 3s_k \quad \textrm{(since $s_k \! >\! |h\! -\! 
h_{n_k}|$ ).} \\
& \leq & \era_h \mu (\theta_{\sigma^{\prime \prime}}T) + 3s_k .
\end{eqnarray*} 
Next observe that $\delta (\sigma, \sigma^{\prime \prime}) \! \leq \! 8 s_k$ and that 
$\theta_{\sigma^{\prime \prime}} T \! \subseteq \! \theta_{\sigma (8s_k)} T \! \subseteq \! \theta_{\sigma (r)} T \!  $. 
Thus, we have proved that $\pi (\theta_\sigma T) \leq   \era_h \mu ( \theta_{\sigma (r)} T) + 3s_k $, which completes 
the proof of (\ref{pitoera}). \cq

\smallskip

Next note that $r\! \mapsto\!  \theta_{\sigma (r)} T$ is non-increasing with respect to inclusion and that 
$\bigcap_{r\in \bbR^*_+} \theta_{\sigma (r)} T\! = \! \theta_\sigma T$. Thus, 
$\era_h \mu( \theta_{\sigma (r)} T) \! \! \rightarrow \! \era_h (\theta_\sigma T)$ and 
$\pi( \theta_{\sigma (r)} T) \!\!  \rightarrow \! \pi(\theta_\sigma T)$ as $r\! \! \rightarrow \! 0$. 
Then, (\ref{eratopi}) and (\ref{pitoera}) show the following: for all $\sigma \! \in \! T_{\! \mu, h}$, 
$\pi (\theta_\sigma T)\! >\! 0$ iff $\era_h \mu (\theta_\sigma T)\! >\! 0$, 
and in that case, $\pi (\theta_\sigma T)\! =\! \era_h \mu (\theta_\sigma T)$. Combined with (\ref{subsupp}), 
this proves that for all $\sigma \! \in \! T$, $\pi (\theta_\sigma T)\! =\! \era_h \mu (\theta_\sigma T)$, 
which entails $\pi\! = \! \era_h \mu$ by Proposition \ref{AboTrprop} $(vi)$ combined with standard arguments. 
We have thus proved that 
$\era_h \mu$ is the only possible weak limit of the tight measures $(\era_{h_n} \mu_n)_{n\in\bbN}$. 
This implies that $\era_{h_n} \mu_n \! \rightarrow \! \era_h \mu$ weakly as $n\! \rightarrow \! \infty$ and the proof 
of Theorem \ref{eraGcont} $(ii)$ is complete.

 Theorem \ref{eraGcont} $(i)$ can now be deduced by use of the uniform bound of Proposition \ref{properased} $(vi)$ and the fact that for all $h' \! \geq \! c\! := \! \mu (T)+ \sup_{n\in \bN} \mu_n (T_n)$, $\era_{h'} \mu\! =\! \era_{h'} \mu_n\! = \! 0$. \emph{Indeed}, let $\eta \ino (0, 1)$. Since the $\deHaus$-discontinuities of $h'\ino (0, c] \! \mapsto T_{\! \mu, h'}$ are at most countable, we can find positive real numbers 
 $h^\eta_1\leko h^\eta_2 \leko \cdots \leko h^\eta_N\eqo c+1$, 
 such that $N \leko 2c/\eta$ and such that 
 $h^\eta_{j+1} \! -\! h^\eta_j \leqo \eta$  and  
 $h'\ino (0, c] \! \mapsto T_{\! \mu, h'}$ is $\deHaus$-continuous at $h^\eta_j$, for all $j\ino \{ 1, \ldots, N\! -\! 1\}$. 
 To avoid trivialities, let us suppose that $h, h'\ino (0, c]$ such that $h'\geqo h$. 
 Then there is $2\leko j  \leko N\! -\! 1$, such that 
 $h'\ino (h^\eta_{j-1}, h^\eta_{j}]$ and  Proposition \ref{properased} $(vi)$ implies $\dePro (\era_{h'} \mu_n ,   \era_{h^\eta_j} \mu_n)\leqo 3\eta \mu_n (T_n)/h'\leqo 3\eta c/h$ and  similarly $\dePro (\era_{h'} \mu ,   \era_{h^\eta_j} \mu)\leqo 3\eta 
 \mu (T)/h'\leqo 3\eta c/h$. Therefore 
 $\sup_{h'\in \bbR_+} \dePro \big( \era_{h'} \mu_n ,   \era_{h'} \mu \big) \leq 6\eta c h^{-1}+ \max_{1\leq j< N}  
 \dePro ( \era_{h^\eta_j} \mu_n, \era_{h^\eta_j} \mu) $. Since $\lim_{n\to \infty} \max_{1\leq j< N}  
 \dePro ( \era_{h^\eta_j} \mu_n, \era_{h^\eta_j} \mu)\eqo 0$ by (\ref{projC0GHP}), we get $\limsup_{n\to \infty}\sup_{h'\in [h, \infty)} \dePro \big( \era_{h'} \mu_n ,   \era_{h'} \mu \big) \leqo \eta$, which implies Theorem \ref{eraGcont} $(i)$ since $\eta$ is arbitrarily small and since $ \dGP (\era_{h'} \fmu_n ,   \era_{h'} \fmu ) \leqo  \dePro ( \era_{h'} \mu_n ,   \era_{h'} \mu )$.

 To prove Theorem \ref{eraGcont} $(iii)$ observe that for all $0\! <\! h_0 \! <\! h_1$, 
\begin{eqnarray*}
\sup_{h, h^\prime \in [h_0, h_1]} \Big( 
\deHaus (T_{\! \mu_n, h}, T_{\! \mu_n, h^\prime} ) \!\!\! &+ &\!\!\!   \dePro 
(\ttP_{T_{\! \mu_n, h}}\mu_n , \ttP_{T_{\! \mu_n, h^\prime}} \mu_n ) \Big) \\
 \leq  \!\!\!  & &  \!\!\!  \!\!\!  \!\!\! 
\deHaus (T_{\! \mu_n, h_0}, T_{\! \mu_n, h_1} )+ \dePro 
(\ttP_{T_{\! \mu_n, h_0}}\mu_n , \ttP_{T_{\! \mu_n, h_1}} \mu_n )
\end{eqnarray*}
by definition of $\deHaus$ and by (\ref{projprop}) in Lemma \ref{projprok}. An 
argument of uniformity, similar to the one used in the proof of Theorem \ref{eraGcont} $(i)$,  
completes the proof Theorem \ref{eraGcont} $(iii)$.\cqfd     

\subsection{Convergence of $\bbR$-trees in the sense of mass erasure}
\label{Gromerasense}
In this section, we introduce the space $\bbT^*\!$ of measured $\bbR$-trees 
with boundary and we define on $\bbT^*$ 
the convergence in the sense of mass erasure, which is weaker than Gromov--Prokhorov convergence. 
We prove that $\bbT^*$ is the completion of $\bbT$ with respect to mass-erasure convergence, whose main properties are studied in Theorem \ref{cveracomplet}. 

\smallskip

\noi
\textbf{The space $\bbT^*\! $}. Before introducing the space of measured $\bbR$-trees with boundary, let us recall the following definitions and results. Let $(T,d,\rho)$ be a pointed Polish $\bbR$-tree. 

\smallskip

\noi
$-$ We recall from Theorem \ref{bordif} the definition of the bordification $(T^*\!, d^*\! , \rho)$ of $(T,d,\rho)$ which is the pointed Polish $\bbR$-tree obtained by the $g$-deformation of the metric $d$ (see Notation \ref{dgtreedef}) corresponding to the function $g(x)\eqo 1\! -\! e^{-x}\ino [0, 1)$, $x\ino \bbR_+$. The boundary of $T$ is then given by 
$\partial T\eqo \{ \mathbf s\ino T^*  :  d^*(\rho, \mathbf s)\eqo 1 \}$.

\smallskip

\noi
$-$ Let $\mu\ino \cM_f(T)$ be minimal, i.e.~$\sspaan \, (\supp \mu)\eqo T$. Thus $(T,d,\rho, \mu) \! \equiv :\! \fmu \ino \bbT$. 
We recall Remark \ref{toporem} $\textbf{(c)}$ and Lemma \ref{embedmeas}, which allow $\cM_{\! f }(T)$ equipped with the weak topology on $(T,d)$ to be viewed as a topologically embedded subspace of $\cM_{\! f }(T^*)$ equipped with the weak topology on $(T^*\! , d^*)$. Therefore, the isometry class of $(T^*\!, d^*\! , \rho, \mu)$, which only depends on $\fmu$ and which is denoted by $\jmath_*(\fmu)$  belongs to the subset $\bbT(1)$ of $\bbT$ defined by 
$$\bbT(1)\eqo \Big\{ (T'\!, d'\!, \rho'\!, \mu')\! \equiv\fmu'\ino \bbT:\,  \mathtt{Ht}_{d'\! , \rho'} (T')\leqo 1 \; \textrm{and} \; \mu \big( \{ \sigma\ino T'\!\! : d'(\rho', \sigma)\eqo 1\} \big)\eqo 0\Big\} $$ 
(see Definition \ref{Tadef}). By Lemma \ref{jgtopoembed}, 
$j_*$ is a $\dGP$-homeomorphism from $\bbT$ onto $\bbT(1)$.

\smallskip

\noi
$-$ We want to equip $T^*$ with finite Borel measures that possibly charge the boundary. More precisely, we 
recall that $\cM^{\mathtt{era}}_{\! f} (T)\eqo \{ \nu \ino \cM_f(T^*): \forall \mathbf s\ino \partial T, \, \nu (\{ \mathbf s\})\eqo 0 \} $ (see  (\ref{meraspace})). 
\begin{definition}
\label{Treestar} We keep the previous notation and we intoduce the following. 
\begin{compactenum}

\smallskip

\item[$(a)$] $(T,d,\rho, \mu)$ is a $\ast$-$\bbR$-\emph{tree} if $(T,d,\rho)$ is a rooted Polish $\bbR$-tree and if $\mu\ino \cM^{\mathtt{era}}_{\! f} (T)$ such that $(T^*\! , d^*\! , \rho, \mu)$ is minimal, i.e.~$T^*\! \eqo \sspaan_* (\mu) $.

\smallskip

\item[$(b)$] Two $\ast$-$\bbR$-trees $(T_i,d_i,\rho_i, \mu_i)$, $i\ino \{ 1,2\}$ 
are \emph{$\ast$-equivalent} if their respective bordifications $(T^*_i,d^*_i,\rho_i, \mu_i)$, $i\ino \{ 1,2\}$ are isometric.

\smallskip

\item[$(c)$] We denote by $\bbT^*\! $ the space of $\ast$-equivalence classes of $\ast$-$\bbR$-trees. We shall 
view $\bbT$ as a subset of $\bbT^*\! $ as explained in Remark \ref{TinTstar} below. \cq 
\end{compactenum}
\end{definition}

\begin{remark} 
\label{TinTstar} Let us justify why $\bbT$ can be considered as a subset of $\bbT^*$. 
Let $(T,d,\rho, \mu) \! \equiv \! \fmu\ino \bbT$. Then $(T,d,\rho, \mu)$ is a $\ast$-$\bbR$-tree and we denote by $\imath (\fmu)$ its $\ast$-equivalence class in $\bbT^*$ (\emph{indeed}, it is easy to check that it only depends on the isometry class of $(T, d,\rho, \mu)$). 
If $(T'\! ,d'\!,\rho'\!, \mu') \! \equiv \! \fmu'\ino \bbT$ is such that $\imath (\fmu')\eqo \imath (\fmu)$, then by definition of $\ast$-equivalence $\jmath_*(\fmu')\eqo \jmath_* (\fmu)$ and we get $\fmu'\eqo \fmu$ since $\jmath_*$ is injective. Thus $\imath \colon  \bbT \hookrightarrow \! \bbT^*$ is injective and we identify $\bbT$ with $\imath (\bbT)$. \cq
\end{remark}
We shall write $(T,d, \rho, \mu)\! \equiv\! \fmu \ino \bbT^*\! $ to mean that $(T,d, \rho, \mu)$ is a representative of the $\ast$-equivalence class $\fmu$. We next recall from (\ref{massheight}) the definition of $\langle \fmu \rangle$ and $\mathtt{Ht} (\fmu) $ and we extend $\jmath_*$ on $\bbT^*$ as follows. 
\begin{definition}
\label{masserastardef} We keep the previous notation and we intoduce the following. 
\begin{compactenum}

\smallskip

\item[(a)]  We denote by $\bbT^*(1)$ the set of $(T'\! ,d'\!,\rho'\!, \mu') \! \equiv \! \fmu'\ino \bbT$ such that $\mathtt{Ht}_{d'\! , \rho'\! } (T') \leqo 1$ and $\mu'(\{\sigma \})\eqo 0$ for all $\sigma \ino T'$ such that $d'(\rho', \sigma) \eqo 1$. Note that $\bbT(1) \subo \bbT^*(1)\subo \bbT$.   

\smallskip

\item[(b)] We extend $\jmath_*$ to $\bbT^*$ as follows: for all $(T,d,\rho, \mu)\! \equiv \! \fmu \ino \bbT^*$, we denote by $\jmath_* (\fmu)$ the isometry class of $(T^*\! , d^*\! , \rho, \mu)$ (by definition, it only depends of the $\ast$-equivalence class of $(T,d,\rho, \mu)$). We see that $\jmath_* $ is a bijection from $\bbT^*$ onto $\bbT^* (1)$.  
\end{compactenum}
\end{definition}
\begin{remark}
\label{noGPonTstar} We do not equip $\bbT^*$ with the topology induced $\dGP$ on $\bbT^* (1)$ via $\jmath_*$. 
We rather endow $\bbT^*$ with the topology of mass erasure which is introduced below. \cq
\end{remark}

\noi
\textbf{Mass erasure on $\bbT^*\!$, projective properties.}  Let us extend mass erasure to $\bbT^*$ as follows. 
\begin{definition}
\label{masserastardef} 
Let $h\ino \bbR_+^*$ and let $(T,d,\rho, \mu)\! \equiv \! \fmu \ino \bbT^*$. By Lemma \ref{descentera} $(i)$, 
$(T_{\! \mu, h+}, d, \rho, \era_h \mu)$ is a rooted  finitely measured minimal compact $\bbR$-tree whose isometry class 
only depends on the $\ast$-equivalence class of $(T,d,\rho, \mu)$ (see Remark \ref{justiferastar} below). 
We denote it by $\cE_h \fmu $ and 
it  belongs to $\bbT \subo \bbT^*$. The semigroup property (\ref{GPmerasemi}) on $\bbT$ trivally 
extends to $\bbT^*$. \cq 

\end{definition}
\begin{remark}
\label{justiferastar} Let $h\ino \bbR^*_+$ and $(T,d,\rho, \mu)\! \equiv \! \fmu \ino \bbT^*$. Then $\jmath_*(\fmu) \ino \bbT^*(1) \subo \bbT$, and $\cE_h (\jmath_*(\fmu)) \! \equiv \! \big( (T^* )_{ \mu, h+}, d^*\! , \rho, \cE_h \mu\big)$ is well-defined as an element of $\bbT$ (and more precisely as an element of $\bbT(1) $ or $\bbT^1_{\! c}(1)$). Indeed, $ (T^*)_{\mu, h+}$ is a  $d^*$-compact subtree such that $\mathtt{Ht}_{d^*\! , \, \rho} ( (T^*)_{\mu, h+} ) \leko 1$. Thus 
$ (T^*)_{ \mu, h+} \!\subseteq\! T$ and by Theorem \ref{bordif}, it is a $d$-compact subtree of $T$. According to Lemma \ref{descentera} $(i)$, we simply denote $ (T^*)_{\mu, h+}$ by $ T_{\! \mu, h+}$. Namely $\cE_h (\jmath_*(\fmu))$ is the isometry class of $\big( T_{\! \mu, h+}, d^*\! , \rho, \cE_h \mu\big)$, which belongs to $\bbT(1)$. Moreover, by Lemma \ref{jgtopoembed}, there is a unique $\fnu\ino \bbT$ such that $\jmath_* (\fnu)\eqo \cE_h (\jmath_*(\fmu))$ and $\fnu$ is the isometry class of 
$(T_{\! \mu, h+}, d , \rho, \cE_h \mu\big)$. This shows that 
\begin{equation}
\label{starinj}
\forall \fmu\ino \bbT^*, \forall h\ino \bbR^*_+, \quad \textrm{$\cE_h \fmu$ is the unique $\fnu\ino \bbT$ such that $\cE_h (\jmath_*(\fmu))\eqo \jmath_* (\fnu)$, }
\end{equation}
which justifies Definition \ref{masserastardef}. \cq
\end{remark}

\begin{lemma}
\label{sepera} 
Let $\fmu, \fmu' \ino \bT^*$ be such that $\era_{h} \fmu\eqo \era_{h} \fmu'$ for all $h\ino \bbR^*_+$. Then, $\fmu\eqo \fmu'$. 
\end{lemma}
\noi
\textbf{Proof.} As a consequence of Proposition \ref{weakcvera}, we first see that for all $\fnu\ino \bbT$, $\lim_{h\to 0}\dGP(\cE_h \fnu, \fnu)\eqo 0$. Then, by (\ref{starinj}), we see 
$\cE_h (\jmath_*(\fmu))\eqo \cE_h (\jmath_*(\fmu'))\ino \bbT^*(1) \subo \bbT$, 
for all $h\ino \bbR^*_+$. Therefore we get $\jmath_*(\fmu)\eqo \jmath_*(\fmu')$, which in turns implies $\fmu\eqo \fmu'$ since $\jmath_*$ is a 
bijection from $\bbT^*$ onto $\bbT^*(1)$. \cqfd 
\begin{lemma}  
\label{Geraconsist} Let $((T_p, d_p , \rho_p , \nu_p) \!\equiv \! \fnu_p)_{p\in \bbN}$ be $\bbT^*$-valued, and let $(h_p)_{p\in \bbN}$ be $\bbR_+^*$-valued and strictly decreasing to $0$. 
We assume that $\era_{h_p -h_{p+1}} \fnu_{p+1} \! = \! \fnu_p$, $p\ino \bbN$. 
Then, there exists a unique $\fmu \! \in \! \bbT^*$ such that $\era_{h_p} \fmu \! = \! \fnu_p$.   
\end{lemma} 
\noi
\textbf{Proof.} Let us prove the uniqueness first: if $\fmu, \fmu'  \! \in \! \bbT^*$ are such that $\era_{h_p} \fmu \! = \! \fnu_p\eqo \era_{h_p} \fmu' $, $p\ino\bbN$, then the semigroup property of mass erasure entails that $\era_{h} \fmu \! = \! \eqo \era_{h} \fmu' $ for all $h\ino \bbR^*_+$ and Lemma \ref{sepera} implies that $\fmu\eqo \fmu'$. 
We now prove existence. For all $p\ino\bN$, there is a bijective isometry 
$$\phi_p\colon  (T_p, d_p , \rho_p , \nu_p)  \longrightarrow   ((T_{p+1})_{\nu_{p+1}, (h_p -h_{p+1})+ }, d_{p+1}, \rho_{p+1} , \era_{h_p -h_{p+1}} \nu_{p+1} ) \; .$$   
Then, we set $\rho\! := \! (\rho_p)_{p\in \bN} $, $M_p \! = \! \big\{ (\sigma_n)_{n\geq p}\! :\!  \sigma_n \ino T_n \; \textrm{and} \; \phi_{n} (\sigma_n)\! = \! \sigma_{n+1} \; \textrm{for all  $n \! \geq \! p$}\big\}$ and $M\! = \! \bigsqcup_{p \in \bN} M_p$. 
Let $\sigma\! = \!  (\sigma_n)_{n\geq p} \! \in \! M_{p}$ and $\sigma^\prime \! = \!  (\sigma^\prime_n)_{n\geq p^\prime} \! \in \! M_{p^\prime}$. 
Then  
$d_n (\sigma_n, \sigma^\prime_n)\! = \! d_{p\vee p^\prime} (\sigma_{p\vee p^\prime} , 
\sigma^\prime_{p\vee p^\prime})$ for all $n\geqo p\! \vee \! p'$. 
We denote this quantity by 
$d(\sigma, \sigma^\prime)$. Clearly, $d$ is a pseudometric on $M$ which satisfies the 
four-point condition (\ref{4ptscondi}) and $d(\sigma, \sigma^\prime)\! = \! 0$ iff one of the two 
sequences $\sigma$ or $\sigma^\prime$ extends the other: let us denote by $\sim$ this 
equivalence relation. We then set $T^o\! = \! M/\! \sim$ and we keep denoting by $d$ the 
quotient metric on $T^o$. We slightly abuse notation by confusing $\rho$ with its $\sim$-equivalence class. 
For all $\gamma\! \in \! T_p$, denote by $\psi_p (\gamma)$ the $\sim$-equivalence class of 
$(\gamma_{n})_{n\geq p}$ where $\gamma_p \! = \! \gamma$ and where $\phi_n (\gamma_n)\! = \! \gamma_{n+1}$ 
for all $n\! \geq \! p$. Then, $\psi_p\colon T_p \! \hookrightarrow \! T^o$ is an isometry such that $\psi_p (\rho_p)\! = \! \rho$. 
We set $T^\prime_p \! = \! \psi_p (T_p)$, which is (compact and) connected. 
Therefore $T^o\eqo \bigcup_{p\in \bbN} T'_p$ is connected and Theorem \ref{4ptsth} implies that 
$(T^o\! , d, \rho)$ is a rooted separable $\bR$-tree. 
Let $(T, d)$ be a completion of $(T^o\! , d)$ such that  
$T^o\!\subseteq\!  T$, which is a Polish $\bR$-tree (by Proposition \ref{spanreafo}). 
We set $\nu_p^\prime\! = \! \nu_p \circ \psi_p^{-1}$. 
Then $\era_{h_p-h_{p+1}} \nu_{p+1}^\prime\! = \! \nu_p^{\prime}$. By Lemma \ref{descentera} $(ii)$, there is a unique 
$\mu \ino \cM_f^{\mathtt{era}} (T)$ such that $\era_{h_p} \mu \! = \! \nu^\prime_p$. If $\fmu$ is the $\ast$-isometry class of 
$(\sspaan \, (\mathtt{\supp} \mu),d, \rho, \mu)$, then $\era_{h_p} \fmu \! = \! \fnu_p$, for all $p\ino \bbN$. \cqfd 

\medskip

\noi
\textbf{Convergence in the sense of mass erasure on $\bbT^*$.} As specified in Remark \ref{noGPonTstar}, 
rather than endowing $\bbT^*$ with $\dGP$-topology, we equip this space with the mass-erasure convergence that is defined as follows. 
\begin{definition}
\label{masseracvdef} Let $(\fmu_n)_{n\in \bN}$ be a $\bbT^*$-valued sequence.  We say that it converges \textit{in the sense of mass erasure} to $\fmu\! \in \! \bbT^*$ if for all $h\ino \bbR^*_+$, 
$\fdelta_{{\mathtt{GP}}} \big( \era_{h} \fmu, \era_{h} \fmu_n \big)\! \to \! 0$. \cq 
\end{definition} 

\begin{lemma}
\label{eracvprop} Let $(T,d,\rho, \mu)\! \equiv \! \fmu$ and $ (T_n,d_n,\rho_n, \mu_n) \! \equiv \! \fmu_n \ino \bT^*$, 
$n\ino \bbN$.  
\begin{compactenum}

\smallskip

\item[$(i)$] If $\fmu,\fmu_n\ino\bT$, $n\ino\bbN$, and $\dGP(\fmu, \fmu_n)\! \to \! 0$, then $\fmu_n \! \to \! \fmu$ in the sense of mass erasure.

\smallskip

\item[$(ii)$] Suppose that $\fmu_n \! \to \! \fmu$ in the sense of mass erasure. 
Let $h\ino \bbR^*_+$ and 
let $h_n \! \to h$. Then $\lim_{n\to \infty} \dGP (\era_{h_n} \fmu_n, \era_{h} \fmu )\eqo 0$.  
If we furthermore assume that 
$h'\! \mapsto T_{\! \mu, h'}$ is $\dHaus$-continuous at $h$. Then 
(\ref{projC0GHP}) in Theorem \ref{eraGcont} $(ii)$ holds true.  
\end{compactenum}
\end{lemma}
\noi
\textbf{Proof.} $(i)$ is a direct consequence of Theorem \ref{eraGcont} $(i)$. To prove $(ii)$, we fix $\epp \ino (0, \inf_{n\in \bbN} h_n)$. As $\fmu_n \! \to \! \fmu$ in the sense of mass erasure, $\lim_{n\to \infty}\dGP (\era_\epp \fmu_n, \era_\epp \fmu)$. We apply Theorem \ref{eraGcont} $(i)$ to $\era_\epp \fmu_n$, $\era_\epp \fmu$, $h_n\! -\epp$ and $h\! - \epp$: combined with the semigroup property, this entails $\dGP (\era_{h_n} \fmu_n, \era_{h} \fmu )\! \to \!  0$.
Since $T_{\! \era_\epp \mu, h-\epp}\eqo T_{\! \mu, h}$ (by the semigroup property), $h'\! \mapsto \!   T_{\! \era_\epp \mu, h'}$ is $d_{\mathtt{Haus}}$ continuous at $h\! -\! \epp$. We then apply Theorem \ref{eraGcont} $(ii)$ to $\era_\epp \fmu_n$, $\era_\epp \fmu$, $h_n\! -\! \epp$ and $h\! -\! \epp$, combined with the semigroup property, this implies (\ref{projC0GHP}).  \cqfd 
\begin{theorem}  
\label{cveracomplet} Let  $(h_p)_{p\in \bN}$ be $\bbR^*_+$-valued and strictly decreasing to $0$.  
We recall from (\ref{dderadef}) the definition of $\dera$. 
Then, the following holds true. 

\smallskip

\begin{compactenum}
\item[$(i)$] $ \fdelta_{{\mathtt{era}}}$ metrizes  the convergence in the sense of mass erasure on $\bbT^*$.

\smallskip

\item[$(ii)$] We recall $\jmath_*$ and $\bbT^*(1)$ from Definition \ref{masserastardef}. Then, $\jmath_*$ is a homeomorphism from $(\bbT^*\!,  \dera)$ onto $(\bbT^*(1), \dera)$.

\smallskip

\item[$(iii)$] The metric space $(\bbT^*\! , \dera)$ is Polish and $\bT$ is $\dera$-dense in $\bT^*$. Namely, $\bT^*\! $ is the $\dera$-completion of $(\bT, \dera)$.

\smallskip

\item[$(iv)$] Let us recall from (\ref{massheight}) that $ \langle \fmu \rangle $ stands for the total mass. Then $\fmu\ino \bbT^*\! \mapsto \! \langle\fmu\rangle\ino\bbR_+$ is $\dera$-continuous.

\item[$(v)$] Let us recall $\mathtt{Ht} ( \fmu )$ from 
(\ref{massheight}). Let $\emptyset \! \neq \! \cC\! \subseteq\!  \bbT^*$. Then, $\cC$ is precompact in $(\bbT^*, \dera)$ iff 
\begin{equation}
\label{totbound}
\sup \big\{  \langle \fmu \rangle  \, ; \, \fmu\ino \cC \big\} <  \infty \quad \;\,\textrm{and} \quad  \; \,   \sup \big\{ \mathtt{Ht} (\era_{h_p} \fmu) \, ; \, \fmu\ino \cC \big\} <  \infty \quad   \forall p\! \in \! \bN. 
\end{equation}
\item[$(vi)$]  Let $(T, d , \rho, \mu) \! \equiv \! \fmu \ino \bbT^*$ and 
$(T_n, d_n , \rho_n, \mu_n) \! \equiv \! \fmu_n\ino \bT$, $n\ino \bbN$ such that $\dera(\fmu_n, \fmu) \!\! \to \! 0$. 
Then $\fmu\ino \bbT$ (without necessarily implying $\dGP(\fmu_n, \fmu)\! \to \! 0$)  iff 
\begin{equation}
\label{inTcondi0}
\lim_{r\to \infty} \, \limsup_{p\to \infty} \, \sup_{n\in \bbN} \era_{h_p} \mu_n \big(T_n \backslash B_{T_n, d_n} (\rho_n, r)\big) = 0.  
\end{equation}
\item[$(vii)$] The sequence $(\fmu_n)_{n\in \bN}$ is $\dGP$-convergent iff
the set $\{ \fmu_n ; n\! \in \! \bN\}$ is $\dGP$-precompact and 
$(\fmu_n)_{n\in \bN}$ converges in the sense of mass erasure. In that case, the limits are the same. 
\end{compactenum}
\end{theorem}
\noi
\textbf{Proof.} Let us prove $(i)$: $\dera$ is clearly nonnegative symmetric and it satisfies the triangle inequality. 
If $\dera (\fmu, \fnu)\eqo 0$, then $\era_{h_p} \fmu\eqo \era_{h_p} \fnu$ for all $p\ino \bbN$ and thus 
$\era_h \fmu\eqo \era_h \fnu$ for all $h\ino \bbR^*_+$, by the semigroup property of mass erasure, and Lemma \ref{sepera} implies that $\fmu\eqo \fnu$. This proves that $\dera$ is a distance. 
If $\fmu_n \!\! \to \! \fmu$ in $\bbT^*$ in the sense of mass erasure, then we immediately get $\dera(\fmu_n, \fmu)\! \to \! 0$. Conversely, let us assume that $\dera(\fmu_n, \fmu)\! \to \! 0$. Let $h\ino \bbR^*_+$ and $p\ino \bbN$ such that $h_p\leko h$. We set $h'\eqo h\! -\! h_p$. 
Since $\lim_{n\to \infty}\dGP(\era_{h_p} \fmu_n, \era_{h_p}\fmu)\eqo 0$, Lemma \ref{eracvprop} $(i)$ 
implies that $\lim_{n\to \infty}\dGP(\era_{h'}(\era_{h_p} \fmu_n), \era_{h'} (\era_{h_p}\fmu))\eqo 0$, and thus 
$\lim_{n\to \infty}\dGP(\era_{h} \fmu_n,  \era_{h}\fmu)\eqo 0$ by the semigroup property, which completes the proof $(i)$. 

Let us prove $(ii)$. By definition $\dera(\fmu_n, \fmu)\!\! \to 0$ iff $\dGP(\cE_h\fmu_n, \cE_h\fmu)\!\! \to 0$ 
for all $h\ino \bbR^*_+$. By Lemma \ref{jgtopoembed}, this is equivalent to $\dGP\big( \jmath_* (\cE_h \fmu_n),
\jmath_* (\cE_h \fmu) \big) \!\! \to \! 0$ for all $h\ino \bbR^*_+$. Then we observe that 
$\jmath_* (\cE_h \fmu_n)\eqo \cE_{h} \jmath_* (\fmu_n)$ and  
$\jmath_* (\cE_h \fmu)\eqo \cE_{h} \jmath_* (\fmu)$ by (\ref{starinj}). Namely, $\dera(\fmu_n, \fmu)\!\! \to 0$ iff $\dera \big(\jmath_* (\fmu_n), \jmath_* (\fmu) \big)\!\! \to 0$, which entails $(ii)$. 

We next prove $(iii)$. 
Let $\fmu\ino \bbT^*$. Then $\era_{h^\prime}\fmu\in\bT$ for all $h^\prime\ino\bbR^*_+$ and Proposition \ref{weakcvera} entails that $\lim_{h\to 0}\dGP(\era_{h^\prime+h}\fmu,\era_{h^\prime}\fmu)\eqo 0$. 
Hence 
$\dera (\era_h \fmu,  \fmu) \!\! \to \! 0$. 
This proves that $\bT$ is $\dera$-dense in $\bT^*$. The separability of $(\bT, \dGP)$, combined with 
Lemma \ref{eracvprop} $(i)$, implies that $(\bT, \dera)$ is separable, and so is $( \bT^*, \dera)$. 
Let $(\fmu_n)_{n\in \bbN}$ be a $\dera$-Cauchy 
sequence in $\bT^*$. For all $p\ino \bbN$, $(\era_{h_p} \fmu_n)_{n\in \bbN}$ is a Cauchy 
sequence in the $\dGP$-complete space $\bT$ and there is $\fnu_p\ino \bT$ such that $\lim_{n\to \infty} 
\dGP( \era_{h_p} \fmu_n,  \fnu_p)\eqo 0$. 
The semigroup property and Theorem \ref{eraGcont} imply the $\dGP$-convergence 
$\era_{h_p} \fmu_n \eqo \era_{h_{p}-h_{p+1}} (\era_{h_{p+1}} \fmu_n) \!\!  \to \! \era_{h_p-h_{p+1}} \fnu_{p+1}$. 
Thus $ \era_{h_p-h_{p+1}} \fnu_{p+1}\eqo \fnu_p$ for all $p\ino \bbN$ and 
by Lemma \ref{Geraconsist} there is $\fmu\ino  \bbT^*$ such that $\era_{h_p} \fmu\eqo \fnu_p$. Namely, 
$\lim_{n\to \infty} \dGP( \era_{h_p} \fmu_n, \era_{h_p} \fmu)\eqo 0$ 
which implies $\lim_{n\to \infty}\dera(\fmu_n, \fmu)\! \to \! 0$ and this completes the proof of $(iii)$.
 
To prove $(iv)$, we first note that $\fnu \ino \bbT\! \mapsto \! \langle \fnu \rangle$ is $\dGP$-continuous and then we observe that 
$|\langle \fmu \rangle \! -\! \langle \era_h \fmu \rangle |\leqo h$ for all $\fmu \ino \bbT^*$. 
Let $(\fmu_n)_{n\in \bbN}$ be a $\bbT^*$-valued sequence that $\dera$-converges to $\fmu\ino \bbT^*$. Then 
$|\langle \fmu \rangle  \! -\! \langle  \fmu_n \rangle| \leqo 2h +|\langle \era_h \fmu \rangle  \! -\! \langle  \era_h \fmu_n \rangle |$  and thus $\limsup_{n\to \infty} |\langle \fmu \rangle  \! -\! \langle  \fmu_n \rangle |\leqo 2h$, which completes the proof of $(iv)$ since $h$ can be arbitrarily close to $0$. 

 Let us prove $(v)$. We first suppose (\ref{totbound}). We set $x\! = \! \sup_{\fmu\in \cC}  \langle \fmu \rangle $ and $a_p\! = \! \sup_{\fmu\in \cC}  \mathtt{Ht} (\era_{h_p} \fmu) $. Let $(T,d, \rho, \mu)\! \equiv \! \fmu \ino \cC$. We recall notation (\ref{entropie})
 and we note, for all $\epp\ino \bbR^*_+$ and all $p\ino \bbN$, that $N( T_{\! \mu , h_p} , \varepsilon) 
\leqo  a_p\varepsilon^{-1}  \# \mathtt{Lf} \big( T_{\! \mu , h_p} \big) \leqo x a_p/(\varepsilon h_p). $ We recall the definition of $\esstree \, ( \fmu, \varepsilon)$ from (\ref{spanescov}). We have proved, for all $\fmu \ino \mathscr C$, that 
$\esstree \, ( \era_{h_p} \fmu, \varepsilon) \leqo  N( T_{\! \mu , h_p} , \varepsilon) + \mathtt{Ht} (\era_{h_p} \fmu)\leqo  
\frac{x a_p}{\varepsilon h_p} +a_p$.
Since  $\sup_{\fmu \in \cC} \langle \era_{h_p} \fmu \rangle \! \leq \! x$, this implies that 
$\{  \era_{h_p} \fmu ; \fmu \! \in \! \cC \}$ is $\dGP$-precompact by Proposition \ref{spanprecp}. 
By a diagonal extraction argument, from any $\cC$-valued sequence we extract a sequence $(\fmu_n)_{n\in \bN}$ such that for all $p\! \in \! \bN$, $(\era_{h_p} \fmu_n)_{n\in \bN}$ is $\dGP$-convergent: $(\fmu_n)_{n\in \bN}$ is thus a $\dera$-Cauchy sequence, which in turn entails that $\cC$ is $\dera$-precompact since $(\bbT^*, \dera)$ is complete.

Conversely, let us assume that $\cC$ is $\dera$-precompact. This easily implies, for all $p\ino \bbN$, that  
$\{ \era_{h_p} \fmu\, ; \, \fmu \! \in \! \cC\} \!\subseteq\! \bbT$ is $\dGP$-precompact, which entails $\sup_{\fmu\in \cC}  
\langle\cE_{h_p} \fmu \rangle   \! <\!   \infty$ by Proposition \ref{spanprecp}. 
Consequently, we get $\sup_{\fmu\in \cC}  \langle \fmu \rangle   \! <\!   \infty$ because  
$\langle \fmu \rangle \! \leq \!  \langle \era_{h_p} \fmu \rangle + h_p$. 
Next, we fix $p\! \in \! \bN$ and set $\varepsilon\! = \! \frac{1}{2} (h_{p}\! -\! h_{p+1})$.   
By Proposition \ref{spanprecp} again, 
since $\{ \era_{h_{p+1}} \fmu; \fmu \! \in \! \cC\}$ is 
$\dGP$-precompact, we get $b_p\! =\! \sup_{\fmu\in \cC} \esstree\,  (\era_{h_{p+1}}\fmu, \varepsilon) \! <\!  \infty$. 
We fix $(T,d,\rho, \mu) \! \equiv \! \fmu \ino \cC$. Let $\rho \! \in \! \tau \! \subseteq \! T_{\! \mu, h_{p+1}}$ be such that $\era_{h_{p+1}} \mu (  T_{\! \mu, h_{p+1}}\backslash \tau) \! < \! \varepsilon$ and 
$N(\tau, \varepsilon) + \mathtt{Ht} (\tau) \! \leq \! 2b_p$. 
We observe that $T_{\! \mu,h_p+} \!\subseteq\! T_{\! \mu, h_{p+1}+ \varepsilon} \!\subseteq\!\tau$. 
Thus, $\mathtt{Ht} (\era_{h_p } \fmu) \! \leq \! \mathtt{Ht} (\tau)$. This implies $\mathtt{Ht} (\era_{h_p} \fmu) \! \leq \! 2b_p$, for all $\fmu \! \in \! \cC$. Namely, $\sup_{\fmu \in \cC} \mathtt{Ht} (\era_{h_p} \fmu)\! < \! \infty$, which entails (\ref{totbound}) and completes the proof of $(v)$.

We now prove $(vi)$. We first suppose (\ref{inTcondi0}). As a consequence of 
Theorem \ref{eraGcont} $(v)$, we note that $\fmu'\ino \bbT \mapsto \flambda_{\fmu'} \ino \cM_{\! f} (\bbR_+)$ is 
continuous when $\bbT$ is equipped with $\dGP$ and $ \cM_{\! f} (\bbR_+)$ with weak convergence. 
Therefore, for all $p\ino \bbN$, $\lim_{n\to \infty}\flambda_{\era_{h_p} \fmu_n} \eqo \flambda_{\era_{h_p} \fmu} $ 
weakly on $\cM_{\! f}(\bbR_+)$ and by the Portmanteau theorem, $\era_{h_p} \mu(T\backslash B_{T, d}(\rho, r))$ 
$\eqo$ $\lambda_{\era_{h_p} \fmu} ((r, \infty))$ $\leqo$ $ \liminf_{n\to \infty} \lambda_{\era_{h_p} \fmu_n} ((r, \infty))$
$ \eqo  $ $\liminf_{n\to \infty}  \era_{h_p} \mu_n(T_n\backslash B_{T_n, d_n}(\rho_n , r)) $ $\leqo$  $\sup_{n\in \bbN}  
\era_{h_p} \mu_n(T_n\backslash B_{T_n, d_n}(\rho_n , r))$. We denote by $(C^*_i)_{i\in I}$ the connected components of 
$T^*\backslash B_{T, d} (\rho, r)$. Then  
$\era_{h_p} \mu (T\backslash B_{T, d}(\rho, r))$ 
$ \eqo $  $\sum_{i\in I} ( \mu(C^*_i)\! -\! h_p)_+$, and 
by monotone convergence for sums, 
$\lim_{p\to \infty} \era_{h_p} \mu(T\backslash B_{T, d}(\rho, r))\eqo \mu(T^*\backslash B_{T, d} (\rho, r))$. Therefore, (\ref{inTcondi0}) and the previous bound imply that $\mu(\partial T)\eqo 0$. Namely $\fmu\ino \bbT$. 

  Conversely, let $\fmu \ino \bbT$. The previous arguments imply for Lebesgue-almost all $r\ino \bbR^*_+$ that 
\begin{eqnarray*}
\lim_{n\to \infty}   \era_{h_p} \mu_n \big( T_n\backslash B_{T_n, d_n}(\rho_n, r) \big)  \!\!\! &=& \!\!\!  \lim_{n\to \infty} \lambda_{\era_{h_p} \mu_n} \big( (r,\infty) \big) =\lambda_{\era_{h_p} \mu} \big( (r,\infty) \big) \\
\!\!\!  & = &\!\!\!  \era_{h_p} \mu \big( T\backslash B_{T, d})(\rho, r)\big) \leq \mu \big( T^*\backslash B_{T, d}(\rho, r)\big) \xrightarrow[r\to 0]{\;} 0 , 
\end{eqnarray*}
which in turn easily entails (\ref{inTcondi0}).

Let us prove the non-trivial implication of $(vii)$: suppose that $\fmu\ino \bbT^*$ is the $\dera$-limit of the $\bbT$-valued sequence $(\fmu_n)_{n\in \bbN}$ and let $\fmu^\prime\ino \bbT$ be a $\dGP$-limit point 
of $(\fmu_n)_{n\in \bbN}$. Theorem \ref{eraGcont} $(i)$ entails that $\era_h \fmu^\prime\! = \! \era_h \fmu$ for all $h\! \in \! \bbR^*_+$ and Lemma \ref{sepera} implies that $\fmu\eqo \fmu'$. 
Then $\fmu\ino \bbT$ and it is the only 
$\dGP$-limit point of $(\fmu_n)_{n\in \bN}$, which entails the desired result. \cqfd

\begin{remark}
\label{musequence} Theorem \ref{cveracomplet} $(i)$ and $(iii)$ can be rephrased as follows: from the 
metric point of view, an element $\fmu$ of $\bT^*$ \emph{is actually} the sequence 
$(\era_{h_p} \fmu)_{p\in \bbN}$. More precisely, let us equip $\bT$ with the topology induced by $\dGP$. 
On the space $\bbT^{\bbN}$ of $\bbT$-valued sequences, we define the distance 
$\dGP^{\bbN} \big( (\fmu_p)_{p\in \bbN}, (\fnu_p)_{p\in \bbN}\big) \eqo \sum_{p\in \bbN}2^{-p}\big(1\! \wedge\! \dGP (\fmu_p,  \fnu_p )  \big)$. Standard arguments imply that 
$(\bbT^{\bbN}, \dGP^{\bbN})$ is Polish and its topology is the product topology. We next introduce the subspace $\bT^{\oast}\eqo \big\{ (\fnu_p)_{p\in \bbN} \ino \bbT^{\bbN} \! : \era_{h_p-h_{p+1}} \fnu_{p+1}\eqo \fnu_p ,\,  p\ino \bbN\}$. Since $\era_h$ is $\dGP$-continuous (Theorem \ref{eraGcont} $(i)$), $\bT^\oast$ is a $\dGP^{\bbN}$-closed subset of $\bbT^{\bbN}$. We now define 
$\jmath\colon \bbT^*\!\!  \to \! \bbT^{\bbN}$, by setting $\jmath (\fmu)\eqo (\era_{h_p} \fmu)_{p\in \bbN}$. Then, Lemma \ref{Geraconsist} yields that $\jmath$ is a bijective isometry from $(\bbT^*\! , \dera)$ onto $(\bbT^\oast\! , \dGP^{\bbN})$. \cq 
\end{remark}

\subsection{Connections between the topologies induced by $\fdelta_{\mathtt{era}}$  $\dGP, \dGHP$ and $\dGH$}
\label{connectsubsec}
In Theorem \ref{GPcvvsera} we discuss conditions that imply $\dGP$-convergence from $\dera$-convergence. 
We recall Notation \ref{GPmeradef} $(b)$ of the $\bhh$-projected height measure $\fLambda_{ \bhh, \fmu}$ and we prove the following lemma which is an analogue of Proposition \ref{Lambdaprop} $(iv)$. 
\begin{lemma}
\label{LambdaGprop} Let $(T,d,\rho, \mu) \equiv \fmu\ino \bbT^*$ and $\fmu_n \ino \bbT^*$, $n\ino \bbN$, be such that 
$\dera(\fmu_n, \fmu) \!\! \to \! 0$. Let $\bhh\eqo (h_p)_{p\in \bbN}$ be $\bbR^*_+$-valued and strictly decreasing to $0$. We assume that $h\mapsto T_{\! \mu, h}$ is $\dHaus$-continuous at each of the $h_p$, $p\ino \bbN$. Then the following holds true.
\begin{compactenum}

\smallskip

\item[$(i)$] $\fLambda^{\! o}_{ \bhh, \fmu_n} \!\! \! \! \to  \fLambda^{\! o}_{ \bhh, \fmu}$ weakly on $\cM_{\! f} (\fell^{\uparrow}(\bbN))$. 

\smallskip

\item[$(ii)$] The measures $(\fLambda_{ \bhh, \fmu_n})_{n\in \bbN}$ are tight on 
$\cM_{\! f} (\fell_\bullet^{\uparrow}(\bbN))$.$\, $Moreover, if $\fLambda$ is a weak limit point, then $\fLambda^{\! o}\eqo  \fLambda^{\! o}_{\bhh, \fmu} $ and $\fLambda$-a.e.~$\sup_{p\in \bbN} L_p \leqo L_\infty$.   
\end{compactenum}
\end{lemma}
\noi
\textbf{Proof.} We denote by $\fLambda^{{\! o, p}}_{{\bhh, \fmu_n}}$ and $\fLambda^{{\! o, p}}_{{\bhh, \fmu}}$ 
the push-forward measures on $\bbR_+^{{p+1}}$ of resp.~$\fLambda^{{\! o}}_{{\bhh, \fmu_n}} $ and $\fLambda^{{\! o}}_{{\bhh, \fmu}} $ via the function $(L_q )_{0\leq q \leq p}$. To prove $(i)$, we only need to prove that $\fLambda^{{\! o, p}}_{{\bhh, \fmu_n}} \! \!\!\!  \to \! \fLambda^{{\! o, p}}_{{\bhh, \fmu}}$ weakly in $\cM_{\! f} (\bbR_+^{{p+1}})$. To this end, for all $n\ino \bbN$, we set 
$$\mu'_n \eqo \mathtt{P}_{T^{n}_{\! \mu_n, h_{p+1}}}  \mu_n \quad \textrm{and} \quad \mu'\eqo \mathtt{P}_{T_{\! \mu, h_{p+1}}} \mu $$ 
and for all $q\ino \{ 0, \ldots, p\}$, we also set $h'_q\eqo h_q\! -\! h_{p+1}$. We complete $(h'_q)_{0\leq q\leq p}$ into a the sequence $\bhh'$ which decreases to $0$ and such that $h\! \mapsto \! T_{\! \mu, h}$ is $\dHaus$-continuous at 
$h_{p+1}+ h'_q$, for all $q\geko p$.  
 Then (\ref{eraC0GHP}) in Theorem \ref{eraGcont} $(ii)$ applies and we get $\dGP(\fmu_n', \fmu') \! \to \! 0$. We then apply Theorem \ref{eraGcont} $(iv)$ to the $\fmu_n'$ to see that 
$\fLambda^{{\! o, p}}_{{\bhh, \fmu_n}}\!  \eqo \fLambda^{{\! o, p}}_{{\bhh', \fmu_n'}}\! $ converges to 
$\fLambda^{{\! o, p}}_{ ^{\bhh, \fmu}}\!  \eqo \fLambda^{{\! o, p}}_{{\bhh', \fmu'}}\! $ weakly on $\cM_{\! f}(\bbR_+^{{p+1}})$, which completes the proof of $(i)$. 
The proof of $(ii)$ is easily adapted from that of Proposition \ref{Lambdaprop} $(iv\textrm{-}b)$: we leave the details to the reader. \cqfd

\begin{theorem}
\label{GPcvvsera} Let $(T,d,\rho, \mu) \! \equiv\! \fmu\ino \bbT^*$. Let $(h_p)_{p\in \bbN}$ be any $\bbR_+^*$-valued sequence which strictly decreases to $0$ and such that 
$h\! \mapsto \! T_{\! \mu, h}$ is $\dHaus$-continuous at each $h_p$, $p\ino \bbN$ (there is always one).  
Let $(\fmu_n)_{n\in \bbN}$ be $\bbT$-valued and such that $\dera (\fmu_n,   \fmu)\!\! \to \! 0$. Then the following assertions are equivalent.
\begin{compactenum}

\smallskip

\item[$(a)$] $\fmu\ino \bbT$ and $\fdelta_{{\mathtt{GP}}} (\fmu_n, \fmu) \! \rightarrow \! 0$,

\smallskip

\item[$(b)$] $\flambda_{\fmu} \ino \cM_{\! f} (\bbR_+)$ and $\flambda_{\fmu_n} \!\!  \to \! \flambda_{\fmu}$ weakly in $\cM_{\! f} (\bbR_+)$.

\smallskip

\item[$(c)$] For all $\epp \ino \bbR^*_+$, $\lim_{p\to \infty} \sup_{n\in \bbN} \fLambda_{\bhh, \fmu_n} (\{ L_\infty \! -\! L_p \geqo \epp \}) \eqo 0$. 
\end{compactenum}  
\end{theorem}
\noi
\textbf{Proof.} 
Theorem \ref{eraGcont} $(v)$ 
show $(a)\Rightarrow (b)$. Let us suppose $(b)$ and prove $(b)\! \Rightarrow \! (c)$. 
We first prove that 
$\fLambda_{{\bhh, \fmu_n}} \!\!\! \to \! \fLambda_{{\bhh, \fmu}}$ weakly on $\cM_f(\fell_\bullet^{\uparrow}(\bbN))$. 
\emph{Indeed}, let $\fLambda$ be any weak limit point of 
$(\fLambda_{{\bhh, \fmu_n}})_{n\in \bbN}$ (there is at least one by Lemma \ref{LambdaGprop} $(ii)$). It must 
satisfy 
$\fLambda^{\! o}\! \eqo  \fLambda^{\! o}_{{\bhh, \fmu}}$ and $\fLambda$-a.e.~$\mathcal L\! :=\! \sup_{p\in \bbN} 
L_p \leqo L_\infty$ (also by Lemma \ref{LambdaGprop} $(ii)$). On one hand, since $\flambda_{\fmu_n} \! \eqo  
\fLambda_{{\bhh, \fmu_n}} \! \circ\!  L_\infty^{-1} $ and since 
$L_\infty$ is continuous, we get 
$\flambda_\fmu \eqo \fLambda \! \circ \! L_\infty^{-1}$. On the other hand, since $\fLambda^{\! o}\eqo 
\fLambda^{\! o}_{{\bhh, \fmu}}$, we get that 
$\fLambda \! \circ\!  \mathcal L^{{-1}}\eqo \flambda_\fmu$ and Lemma \ref{adhoclem} (for a deterministic 
$\fLambda$) 
implies $\fLambda$-a.e.~that $L_\infty\eqo \mathcal L \eqo \sup_{p\in \bbN} L_p$ and thus $\fLambda\eqo  \fLambda_{{\bhh, \fmu}}$. This shows that  $\fLambda_{{\bhh, \fmu}}$ is the only weal limit point of $(\fLambda_{{\bhh, \fmu_n}})_{n\in \bbN}$, which entails the desired result. 
We complete the proof of $(c)$ exactly as in Theorem \ref{eravsweakfix}. 

  Finally, we assume $(c)$ and we prove $(c)\Rightarrow (a)$. We show that $\fLambda_{{\bhh, \fmu_n}} \!\! \to \! \fLambda_{{\bhh, \fmu}}$ weakly on $\cM_f(\fell_\bullet^{\uparrow}(\bbN))$ exactly as in Theorem \ref{eravsweakfix}. This implies that $\fLambda_{{\bhh, \fmu}}$-a.e.~$L_\infty \leko \infty$. Thus, $\flambda_\fmu\ino \cM_f(\bbR_+)$ and $\fmu \ino \bbT$.  
By Proposition \ref{Lambdaprop} $(ii)$, we next observe for all $p, n\ino \bbN$, and all $\epp\ino \bbR^*_+$ that 
$$  \dGP( \fmu_n, \fmu) \leqo \epp \! \vee \! \fLambda_{\bhh, \fmu_n} (\{ L_\infty\! -\! L_p \geqo \epp \}) +  \dGP( \ttP_{T^n_{\! \mu_n, h_p}} \mu_n , \ttP_{T_{\! \mu, h_p}} \mu) +  \epp \! \vee \! \fLambda_{\bhh, \fmu} (\{ L_\infty\! -\! L_p \geqo \epp \}) $$
and we complete the proof exactly as in Theorem \ref{eravsweakfix}, using \eqref{eraC0GHP} in Theorem \ref{eraGcont} $(ii)$ instead of Lemma \ref{weakvsmera1} $(ii)$. \cqfd  

\begin{corollary}
\label{GPprecpera} 
Let $\bhh\eqo (h_p)_{p\in \bN}$ be an $\bbR^*_+$-valued sequence decreasing to $0$.  
Let $\emptyset \! \neq \! \cC \! \subseteq \! \bbT$. Then, $\cC$ is $\dGP$-precompact iff it satisfies (\ref{totbound}) and 
\begin{equation}
\label{precpaddi}
\forall \epp \ino \bbR^*_+, \qquad \sup  \big\{\,  \mu \big( T\backslash (T_{\! \mu, h_p})^{(\varepsilon)} \big) \,  ; \; (T,d, \rho, \mu) \! \equiv \! \fmu \ino  \cC \, \big\} \underset{p\rightarrow \infty}{-\!\!\! -\!\!\! -\!\!\! \longrightarrow} \;  0.
\end{equation} 
\end{corollary}
\noi
\textbf{Proof.} For all $(T,d, \rho, \mu) \! \equiv \! \fmu \ino \bbT$, we observe that 
$ \mu ( T\backslash (T_{\! \mu, h_p})^{(\varepsilon)} ) \eqo \fLambda_{{\bhh, \fmu}}(\{ L_\infty \! -\! L_p \geko \epp\})$. The proof is then a simple consequence of Theorem \ref{GPcvvsera}: we leave the details to the reader. \cqfd 

\medskip

\noi
\textbf{Connections with Gromov--Haussdorf--Prokhorov convergence.}
We conclude this section by discussing necessary and sufficient conditions which allow to derive $\dGHP$-convergence from joint $\dGP$- and $\dGH$-convergence. 
Indeed, recall that a joint $\dGP$- and $\dGH$-convergence does not necessarily entail a $\dGHP$-convergence as illustrated by Example \ref{GHGPnotGHPex}.

Let us recall from Definition \ref{mT01def} the spaces $(\bbT^0_{\! c}, \dGH)$ and $(\bbT^1_{\! c}, \dGHP)$ and from Definition \ref{mdef} the space $(\bbT, \dGP)$. 
We also recall the following from Definition \ref{TcTmindef}: the definition of the $1$-($\dGHP, \dGH$)-Lipschitz function $\Phi^1_0\colon \bbT^1_{\! c} \! \to \! \bbT^0_{\! c}$, which consists in forgetting the measure; the function $\spaan \colon \bbT^1_{\! c} \! \to \! \bbT^1_{\! c}$ such that $\spaan \, (\fmu)$ is the isometry class of $(\spaan \, (\supp \mu) , d, \rho, \mu)$;   
the set $\bbT^1_{\! c, \mathtt{min}}$ of the minimal elements of $\bbT^1_{\! c}$, namely the set of $\fmu\ino \bbT^1_{\! c}$ such that $\fmu \eqo \spaan\, (\fmu)$; 
we also recall that $\bbT_{\! c}$ stands for the set of $(T,d,\rho, \mu)  \! \equiv\!  \fmu\ino \bbT$ such that $(T,d)$ is compact and we recall for all $(T,d,\rho, \mu)  \! \equiv\!  \fmu\ino \bbT_{\! c}$ that $\Phi_1 (\fmu)\ino\bbT^1_{\! c} $ is the isometry class in $\bbT^1_{\! c}$ of $(T,d, \rho, \mu)$ and that $\Phi_0 (\fmu)\ino \bbT^0_{\! c}$ is the isometry class of $(T,d, \rho)$. 
We first need to prove the following. 
\begin{lemma}
\label{cpctsubsets} We keep the previous notation. Then 
$\bbT_{\! c}$ is a Borel-measurable subset of $(\bbT, \dGP)$, $\bbT^1_{\! c, \mathtt{min}}$ is a Borel-measurable subset of $(\bbT^1_{\! c}, \dGHP)$, $\spaan\colon \bbT^1_{\! c} \! \to \! \bbT^1_{\! c, \mathtt{min}}$ is a Borel-measurable surjection, 
$\Phi_1\colon  \bbT_{\! c} \! \to \! \bbT^1_{\! c, \mathtt{min}}$ is a Borel-bi-measurable bijection and $\Phi_0\colon  \bbT_{\! c} \! \to \! \bbT^0_{\! c}$ is a Borel-measurable surjection.
\end{lemma}
\noi
\textbf{Proof.} Let us first recall Kuratowski's theorem which is used several times in the proof: \emph{Let $(E_1,d_1)$ and $(E_2, d_2)$ be two Polish spaces, let $B_1$ be a Borel subset of $E_1$ and let $\varphi\colon B_1\! \to \! E_2$ be injective and Borel-measurable. Then $B_2\! :=\! \varphi (B_1)$ is a Borel subset of $E_2$ and $\varphi^{-1}\colon  B_2 \! \to \! B_1$ is also Borel-measurable} (see e.g.~Ethier and Kurtz \cite[Thm.~10.5]{EK}).

We first prove that $\spaan\colon \bbT^1_{\! c} \! \to \! \bbT^1_{\! c}$ is Borel-measurable. 
To this end, we fix $(T,d,\rho, \mu)\! \equiv \! \fmu \ino \bbT^1_{\! c}$, $h\ino \bbR^*_+$ and $s\ino (0, 1)$.
By Proposition \ref{properased} $(iv)$, $T_{\! \mu, sh+}\! \! \eqo \mathtt{Supp}\,  \era_{sh} \mu$, and we 
observe that $s\ino (0, 1) \! \mapsto \! T_{\! \mu, sh+}$ is non-increasing with respect to inclusion and thus $d_{\mathtt{Haus}}$-c\`adl\`ag by the standard results $\textbf{Haus}$-$(iv)$ and $\textbf{Haus}$-$(v)$ that are recalled 
at the beginning of Section \ref{Merafixsec}. 
We then also recall from Proposition \ref{properased} $(vi)$ that $s\ino (0, 1)\! \mapsto \era_{sh}\mu$ is $d_{\mathtt{Pro}}$-continuous. Let 
$F(sh,\fmu)$ stand for the isometry class of $(T_{\! \mu, sh+}, d, \rho, \era_{sh} \mu)$. 
Therefore $s\ino (0, 1) \! \mapsto\!  F(sh,\fmu) \ino \bbT_{\! c}^1$ is $\dGHP$-c\`adl\`ag and thus Borel-measurable. 
Let $U\colon \Omega \! \to \! (0, 1)$ be a uniformly distributed r.v. Then $F(Uh,\fmu)$ is a r.v., by the 
measurability just noted and we denote its law by $L (h,\fmu)\ino \cM_1(\bbT^1_{\! c})$.
 By Theorem \ref{eraGcont} $(ii)$, $\bP$-a.s.~$\fmu' \ino \bbT^1_{\! c} \! \mapsto \! F(Uh,\fmu')\ino \bbT^1_{\! c}$ 
 is continuous at $\fmu$, which implies that $\fmu \! \mapsto \! L(h,\fmu)$ is continuous when 
 $\cM_1(\bbT^1_{\! c})$ is equipped with the topology of weak convergence. 

By Proposition \ref{massTprop} $(vi)$, we next get $\lim_{h\to 0^+} 
d_{\mathtt{Haus}} (T_{\! \mu, h}, \spaan \, (\supp \mu))\eqo 0$ and by Proposition \ref{weakcvera}, 
$\lim_{h\to 0^+} d_{\mathtt{Pro}} (\era_h \mu, \mu)\eqo 0$. Therefore, a.s.
$ \lim_{h\to 0+}$ $\dGHP(F(Uh,\fmu), \spaan \, (\fmu)) $ $\eqo$ $ 0$ and 
$ \lim_{h\to 0+}L(h,\fmu)\! = \! \delta\mathtt{irac} (\spaan (\fmu))$ weakly in $\cM_1(\bbT^1_{\! c})$, 
where $\delta\mathtt{irac}(x)$ is a generic 
notation for the Dirac mass at $x$. 
This proves that the function 
$\mathtt{q} \colon  \fmu \ino \bbT^1_{\! c} \! \mapsto \!\delta\mathtt{irac} (\spaan \, (\fmu))\ino \cM_1(\bbT^1_{\! c})$ is 
Borel-measurable since it is a pointwise limit of continuous functions. 

We next observe that the function $\varphi \colon  \fmu \ino \bbT^1_{\! c} \! \mapsto\! \delta\mathtt{irac}(\fmu )  \ino  \cM_1(\bbT^1_{\! c})$ 
is continuous and injective. By Kuratowski's theorem the range 
$A\! := \! \{ \delta\mathtt{irac} ({\mathtt t})\,  ; \, \mathtt t\ino \bbT^1_{\! c} \}$ of $\varphi$ 
is a Borel subset of $\cM_1(\bbT^1_{\! c})$ and $\varphi^{-1} \colon A\! \to \! \bbT^1_{\! c}$ is Borel-measurable. This implies that $\fmu \ino \bbT^1_{\! c} \! \mapsto \! \spaan \, (\fmu) \ino \bbT^1_{\! c}$ is Borel-measurable since $\spaan \eqo \varphi^{-1} \! \circ\!  q $.   
 
 We now complete the proof of the lemma. To this end, we denote by $\Delta$ the diagonal of the space $(\bbT^0_{\! c})^2$  that is endowed with the product topology: $\Delta $ is therefore a closed subset of $(\bbT^0_{\! c})^2$. We then set 
 $g \colon\fmu \ino \bbT^1_{\! c}\! \mapsto \! \big( \Phi^1_0 ( \spaan\, (\fmu)), \Phi^1_0 (\fmu)\big) \ino (\bbT^0_{\! c})^2$, which is Borel-measurable (since $\spaan$ is Borel-measurable and $\Phi^1_0$ continuous) and we 
observe that $\bbT^1_{c,\mathtt{min}}\eqo g^{-1} (\Delta)$, which is therefore a Borel subset of $ \bbT^1_{\! c}$.

For all $(T,d,\rho, \mu) \! \equiv \! \fmu  \ino \bbT^1_{c, \mathtt{min}}$, we denote by $\phi(\fmu)$ 
the corresponding element of $\bbT_{\! c}$. The reader may note that $\phi$ is actually the identity map, but the domain and the co-domain have different topologies. We note that $\phi$ is injective $1$-Lipschitz from $(\bbT^1_{c, \mathtt{min}}, \dGHP)$ to $(\bbT, \dGP)$ and that 
$\bbT_{\! c}\eqo \phi (\bbT^1_{c, \mathtt{min}})$. Kuratowski's theorem applies and yields that $\bbT_{\! c}$ is a Borel subset of $\bbT$ and that $\Phi_1\eqo \phi^{-1}\colon\bbT_c\rightarrow\bbT_{c,\mathtt{min}}^1$ is Borel-measurable. 
Therefore,  $\Phi_0\eqo \Phi^1_0 \circ \Phi_1 \colon \bbT_{\! c} \! \to \! \bbT^0_{\! c}$ is also Borel-measurable, which completes the proof of the lemma.   \cqfd 

 \smallskip
 
The following lemma provides a simple criterion for elements of $\bbT$ to belong to $\bbT_{\! c}$. 
\begin{lemma} 
\label{Tccriterion} Let $\fmu \ino \bbT$. Then $\fmu \ino \bbT_{\! c}$ iff the trees $(\Phi_0 (\era_h \fmu))_{h\in \bbR_+^*}$ are precompact in $(\bbT^0_{\! c}, \dGH)$. In that case, we furthermore get $\lim_{h\to 0+} \dGHP( \era_h \fmu, \Phi_1 (\fmu))\eqo 0$. 
\end{lemma}
\noi
\textbf{Proof.} Let $(T,d, \rho, \mu) \! \equiv\! \fmu \ino \bbT$. Proposition \ref{properased} $(iv)$ implies that $T_{\! \mu, h+}\eqo \mathtt{Supp}\,  \era_{h} \mu$, which is compact. Therefore $(T_{\! \mu, h+}, d, \rho, \era_h\mu) $ is minimal, it belongs to $\bbT_{\! c}$ and $\Phi_0 (\era_h \fmu)\equiv (T_{\! \mu, h+}, d, \rho)$. We next observe that $h\! \mapsto \! T_{\! \mu, h+}$ 
is non-increasing and by Proposition \ref{massTprop} $(vi)$, the closure of $ \bigcup_{h\in \bbR^*_+} T_{\! \mu , h+}$ is 
$\sspaan \, (\supp  \mu )\eqo T$, since $\fmu$ is by definition minimal. Thanks to 
$\textbf{Haus}$-$(iv)$, we see that $T$ is compact iff $(T_{\! \mu, h+})_{h\in \bbR_+^*} $ are precompact in $(\mathcal K_T, d_{\mathtt{Haus}})$, which is equivalent to the fact that $(\Phi_0 (\era_h \fmu))_{h\in \bbR_+^*}$ are precompact in $(\bbT^0_{\! c}, \dGH)$ since $h\! \mapsto \! T_{\! \mu, h+}$ is non-increasing. Moreover, in this case, $\textbf{Haus}$-$(iv)$ also asserts that $\lim_{h\to 0}d_{\mathtt{Haus}}(T_{\! \mu, h+}, T)\eqo 0$. By Proposition \ref{weakcvera}, we also see that $\lim_{h\to 0+} d_{\mathtt{Pro}} (\era_h \mu , \mu) \eqo 0$, which completes the proof of the lemma. \cqfd 

\smallskip

Before proving the main result of this part of the section, which is Theorem \ref{GPGH+GHPdeter}, we need to prove 
an auxiliary result on embedded rooted compact $\bbR$-trees. 
\begin{definition}
\label{order} For all $\widetilde{T}_1, \widetilde{T}_2 \ino \bbT^0_{\! c}$, we write $\widetilde{T}_1\! \preceq \! \widetilde{T}_2$ if we can find an isometry $f\colon T_1\! \to \! T_2$ such that $f(\rho_1)\eqo\rho_2$, where $(T_i, d_i, \rho_i)$ stands for a representative of $\widetilde{T}_i$, $i\ino \{1, 2\}$. \cq 
\end{definition}
Only the first statement of the following lemma is used in the proof of Theorem \ref{GPGH+GHPdeter}. However, it is convenient to also prove $(ii)$ and $(iii)$: they are used in Theorem \ref{GPGH+GHP}, which is a random version of Theorem \ref{GPGH+GHPdeter}. 
\begin{lemma}
\label{orderdeter} We keep the previous notation. Then the following holds true. 
\begin{compactenum}

\smallskip

\item[$(i)$] $\preceq$ is a partial order on $\bbT^0_{\! c}$.

\smallskip

\item[$(ii)$] $\big\{ (\widetilde{T},\widetilde{T}')\ino (\bbT^0_{\! c})^2: \widetilde{T}\! \preceq \! \widetilde{T}'\big\}$ is a closed subset of $(\bbT^0_{\! c})^2$, equipped with the product topology.  

\smallskip

\item[$(iii)$] Let $\widetilde{\mathbf T}, \widetilde{\mathbf T}'\! $ be two $\bbT^0_{\! c}\! $-valued and identically distributed r.v.s 
such that $\bP$-a.s.~$\widetilde{\mathbf T}\! 
\preceq \! \widetilde{\mathbf T}'\! $. Then $\bP$-a.s.~$\widetilde{\mathbf T}\eqo \widetilde{\mathbf T}'$. 

\end{compactenum}
\end{lemma}
\noi
\textbf{Proof.} In $(i)$ antisymmetry, which is the only non-trivial point to check, follows from the fact that a compact metric space cannot be isometric to a proper subspace: see e.g.~Burago, Burago and Ivanov \cite[Thm 1.6.14]{BuBuIv}. 

Let us prove $(ii)$. Let $(T,d, \rho) \! \equiv \! \widetilde{T}$ and $(T'\! , d'\!, \rho') \! \equiv \! \widetilde{T}'$ 
belong to $\bbT^0_{\! c}$ and suppose $\dGH(\widetilde{T},\widetilde{T}')\leko \epp$. Then there exists 
$f\colon  T\! \to \! T'$ such that 
$\mathtt{dis} (f)\! :=\! \sup_{x,y\in T} |d'(f(x),f(y)) \!- \! d(x,y)|\leko 2\epp$ and $\sup_{y\in T'}d'(y,f(T))\leqo 2\epp$ 
(see e.g.~\cite[Cor.~7.3.28]{BuBuIv}). 
Let $(T_n, d_n, \rho_n) \! \equiv \! \widetilde{T}_n \! \preceq \! \widetilde{T}'_n \! \equiv \!(T'_n, d'_n ,\rho'_n)$, 
$n\ino \bbN$, and suppose 
$\epp_n \! :=\! \dGH(\widetilde{T}_n,\widetilde{T})+  \dGH(\widetilde{T}'_n,\widetilde{T}')\! \to \! 0$. 
Then for all $n\ino \bbN$, there are an isometry $\phi_n \colon  T_n \! \to \! T_n'$ and 
two functions $f_n \colon  T\! \to \! T_n$ and $f'_n\colon T_n'\! \to \! T'$ such that $\mathtt{dis} (f_n)\vee \mathtt{dis} (f'_n) \leko 2\epp_n$. 
We set $\psi_n\eqo f'_n\circ \phi_n \circ f_n$ which is a 
function from $T$ to $T'$ with $\mathtt{dis} (\psi_n) \leko 4\epp_n$. Let $Q\subo T$ be countable and dense. 
Since $T'$ is compact, a diagonal extraction argument implies that there is an increasing sequence of integers 
$(n_k)_{k\in \bbN}$ such that $d'$-$\lim_{k\to \infty} \psi_{n_k} (\sigma) \! =:\! \psi(\sigma)$ exists for all $\sigma \ino Q$. 
Clearly $\psi\colon Q\! \to \! T'$ is an isometry which extends uniquely to $T$ by standard arguments since $Q$ is dense. 
This completes the proof of $(ii)$.

Let us prove $(iii)$. To this end, for all $\epp \ino \bbR^*_+$ and all $(T,d,\rho) \! \equiv \! \widetilde{T}\ino \bbT^0_{\! c}$, we set 
 $V(\widetilde{T}, \epp)\eqo \max \{ \, p\ino \bbN^*\!  : \!  \exists 
 \sigma_1, \ldots, \sigma_p \ino T\! :\!  \min_{1\leq i< j\leq p}d(\sigma_i, \sigma_j) \geko \epp\}$, 
 which is a well-defined integer since $T$ is compact. First note that 
$\widetilde{T} \! \mapsto \! V(\widetilde{T} , \epp)$ is $\preceq$-non-decreasing. 
It is also easy to check that it is $\dGH$-lower semicontinuous (we leave the details to the reader), and therefore Borel-measurable. 

We next observe that if $\widetilde{T} \! \preceq \! \widetilde{T}'$ and 
$\widetilde{T} \! \neq \! \widetilde{T}'$, then there must be $\epp_0\ino \bbR^*_+$ such that $V(\widetilde{T}, \epp) \leko V(\widetilde{T}', \epp)$ for all $\epp \ino (0, \epp_0)$. Our assumptions imply for all $\epp\ino \bbR^*_+$ that 
$V(  \widetilde{\mathbf T}, \epp)$ and $V(\widetilde{\mathbf T}', \epp)$, which are well-defined r.v.s, 
have the same law and that $\bP$-a.s.~$V(  \widetilde{\mathbf T}, \epp)\leqo V(\widetilde{\mathbf T}', \epp)$. 
By Lemma \ref{adhoclem}, 
$\bP \big( V(  \widetilde{\mathbf T}, \epp)\eqo V(\widetilde{\mathbf T}', \epp) \big)\eqo 1$, for all $\epp \ino \bbR^*_+$, which implies that $\bP$-a.s.~$\widetilde{\mathbf T}\eqo \widetilde{\mathbf T}'$ by the previous arguments.  \cqfd

\begin{theorem}
\label{GPGH+GHPdeter} Let $\fmu,\fmu_n\ino \bbT_{\! c}$, $n\ino \bbN$. Then, the following holds true. 
\begin{compactenum}

\smallskip

\item[$(i)$] Let us assume that  $\fmu_n \!\!\!  \overset{\dGP}{\, \longrightarrow}\!  \fmu$ and that $(\Phi_1 (\fmu_n))_{n\in \bbN}$ is $\dGHP$-precompact. 
Then, any $\dGHP$-limit point $\fmu'$  of $(\Phi_1 (\fmu_n))_{n\in \bbN}$ is such that $\Phi_1 (\fmu) \eqo \spaan \, (\fmu')$.

\smallskip

\item[$(ii)$] Let us assume that  $\fmu_n \! \! \! \overset{\dGP}{\, \longrightarrow}\!  \fmu$ and that 
$(\Phi_0 (\fmu_n))_{n\in \bbN}$ is $\dGH$-precompact. Then, any $\dGH$-limit point $\widetilde{T}'$ of the sequence $(\Phi_0 (\fmu_n))_{n\in \bbN}$ is such that $\Phi_0 (\fmu) \preceq \widetilde{T}'$. 

\smallskip

\item[$(iii)$] $\Phi_1 (\fmu_n) \! \! \overset{\dGHP}{\, \longrightarrow}\!  \Phi_1 (\fmu)\; $ iff $\; \fmu_n \!\!  \overset{\dGP}{\, \longrightarrow}\!  \fmu$ and $\Phi_0 (\fmu_n) \!\!  \overset{\dGH}{\, \longrightarrow}\!  \Phi_0 (\fmu)$. 
\end{compactenum}
\end{theorem}
\noi
\textbf{Proof.} 
Let us prove $(i)$. Let  $(T'\! , d'\! , \rho'\!, \mu')\! \equiv \fmu'\ino \bbT^1_{\! c}$ be a $\dGHP$-limit point of  $(\Phi_1 (\fmu_n))_{n\in \bbN}$, i.e.~there is an increasing sequence of integers 
$(n_k)_{k\in \bbN}$ such that $\lim_{k\to \infty} $ 
$ \dGHP ( \Phi_1(\fmu_{n_k} ), \fmu')$ $\eqo 0$. We recall the notation 
$(\spaan\, (\supp \mu'), d'\! , \rho'\! , \mu')\! \equiv \! \spaan \,  (\fmu')  \ino \bbT^1_{\! c, \mathtt{min}}$. 
We then set $\fmu^\bullet\! :=\! \Phi_1^{-1} (\spaan \,  (\fmu') ) \ino \bbT_{\! c}$ (which is also represented by 
$(\spaan\, (\supp \mu'), d'\! , \rho'\! , \mu')$).  
Since $h\ino \bbR^*_+\! \mapsto \! T'_{\! \mu'\! ,\,  h}$ has at most countably many 
$d'_{\mathtt{Haus}}$-discontinuities, Theorem \ref{eraGcont} $(ii)$ implies that 
there is $S\!\subseteq\! \bbR_+^*$, such that $\bbR^*_+ \backslash S$ is countable and such that    
$\lim_{k\to \infty} \dGHP(\era_h \Phi_1 (\fmu_{n_k}) , \era_h \fmu' )\eqo 0$, for all $h\ino S$, 
which entails $\lim_{k\to \infty}$ $\dGHP(\era_h \Phi_1 (\fmu_{n_k}) , \era_h \Phi_1 (\fmu^\bullet) )$ $\eqo 0$, 
since $\era_h \fmu'\eqo \era_h \spaan\, (\fmu') \eqo \era_h \Phi_1 (\fmu^\bullet)$. This implies 
$\lim_{k\to \infty}$ $\dGP(\era_h \fmu_{n_k} , \era_h \fmu^\bullet )$ $\eqo 0$. Our assumption 
$\dGP (\fmu_n, \fmu) \! \to \! 0$, combined with the $\dGP$-continuity of mass erasure proved in 
Theorem \ref{eraGcont}, yield $\lim_{k\to \infty}$ $\dGP(\era_h \fmu_{n_k} , \era_h \fmu )\eqo 0$ 
for all $h\ino \bbR^*_+$. So we get $\era_h \fmu\eqo \era_h \mu^\bullet$, for all $h\ino S$. Since $h$ can be arbitrarily close to $0$, Lemma \ref{sepera} shows that $\fmu\eqo \fmu^\bullet$. Thus it shows $\Phi_1 (\fmu) \eqo \Phi_1 (\fmu^\bullet) \eqo \spaan \, (\fmu')$, which completes the proof of $(i)$.

   Let us prove $(ii)$. Let $(n_k)_{k\in \bbN}$ be an increasing sequence of integers and $(T'\! , d'\! , \rho')\! \equiv \widetilde{T}'\ino \bbT^0_{\! c}$ such that $\lim_{k\to \infty} $ 
$ \dGH ( \Phi_0(\fmu_{n_k} ), \widetilde{T}')$ $\eqo 0$. By the $\dGHP$-precompactness criterion 
$\textbf{Grom}$-$(v)$, which is recalled at the beginning of Section \ref{CVspacemerasec}, $(\Phi_1 (\fmu_{n}))_{n\in \bbN}$ and thus $(\Phi_1 (\fmu_{n_k}))_{k\in \bbN}$ are $\dGHP$-precompact. If $\fmu'\ino \bbT^1_{\!c}$ is a $\dGHP$-limit point of 
$(\Phi_1 (\fmu_{n_k}))_{k\in \bbN}$ then $\Phi^1_0 (\fmu') \eqo \widetilde{T}'$ since $\Phi^1_0$ is ($\dGHP, \dGH$)-Lipschitz. 
Now observe by $(i)$ that $\Phi_1(\fmu)\eqo \spaan \, (\fmu') $. Thus $\Phi_0 (\fmu)\eqo \Phi^1_0 (\Phi_1 (\fmu)) \! \equiv \!  (\spaan \, (\supp \mu') ,d'\! , \rho' ) \preceq (T'\! ,d'\! ,\rho')\! \equiv \! \widetilde{T}'$, which completes the proof of $(ii)$.

To prove $(iii)$, we first observe that  $\dGHP$-convergence implies a joint $\dGH$- and $\dGP$-convergence because here the $\dGHP$-limit is assumed from the beginning to be minimal. Let us prove the non-trivial implication. 
Namely, we assume that $\lim_{n\to \infty}$ $\dGP(\fmu_n, \fmu)\eqo $
$ \lim_{n\to \infty} $ $\dGH(\Phi_0 (\fmu_n), \Phi_0 (\fmu)) $ $\eqo 0$. By the $\dGHP$-precompactness criterion 
$\textbf{Grom}$-$(v)$, $(\Phi_1 (\fmu_n))_{n\in \bbN}$ is $\dGHP$-precompact sequence. 
Let $\fmu'$ be a $\dGHP$-limit point. By $(i)$, $\spaan \, (\fmu') \eqo \Phi_1 (\fmu)$ and thus 
$\Phi^1_0 (\spaan \, (\fmu'))\eqo \Phi^1_0 ( \Phi_1(\fmu))\eqo \Phi_0 (\fmu)$. 
Since $\Phi^1_0$ is ($\dGHP, \dGH$)-Lipschitz, $\Phi^1_0 (\fmu') $ is also the $\dGH$-limit of 
$\Phi_0 (\fmu_n)\eqo \Phi^1_0 (\Phi_1(\fmu_n))$, $n\ino \bbN$. Thus $\Phi^1_0 (\fmu')\eqo \Phi_0 (\fmu)$. This shows that 
$\Phi^1_0 (\spaan \, (\fmu'))  \eqo \Phi^1_0 (\fmu')$, i.e.~$(\spaan\, (\supp \mu'), d'\! , \rho')$ is isometric to 
$(T'\! , d'\! , \rho'\!)$, where $(T'\! , d'\! , \rho'\! , \mu')$ stands for a representative of $\fmu'$. Since a compact metric 
space cannot be isometric to a proper subspace (see e.g.~\cite[Thm~1.6.14]{BuBuIv}), we get $T'\eqo  \spaan \, (\supp \mu')$, which implies $\fmu'\eqo \Phi_1 (\fmu)$. Thus $ \Phi_1 (\fmu)$ is the only $\dGHP$-limit point of $(\Phi_1 (\fmu_n))_{n\in \bbN}$, which completes the proof of $(iii)$. \cqfd

\subsection{Convergence in distribution of random measured $\bbR$-trees}
\label{cvrandomtrsec}
For any Polish $(E,d)$, we denote by $\mathcal M_1 (E)$ the space of its Borel probability measures. We always endow $\mathcal M_1 (E)$ with the topology of the convergence in distribution, which is Polish too. 
We generically denote $\bbT^*$-valued r.v.s by $\mathbf m\! \equiv \!  (\mathbf T, \mathbf d, \rho,  \mathtt{m})$. 
\begin{lemma}
\label{cvlawmeasR+} Let $\flambda, \flambda_n $, be $\cM_f(\bbR_+)$-valued (resp.~$\cM_f([0,\infty])$-valued) r.v.s, $n\ino \bbN$. Then the following holds true. 
\begin{compactenum}

\smallskip

\item[$(i)$] There is a (deterministic) subset $C$ of $\bbR_+$ such that $\bbR_+ \backslash C$ is countable and for all $r\ino C$, $\bP$-a.s.~$\flambda(\{ r\})\eqo 0$.  

\smallskip

\item[$(ii)$] $\flambda_n \! \to \! \flambda$ in law in $\cM_f(\bbR_+)$ (resp.~in $\cM_f([0,\infty])$) iff for all $p\ino \bbN^*$ and all $r_1\leko \cdots \leko r_p $ in $C$,  $\big( \flambda_n((r_j, \infty))\big)_{1\leq j\leq p} \!\!\! \to \! \big( \flambda ((r_j, \infty)) \big)_{1\leq j\leq p} $ 
(resp.~$\big( \flambda_n((r_j, \infty])\big)_{1\leq j\leq p} \!\! \! \to \! \big( \flambda ((r_j, \infty]) \big)_{1\leq j\leq p} $) in law in $\bbR_+^{p}$. 
\end{compactenum}
\end{lemma}
\noi
\textbf{Proof.} This is left to the reader. \cqfd 

\begin{lemma}
\label{consistlaw} 
Let $(h_p)_{p\in \bN}$ be $\bbR^*_+$-valued and strictly  decreasing to $0$. 
\begin{compactenum}

\smallskip

\item[$(i)$] Let $\mathbf m$ and $\mathbf m'$ be two $\bbT^*\!$-valued r.v.s such that for all $p\ino \bbN$, 
$\era_{h_p} \mathbf m$ has the same law as $\era_{h_p} \mathbf m'$. Then $\mathbf m$ and $\mathbf m'$ have the same law.  

\smallskip

\item[$(ii)$] Let 
$\mathbf m'_p
  $, $p\ino \bbN$, be a sequence of $\bbT^*\! $-valued r.v.s such that $\era_{h_p-h_{p+1}}\mathbf m_{p+1}'$ has the same law as $\mathbf m'_p$. Then, there exists a $\bbT^*$-valued r.v.~$\mathbf m
 $, unique in law and such that $\era_{h_p}  \mathbf m$ has the same law as $\mathbf m_p'$, for all $p\ino \bbN$,
\end{compactenum}
\end{lemma}
\noi
\textbf{Proof.} We first prove $(i)$. To this end, we ecall the notation in Remark \ref{musequence}. Let $p\ino \bbN$. The semigroup property of mass erasure implies that $(\era_{h_q} \mathbf m )_{0\leq q\leq p}\eqo (\era_{h_q-h_p} (\era_{h_p} \mathbf m))_{0\leq q\leq p}$ and that $(\era_{h_q} \mathbf m')_{0\leq q\leq p}\eqo (\era_{h_q-h_p} (\era_{h_p} \mathbf m'))_{0\leq q\leq p}$. 
Since $\era_{h_p} \mathbf m$ and $\era_{h_p} \mathbf m'$ have the same law, this implies that $(\era_{h_q} \mathbf m)_{0\leq q\leq p}$ and $(\era_{h_q} \mathbf m')_{0\leq q\leq p} $ have the same law in $\bbT^{p+1}$. 
Since this holds true for arbitrarily large $p$, Kolmogorov's extension theorem implies that 
$\jmath (\mathbf m)$ and  $\jmath (\mathbf m')$ have the same law in $\bbT^{\bbN}$. Therefore $\mathbf m\eqo \jmath^{-1} (\jmath (\mathbf m)) $ and  $\mathbf m'\eqo \jmath^{-1} (\jmath (\mathbf m')) $ have the same law in $\bT^*$. 

Let us prove $(ii)$. For all integers $p\geqo q \geqo 0$, we set $\mathbf m^{{(p)}}_q\eqo \era_{h_q-h_p} \mathbf m'_{p}$. Our assumptions imply that $(\mathbf m^{{(p)}}_q)_{0\leq q \leq p} $ has the same law as 
$( \era_{h_q-h_p} (\era_{h_p-h_{p+1}} \mathbf m'_{p+1}))_{0\leq q \leq p}$ which is equal to $(\mathbf m^{{(p+1)}}_q)_{0\leq q \leq p}  $ by the semigroup property of mass erasure. By Kolmorogov's extension theorem there is a 
$\bT^{\bbN} $-valued r.v.~$(\mathbf m_p)_{p\in \bbN}$ such that for all $p\ino \bbN$, $(\mathbf m_q)_{0\leq q\leq p}$ and $(\mathbf m^{{(p)}}_q)_{0\leq q \leq p} $ have the same law. In particular, this implies that $\mathbf m_p$ and $\mathbf m_p'$ have the same law. Moreover, for all $p\ino \bN$, $\era_{h_p-h_{p+1}} \mathbf m_{p+1} \eqo \mathbf m_p$, i.e.~$\bP$-a.s.~$(\mathbf m_p)_{ p\in \bbN }\ino \bbT^\oast \eqo \jmath(\bbT^*)$. We then set 
$\mathbf m\eqo \jmath^{-1}( (\mathbf m_p)_{ p\in \bbN })$ and, by definition of $\jmath$, $\era_{h_p} \mathbf m\eqo \mathbf m_p $, which has therefore the same law as $\mathbf m'_p$. Uniqueness in law is an immediate consequence of $(i)$.  \cqfd 

\begin{proposition}
\label{traducvlawera} Let $(h_p)_{p\in \bN}$ be $\bbR^*_+$-valued and strictly  decreasing to $0$. Let 
$\mathbf m\! \equiv\!  (\mathbf T, \mathbf d, \rho, \mathtt m)$ and $\mathbf m_n\! \equiv\!  (\mathbf T_n, \mathbf d_n, \rho_n, \mathtt m_n)$, $n\ino \bbN$, be $\bbT^*$-valued r.v.s.
\begin{compactenum}

\smallskip

\item[$(i)$] The following assertions are equivalent. 

\begin{compactenum}

\smallskip

\item[$(a)$] $\mathbf m_n \! \to \! \mathbf m$ in law in $(\bT^*, \dera)$.

\smallskip

\item[$(b)$] For all $h\ino \bbR^*_+$, $\era_{h} \mathbf m_n \! \to \! \era_{h} \mathbf m$ in law in $(\bT, \dGP)$. 

\smallskip

\item[$(c)$] For all $p\ino \bbN$, $\era_{h_p} \mathbf m_n \! \to \! \era_{h_p} \mathbf m$ in law in $(\bT, \dGP)$.

\end{compactenum}

\smallskip

\item[$(ii)$] The laws of $(\mathbf m_n)_{n\in \bbN}$ are tight in $(\bbT^*, \dera)$ iff the laws of 
$(\langle \mathbf m_n \rangle)_{n\in \bbN} $ are tight in $\bbR_+$ and for all fixed $p\ino \bbN$, the laws of $(\mathtt{Ht} (\era_{h_p} \mathbf m_n))_{n\in \bbN}$  are tight in $\bbR_+$.

\smallskip

\item[$(iii)$] Let us suppose $(\mathbf m_n)_{n\in \bN}$ to be $\bbT$-valued and $\mathbf m_n \! \to \! \mathbf m$ in law in $(\bbT^*, \dera)$. 
Then a.s.~$\mathbf m\ino \bbT$ (without necessarily implying that $\mathbf m_n \! \to \! \mathbf m $ in law in $(\bbT, \dGP)$) 
iff
\begin{equation}
\label{inTcondi}
\forall \epp \ino \bbR^*_+, \quad \lim_{r\to \infty} \, \limsup_{p\to \infty} \, \sup_{n\in \bbN}  \bP \big(\era_{h_p} \mathbf m_n \big(T_n \backslash B_{T_n} (\rho_n, r)\big) \geko \epp  \big) = 0.  
\end{equation}

\item[$(iv)$]  Let us suppose that $(\mathbf m_n)_{n\in \bN}$ is a $\bbT$-valued sequence. Then $(\mathbf m_n)_{n\in\bbN}$ converges in law in $(\bT, \dGP)$ iff the laws of 
$(\mathbf m_n)_{n\in \bbN}$ are tight in $(\bbT, \dGP)$ and $(\mathbf m_n)_{n\in\bbN}$ converges in law in $(\bT^*, \dera)$. In this case the limiting laws are concentrated on $\bT$ and coincide. 
\end{compactenum}
\end{proposition}
\noi
\textbf{Proof.} Let us prove $(i)$. By the very definition of convergence in the sense of mass erasure, $\era_h$ 
is continuous on $(\bbT^*\! , \dera)$ (see Lemma \ref{eracvprop}). So $(a)$ implies $(b)$ which in turns entails $(c)$, 
trivially. Let us prove that $(c)\! \Rightarrow \! (a)$. To this end, recall the notation in Remark \ref{musequence}. Let 
$p\ino \bbN$. The semigroup property combined with $(c)$ and the 
$\dGP$-continuity of $\era_h$ (Theorem \ref{eraGcont} $(i)$) imply that 
$(\era_{h_q} \mathbf m_n)_{0\leq q \leq p}\eqo (\era_{h_q-h_p} (\era_{h_p}\mathbf m_n))_{0\leq q \leq p}$ converges to 
$ (\era_{h_q-h_p} (\era_{h_p}\mathbf m))_{0\leq q \leq p}\eqo (\era_{h_q} \mathbf m)_{0\leq q \leq p}$ in law in $\bbT^{p+1}$.  
This implies that $(\era_{h_p} \mathbf m_n)_{p\in \bbN} \! \to \! (\era_{h_p} \mathbf m)_{p\in \bbN}$ 
in law in 
$\bbT^{\bbN}$ (\emph{indeed}, since $\bbT^{\bbN}$ is equipped with the product topology, the convergence in 
law in $\bbT^{\bbN}$ coincides with convergence in law of finite dimensional marginal laws). 
Therefore $\mathbf m_n \eqo \jmath^{-1} \big( (\era_{h_p} \mathbf m_n)_{p\in \bbN}\big) \! \to \!
  \jmath^{-1} \big( (\era_{h_p} \mathbf m_n)_{p\in \bbN}\big) \eqo \mathbf m$ in law in $(\bbT^*\! , \dera)$.

Let us prove $(ii)$. We first suppose that the laws of $(\mathbf m_n)_{n\in \bbN}$ are tight in $(\bbT^*, \dera)$. Let $\epp \ino (0, 1)$. 
There exists a $\dera$-compact subset $K$ of $\bbT^*$ such that $\sup_{n\in \bbN} \bP (\mathbf m_n \! \notin \!  K) \leqo\epp $. By Theorem \ref{cveracomplet} $(iii)$, $x\! :=\! \sup \big\{  \langle \fmu \rangle  \, ; \, \fmu\ino K \big\} \leko  \infty$ and for all $p\ino \bbN$, 
$y_p\! :=\! \sup \big\{ \mathtt{Ht} (\era_{h_p} \fmu) \, ; \, \fmu\ino K \big\} \leko  \infty$. Thus $\bP (\langle \mathbf m_n \rangle \geko x)\vee \max_{p\in \bbN} \bP (\mathtt{Ht} (\era_{h_p} \mathbf m_n) \geko y_p)   \leqo \bP (\mathbf m_n \! \notin \! K) \leqo \epp$. Since $\epp$ can be chosen arbitrarily small, this implies the tightness of the laws of  
$(\langle \mathbf m_n \rangle)_{n\in \bbN} $ and $(\mathtt{Ht} (\era_{h_p} \mathbf m_n))_{n\in \bbN}$, $p\ino \bbN$. 

 Conversely, we assume the tightness of the laws of  
$(\langle \mathbf m_n \rangle)_{n\in \bbN} $ and of $(\mathtt{Ht} (\era_{h_p} \mathbf m_n))_{n\in \bbN}$, $p\ino \bbN$. 
Let $\epp \ino (0, 1)$. There are $x, y_{p}\ino \bbR^*_+$ such that we have $\sup_{n\in \bbN} \bP(\langle \mathbf m_n \rangle\geko x) \leko \frac{1}{2}\epp$ and $\sup_{n\in \bbN} \bP(\mathtt{Ht} (\era_{h_p} \mathbf m_n) \geko y_p) \leko 2^{-p-2}\epp $, $p\ino\bbN$.  
We then define $K_\epp$ as the set of $\fmu\ino \bbT^*$ such that $ \langle \fmu \rangle \leqo x$ and $\mathtt{Ht} (\era_{h_p} \fmu)\leqo y_p$ for all $p\ino \bbN$. Theorem \ref{cveracomplet} $(iii)$ implies that $K_\epp$ has a $\dera$-compact closure and observe that for all $n\ino \bbN$, $\bP ( \mathbf m_n \! \notin \! K_\epp )\leqo \bP(\langle \mathbf m_n \rangle\geko x) + \sum_{p\in \bbN}  \bP(\mathtt{Ht} (\era_{h_p} \mathbf m_n) \geko y_p) \leko \epp$, which implies the desired tightness.  
 
Let us prove $(iii)$. Since $\era_{h_p} \mathbf m_n \! \to \! \era_{h_p} \mathbf m$ in law in $(\bbT, \dGP)$, Theorem \ref{eraGcont} $(v)$ easily implies that $\flambda_{\era_{h_p} \mathbf m_n } \! \!\! \to \!  \flambda_{\era_{h_p} \mathbf m } $ in law in $\cM_f(\bbR_+)$. 
By Lemma \ref{cvlawmeasR+} $(ii)$ there is a dense subset $C\!\subseteq\! \bbR_+$ such that for all $r\ino C$, $\flambda_{\era_{h_p} \mathbf m_n } ((r, \infty)) \! \to \!  \flambda_{\era_{h_p} \mathbf m } ((r, \infty))$ in law in $\bbR_+$ and the Portmanteau theorem implies the following. 
\begin{eqnarray*}
\bP \big(  \flambda_{\era_{h_p}\!  \mathbf m } ((r, \infty)) \geko \epp \big) &\leqo & \liminf_{n\to \infty}  \bP \big(\era_{h_p} \mathbf m_n \big( \mathbf T_n \backslash B_{\mathbf T_n,d_n} (\rho_n, r)\big) \geko \epp  \big)  \\
& \leqo&  \limsup_{n\to \infty}  \bP \big(\era_{h_p} \! \mathbf m_n \big( \mathbf T_n \backslash B_{\mathbf T_n, d_n} (\rho_n, r)\big) \geqo \epp  \big)\leq \bP \big(  \flambda_{\era_{h_p}\!  \mathbf m } ((r, \infty))\geqo \epp \big) \, .
\end{eqnarray*}
We denote by $(C_i)_{i\in I}$ the connected components of $\mathbf T\backslash B_{\mathbf T,d}(\rho, r)$ and we recall that 
$\lambda_{\era_{h_p}\!  \mathbf m } ((r, \infty))$ $ \eqo$ $ \era_{h_p}\!  \mathbf m (\mathbf T\backslash B_{\mathbf T,d}(\rho, r))\eqo \sum_{i\in I} (\mathbf m (C_i)\!-\! h_p )_+ \uparrow \sum_{i\in I} \mathbf m (C_i)=\mathbf m (\mathbf T\backslash B_{\mathbf T,d}(\rho, r))\eqo  \flambda_{ \mathbf m } ((r, \infty])) $ by monotone convergence for sums. 
Thus, $\lim_{p\to \infty}  \bP \big(\flambda_{\era_{h_p} \mathbf m } ((r, \infty)) \geko \epp  \big)\eqo \bP \big(\flambda_{\mathbf m} ((r, \infty]) \geko \epp  \big) $, which easily completes the proof of $(iii)$. 

The proof of $(iv)$ is quite similar to the proof of Theorem \ref{cveracomplet} $(v)$. Let us prove the non-trivial implication. To this end,  
let us denote by $\mathbf m$ a $\bT^*$-valued r.v.~such that $\mathbf m_n \! \to \! \mathbf m$ in law in $(\bT^*, \dera)$. Let $\mathbf m'$ be a $\bT$-valued r.v.~such that $\mathbf m_{n_k} \! \to \! \mathbf m'$ in law in $(\bT, \dGP)$, where $(n_k)_{k\in \bbN}$ is some increasing $\bbN$-valued sequence. We only need to prove that $\mathbf m$ and $\mathbf m'$ have the same law on $(\bT^*, \dera)$. 
Let $h\ino \bbR^*_+$. On one hand, Theorem \ref{eraGcont} $(i)$ implies $\era_{h} \mathbf m_{n_k}\! \to \! \era_h \mathbf m'$ in law in $(\bT, \dGP)$. On the other hand, by definition, $\era_{h} \mathbf m_{n_k}\! \to \! \era_h \mathbf m$ in law in $(\bT, \dGP)$. Thus for all $h\ino \bbR^*_+$, $\era_h \mathbf m$ and $\era_h \mathbf m'$ have the same law. 
Lemma \ref{consistlaw} $(i)$ implies the desired result. \cqfd

\medskip

\noi
\textbf{Mass erasure and Gromov--Prokhorov convergence in distribution.} Let us start with the following technical remark. 
\begin{remark}
\label{stochaconti} Let $(\mathbf T,d, \rho, \mathtt m)\! \equiv\! \mathbf m$ be a $\bbT^*\! $-valued r.v. We recall the following: if $(X_{\! s})_{s\in \bbR_+}\! $ is a c\`adl\`ag process which takes its values in a Polish space, 
then the set of times $\{s \ino \bbR_+\!  :\!  \bP (X_{s-}\! \neq \! X_s) \geko 0\}$ is countable. 
Accordingly, there exists an $\bbR^*_+$-valued sequence $\bhh\eqo (h_p)_{p\in \bbN}$ strictly decreasing to $0$ such that 
$\bP$-a.s.~$h\! \mapsto \! \mathbf T_{\! \mathtt m, h}$ is $\dHaus$-continuous at each $h_p$, $p\ino \bbN$. \cq
\end{remark}

\begin{theorem}
\label{condilawGP} Let $(\mathbf T,d, \rho, \mathtt m)\! \equiv\! \mathbf m$ be a $\bbT^*\! $-valued r.v. 
Let $(h_p)_{p\in \bbN}$ be any $\bbR_+^*$-valued sequence which strictly decreases to $0$ and is such that 
$\bP$-a.s.~$h\! \mapsto \! \mathbf T_{\! \mathtt m, h}$ is $\dHaus$-continuous at each $h_p$, $p\ino \bbN$ 
(there is always one by Remark \ref{stochaconti}).  
Let $(\mathbf m_n)_{n\in \bbN}$ be $\bbT$-valued r.v.s such that $\mathbf m_n \! \to \! \mathbf m$ in law in 
$(\bT^*\! , \dera)$.  Then the following assertions are equivalent.
\begin{compactenum}

\smallskip

\item[$(a)$] $\bP$-a.s.~$\mathbf m\ino \bbT$, and $\mathbf m_n \! \to \! \mathbf m$ in law in $(\bT, \dGP)$. 

\smallskip

\item[$(b)$] $\bP$-a.s.~$\flambda_{\mathbf m} \ino \cM_f(\bbR_+)$, and the random finite measures 
$(\flambda_{\mathbf m_n})_{n\in \bbN}$  on $\bbR_+$ converge to $ \flambda_{\mathbf m}$ in law in $\cM_f(\bbR_+)$.

\smallskip

\item[$(c)$] For all $\epp, \eta\ino \bbR_+^*$, $\lim_{p\to \infty} \sup_{n\in \bbN} \bP \big( \fLambda_{\bhh, \mathbf m_n} (\{ L_\infty \! -\! L_p \geqo \epp\} ) \geqo \eta \big)\eqo 0$. 
\end{compactenum}  
\end{theorem}
\noi
\textbf{Proof.} 
Theorem \ref{eraGcont} $(v)$ entails that $(a)\! \Rightarrow \! (b)$. 
We next prove that $(b)$ or $(c)$ imply the following 
convergence in law in $\bbT^* \! \times \! \cM_f(\fell^\uparrow_\bullet (\bbN))$.
\begin{equation}
\label{jontcvmLambda}
(\mathbf m_{n} , \fLambda_{\bhh, \mathbf m_{n} })\xrightarrow[n\to \infty]{\;} (\mathbf m,  \fLambda_{\bhh, \mathbf m }) \quad \textrm{and} \quad \flambda_{\mathbf m} \eqo \fLambda_{\bhh, \mathbf m }\circ L_\infty^{-1} \in \cM_f(\bbR_+).
\end{equation}
By Lemma \ref{LambdaGprop} $(i)$,  
$\fnu \ino \bT^* \! \mapsto \! \fLambda^{\! o}_{{\bhh, \fnu}}$ is a.s.~continuous at $\mathbf m$. Thus, 
$(\mathbf m_n, \fLambda^{\! o}_{{\bhh, \mathbf m_n}}) \! \to \! (\mathbf m , \fLambda^{\! o}_{{\bhh, \mathbf m}})$ in law in $ \bbT^* \! \times \!\cM_f(\fell^\uparrow (\bbN))$. By Remark \ref{tightLambda}, the laws 
of $(\fLambda_{{\bhh, \mathbf m_n}})_{n\in \bbN}$ 
are tight in $\cM_f(\fell^\uparrow_\bullet (\bbN))$ and thus the laws of 
$((\mathbf m_n, \fLambda_{{\bhh, \mathbf m_n}}))_{n\in \bbN}$ are tight in 
$\bbT^* \! \times \! \cM_f(\fell^\uparrow_\bullet (\bbN))$. Let $(n_k)_{k\in \bbN}$ be any increasing sequence of 
integers such that $(\mathbf m_{n_k} , \fLambda_{{\bhh, \mathbf m_{n_k} }}\! )\!  \to \! (\mathbf m, \fLambda)$ in 
distribution in $\bbT^* \! \times \! \cM_f(\fell^\uparrow_\bullet (\bbN))$. By the previous arguments, we can 
assume that 
a.s.~$\fLambda^{\! o}\eqo \fLambda^{\! o}_{{\bhh, \mathbf m }}$. By Skorokhod's representation theorem 
(and a slight abuse of notation), we can also assume that the previous convergence in law actually holds 
$\bP$-almost surely on $\bbT^* \! \times \! \cM_f(\fell^\uparrow_\bullet (\bbN))$.  

Under $(b)$ and arguing as in Theorem \ref{GPcvvsera}, we see that $\fLambda \circ \mathcal L^{-1}$ has the same law as 
$\fLambda \circ L_\infty^{-1}$, which is the law of $\flambda_{\mathbf m}$ and that 
$\bP$-a.s $\fLambda$-a.e.~$\mathcal L \leqo L_\infty \leko \infty$. Lemma \ref{adhoclem} applies and we get that 
$\bP$-a.s.~$\fLambda$-a.e.~$\mathcal L\eqo L_\infty$ and thus a.s.~$\fLambda\eqo \fLambda_{{\bhh, \mathbf m}}$.
Consequently the law of $(\mathbf m,  \Lambda_{\bhh, \mathbf m }) $ is the only limit point of the laws of $((\mathbf m_n, \fLambda_{{\bhh, \mathbf m_n}}))_{n\in \bbN}$, which entails (\ref{jontcvmLambda}). 

We assume $(c)$. As in the proof of Theorem \ref{GPcvvsera}, the Portmanteau theorem implies a.s.~that 
$\fLambda (L_\infty\eqo \infty) \leqo \fLambda (\{ L_\infty \! -\! L_p \geko \epp \}) \leq \liminf_{k\to \infty} 
\fLambda_{{\bhh, \mathbf m_{n_k}}} (\{ L_\infty \! -\! L_p \geko \epp \}) $. Then $(c)$ easily entails that $\bP$-a.s.~$\fLambda$-a.e.~$\mathcal L\eqo L_\infty \leko \infty$ and thus a.s.~$\fLambda\eqo \fLambda_{{\bhh, \mathbf m}}$, which implies (\ref{jontcvmLambda}) as already explained. 

We next observe that (\ref{jontcvmLambda}) implies $(c)$. \emph{Indeed}, by Skorokhod's representation theorem, we can assume that (\ref{jontcvmLambda}) holds 
$\bP$-a.s.~on $\bbT^* \! \times \! \cM_f(\fell^\uparrow_\bullet (\bbN))$. Arguing deterministically as in Theorem \ref{GPcvvsera}, we get a.s.~for all $\epp\ino \bbR_+^*$, $\lim_{p\to \infty} \sup_{n\in \bbN} \fLambda_{{\bhh, \mathbf m_{n}}} (\{ L_\infty \! -\! L_p \geko \epp \})\eqo 0$, which implies $(c)$. 

We have proved that $(b) \! \Rightarrow \! (\ref{jontcvmLambda}) \! \Rightarrow (c)$. Now let us assume $(c)$, which implies  (\ref{jontcvmLambda}). By Skorokhod's representation theorem, we can assume that (\ref{jontcvmLambda}) holds 
$\bP$-a.s.~on $\bbT^* \! \times \! \cM_f(\fell^\uparrow_\bullet (\bbN))$. Then, we deterministically 
argue exactly as in the proof of $(c)\! \Rightarrow \! (a)$ in Theorem \ref{GPcvvsera} to get $(a)$.  \cqfd

\begin{corollary}
\label{critightGP}   
For all $n\ino \bbN$, let 
$ (\mathbf T^{ n}, \mathbf d_n, \rho_n, \mathtt m_n)\! \equiv\! \mathbf m_n $ be a $\bbT$-valued r.v. 
Then, the laws of $(\mathbf m_n)_{n\in \bbN}$ are tight in $(\bbT, \dGP)$ iff they are tight in $(\bbT^*\! , \dera)$ and 
\begin{equation}
\label{critightGP1}
\forall \epp, \eta \ino \bbR^*_+, \quad \lim_{h\to 0} \, \sup_{n\in \bbN} \bP \Big( \mathtt m_n \big( \mathbf T^{ n} \backslash (\mathbf T^n_{\! \mathtt m_n , h})^{(\epp)}\big) \geko \eta \Big) = 0\; .
\end{equation}
\end{corollary}
\noi
\textbf{Proof.} This is a simple consequence of Theorem \ref{condilawGP}: we leave the details to the reader. \cqfd

\medskip

\noi
\textbf{Connections with Gromov--Haussdorf--Prokhorov convergence in law.}
We recall Example \ref{GHGPnotGHPex}, Lemma \ref{Tccriterion}, 
which provides a criterion 
for an element of $\bbT$ to belong to $\bbT_{\! c}$, Lemma \ref{orderdeter} about embedded trees and Theorem \ref{GPGH+GHPdeter}, which discusses $\dGHP$-convergence in connection with joint $\dGP$- and $\dGH$-convergences.  
We first prove the analogue of Lemma \ref{Tccriterion}. 
\begin{lemma} 
\label{Tcrandomcriterion} 
Let $\mathbf m $ be a $\bbT$-valued r.v. Then $\bP (\mathbf m \ino \bbT_{\! c})\eqo 1$ iff 
the laws of $(\Phi_0 (\era_h \mathbf m))_{h\in \bbR_+^*}$ are tight in $(\bbT^{0}_{{\! c}}, \dGH)$. 
In this case, we furthermore get $\lim_{h\to 0+} \dGHP( \era_h \mathbf m, \Phi_1 (\mathbf m))\eqo 0$, almost surely. 
\end{lemma}
\noi
\textbf{Proof.} If $\bP (\mathbf m \ino \bbT_{\! c})\eqo 1$, then we note that 
a.s.~for all $h\ino \bbR_+^*$, $\Phi_0 (\era_h \mathbf m) \preceq \Phi_0 (\mathbf m)$ and we easily check that the laws of  $(\Phi_0 (\era_h \mathbf m))_{h\in \bbR_+^*}$ are tight in $(\bbT^{0}_{\! c}, \dGH)$. 

Conversely, let us assume that the laws of  $(\Phi_0 (\era_h \mathbf m))_{h\in \bbR_+^*}$ are tight in 
$(\bbT^{0}_{\! c}, \dGH)$. For all $\epp \ino (0, 1)$, there is a compact subset $K_\epp $ of $(\bbT^{0}_{\! c} , \dGH)$  
such that 
$\bP (A_\epp)\geqo 1\! -\! \epp$, where $A_\epp$ is the event $ \{ \forall h\ino \bbR_+^*\! : \Phi_0 (\era_h \mathbf m) \ino K_\epp \}$. By Lemma \ref{Tccriterion}, on $A_\epp$, $\mathbf m \ino \bbT_{\! c}$ and thus 
$\bP (\mathbf m \ino \bbT_{\! c}) \geqo 1\! -\! \epp$ 
(we recall from Lemma \ref{cpctsubsets} that $\bbT_{\! c}$ is a Borel subset of $(\bbT, \dGP)$, which implies that 
$\{ \mathbf m \ino \bbT_{\! c}\}$ is a measurable event). Since $\epp$ can be arbitrarily small, we get 
$\bP (\mathbf m \ino \bbT_{\! c} )\eqo 1$. 
To complete the proof of the lemma, we apply
Lemma \ref{Tccriterion} again on the event $\{ \mathbf m \ino \bbT_{\! c}\} $ to get $\lim_{h\to 0+} \dGHP( \era_h \mathbf m, \Phi_1 (\mathbf m))\eqo 0$. \cqfd 

\smallskip

\smallskip

We now prove the analogue of Theorem \ref{GPGH+GHPdeter}. 

\begin{theorem}
\label{GPGH+GHP} Let $ \mathbf m,\mathbf m_n$, $n\ino \bbN$, be $\bbT$-valued r.v.s such that $\bP (\mathbf m \ino \bbT_{\! c} )\eqo \bP (\mathbf m_n \ino \bbT_{\! c} )\eqo 1$, $n\ino \bbN$. Then the following holds true. 
\begin{compactenum}

\smallskip

\item[$(i)$] If $\mathbf m_n \! \to \! \mathbf m$ in law in $(\bbT, \dGP)$, if the laws of $(\Phi_1 (\mathbf m_n))_{n\in \bbN}$ are tight on $(\bbT_{\! c}^1, \dGHP)$ and if $\mathbf m'$ is distributed according to a weak limit point of the laws of $(\Phi_1 (\mathbf m_n))_{n\in \bbN}$, then $\Phi_1 (\mathbf m) $ has the same law as $\spaan \, (\mathbf m')$. 

\smallskip

\smallskip

\item[$(ii)$] $\Phi_1(\mathbf m_n) \! \to \! \Phi_1(\mathbf m)$ in law in $(\bbT^1_{\! c}, \dGHP)$ iff 
$\mathbf m_n \! \to \! \mathbf m$ in law in $(\bbT, \dGP)$ and $\Phi_0(\mathbf m_n) \! \to \! \Phi_0(\mathbf m)$ in law in $(\bbT^0_{\! c}, \dGH)$

\end{compactenum}
\end{theorem}
\noi
\textbf{Proof.} We first prove $(i)$. Let $\mathbf m'$ be distributed 
according to a weak limit point of the laws of $(\Phi_1 (\mathbf m_n))_{n\in \bbN}$. 
We note that the laws of $(\mathbf m_n , \Phi_1(\mathbf m_n))$, $n\ino \bbN$, are tight in $\bbT \! \times \! \bbT^{1}_{\! c}$. 
Thus there exists an increasing sequence of integers $(n_k)_{k\in \bbN}$ such that 
$(\mathbf m_{n_k}, \Phi_1(\mathbf m_{n_k}))\! \to \! (\mathbf m, \mathbf m'')$ in law in 
$\bbT \! \times \! \bbT^{1}_{\! c}$, where $\mathbf m'$ and $\mathbf m''$ have the same distribution. 
By Skorokhod 
Representation Theorem (and a slight abuse of notation), we suppose that the latter convergence holds $\bP$-almost 
surely. Theorem \ref{GPGH+GHPdeter} $(i)$ then implies $\bP$-a.s.~that 
$\Phi_1(\mathbf m) \eqo \spaan \, (\mathbf m'')$, which entails that $\Phi_1 (\mathbf m) $ has the same law as 
$\spaan \, (\mathbf m')$.

To prove $(ii)$, we first observe that $\dGHP$-convergence implies $\dGH$- and $\dGP$-convergences in law because, here,  the $\dGHP$-limit is assumed to be minimal. Let us prove the non-trivial implication. 
To this end, we first observe that $\dGP$- and $\dGH$-convergences in law imply, by \textbf{Gromov}-$(v)$, 
that the laws of $(\Phi_1(\mathbf m_n))_{n\in \bbN}$ are tight on $(\bbT^{1}_{\! c}, \dGHP)$. Let $\mathbf m'$ be a 
$\bbT^{1}_{\! c}$-valued r.v.~distributed according to a weak limit point of the laws of $(\Phi_1(\mathbf m_n))_{n\in \bbN}$. 
By $(i)$, $\Phi_1(\mathbf m)$ and $\spaan\,  (\mathbf m')$ have the same law. 
Since $\Phi^1_0$ is ($\dGHP, \dGH$)-Lipschitz, 
the law of $\Phi^1_0(\mathbf m')$ is also a weak limit point of the law of $\Phi^1_0 (\Phi_1 (\mathbf m_n))\eqo \Phi_0 (\mathbf m_n)$, $n\ino\bbN$ and since $\Phi_0(\mathbf m_n)\! \to \! \Phi_0(\mathbf m)$ in law in $(\bbT^{0}_{\! c}, \dGH)$, this implies that $\Phi^1_0(\mathbf m')$ has the same law as $\Phi_0 (\mathbf m)$ and thus as $\Phi^1_0 (\spaan \, (\mathbf m'))$. 
Since clearly $\Phi^1_0 (\spaan \, (\mathbf m')) \preceq \Phi^1_0(\mathbf m')$, Lemma \ref{orderdeter} $(iii)$ implies that 
$\Phi^1_0 (\spaan \, (\mathbf m'))\eqo \Phi^1_0(\mathbf m')$ almost surely, which entails that $\spaan \, (\mathbf m')\eqo \mathbf m'$. Thus $\mathbf m'$ and $\Phi_1(\mathbf m)$ have the same law and the law of $\Phi_1(\mathbf m)$ is therefore the unique $\dGHP$-weak limit point of the laws of $(\Phi_1(\mathbf m_n))_{n\in \bbN}$, which yields the desired result. \cqfd 

\section{Application to mass erasure of Galton--Watson and Lévy trees}
\label{secappl}

\subsection{Measured Galton--Watson $\bbR$-trees with regenerative branching property} \label{applGW}
We introduce the class of measured finite-type random trees enjoying the \emph{regenerative branching property} (roughly speaking, trees such that conditionally given the subtree below height $r$, the subtrees above height $r$ are i.i.d.~copies of the original tree). When there is no measure, it has been proved by Weill \cite{Weill} (see also \cite{DuWi2}) that the only finite-type trees enjoying the (purely metric) regenerative property are Galton--Watson trees with i.i.d.~exponential edge lengths. In the context of trees equipped with a measure, regenerative finite-type trees are Galton--Watson trees with i.i.d.~exponential edge lengths equipped with specific measures on edges and specific weights at branch points and leaves, which possibly adds 
extra randomness (see Definition \ref{meaGWregedef}).  
To introduce these laws on the space $(\bT_c^1, \dGHP)$ of isometry classes of compact rooted and 
measured $\bbR$-trees, it is convenient to construct first a model in a specific $\boldsymbol{\ell}_1$-space 
as explained below (here we follow \cite{RW1} and a variant of Aldous's \cite{Al1,Al91,Al2} idea of studying trees as subsets of $\boldsymbol{\ell}_1(\bN)$).

We first briefly recall Ulam's formalism on \emph{rooted ordered trees}. 
We denote by $\mathbb{U}\eqo\bigcup_{n\in\bN}(\bN^*)^n$ the set of finite integer-words, with the convention that 
$(\bN^*)^0\eqo  \{ \varnothing \}$ ($\varnothing $ is called the empty word). For all $u\eqo (u_1, \ldots, u_n)\ino \mathbb{U}\backslash \{ \varnothing \}$, we set $\overleftarrow{u} \! :=\!  (u_1, \ldots, u_{n-1})$ for the \emph{direct parent} of $u$
($\overleftarrow{u} \eqo \varnothing$ if $n\eqo 1$). For any $v\eqo  (v_1, \ldots, v_m)\in \mathbb U $, $u\ast v\! :=\!   (u_1, \ldots, u_n, v_1, \ldots, v_m) $ stands for the concatenation of $u$ and $v$. 
\begin{definition}
\label{ellonemodel} 
\begin{compactenum}
\item[$(a)$] We denote by $\mathbb{T}_{\! \mathtt{dis}}$ the countable set of \emph{finite rooted ordered trees}, namely finite subsets $\mathbf{t}\!\subset\!\mathbb{U}$ which satisfy the following.
\begin{compactenum}

\smallskip

\item[$(a1)$] For all $u\ino \mathbf{t} \backslash \{ \varnothing \} $, $\overleftarrow{u} \ino \mathbf t$.

\smallskip

\item[$(a2)$] For all $u\ino \mathbf t$, there is $\mathbf k_u(\mathbf t) \ino \bN$ such that $u\ast (j) \ino \mathbf t$ iff $1\leqo j \leqo \mathbf k_u(\mathbf t)$ (if $\mathbf k_u(\mathbf t) \! \neq \! 0$). Here $\mathbf k_u(\mathbf t)$ is interpreted 
as the \emph{number of children} of $u$ in $\mathbf{t}$.
\end{compactenum}

\smallskip

\item[$(b)$] We equip $\boldsymbol{\ell}_1(\mathbb{U}) \! :=\!  \big\{ \mathbf{s} \eqo (s(u))_{u\in\mathbb{U}}\ino\bbR^{\mathbb{U}}\! 
\colon \lVert \mathbf s \rVert_1\! :=\! \sum_{u\in\mathbb{U}} |s (u)|\! <\! \infty\big \}$ with the metric $d_{\ell_1}$ induced by the $\ell_1$-norm $ \lVert \cdot \rVert_1$. 
We view the origin $\mathbf{0}$ as a root. Let $\mathbf{e}_u\ino \boldsymbol{\ell}_1(\mathbb{U})$ be such that 
$\mathbf{e}_u(u)\eqo 1$ and $\mathbf{e}_u(v)\eqo 0$ 
for all $u,v\ino\mathbb{U}$, $u\!  \neq \! v$. For all $\sigma, \sigma'\ino \boldsymbol{\ell}_1(\mathbb{U})$, we 
set $\lgeo \sigma, \sigma'\rgeo_1\eqo  \big\{ \sigma \! +\!  s\! \cdot \! (\sigma'\! -\! \sigma)\, ; \, s\ino [0, 1] \big\}$.

\smallskip

\item[$(c)$] Let $\mathbf t\ino \mathbb{T}_{\! \mathtt{dis}}$. For $u\ino \mathbf t$, we fix an \emph{edge length} $\ell (u)\ino \bbR_+$, a \emph{weight} $w_u\ino \bbR_+$ and an \emph{edge measure} $\mu_u \ino \cM_{\! f} (\bbR_+)$ such that $\mathtt{Supp} \, \mu_u \subo [0, \ell (u)]$.  
\begin{compactenum}

\smallskip

\item[$(c1)$] We define the compact $\bbR$-tree $T \! := \! \lgeo \mathbf 0, \sigma_\varnothing \rgeo_1 \! \cup  \bigcup_{u\in \mathbf t \backslash \{ \varnothing \}}  \lgeo   \sigma_{\overleftarrow{u}} , \sigma_u \rgeo_1 $, where 
the $\sigma_{\! u} \ino \boldsymbol{\ell}_1(\mathbb{U})$ are defined by 
$\sigma_\varnothing \eqo \ell(\varnothing) \! \cdot \!  \mathbf{e}_\varnothing $ and $\sigma_u \eqo \sigma_{\overleftarrow{u}} + \ell (u) \! \cdot \! \mathbf e_u $, for all $u\ino \mathbf t \backslash \{ \varnothing \}$ (with the convention  $ \sigma_{\overleftarrow{\varnothing}}\eqo \mathbf 0$). 

\smallskip

\item[$(c2)$] $T$ is equipped with the finite measure $\mu$ given by  
$$ \int_T\!   f(\sigma) \, \mu(d\sigma)= \! 
\sum_{u\in \mathbf t }  \int_{\bbR_+} \!\!\!  f(\sigma_{\overleftarrow{u}}+x\! \cdot \! \mathbf e_u) \, \mu_u (dx)+ \sum_{u\in \mathbf t} w_u f(\sigma_u) , $$
for all bounded measurable $f\colon T\!  \to\!  \bbR$. 

\smallskip

\end{compactenum}
We denote by $\widetilde{\mathtt T}\mathtt{ree} \big( \mathbf t \, ; \, (\ell(\cdot ), w_\cdot , \mu_\cdot) \big)$ the isometry class 
in $\bbT^1_{\! c}$ of the rooted compact measured $\bbR$-tree $\mathtt{Tree} \big( \mathbf t \, ;\,  (\ell(\cdot ), w_\cdot , \mu_\cdot) \big) \! :=\!  (T,d_{\boldsymbol{\ell}_1}, \mathbf 0, \mu)$.  \cq 
\end{compactenum}
\end{definition}

\noi
It is convenient to introduce 
\begin{align} 
\label{T1discr}
\mathbb{T}^1_{\! \mathtt{dis}}=\!\!\!  \bigsqcup_{\mathbf t \in \mathbb{T}_{\mathtt{dis}}} \!\!\! 
(\{ \mathbf t \} \! \times\!  \mathscr P_{\mathbf t} ) \quad  &  \textrm{where} \\[-3mm] \mathscr P_{\mathbf t} = & \Big\{  \big(\ell(u), w_u, \mu_u\big)_{\! u\in\mathbf t} \! \in \!\big( \bbR_+^2\! \times \! \mathcal M_f (\bbR_+)\big)^{\! \mathbf t}\! : \mathtt{Supp}\,  \mu_u \! \subseteq \! [0, \ell(u)], u \! \in \!  \mathbf t \Big\}.\nonumber 
 \end{align}
We equip $ \mathbb{T}^1_{\! \mathtt{dis}}$ with the Polish topology where 
$(\mathbf t_p ;(\ell_p(u), w_{p,u}, \mu_{p,u})_{u\in \mathbf t_p})$ $\! \to\! $ $(\mathbf t;   $ $(\ell(u), w_{u}, \mu_{u})_{u\in \mathbf t})$ iff for all sufficiently large $p$, 
$\mathbf{t}_p\eqo  \mathbf t$, and for all $u\ino \mathbf t$, $\lim_{p\to \infty} (\ell_p(u), w_{p,u})\eqo (\ell(u), w_{u})$ and $\mu_{p,u} \! \to \! \mu_u$ weakly in $\mathcal M_f(\bbR_+)$. \emph{Indeed,} since  $ \mathtt{Supp}\,  \mu_{p,u} \! \subseteq \! [0, \ell_p(u)]$ for all $p\ino \bbN$, and since $\ell_p (u) \! \to \! \ell (u)$ and 
$\mu_{p, u} \! \to \! \mu_u$ weakly on $\bbR_+$, the Portmanteau theorem shows that 
$\mu_u ((\varepsilon+ \ell (u), \infty) ) \leqo \liminf_{p\to \infty} \mu_{p,u} ((\varepsilon+ \ell (u), \infty) ) \eqo 0$ for all $\varepsilon \ino \bbR^*_+$, which implies that $\mathtt{Supp}\,  \mu_u \! \subseteq \! [0, \ell(u)]$. We easily check that 
$\mathtt{Tree}$ is continuous with respect to the Hausdorff distance on compact subsets of $\boldsymbol{\ell}_1(\mathbb U)$ and weak convergence of finite Borel measures on $\boldsymbol{\ell}_1(\mathbb U)$, which implies that 
\begin{equation}
\label{contTree}
\textrm{$\widetilde{\mathtt{T}}\mathtt{ree}\colon    \mathbb{T}^1_{\! \mathtt{dis}} \!\! \to \! \mathbb{T}^1_{\! c}\ $ is continuous, and hence so is $\ \widetilde{\mathtt{T}}\mathtt{ree}_0\! :=\!  \Phi^1_0 \! \circ\!  \widetilde{\mathtt{T}}\mathtt{ree} \colon   \mathbb{T}^1_{\! \mathtt{dis}} \!\!\to\!  \mathbb{T}^0_{\! c}$,}
\end{equation}
where we recall from Definition \ref{TcTmindef} $(a)$ that $\Phi^1_0\!: \! \bbT^1_{\! c} \!\! \to \! \bbT^0_{\! c} $ is the $1$-($\dGHP, \dGH$)-Lipschitz function consisting in forgetting the measure. 

 We now introduce the main class of measured GW-trees and forests. For all $\mathtt m \ino \mathcal M_f (\bbR_+)$ and all bounded  measurable $f\colon   \bbR_+ \! \to \! \bbR$, we adopt the notation $\langle\mathtt m,f\rangle\eqo\int_{\bbR_+}\!f(x)\,\mathtt m(dx)$.
\begin{definition}
\label{meaGWregedef} We fix $\xi$ and $\varrho$, two probability distributions on $\bN$. 
We assume that $\xi$ is (sub)critical: $\sum_{k\in \bbN} \xi (k) \leqo 1$. 
We fix $c\ino \bbR^*_+$ and a (conservative) c{\`a}dl{\`a}g subordinator $S\eqo(S_s)_{s\in \bbR_+}$ 
whose Laplace exponent $\phi (y)\eqo  -\! \log \bE [e^{-yS_1}]$ is given by 
\begin{equation}
\label{Khint}
\forall y\ino \bbR_+, \quad \phi(y)= \kappa y + \int_{\bbR^*_+} \!\! \! \! \big( 1\! -\! e^{-yz} \big) \, \Gamma (dz) , 
\end{equation} 
where $\kappa\ino \bbR_+$ and $\int_{\bbR^*_+}\!  \min (1, z)\,  \Gamma (dz) \! <\!  \infty$. For all $n\ino \bN$, we fix $\Lambda_n, \Lambda_{\varnothing, n} \ino \cM_{ 1} (\bbR_+)$ and we write $\Lambda\eqo(\Lambda_n)_{n\in\mathbb{N}}$ and $\Lambda_\varnothing\eqo(\Lambda_{\varnothing,n})_{n\in\mathbb{N}}$.
\begin{compactenum}

\smallskip

\item[$(a)$] A random variable $\tau\colon \Omega\!  \to\!  \mathbb{T}_{\!  \mathtt{dis}}$ is a ($\xi, \varrho$)-Galton--Watson forest (a GW($\xi, \varrho$)-forest for short) if $\bP (\tau \eqo \mathbf t)\eqo \varrho (\mathbf k_\varnothing(\mathbf t))\prod_{u\in \mathbf t\backslash \{ \varnothing\}} \xi (\mathbf k_u (\mathbf t))$ for all $\mathbf t \ino  \mathbb{T}_{\! \mathtt{dis}}$. If $\varrho \eqo \xi$, $\tau$ is a $\xi$-Galton--Watson tree  (a GW($\xi$)-tree for short). 

\smallskip

\item[$(b)$] We denote by $Q_{\xi, c, \phi, \Lambda} $ the law of $\widetilde{\mathtt{T}}\mathtt{ree}  (\tau\, ; (\ell (\cdot), w_\cdot, \mu_\cdot))$ where $\tau$ is a GW($\xi$)-tree and conditionally given $\tau$, the r.v.s
$(\ell (u),w_u,\mu_u)$, $u\ino \tau$, are independent and distributed as follows: for all $u\ino \tau$, all measurable and bounded $f_1, f_2 \colon   \bbR_+ \! \to\!  \bbR$, $f_3\colon \mathcal M_{\! f}(\bbR_+) \! \to\!  \bbR$, $\bP$-a.s.\vspace{-1mm}
\begin{equation}
\label{lawatu} 
\bE \big[ f_1(w_u)f_2(\ell(u))f_3 (\mu_u) \,\big| \, \tau \big]\!=\! \langle \Lambda_{\mathbf k_u (\tau)} ,  f_1 \rangle   \! \int_0^\infty \!\! \! dt \, ce^{-ct} f_2 (t) \, \bE \big[ f_3 \big(\un_{[0, t]} (s) dS_s  \big)\big].\vspace{-1mm}
\end{equation}
\item[$(c)$] We denote by $P_{\xi, c, \phi, \Lambda, \varrho, \Lambda_{\varnothing}} $ the law of $\widetilde{\mathtt{T}}\mathtt{ree}  (\tau; (\ell (\cdot), w_\cdot, \mu_\cdot))$ where $\tau$ is a GW($\xi, \varrho$)-forest and conditionally given $\tau$ the r.v.s
$(\ell (u),w_u,\mu_u)$, $u\ino \tau$, are independent and distributed as follows: $\bP$-a.s.~$\mu_\varnothing\eqo  0$, $\ell (\varnothing)\eqo 0$ and $\bE [ f_1 (w_\varnothing) \, | \, \tau]\eqo  \langle \Lambda_{\varnothing , \mathbf k_{\varnothing} (\tau)} , f_1   
\rangle $ and for all $u\ino \tau \backslash \{\varnothing \}$, (\ref{lawatu}) holds true.

\smallskip

\item[$(d)$] We next denote by $\texttt Q_{\xi, c}$ the law on $(\bbT^0_{\! c} , \dGH)$ 
of $\Phi^1_0$ under $Q_{\xi, c, \phi, \Lambda}$, which is the law of a purely metric GW($\xi, c$)-tree (we simply forget about the measure component under $Q_{\xi, c, \phi, \Lambda}$). Similarly we denote by $\texttt P_{\! \xi, c, \varrho}$ the law on $(\bbT^0_{\! c} , \dGH)$ 
of $\Phi^1_0$ under $P_{\xi, c, \phi, \Lambda, \varrho, \Lambda_{\varnothing}} $. \cq 
\end{compactenum}
\end{definition}

\begin{remark}
\label{minimalGW} Using the notation of Definition \ref{meaGWregedef}, we observe that 
$Q_{\xi, c, \phi, \Lambda} (\bbT^1_{\!c , \texttt{min}})\eqo 1$ or equivalently 
$P_{\xi, c, \phi, \Lambda, \varrho, \Lambda_{\varnothing}} (\bbT^1_{\! c, \texttt{min}})\eqo 1$ iff either $\kappa \geko 0$, or $\Gamma ((0, 1))\eqo \infty$ or $\Lambda_{0} (\{ 0\}) \eqo 0$. \cq 
\end{remark}
\begin{lemma}\label{poublm2} For all $p\ino \bbN \! \cup \! \{ \infty\}$ let us fix parameters 
$\xi_p$, $c_p$, $\phi_p$, $(\Lambda^{{\! (p)}}_{k})_{k\in \bbN}$, $\varrho_p$ and 
$(\Lambda^{{\! (p)}}_{\varnothing, k})_{k\in \bbN}$ as in Definition \ref{meaGWregedef}. We assume that 
$c_p\! \to \! c_\infty$, that 
$\xi_p\! \to \xi_\infty$ and $\varrho_p\! \to \varrho_\infty$ weakly on $\bbN$, that for all $k\ino \bbN$ such that $\xi_{\infty}(k)\! \neq \! 0$, $\Lambda^{{\! (p)}}_{k}\! \! \to \! \Lambda^{{\! (\infty)}}_{k}$ weakly on $\bbR_+$, 
that for all $k\ino \bbN$ such that $\varrho_{\infty}(k)\! \neq \! 0$,
$\Lambda^{{\! (p)}}_{\varnothing, k}\!\!  \to \! \Lambda^{{\! (\infty)}}_{\varnothing, n}$ weakly on $\bbR_+$,  
and that $\phi_p \! \to \! \phi_\infty$ pointwise on $\bbR_+$. Then,  
$$ Q_p\! :=\! Q_{\xi_p, c_p, \phi_p, \Lambda^{\! (p)}_\cdot}\! \! \longrightarrow \! Q_\infty \! :=\! Q_{\xi_\infty, c_\infty, \phi_\infty, \Lambda^{\! (\infty)}_\cdot}\; \textrm{and} \;  
P_{\! \xi_p, c_p, \phi_p, \Lambda^{\! (p)}_\cdot \!\!  , \, \varrho_p,\Lambda^{{\! (p)}}_{\varnothing, \cdot} } \!\!  \longrightarrow \! 
P_{\! \xi_\infty, c_\infty, \phi_\infty, \Lambda^{\! (\infty)}_\cdot\!\!  ,\,   \varrho_\infty,\Lambda^{{\! (\infty)}}_{\varnothing, \cdot}}, $$ 
weakly in 
$(\bbT^1_{\! c}, \dGHP)$.  
\end{lemma}
\noi 
\textbf{Proof.} We recall Definition \ref{ellonemodel} $(a)$ of $\bbT_{\! \mathtt{dis}}$ and we recall from (\ref{T1discr}) the definition of the Polish space $\bbT_{\! \mathtt{dis}}^1$. For all $p\ino \bbN\cup \{ \infty\}$, let $\tau_p$ be a GW($\xi_p$)-discrete tree. Namely, for all $\mathtt t\ino \bbT_{\! \mathtt{dis}}$,  $\bP (\tau_p \eqo \mathtt t)\eqo \prod_{u\in \mathtt t} 
\xi_p (\mathbf k_u (\mathtt t))$. For all $p\ino \bbN \cup \{ \infty\}$ and all $\mathbf t\ino \bbT_{\! \mathtt{dis}}$, we denote by 
$\mathscr Q_{p,\mathtt t}$ the law on $\bbT_{\! \mathtt{dis}}^1$ of 
$(\mathtt t; (\ell^p_u , w^p_u , \mu^p_u)_{u\in \mathtt t})$, where the $(\ell^p_u , w^p_u , \mu^p_u)_{u\in \mathtt t}$ 
are independent with respective laws $c_p e^{-c_p \ell} d\ell \Lambda^{{\! (p)}}_{\mathbf k_u (\mathtt t) } (dw) 
\un_{[0, \ell^p_u]} (s) dS_s^p$, where $(S^p_s)_{s\in \bbR_+}$ is a subordinator with Laplace exponent $\phi_p$. 
If $\bP(\tau_\infty\eqo \mathtt t)\! \neq \! 0$, our assumptions imply that 
$\mathscr Q_{p,\mathtt t} \!\!  \to \! \! 
\mathscr Q_{\infty, \mathtt t}$ in law  on 
$\bbT_{\! \mathtt{dis}}^1$. Since $\lim_{p\to \infty} $ $\bP (\tau_p \eqo \mathtt t) $ $\eqo $ $
\prod_{u\in \mathtt t} \xi_\infty (\mathbf k_u (\mathtt t)) $ $\eqo $ $ \bP (\tau_\infty \eqo \mathtt t)$, we get  
$ \sum_{\mathbf t\in\bbT_{\! \mathtt{dis}}^1 } \! 
|  \bP (\tau_p \eqo \mathtt t)\! -\!  \bP (\tau_\infty \eqo \mathtt t)| \! \to \!  0$, by standard arguments. Thus, 
$\sum_{\mathbf t\in\bbT_{\! \mathtt{dis}}^1}\!  \bP (\tau_p \eqo \mathtt t)
\mathscr Q_{p,\mathtt t}  \! \! \to \! \! 
\sum_{\mathbf t\in\bbT_{\! \mathtt{dis}}^1} \! \bP (\tau_\infty \eqo \mathtt t)
\mathscr Q_{\infty, \mathtt t}$ in law too and we obtain the desired result regarding the weak 
convergence of the $Q_{p}$ since they are the push-forward measure of 
$\sum_{\mathbf t\in\bbT_{\! \mathtt{dis}}^1} \bP (\tau_p \eqo \mathtt t)
\mathscr Q_{p, \mathtt t} $ via the continuous function 
$\widetilde{\mathtt T}\mathtt{ree}$. The proof concerning GW forests is similar: we leave the details to 
the reader. \cqfd

 \smallskip

One of the main goals of this section is to reformulate Definition \ref{meaGWregedef} in a more intrinsic way 
and to characterize it by the regenerative property. To this end, we need next to introduce the following notation.

\smallskip

\noi
$\bullet$ For $i\ino \{ 1,2\}$, let $(T_i, d_i, \rho_i, \mu_i) \! \equiv\! \fmu_i\ino \mathbb{T}^1_c$. 
We denote by $(T^\circ, d^\circ, \rho^\circ, \mu^\circ)$ the trees $T_i$ \emph{fused at their roots}: 
$T^\circ\!\! =\!  \{ \rho^\circ\} \sqcup (T_1\backslash \{ \rho_1\}) \sqcup (T_2\backslash \{ \rho_2\})$ equipped with the metric 
$d^\circ\!$ which extends the $d_i$ and which satisfies
$d^\circ (\rho^\circ\! , \sigma_i) \! =\!  d_i(\rho, \sigma_i)$ and $d^\circ (\sigma_1, \sigma_2)\! =\!  d_1(\rho_1, \sigma_1)\! +\! d_2(\rho_2, \sigma_2)$ for all 
$\sigma_i \!  \in \! T_i\backslash \{ \rho_i\}$, $i\ino \{ 1,2\}$, and with the measure 
$\mu^\circ \!\! =\!  \big(  \mu_1 (\{\rho_1 \})\! +\! \mu_2(\{\rho_2 \}) \big)\delta_{\! \rho^\circ}\! + \! \sum_{i\in \{ 1,2\}} 
\mu_i \big(  \cdot  \cap (T_i\backslash \{ \rho_i\}) \big)$. 

We denote its isometry class by $\fmu_1\circledast \fmu_2$, which makes sense since it only depends on the 
$\fmu_i$. We check easily that $\circledast$ is 
$\dGHP$-continuous from $(\mathbb{T}^1_c)^2$ to $\mathbb{T}^1_c$ and that it is an associative relation.

\medskip

\noi
$\bullet$ Let $Q_1, Q_2\ino \mathcal M_1(\mathbb{T}^1_c)$: we denote by $Q_1\circledast Q_2$ the law of 
$\fmu_1\circledast \fmu_2$ under $Q_1(d\fmu_1) \otimes Q_2( d\fmu_2)$. For all $n\ino \bN$, we recursively define $Q^{\circledast n}$ by setting $Q^{\circledast 0}\eqo  \delta_{\Upsilon_{\! 0}}$, where $\Upsilon_{\! 0}$ istands for \emph{the tree reduced to a point, equipped with the null measure}, 
and $Q^{\circledast n+1}\eqo Q^{\circledast n}\circledast Q\eqo  Q \circledast Q^{\circledast n}$. 
For all $w\ino \bbR_+$, we denote by $\Upsilon_{\! w}$ the isometry class of $(\{ \rho\}, d, \rho, w\delta_{\rho})$ and we check that $w \! \mapsto \! \Upsilon_{\! w}$ is $\dGHP$-continuous. 
Definition \ref{meaGWregedef} $(b)$ then entails 
\begin{equation}
\label{forestintrin}
P_{\xi, c, \phi, \Lambda , \varrho, \Lambda_{\varnothing}} = \sum_{k\in \bbN} \varrho (k) \!  \int_{\bbR_+} \!\!\!\!\Lambda_{\varnothing, k} (dw) \, \delta_{\Upsilon_{\! w}} \!\!  \circledast   Q^{\circledast k}_{\xi, c, \phi, \Lambda}.
\end{equation}

\noi
$\bullet$ Let $(T,d,\rho, \mu) \! \equiv\!\fmu\ino \mathbb{T}^1_c$. To simplify the notation, for all 
$r\ino \bbR_+$, we set $B_r \eqo B_{T,d} (\rho, r)$ and $B^c_r\eqo T\backslash 
B_r$. We denote by $Z^+_r (\fmu) \ino \bN\cup \{ \infty\}$ the number of connected components of $B_r^c$. 
We then define $\mathtt{Blw}_r (\fmu)$ as the isometry class in $\mathbb{T}^1_c$ of 
$(B_r, d, \rho, \mu_{\leq r})$, where $\mu_{\leq r}\! :=\!  \mu (\, \cdot \, \cap B_r))$. This makes sense since it only 
depends on $r$ and $\widetilde{T}$. 

Similarly, we also introduce the tree $\mathtt{Abv}_r (\fmu)$ which is obtained by grafting at their roots the 
connected components of $B_r^c$. More precisely, for all $\sigma \ino T$, we set 
$f_r (\sigma) \eqo (d(\rho, \sigma)\! -\! r)_+$. 
On $B_r^c$ we define the $0$-hyperbolic distance 
$d_r (\sigma, \sigma')\eqo  f_r(\sigma) + f_r(\sigma) \! -\! 2\min_{\gamma \in \lgeo
 \sigma, \sigma'\rgeo_T}  f_r(\gamma) $. We observe that the $d_r$-completion $\overline{B}_r^c$ of 
 $B^c_r$ is a compact $\bbR$-tree which differs from $B_r^c$ by a single point that we denote by $\rho_r$. 
We then define $\mathtt{Abv}_r (\fmu)$ as the isometry class of $(\overline{B}_r^c, d_r, \rho_r, \mu_{>r} )$ 
where $\mu_{>r}\! =\! \mu (\, \cdot \, \cap B_r^c))$. This makes sense since it only depends on $r$ and $\fmu$.
It is not difficult to prove (e.g.~using Lemma \ref{extendist}) that for all $a,b\ino\bbR_+$ such that 
$a\leko b$, 
$$  \fdelta_{\mathtt{GHP}} \big( {\mathtt{Abv}}_a(\fmu),{\mathtt{Abv}}_b(\fmu) \big)
\leq 2(b\! -\! a)+\flambda_\fmu((a,b]),  $$
which entails that $r\! \mapsto \! {\mathtt{Abv}}_r(\fmu)$ is right-continuous. Here, we recall that 
$\flambda_{\fmu} \ino \mathcal M_{\! f} (\bbR_+)$ is given by $\langle \flambda_{\fmu}, f \rangle \eqo \int_{T}\!  f(d(\rho, \sigma))\,  \mu (d\sigma) $ for all bounded measurable $f\colon\bbR_+ \! \to \! \bbR$ 
(see Notation \ref{GPmeradef} $(b)$).

\medskip

\noi
$\bullet$ We  set $D(\fmu)  \eqo \inf\{d(\rho,\gamma);\, \sigma\ino ({\mathtt{Lf}}(T)\cup{\mathtt{Br}}(T))\backslash \{ \rho \}\}$, with the convention that $\inf \emptyset \eqo \infty$ (this quantity $D(\fmu)$ may also be equal to zero). 
It is the height of the first leaf or branch point nearest to the root. 
We also denote by $\mathbf{k}(\fmu)\ino \bN\cup\{\infty\}$ the number of connected components of 
the open subset $\{ \sigma \ino T \! : d(\rho, \sigma ) \! >\! D(\fmu)\}$. If $D(\fmu)\! <\! \infty$, we observe that 
$\mathbf{k}(\fmu) \eqo  Z^+_{{D(\fmu)} } (\fmu)$. We also recall that $\langle \fmu \rangle \eqo \mu(T)$ is the total mass of $\fmu$. 

 \smallskip

The following lemma allows to provide an intrinsic definition of $Q_{\xi, c, \phi, \Lambda}$
 \begin{lemma} 
\label{measwithmeas} The following holds true. 
\begin{compactenum}

\smallskip

\item[$(i)$] 
$D\colon   \mathbb{T}^1_c\! \to \! [0, \infty]$, $\mathbf k\colon   \mathbb{T}^1_c\! \to \! \bN\cup \{ \infty\}$, 
$Z^+_r  \colon   \mathbb{T}^1_c\! \to \! \bN\cup \{ \infty\}$, $\flambda \colon   \mathbb{T}^1_c\! \to \! \mathcal M_{\! f} (\bbR_+)$, 
$(r, \fmu)\ino \bbR_+ \! \times \! \mathbb{T}^1_c \! \mapsto \!   \mathtt{Abv}_r (\fmu) \ino  \mathbb{T}^1_c$ and $(r, \fmu)\ino \bbR_+ \! \times \! \mathbb{T}^1_c \! \mapsto \!   \mathtt{Blw}_r (\fmu) \ino  \mathbb{T}^1_c$  are Borel-measurable.

\smallskip

\item[$(ii)$] Let $\xi, c, \phi, \Lambda\eqo(\Lambda_n)_{n\in \bbN}$ and $S\eqo(S_s)_{s\in\bbR_+}$ be as in Definition \ref{meaGWregedef}.  
In addition, we assume that $\xi$ is proper, i.e.~that $\xi(1)\eqo 0$. Let $Q\in \mathcal M_1(\mathbb{T}^1_c)$. Then $Q=Q_{\xi, c, \phi, \Lambda}$ iff 
\begin{eqnarray}\label{df:GWweight}
Q \big[ f_1 \big( D, \flambda ([0, D)) \big)\!\!\!\!\!\! \!\! & &\!\!\!  \!\!\!\!  f_2 (\flambda (\{ D\}) ) \un_{\{ \mathbf k=n\}} f_3 \big( \mathtt{Abv}_{D}\big) \big] \\
\!\!\!\!\!\!&=& \!\!\!\! \int_0^\infty \!\!  dt \, ce^{-ct} \bE \big[ f_1 \big( t, \un_{[0, t]} (s) dS_s  \big) \big] \xi (n) \, \langle \Lambda_n , f_2\rangle\,  Q^{\circledast n} [ f_3] \nonumber 
\end{eqnarray}
for all bounded measurable $f_1, f_2, f_3$ and all $n\ino \bbN$.  
\end{compactenum}
\end{lemma}
\noi
\textbf{Proof.} See Appendix \ref{pfsecmeaswithmeas}. \cqfd 

\smallskip

We next characterise the law of the masses of balls of a measured GW tree, as a process. To this end, we introduce the following notation. 
We denote Laplace transforms as follows. 
\begin{equation}
\label{notaLapla}
 \forall \mu \ino \mathcal M_{\! f} (\bbR_+), \, \forall y \ino\bbR_+, \quad \widehat{\mu} (y)= \!\! \int_{\bbR_+} \!\!\! e^{-yz} \, \mu(dz) \; .
 \end{equation}
Let $\xi , c, \phi $, $(\Lambda_n)_{n\in \bN}$, $\varrho$ and  $(\Lambda_{ \varnothing,n})_{n\in \bN}$ be as in Definition \ref{meaGWregedef}. 
For all $s\ino [0, 1]$ and 
$y\ino \bbR_+$, we set 
\begin{equation}
\label{masterfunc}
F_{\xi,c, \Lambda } (s,y) \eqo c \sum_{n\in \bN}
 \xi (n) s^n\widehat{\Lambda}_n(y) \quad \textrm{and} \quad F^{\varnothing}_{\! \varrho, \Lambda_{\varnothing }} (s,y)\eqo 
\sum_{n\in \bN} \varrho (n) s^n\widehat{\Lambda}_{\varnothing, n}(y).
\end{equation}
We note that $F_{\xi,c, \Lambda}$  (resp.~$F^{\varnothing}_{\! \varrho, \Lambda_{\varnothing}}$) 
characterises $\xi,c$ and the 
$\Lambda_{k}$ such that $\xi(k) \! >\! 0$ (resp.~$\varrho$ and the $\Lambda_{\varnothing , k}$ 
such that $\varrho(k) \! >\! 0$). Moreover, we see that the laws of measured GW trees and forests depend continuously on these transforms as shown by the following proposition used in Section \ref{pfsecthmIP}.
\begin{proposition}
\label{proppoubelle} 
For all $p\ino \bbN $ let $\xi_p$, $c_p$, $\phi_p$, $\Lambda^{\!(p)}\eqo(\Lambda^{{\! (p)}}_{n})_{n\in \bbN}$, $\varrho_p$ and 
$\Lambda_{\varnothing}^{\!(p)}\eqo(\Lambda^{{\! (p)}}_{\varnothing, n})_{n\in \bbN}$ be as in Definition \ref{meaGWregedef}. Further suppose that $\xi_p$ is proper for all $p\in\bbN$. 
We assume that, as $p \!\rightarrow \! \infty$,
\begin{compactenum}
\smallskip
\item[$(a)\!$] $F_{\xi_p,c_p,\Lambda^{{\! (p)}} } \!\!  \rightarrow \! F$ pointwise on $[0,1]\! \times\! \mathbb{R}_+$, for 
a limit $F\! \not \equiv\!  0 $ that is continuous at the point $(1,0)$;%
\smallskip
\item[$(b)\!$] $\phi_p \! \rightarrow\! \phi$ pointwise on $\mathbb{R}_+$, for a limit $\phi$ that is continuous in a neighbourhood of $0$;
\smallskip
\item[$(c)\!$] $F^\varnothing_{\!\! \varrho_p,\Lambda^{{\! (p)}}_{\varnothing, \cdot} } \!\! \rightarrow \! F^\varnothing$
pointwise on $[0,1]\times\mathbb{R}_+$, for a limit $F^\varnothing$ that is continuous at the point $(1,0)$.
\smallskip
\end{compactenum}
Then there are admissible GW parameters $(\xi,c,\phi,\Lambda ,\varrho,\Lambda_{\varnothing})$ as in Definition \ref{meaGWregedef} such that 
the convergence $P_{\! \! \xi_p ,c_p,\phi_p,\Lambda^{{\! (p)}} \!\! ,\, \varrho_p,\Lambda^{{\! (p)}}_\varnothing }$ $\! \rightarrow \! 
P_{ \xi,c,\phi,\Lambda ,\varrho,\Lambda_{\varnothing}}$ holds  weakly
in $(\mathbb{T}_{\! c}^1,\dGHP)$. 

Moreover, the same result holds true for unmeasured GW forests: If $c_p\varphi_{\xi_p} \! \to \! F$  pointwise on $[0,1]$ with $F\! \not \equiv\!  0 $ continuous at the point $1$, and if $\varphi_{\varrho_p} \! \to \! f$ pointwise on $[0,1]$ with $f$ continuous at the point $1$, then 
there are admissible GW parameters $(\xi,c,\varrho)$ as in Definition \ref{meaGWregedef} such that the convergence 
$\texttt P_{ \xi_p ,c_p, \varrho_p}$ $\! \rightarrow \! 
\texttt P_{ \xi,c, \varrho}$ holds  weakly
in $(\mathbb{T}_{\! c}^0,\dGH)$.  
\end{proposition}
\noi
\textbf{Proof.} We set $c\! :=\! F( 1, 0)$ and we first prove it is positive. Indeed, by $(a)$, $F$ has to be non-decreasing in its first coordinate and non-increasing in its second one. Thus $c\eqo \max \{ F(s,y); s\ino [0, 1], y\ino \bbR_+ \}$, and this is positive since we assume that $F\! \not \equiv \! 0$. By $(a)$ again, we also see that $c_p\eqo F_{\xi_p,c_p,\Lambda^{{\! (p)}} } (1,0) \! \to \! F(1,0)\eqo c$. To simplify notation we set $F_{\! p} \! :=\! \frac{1}{c_p} F_{\xi_p,c_p,\Lambda^{{\! (p)}} }$ and 
$F^*\! :=\! \frac{1}{c} F$. The previous argument combined with $(a)$ then entails that $F_{\! p}\! \to \! F^*$ pointwise on $[0,1]\! \times\! \mathbb{R}_+$, that $F^*$ is continuous at the point $(1, 0)$ and that $F^*(0,1)\eqo 1$. 

 We fix $y\ino \bbR_+$ and observe that $(\xi_p (k) \widehat{\Lambda}^{{\! (p)}}_k(y))_{k\in \bbN}$ belongs 
 to the compact metrisable product space $[0,1]^{\bbN}$. Thus, it has at least one limit point, say $(g_k(y))_{k\in \bbN}$. 
 By dominated convergence and since $F^*$ is the pointwise limit of the $F_{\! p}$, we get 
 $\sum_{k\in \bbN} s^k g_k(y)\eqo F^*(s,y)$, for all $s\ino [0, 1)$. 
This shows that $(g_k(y))_{k\in \bbN}$ is actually the unique limit point of the 
$(\xi_p (k) \widehat{\Lambda}^{{\! (p)}}_k(y))_{k\in \bbN}$. Namely, $\xi_p (k) \widehat{\Lambda}^{{\! (p)}}_k(y)\! \to \! g_k(y)$ 
for all $k\ino \bbN$ and all $y\ino \bbR_+$. This implies in particular that the $g_k $ are non-increasing. 
We then set $\xi(k)\eqo g_k(0)$, $k\ino \bbN$. 
By monotone convergence we get $\sum_{k\in \bbN} \xi (k)\eqo \lim_{s\uparrow 1^-}\!  \uparrow \! F^*(s,0)\eqo 1$. 
Since $\xi_p(1)\eqo 0$ and $\sum_{k\in \bbN} k\xi_p(k)\leq 1$ for all $p$, we get $\xi (1)\eqo 0$ and 
$\sum_{k\in \bbN} k\xi (k)\leqo 1$, by Fatou's Lemma. This proves that $\xi_p \! \to \! \xi$ and that $\xi$ is a proper and 
(sub)critical offspring distribution. 

For all $y\ino \bbR_+$ and $k\ino \bbN$, we next set $g^*_k(y)\eqo g_k (y)/\xi (k)$ if $\xi (k)\! \neq \!  0$ 
and $g^*_k(y)\eqo 1$ otherwise. Therefore, $\sum_{k\in \bbN} s^k \xi (k) g^*_k(y)\eqo F^*(s,y)$, for all $s\ino [0, 1)$ and $y\ino \bbR_+$. Since the $g^*_k$ are noncreasing and since $g_k^*(0)\eqo 1$, for all $(s,y)\ino [0, 1)\! \times \! \bbR_+$, we get $1\! -\! g^*_k(y) \geqo 0$ and $F^*(s,0) \! -\! F^*(s,y)\eqo \sum_{k\in \bbN} s^k \xi (k) (1\! -\! g_k^*(y)) \! \geqo 0$. Thus monotone convergence applies to show, as $s\! \uparrow \! 1^-$, that we get $1\! -\! F^*(1,y)\eqo \sum_{k\in \bbN} \xi(k) (1\! -\! g_k^*(y))$. 
Note that the right limit $g^*_k(0^+)\eqo \sup_{y\in \bbR_+} g^*_k (y) $ is well-defined and belongs to $[0, 1]$. Since $\sum_{k\in \bbN} \xi (k)\eqo 1$, dominated convergence applies when $y\! \downarrow \! 0^+$ to show that $\sum_{k\in \bbN} \xi (k) (1\! -\! g^*_k (0^+))\eqo 0$. Therefore for all $k\ino \bbN$, $g^*_k(0^+)\eqo 1\eqo g^*_k (0)$ and $g^*_k$ is continuous at $0$. We have proved that for all $k\ino \bbN$ with $\xi (k)\! \neq \! 0$, $\widehat{\Lambda}^{{\! (p)}}_k \! \to \! g^*_k$ pointwise on $\bbR_+$ and that $g^*_k$ is continuous at $0$ and equal to $1$. By 
standard results on Laplace transforms, for all 
$k\ino \bbN$ with $\xi (k)\! \neq \! 0$, there exists $\Lambda_k\ino \cM_1(\bbR_+)$ such that $ \Lambda^{{\! (p)}}_k\! \to \! \Lambda_k$, weakly on $\bbR_+$. 

 Similar arguments show that there is a probability measure $\varrho$ on $\bbN$ such that $\varrho_p\! \to \! \varrho$ weakly on $\bbN$ and that for all $k\ino \bbN$ with $\varrho (k) \! \neq \! 0$, there is $\Lambda_{\varnothing ,k}\ino \cM_1(\bbR_+)$ such that $ \Lambda^{{\! (p)}}_{\varnothing, k}\! \to \! \Lambda_{\varnothing , k}$, weakly on $\bbR_+$. 
Now Lemma \ref{poublm2} applies and entails the desired result.  The last assertion of the lemma, which concerns the unmeasured case, is proved similarly: we leave the details to the reader.  \cqfd

\smallskip

We now prove the following lemma that characterises the law of the masses of balls of measured GW trees as a process.
\begin{lemma}
\label{ballGWmass} We keep the above notation. 
Let $\xi$ be proper (i.e.~$\xi (1)\eqo 0$) and (sub)critical. Let $(S_r)_{r\in \bbR_+}$, 
$(N_r)_{r\in \bbR_+}$ and $(Y_k,W_k)_{k\in \bN^*}$ be independent processes that are distributed as follows: $S$ is a subordinator with Laplace exponent $\phi$ as in (\ref{Khint}), $N$ is a homogeneous Poisson process with unit intensity, the $(Y_k,W_k)$ are independent 
$(\{ -1\} \! \cup \! \bN) \! \times \! \bbR_+$-valued r.v.s whose laws are given by 
$$  \bE \big[s^{Y_0} e^{-yW_0}  \big] \eqo F^\varnothing_{\!\! \varrho, \Lambda_\varnothing} (s,y) \quad \textrm{and} \quad  \bE \big[s^{Y_k} e^{-yW_k}  \big] \eqo (cs)^{-1}F_{\xi, c, \Lambda} (s,y) , \;  k\ino \bN^*.$$
We next define the following processes: 
$$\forall r\ino \bbR_+, \quad (X^+_r,M_r) \! := \! (Y_0, W_0)+ \!\!\!  \sum_{1\leq k \leq N_{cr}} \!\!\! (Y_k,W_k) \; , \quad 
R_0 \! :=\!   \inf \big\{r\ino \bbR_+ \! : \! X^+_r\eqo 0 \big\}, $$
and  $C_r \! :=\!  R_0 \wedge \inf \{ u\ino \bbR_+ \!  : A_u \! >\! r \}$ where $ A_r \eqo 
 \int_0^r du/X^+_u$ if $r\! < \! R_0$ and $A_r \eqo \infty$ otherwise, and with the convention that $\inf \emptyset\eqo \infty$. Then the following holds true. 

\begin{compactenum}

\smallskip

\item[$(i)$]  $\big( Z^+_r\! , \flambda ([0,r])\big)_{\! r\in \bbR_+}$ under $P_{\xi,c, \phi ,\Lambda, \varrho, \Lambda_{ \varnothing}}$ has the same law as $(X^+_{C_r} , S_{C_r}\! +\!  M_{C_r})_{r\in\bbR_+} $.  In particular 
the total mass $\langle \flambda \rangle $ has the same law as $S_{R_0}+ M_{R_0}$. 

\smallskip

\item[$(ii)$] Let us set 
$g_r( s,y) \eqo Q_{\xi, c, \phi, \Lambda} \big[ s^{Z^+_r} e^{-y\flambda ([0, r])} \big] $ and $\gamma (y)\eqo  Q_{\xi, c, \phi, \Lambda} \big[ e^{-y \langle \flambda \rangle}\big] $, for all $r,y\ino \bbR_+$ and all $s\ino [0, 1]$. Then the following holds true for all $s\ino [0, 1]$ and $y\ino \bbR_+$.

\begin{compactenum}

\smallskip

\item[$(ii$-$a)$] If $s \! < \! \gamma (y)$ (resp.~if $s \! > \! \gamma (y)$), then $F_{\xi,c, \Lambda} (s,y)  - s(c+\phi(y)) \! >\! 0$ (resp.~$< \! 0$). Otherwise, if $s \! = \! \gamma (y)$, $F_{\xi,c, \Lambda} (\gamma(y),y) \! =\! \gamma (y)(c+\phi(y))$. 

\smallskip

\item[$(ii$-$b)$] If $s \! \neq \gamma (y)$, then 
\begin{equation}
\label{discrGrey}r= \! \int^{g_r(s,y)}_s \!\!\!\! \frac{du}{F_{\xi,c, \Lambda} (u,y)\! -\! u(c+ \phi(y))}  \; .
\end{equation}
Moreover, $P_{\xi, c, \phi, \Lambda, \varrho, \Lambda_\varnothing} \big[ s^{Z^+_r} e^{-y\flambda ([0, r])} \big] \eqo F^\varnothing_{\!\! \varrho, \Lambda_\varnothing} \big( g_r (s,y), y \big)$. 
\end{compactenum}
\end{compactenum}
\end{lemma}
\noi
\textbf{Proof.} To simplify we set $Q\eqo Q_{\xi, c , \phi, \Lambda}$, $F\eqo F_{\xi, c , \Lambda}$, 
$P\eqo P_{\xi, c ,\phi, \Lambda, \varrho, \Lambda_{\varnothing}}$ and 
$F_\varnothing\eqo F^\varnothing_{\! \varrho, \Lambda_{\varnothing}}$. 
We denote by $\tau$ a GW($\xi, \varrho$)-forest and we assume that 
$(\ell (u), w_u, \mu_u)_{u\in \tau}$ is such that the isometry class of 
$(\mathbf T, d, \mathbf 0, \mathbf m)\! :=\!  
\mathtt{T}\mathtt{ree} (\tau; (\ell (\cdot), w_\cdot, \mu_\cdot ))$ has law $P$. We recall that 
$\mathbf T$ is a random subset of 
$\boldsymbol{\ell}_1 (\bbU)$ as introduced in Definition \ref{ellonemodel} $(c)$ and we observe that $\sigma_\varnothing \eqo \mathbf 0$ and that $\{ \sigma_u \, ; \, u\ino \tau \backslash \{ \varnothing \}\}\eqo (\mathtt{Br} (\mathbf T) 
\cup \mathtt{Lf} (\mathbf T)) \backslash \{ \mathbf 0 \}$. 

We also recall that 
$r\ino \bbR_+ \longmapsto Z^+_r(\widetilde{\mathbf T})\eqo \lim_{\epp \downarrow 0} \# \{ \sigma \ino \mathbf T \! : \! d(\mathbf 0 , \sigma) \eqo r\! +\!  \epp \}$ is a 
GW-branching process whose infinitesimal generator $G\eqo (g_{i,j})_{i,j\in \bN}$ is given by 
$g_{i,j}\eqo c i \xi(j\! -\! i+1)$ if 
$j\geq i\! -\! 1$ and $j\! \neq \! i$, $g_{i,i} \eqo -ci$ and $g_{i,j}\eqo 0$ otherwise. As read off from the generator, 
GW-branching processes are time-changed stopped random walks. More precisely, we set 
$$C_r\eqo \int_0^r \! Z^+_u(\widetilde{\mathbf T}) \, du, \quad C^{-1}_r \eqo \inf \big\{ u\ino \bbR_+ : C_u \! > \! r \big\} \quad \textrm{ and} \quad  X'_r\! :=\! Z^{+}_{{C^{-1}_r}}( \widetilde{\mathbf T}), $$ 
with the convention that $\inf \emptyset \eqo \infty$ and that $Z^+_\infty (\widetilde{\mathbf T}) \! :=\!  \lim_{u\to \infty} Z^+_u(\widetilde{\mathbf T}) \eqo 0$. Then 
$(X'_r)_{r\in \bbR_+}$ has the same law as $(X^+_{{r\wedge R_0}})_{r\in \bbR_+}$. Actually, without loss of generality, $\widetilde{\mathbf T}$ and $X^+$ can be coupled to get $X'\eqo X^+_{{\cdot \wedge R_0}}$, which therefore implies $C^{-1}_r\eqo A_r$, $r\ino \bbR_+$.  

It is furthermore possible to couple $\widetilde{\mathbf T}$ and $(Y_k)_{k\in \bN}$ as follows. We denote by $(u_k)_{0\leq k< \# \tau}$ the vertices of $\tau$ listed in ($\mathbf{P}$-a.s.~strictly) increasing height-order: namely, $d(\mathbf 0, \sigma_{u_{k-1}})  \! < \! d(\mathbf 0, \sigma_{u_{k}})$, $1\leqo k \! < \! \# \tau$. Note that $\{ d(\mathbf 0, \sigma_{u_{k}}); 1\leqo k \! < \! \# \tau\}\eqo 
\{ r\ino \bbR^*_+ \! : \! Z^+_{r-} (\widetilde{\mathbf T})\! \neq \!  Z^+_{r} (\widetilde{\mathbf T})\} $ and by the previous time change, we get $\{ C_{d(\mathbf 0, \sigma_{u_k}) } \, ; 1\leqo k \! < \! \# \tau\}\eqo \{ r\ino \bbR^*_+ \! : \! X^+_{r-} \! \neq \!  
X^+_{r} \} \eqo \{ r \ino \bbR^*_+\! : N_{(cr)-}\eqo  1 +N_{cr} \} $. Here we use the fact that $\xi$ is proper. 
We then set $U\eqo \mathtt{Br} (\mathbf T) \cup \mathtt{Lf} (\mathbf T) \cup \{ \mathbf 0 \} $ and we 
check that it is possible to extend the previous coupling in order to get $\mathbf m (B_{\mathbf{T} } (\mathbf 0, r) \cap  (\mathbf T\backslash U))\eqo S_{C_r}$, for all $r\ino \bbR_+$, and $w_{u_k}\eqo \mathbf m (\{ u_k\})\eqo W_k$, $0\leq k \! < \! \# \tau$. Namely, $ \mathbf m (B_{\mathbf{T} } (\mathbf 0, r) \cap   U)\eqo\sum_{0\leq k\leq N_{cr}}\!  W_k$, for all $r\ino \bbR_+$, which  entails $(i)$. 

To prove $(ii)$, we deduce from (\ref{df:GWweight}) that $F(\gamma (y), y)\eqo (c+ \phi (y))\gamma (y)$. 
We next fix $y$ and set $f(s)\eqo F(s,y)\! -\! s (c+ \phi(y))$. 
Then $f''(s)\eqo c \sum_{n\geq 2} n(n-1) s^{n-2} \xi (n) \widehat{\Lambda}_n (y) \geko 0$ and hence $f$ strictly convex on $[0,1]$ 
if $\xi$ is non-trivial (and $f$ is linear if $\xi(0)=1$). This implies $(ii$-$a)$ because 
$f(0) \eqo \xi (0)  \widehat{\Lambda}_0 (y) \geko 0$, $f(\gamma (y))\eqo 0$ and $f(1)\leqo   c \sum_{n\geq 0} \xi (n) (\widehat{\Lambda}_n (y) \! -\! 1)  \leqo  0$. To prove $(ii$-$b)$ we use (\ref{df:GWweight}) again to get 
$$ g_r (s,y) = se^{-cr-r\phi(y)} + \int_0^r \!\! du \, ce^{-cu -u\phi(y)}\sum_{n\in \bN} \xi (n) \widehat{\Lambda}_n (y) g_{r-u} (s,y)^n, $$
which easily implies (\ref{discrGrey}) by standard arguments. The last point of $(ii$-$b)$ is a simple consequence of (\ref{forestintrin}). \cqfd 
\begin{remark}
\label{heightGW} For all $(T,d,\rho) \! \equiv \! \widetilde{T}\ino \bbT_{\! c}^0$, we recall from Definition \ref{spandex} that 
$\mathtt{Ht} (\widetilde{T})\eqo \max_{\sigma \in T} d(\rho, \sigma)$ stands for the total height of $\widetilde{T}$. 
Clearly $\mathtt{Ht}$ is a $\dGH$-continuous function and we observe that (\ref{discrGrey}) characterizes the law of 
$\mathtt{Ht}$ under $\texttt Q_{\xi, c}$, which is the law of a purely metric GW($\xi, c$)-tree. Namely, observe that 
$F_{\xi, c, \Lambda} (s,0)\eqo c\varphi_\xi (s)$, where $\varphi_{\xi} (s)\eqo \sum_{k\in \bbN} s^k\xi (k)$, $s\ino [0,1]$, 
stands for the generating function of $\xi$. Then we derive from  (\ref{discrGrey}) that 
\begin{equation}
\label{heightGWexpli}
\forall r\ino \bbR_+^*, \quad cr= \int_0^{\texttt Q_{\xi, c} (\mathtt{Ht}  \leq r)} \!\!\!\!\! \frac{du}{\varphi_\xi (u)-u} \, , 
\end{equation}
since $g_r(s,0)\eqo Q_{\xi,c,w,\Lambda} [s^{Z^+_r}]\eqo \texttt Q_{\xi,c}[s^{Z^+_r}] \! \to \! \texttt Q_{\xi, c} (\mathtt{Ht}  \leq r)$ as $s\! \downarrow \! 0^+$. \cq
\end{remark}

  We now prove that GW-trees are characterised by the regenerative property. We rely on the unmeasured cases that have been proved by  
Duquesne and Winkel \cite[Theorem 3.21]{DuWi2}; we recall the statement below as Theorem \ref{recallTzero}. 
We recall that $\Phi^1_0$ stands for the 1-$(\dGHP, \dGH$)-Lipschitz function which consist in forgetting the measure. To simplify the notation, we denote in the same way 
$D$ and $D\! \circ\!  \Phi^1_0$, $Z^+_r$ and $Z^+_r \! \circ\!  \Phi^1_0$, $\mathtt{Abv}_r $ and 
$\mathtt{Abv}_r \! \circ\!  \Phi^1_0$ and we recall from Definition \ref{meaGWregedef} $(d)$ the notation 
$\texttt Q_{\xi, c}\eqo Q_{\xi,c,\phi, \Lambda} \! \circ\!  (\Phi^1_0)^{-1}$.   
\begin{theorem}
\label{recallTzero} Let $\texttt Q_0 \ino \mathcal M_1( \mathbb{T}^0_{\! c})$ be such that 
 $\texttt Q_0\big( D \ino \bbR^*_+ \big)\!>\!0$ and $\texttt Q_0(Z_r^+\! <\! \infty)\! =\! 1$ for all $r \! \in\! \bbR_+$. Then the following assertions are equivalent.
\begin{compactenum}


\item[$(a)$] For every $r\ino \bbR_+$, the conditional distribution of $\mathtt{Abv}_r$ given $Z^+_r$ under $\texttt Q_0$ is $\texttt Q_0^{\oast Z^+_r} $. 

\smallskip

\item[$(b)$] $\texttt Q_0\eqo \texttt Q_{\xi,c}$ for some $c\ino\bbR^*_+$ and some (sub)critical offspring distribution $\xi$ with $\xi (1)\eqo 0$. 

\smallskip
\end{compactenum}  
In particular, if $(a)$ or $(b)$ hold true, then $\texttt Q_0$-a.s.~$r\ino \bbR_+ \! \mapsto \! Z^+_r$ is an $\bbN$-valued c{\`a}dl{\`a}g process.   
\end{theorem}
\noi
\textbf{Proof.} See \cite[Lemma 2.16]{DuWi2} under the assumption that $(Z^+_r)_{r\in \bbR_+}$ is $\texttt Q_0$-a.s.~c{\`a}dl{\`a}g and \cite[Theorem 3.21]{DuWi2} for the remaining part of the proof.  \cqfd 

\smallskip

\begin{theorem}\label{branchingprop} Let $Q \ino \mathcal M_1( \mathbb{T}^1_c)$ be such that 
 $Q\big( D \ino \bbR^*_+ \big)\!>\!0$ and $Q(Z_r^+\! <\! \infty)\! =\! 1$ for all $r \! \in\! \bbR_+$. Then the following assertions are equivalent.
\begin{compactenum}

\smallskip

\item[$(a)$] $\! Q \big[g \big( \flambda ([0,r]) \big)\un_{\{Z_r^+=n\}}H(\mathtt{Abv}_r)\big]\eqo Q \big[ g \big( \flambda ([0,r])\big)\un_{\{Z_r^+=n\}} \big] Q^{\circledast n}[H]$, for all $r\ino\bbR_+$, $n\ino\bN$, all measurable $g\colon \bbR_+ \! \rightarrow \! \bbR_+$ and all $H \colon \bT_c^1 \! \rightarrow \! \bbR_+$.

\smallskip

\item[$(b)$] $Q\big[ G(\mathtt{Blw}_r )\un_{\{Z_r^+=n\}}H(\mathtt{Abv}_r)\big] \eqo Q\big[ G(\mathtt{Blw}_r)\un_{\{Z_r^+=n\}} \big] Q^{\circledast n}[H]$ for all $r\ino\bbR_+$, $n\ino\bN$, all measurable $G  \colon  \bT_c^1 \! \rightarrow \! \bbR_+$ and $H \colon  \bT^1_c \! \rightarrow\! \bbR_+$.
\smallskip

\item[$(c)$] $Q\eqo Q_{\xi,c,\phi,\Lambda}$ for a (sub)critical offspring distribution $\xi$ with $\xi(1)\eqo 0$, $c\ino\bbR^*_+$, the Laplace exponent $\phi$ of a subordinator and an $\cM_{1} (\bbR_+)$-valued sequence $\Lambda\eqo(\Lambda_n)_{n\in\bN}$. 
\end{compactenum}  
\end{theorem} 
\noi
\textbf{Proof.} $(b) \! \!\Rightarrow\!\! (a)$ holds trivially. Let us prove $(c)\!\!  \Rightarrow \! \! (b)$ by use of our representation of measured trees in $\boldsymbol{\ell}_1(\bbU)$. For all $u\ino \bbU$ let 
$\ell_u$, $(S_u (s))_{s\in\bbR_+}$, $(\Sigma_u (s))_{s\in\bbR_+}$, $\mathtt w_u(n)$, $n\ino \bbN$, 
be independent r.v.s such that $\ell_u$ is exponentially distributed with parameter $c$, 
$S_u(\cdot)$ and $\Sigma_u(\cdot)$ are subordinators with exponent $\phi$, and $\mathtt w_u(n)$ has 
law $\Lambda_n$. Let $\tau$ be a GW($\xi$)-tree. We recall the notation $(\sigma_u)_{u\in \tau}$ from Definition \ref{meaGWregedef} and we set $\zeta_u\eqo d_{\ell_1}(\mathbf 0, \sigma_u)$, which is 
the death time of $u$. Then $\zeta_{\overleftarrow{u}}$ is its birth time (with the convention that 
$\zeta_{\overleftarrow{\varnothing}}\eqo 0$) and $\mathcal Z_r\! :=\!  \{ u\ino \tau: \zeta_{\overleftarrow{u}} \leqo r \leko \zeta_u\}$ is the set of individuals alive at time 
$r\ino \bbR_+$. We set $\tau_{\leq r}\eqo \{ u\ino \tau :  \zeta_{\overleftarrow{u}} \leqo r \}$ and for all $u\ino \mathcal Z_r$, we denote by $\theta_u\tau\eqo \{ v\ino \bbU: u\ast v \ino \tau \}$ the subtree of $\tau$ stemming from $u$.  
We next introduce the following r.v.s. 

\smallskip

\noi
$-$ For all $u\ino \mathcal Z_r$, we set $\ell^{\leq r}_u \eqo r\! -\! \zeta_{\overleftarrow{u}}$ and $\ell^{ >r}_u 
\eqo \zeta_u \! -\! r$, $w^{\leq r}_u\eqo 0$ and $w^{>r}_u\eqo w_u \eqo \mathtt w_u(\mathbf k_u (\tau))$, 
$S^{\leq r}_u (s)\eqo S_u ( s\wedge \ell^{\leq r}_u)+ \Sigma_u ( (s\! -\!\ell^{\leq r}_u)_+ )$ and $S^{>r}_u(s) \eqo S_u (r+s)\! -\! S_u(r)$, $s\ino \bbR_+$, and $\mu^{\leq r}_u (ds)\eqo \un_{ [0, \ell^{{\leq r}}_u]} (s) \, dS^{\leq r}_u(s)$, 
$\mu_u(ds)\eqo \un_{[0, \ell_u]} (s) \, dS_u(s)$ and $\mu^{> r}_u (ds)\eqo \un_{[0, \ell^{> r}_u]} (s) \, dS^{> r}_u(s)$. 

\noi
$-$ For all $u\ino \tau $ such that $\zeta_u \leko r $, we set $\ell^{\leq r}_u\eqo \ell_u$, $w^{\leq r}_u\eqo w_u \eqo \mathtt w_u(\mathbf k_u (\tau))$,  $S^{\leq r}_u \eqo S_u$ and $\mu^{\leq r}_u (ds)$ $\eqo$  $\mu_u(ds)$ $\eqo$  $\un_{[0, \ell_u]} (s) \, dS_u(s)$. 

\noi
$-$ For all $u\ino \tau $ such that $r \leko \zeta_{\overleftarrow{u}}$, we set $\ell^{>r}_u\eqo \ell_u$, $w^{>r}_u\eqo w_u \eqo \mathtt w_u(\mathbf k_u(\tau))$, $S^{> r}_u \eqo S_u$ and $\mu^{> r}_u (ds)$ $\eqo$ $\mu_u(ds)$ $\eqo$ $\un_{[0, \ell_u]} (s) \, dS_u(s)$. 

\smallskip

\noi
We next define $\mathbf T\! :=\!  \mathtt{Tree} \big( \tau\, ; \, (\ell_\cdot, w_\cdot, \mu_\cdot)\big)$, 
$\mathbf T_{\! \leq r} \! :=\!  \mathtt{Tree} \big( \tau_{\leq r}\, ; \, (\ell^{\leq r}_v, w^{\leq r}_v, \mu^{\leq r}_v; v\ino \tau_{\leq r} ) \big)$ and for all $u\ino \mathcal Z_r$, $\theta_u\mathbf T\! :=\! \mathtt{Tree} \big( \theta_u \tau \, ; \, (\ell^{> r}_{u\ast v}, w^{> r}_{u\ast v,} \mu^{> r}_{u\ast v}\, ; \,  v\ino \theta_u \tau )\big)$ and we observe that 
\begin{equation}
\label{represellone}\widetilde{\mathbf T}\overset{\textrm{law}}{=} Q_{\xi,c,\phi,\Lambda}, \quad \widetilde{\mathbf T}_{\! \leq r}\eqo 
\mathtt{Blw}_r ( \widetilde{\mathbf T} ) , \quad  \underset{u\in \mathcal Z_r}{\circledast} \! \theta_u \widetilde{\mathbf T}\eqo \mathtt{Abv}_r ( \widetilde{\mathbf T} ) \quad \textrm{and} \quad Z^+_r (\widetilde{\mathbf T})\eqo \# \mathcal Z_r 
\end{equation}

\vspace{-2mm}

\noi
(here, $\widetilde{\mathbf T}$, $ \widetilde{\mathbf T}_{\! \leq r}$ and $\theta_u \widetilde{\mathbf T}$ are the isometry classes in $\bbT^1_{\! c}$ of $\mathbf T$, $\mathbf T_{\!\leq r}$ and 
$\theta_u\mathbf T$). Basic properties of subordinators imply that conditionally given 
$\big( \tau; (\ell_v)_{v\in \tau}\big)$, the subordinators 
$(S^{\leq r}_v(\cdot))_{v\in \tau_{\leq r}}$ and $(S^{>r}_{u*v} (\cdot))_{v\in \theta_u \tau}$, $u\ino \mathcal Z_r$, are independent with the same Laplace exponent $\phi$. By basic properties of GW-trees, and 
since exponential r.v.s are memoryless, conditionally given  
$\big( \tau_{\leq r} ;  (\ell^{\leq r}_v)_{v\in \tau_{\leq r}}\big)$, the subtrees 
$\big(\theta_u \tau ;  (\ell^{> r}_{u\ast v} )_{v\in\theta_u \tau}\big)$, $u\ino \mathcal Z_r$, are independent copies of 
$\big( \tau; (\ell_v)_{v\in \tau}\big)$. We therefore see that conditionally given $ \big( \tau_{\leq r} ; (\ell^{\leq r}_v\! , w^{\leq r}_v \! , \mu^{\leq r}_v)_{v\in \tau_{\leq r} }\big)$, the subtrees   
$\big( \theta_u \tau ; (\ell^{> r}_{u\ast v}, w^{> r}_{u\ast v,} \mu^{> r}_{u\ast v})_{v\in \theta_u \tau }\big)$, $u\ino \mathcal Z_r$, are independent copies of $ \big( \tau\, ; \, (\ell_v, w_v, \mu_v)_{v\in \tau}\big)$, which implies $(b)$ by (\ref{represellone}).  

\smallskip

It remains to prove $(a)\! \Rightarrow \! (c)$. We first note that Theorem \ref{branchingprop} $(a)$ implies Theorem \ref{recallTzero} $(a)$. Therefore 
$r\ino \bbR_+ \! \mapsto \! Z_r^+$ is $Q$-almost surely right-continuous
with $Z_0^+\eqo 1$, and the unweighted tree under $Q$ is a GW-tree for some $c\ino \bbR_+^*$ and some 
offspring distribution $\xi$, which is (sub)critical and proper, i.e.~$\sum_{k\in \bbN} k\xi (k)\leqo 1$ and $\xi(1)\eqo 0$. Thus $D$ under $Q$ is an exponential r.v.~with parameter $c$. 
Also, $r\eqo 0$ in $(a)$ yields $Q(\flambda (\{0\})\eqo 0)\eqo 1$.

We next fix a bounded $\fdelta_{\mathtt{GHP}} $-continuous function $H\colon \bbT^1_{\! c}\! \to \! \bbR_+$, a bounded continuous function $g_1\colon \bbR_+ \! \to \! \bbR_+$ and a 
right-continuous function $h\colon \bbR_+ \! \to \! \bbR_+$ and for all $p\ino \bbN$ and all $s\ino \bbR_+$ 
we use the notation $[s]_p\eqo 2^{-p}\lceil 2^p s\rceil $ for the right dyadic approximation of $s$ of order $p$. For all $k\ino \bbN^*$, 
We observe that $\int_{[0, k2^{-p}]} h([s]_p) \flambda (ds) \eqo \sum_{1\leq j\leq k}  h(j2^{-p})\flambda (I_{p,j} )$, where 
$I_{p,j}\! :=\!  \big(  (j\! -\! 1)2^{-p}\! ,j2^{-p} \big]$ for all $j,p\ino\bN$. It is also convenient to introduce the following notation: 
$$\Pi_{p,h} (k):= \!\! \prod_{1\le j\le k} \!\! Q  \Big[\! \exp \! \big(\!-\! h(j2^{-p})  \flambda (I_{p,1}) \big)
\un_{\! \{Z_{2^{-p}}^+=1\}}\Big]  \quad \textrm{and} \quad  B_{p,k}\! :=\! \bigcap_{1\leq j \leq k} \big\{ Z_{{j2^{-p}}}^{+}\eqo 1 \big\}. $$
By repeatedly applying $(a)$ with $r\eqo 2^{-p}$, we see that the following formula $(*)$ holds for all $k\ino \bbN^*$:
\vspace{-0.1cm}
  \begin{align*}
  &Q\bigg[ \mathrm{e}^{-\int_{[0, k2^{-p}]} h([s]_p) \flambda (ds) } 
\, \un_{\! B_{p,k} \cap \{Z_{{\! (k+1)2^{-p}}}^{+}=n\}}\, 
g_1 \big( \flambda (I_{p, k+1}) \big) H \big( {\mathtt{Abv}}_{(k\! +\! 1)2^{-p}} \big)\bigg] \\
&\qquad \qquad\qquad \qquad = \Pi_{p,h} (k) \, Q
  \big[\,  g_1 \big( \flambda (I_{p,1}) \big) \un_{\! \{Z_{2^{-p}}^+=n\}}\big] \, Q^{\circledast n}[H].\qquad \qquad \qquad \qquad (*)
  \end{align*}
For a first application of $(*)$, let us take $g_1\! \equiv \! 1$ and $H\! \equiv \! 1$.  
We fix $r\ino \bbR^*_+$ and we sum over $n\ino\bN$ to find, for $k\eqo \lfloor 2^pr\rfloor$,
\begin{equation}
 \label{indeplambda}
  Q\Big[
\mathrm{e}^{-\int_{[0, [r]_p-2^{-p} ]} h( [s]_p) \flambda (ds) } \un_{B_{p ,  \lfloor 2^pr\rfloor}}\Big]=
\Pi_{p,h} (  \lfloor 2^pr\rfloor ) 
  \end{equation}
We observe that $B_{p+1,  \lfloor 2^{p+1}r\rfloor} \!\subseteq\! B_{p,  \lfloor 2^pr\rfloor}$.   
We denote by $A$ the event that $r\ino \bbR_+ \! \mapsto \! Z^+_r$ is right-continuous. As already mentionned, $Q(A)\eqo 1$. We then observe that $A \cap \{ D\geko r \} \eqo \bigcap_{p\in \bbN} A\cap B_{p, \lfloor 2^pr\rfloor}$, which then implies that $Q$-a.s.~$\lim_{p\to \infty} \un_{B_{p,  \lfloor 2^pr\rfloor}} \eqo \un_{\{ D > r\}}$. 
Let us fix $m\ino \bbN^*$ and let us take 
$h\eqo \sum_{1\leq i\leq m } y_i\un_{ [(i-1)r/m \, , \, ir/m ) }$ for some $y_i\geqo 0$: we deduce from (\ref{indeplambda}), the previous limit and basic results on Laplace transforms that given $\{D\geko r\}$, the r.v.s $\big(\flambda \big( [\frac{{(i- 1)r}}{m},\frac{{ir}}{m})\big)\big)_{1\leq i\leq m}$ are independent and identically distributed under $Q$. 
Hence
$S_r\! :=\! \flambda ((0,r])$ is infinitely divisible under $Q(\, \cdot \, |\, D\geko r)$, and indeed $(S_r)_{0\leq r<D}$
is distributed like a subordinator up to an independent exponential time with parameter $c$. We denote its
drift coefficient by $\kappa$, its L\'evy measure by $\Gamma$ and its Laplace exponent by
$\phi(y)\eqo\kappa y \!+\!\int_{\bbR^*_+}(1\!-\!e^{-y z})\Gamma(dz)$. 
\begin{equation}
\label{plimrfixed}
 \lim_{p\to \infty} \Pi_{p,h} (  \lfloor 2^pr\rfloor ) = Q\Big[
\mathrm{e}^{-\int_{[0, r)} h(s) \, dS_s } \un_{\{D>r\}}\Big]= \mathrm{e}^{-cr-\int_0^r \! \phi (h(s)) \, ds} .
\end{equation}
Taking $h\! \equiv \! 0$ in (\ref{plimrfixed}), we then get $ \lim_{p\to \infty} Q(Z^{+}_{{2^{-p}}}\eqo 1)^{\lfloor 2^p r\rfloor } \eqo  
\lim_{p\to \infty} \Pi_{p,0} ( \lfloor 2^pr\rfloor ) \eqo \mathrm{e}^{-cr}$, which implies 
\begin{equation}
\label{plimZdiff1}
\lim_{p\to \infty} 2^p \big( 1\! -\! Q( Z^{+}_{{2^{-p}}} \eqo 1)\big)=\lim_{p\to \infty} 2^p Q (Z^{+}_{{2^{-p}}} \! \neq \!  1)= c\; .
\end{equation}
We next introduce the following notation. 
$$ D_p:=2^{-p} \inf \big\{ j\ino \bbN: Z_{j2^{-p}}^+\! \neq \! 1\big\} \quad \textrm{and} \quad  
M_{D_p}^{(p)}:= \flambda\big( (D_p\! -\! 2^{-p},D_p] \big)\; .$$
Let $f\colon \bbR_+ \! \to \! \bbR_+$ be bounded and continuous. We now fix $n\ino \bbN  \backslash \{ 1\}$: 
by multiplying $(*)$ by $f((k+1)2^{-p})$ and summing over $k$, we get the following. 
\begin{align}
 &Q\Big[ \mathrm{e}^{-\int_{[0, D_p-2^{-p})} h([s]_p) \flambda (ds) } 
f(D_p) \, \un_{\! \{Z_{{\! D_p}}^{+}=n\}}\, 
g_1 \big(M_{D_p}^{(p)} \big) H \big( {\mathtt{Abv}}_{D_p} \big)\Big] \nonumber \\
&\qquad = \Big( \!\int_0^\infty \!\!\! f \big( [r]_p \big) \Pi_{p,h} (  \lfloor 2^pr\rfloor ) dr  \Big)  2^{p} Q
\big[\,  g_1 \big( \flambda (I_{p,1}) \big) \un_{\! \{Z_{2^{-p}}^+=n\}}\big] \, Q^{\circledast n}[H]\nonumber \\
& \qquad = \Big( \!\int_0^\infty \!\!\! f \big( [r]_p \big) \Pi_{p,h} (  \lfloor 2^pr\rfloor ) dr  \Big)  
2^{p} Q\big( Z^{+}_{{2^{-}}}\! \neq \! 1 \big)  \,  Q\big[ g_1 \big(M_{D_p}^{(p)} \big)  \un_{\! \{Z_{{\! D_p}}^{+}=n\}} \big]   
 \, Q^{\circledast n}[H]\; . \label{discrdecomp}
\end{align}
\emph{Indeed,} the last equality is derived from the previous ones with $f\! \equiv \! H\! \equiv \! 1$ and $h\! \equiv \! 0$ 
and the fact that $2^pD_p$ is a geometric r.v.~such that $Q(2^pD_n \geqo k)\eqo Q\big( Z^{+}_{{2^{-}}}\! = \! 1 \big)^k$.

We next recall that $Q$-a.s.~$D_p\! \downarrow \! D$ and $D_p\! -\! 2^{-p}$ 
eventually increases to $D$, as $p\! \rightarrow\! \infty$ (since $r\! \mapsto\!  Z_{r}^{+}$ is integer-valued and 
$Q$-a.s.~right-continuous). We next note that $M_{{D_p}}^{{(p)}}$ converges $Q$-a.s.~as 
$p\! \rightarrow\! \infty$, to the mass $M_D\! =\! \flambda(\{D\})$ in the first branch point, while 
$Z_{{D_p}}^{+}\eqo \mathbf{k}$ for all sufficiently large $p$. We recall that $r\ino \bbR_+ \! \mapsto \! \mathtt{Abv}_r $ is 
continuous. Therefore,  
$Q$-a.s.~$\lim_{p\to \infty} f(D_p)\eqo  f(D)$, $\lim_{p\to \infty}g_1(M_{{D_p}}^{{(p)}})\eqo g_1(M_D)$ and 
$\lim_{p\to \infty} H({\mathtt{Abv}}_{D_p})\eqo H({\mathtt{Abv}}_D)$. Let us furthermore assume that $f$ has 
a compact support: by dominated convergence, (\ref{plimrfixed}) and (\ref{plimZdiff1}), we can take the limit in (\ref{discrdecomp}) and get 
\begin{align}
 &Q\Big[ \mathrm{e}^{-\int_{[0, D)} h(s) \flambda (ds) } 
f(D) \, \un_{\! \{\mathbf k =n\}}\, 
g_1 (M_{D} ) H \big( {\mathtt{Abv}}_{D} \big)\Big] \nonumber \\
&\qquad = \Big( \!\int_0^\infty \!\! c\, \mathrm{e}^{-cr-\int_0^r \phi (h(s)) ds } f (r) dr \Big) 
  \,  Q\big[ g_1 (M_{D} )  \un_{\! \{\mathbf k =n\}} \big]   
 \, Q^{\circledast n}[H]\; , \label{limidecomp}
\end{align}
which easily entails  $(c)$ by Lemma \ref{measwithmeas} $(ii)$. This completes the proof of the theorem. \cqfd

\subsection{Mass erasure of Galton--Watson $\bbR$-forests}
In this section we prove that mass erasure preserves the class of measured GW $\bbR$-trees as defined in Definition \ref{meaGWregedef}. 
Before stating this result, we recall a similar, simpler statement from \cite[Thms 2.18 and 2.19]{DuWi2} which asserts that 
GW-trees remain stable under length erasure (and even under erasure by any hereditary property). This result is used in the proof of the invariance principle of Section \ref{sec:IPnew} in Section \ref{pfsecthmIP}. 
In what follows, the generating function of any probability distribution $\nu$ on $\bbN$ is denoted by $\varphi_\nu (s)\eqo \sum_{k\in \bbN} s^k \nu(k)$, $s\ino [0, 1]$.
\begin{proposition}
\label{GWlengthera} Let $\xi, c,\varrho, \texttt Q_{\xi,c}$ and $\texttt P_{\! \xi, c, \varrho}$ be as in Definition \ref{meaGWregedef}. We assume that $\xi(1)\eqo 0$. Let $b\ino \bbR_+^*$. We recall from (\ref{deflengthera}) (and Lemma \ref{lengtheraprop}) the definition of $b$-length erasure $R_b \colon \bbT^{0}_{\! c} \! \to \!   \bbT^{0}_{\! c}$. We also recall that $\mathtt{Ht} \colon\bbT^{0}_{\! c} \! \to \!  \bbR_+$  stands for the height of compact $\bbR$-trees. To simplify, we set $\alpha \! :=\! \texttt Q_{\xi, c} (\mathtt{Ht} \leqo b) $, which belongs to $(0, 1)$ by (\ref{heightGWexpli}). Then $R_b $ under $ \texttt Q_{\xi,c}(\, \cdot \, | \, \mathtt{Ht} \geko b )$ (resp.~under 
$\texttt P_{\! \xi, c, \varrho}$) is a GW($\xi^*\! , c^*$)-tree (resp.~a GW($\xi^*\! , c^*\! , \varrho^*$)-forest) where for all $s\ino [0, 1]$, 
\begin{equation}
\label{lengtheraparam}
c^*\!\!  :=\! c (1\! -\! \varphi'_{\xi} (\alpha)), \; \varphi_{\xi^*} \! (s) \eqo s+ \frac{\varphi_\xi (\alpha + (1\! -\! \alpha) s) \! -\! \alpha \! -\! (1\! -\! \alpha) s }{(1\! -\! \alpha)(1\! -\! \varphi'_{\xi} (\alpha))}, \; \varphi_{\varrho^*}\!  (s) \eqo \varphi_{\varrho} ( \alpha + (1\! -\! \alpha) s). 
\end{equation}
\end{proposition}
\noi
\textbf{Proof.} See \cite[Thm 2.18]{DuWi2} for the statement concerning $\texttt Q_{\xi, c}$ and \cite[Thm 2.19]{DuWi2} for the statement concerning $\texttt P_{\! \xi, c, \varrho}$. \cqfd 

\smallskip

We recall notation (\ref{notaLapla}) for Laplace transforms and we recall the definitions of 
$F_{\xi,c, \Lambda }$ and  $F^{\varnothing}_{\varrho, \Lambda_{\varnothing }}$. 
We use the notation $\partial_1^n F_{\xi,c, \Lambda} (s,y)\eqo 
\frac{d^n}{ds^n} F_{\xi,c, \Lambda}(s,y)$. We also recall that  
if $(T,d,\rho, \mu)\! \equiv \! \fmu \ino \mathbb{T}^1_c$, then $\mu(T) \eqo \langle \fmu \rangle\eqo 
 \langle \flambda_{\fmu} \rangle\eqo  \flambda_{\fmu} (\bbR_+)$ is the total mass of $\fmu$. 
\begin{proposition}
\label{prop:eraGW} Let $\xi, c, \phi, \Lambda, \varrho, \Lambda_{\varnothing}$ be as in 
Definition \ref{meaGWregedef}. Let $h \ino \bbR^*_+$. For all $y\ino \bbR_+$, we set 
$$ \falpha_h(y)\eqo Q_{\xi, c, \phi, \Lambda} \big[ e^{-y \langle \flambda \rangle }
\un_{\{ \langle \flambda \rangle \leq h \}}\big] \quad \textrm{and} \quad \alpha\eqo \falpha_h (0)\eqo Q_{\xi, c, \phi, \Lambda} \big( \langle \flambda \rangle \leqo h\big) \; .$$
Then the following holds true. 
\begin{compactenum}

\smallskip

\item[$(i)$] $\mathcal E_h$ under $Q_{\xi, c, \phi, \Lambda} \big(\, \cdot \, \big| \,  \langle \flambda \rangle \! >\!  h\big)$ has law $Q_{\xi^*\!, c^*\! , \phi^*\! , \Lambda^*}$ where for all $s\ino [0, 1]$ and all $y\in \bbR_+$
\begin{eqnarray}
\label{charaQera}
c^*+ \phi^*(y) \!\!\!\! &\eqo & \!\!\!\! c+ \phi(y) -\partial_1 F_{\xi,c, \Lambda} \big( \falpha_h(y),y\big)  \qquad  \textrm{and} \quad \\
(1\! -\! \alpha) e^{-hy}   F_{\xi^*, c^*, \Lambda^*} (s,y) \!\!\!\! & \eqo&  \!\!\!\!  F_{\xi,c, \Lambda}  \big( (1\! -\! \alpha) e^{-hy}s +  \falpha_h(y)\, , \, y \big)
 - \falpha_h(y)(c+ \phi(y)) \nonumber \\
 & & \quad  - (1\! -\! \alpha) e^{-hy}s \,  \partial_1 F_{\xi,c, \Lambda} \big( \falpha_h(y),y \big) , \nonumber 
\end{eqnarray}  
which characterise $\xi^*\!, c^*\! , \phi^*\!$ and $\Lambda^*\! $. Furthermore, $\xi^*$ is proper whether or not $\xi$ is proper.
  \item[$(ii)$] $\mathcal E_h$ under $P_{\xi, c, \phi, \Lambda , \varrho, \Lambda_{\varnothing}}$ has law 
  $P_{\xi^*\! , c^*\! , \phi^*\! , \Lambda^*\! , \varrho^*\! , \Lambda^*_{\varnothing}}$ where $\xi^*\!, c^*\! , \phi^*\! , \Lambda^*\! $ are as in $(i)$ and   where for all $s\ino [0, 1]$ and all $y\in \bbR_+$
\begin{eqnarray}
 e^{-hy}   F^\varnothing_{\! \varrho^*\! , \Lambda^*_\varnothing} (s,y) \!\!\!\! & \eqo&  \!\!\!\!  F^\varnothing_{\! \varrho, \Lambda_\varnothing}  \big( (1\! -\! \alpha) e^{-hy}s +  \falpha_h(y)\, , \, y \big)-F^\varnothing_{\! \varrho, \Lambda_\varnothing}  \big(   \falpha_h(y) ,  y \big)
 \nonumber \\
 & & \quad  +e^{-hy} \varphi_{\varrho} (\alpha) \bE \big[ e^{-y(M-h)_+}\big]   . \label{charaPera}
\end{eqnarray}  
Here $M$ is a r.v.~such that $\varphi_{\! \varrho} (\alpha) \bE \big[ e^{-yM}\big]\eqo  
F^\varnothing_{\! \varrho, \Lambda_\varnothing}  \big(   \falpha_h(y) ,  y \big) $, $y\ino \bbR_+$. 
\end{compactenum}
\end{proposition}
\noi
\textbf{Proof.} To simplify, we set $F\eqo F_{\xi,c, \Lambda}$, $F_\varnothing \eqo   
F^\varnothing_{\! \varrho, \Lambda_\varnothing}  $, $Q\eqo  Q_{\! \xi, c, \phi, \Lambda} $, 
$P\eqo P_{\! \xi, c, \phi, \Lambda , \varrho, \Lambda_{\varnothing}}$ and $Q_h \eqo Q(\, \cdot \, | \, \langle \flambda \rangle \! > \! h)$. We first prove that $Q_h$ satisfies Proposition \ref{branchingprop} $(a)$. To this end, we fix $n\in \bN^*$ and $r\ino \bbR_+$, and we prove deterministic results on $\mathcal E_h \circ \mathtt{Abv}_r$ and  $\mathtt{Abv}_r\circ \mathcal E_h$. 

Let $(T,d,\rho, \mu)\! \equiv \! \widetilde{T}\ino \mathbb{T}^1_c$ and $r\ino \bbR_+$. We denote by $B_r$ the closed ball with center $\rho$ and radius $r$. Let $(T^o_i)_{i\in I}$ be the connected components of $T\backslash B_r \eqo \{ \sigma \ino T: d(\rho, \sigma) \! > \! r\}$. Then 
$\bigcup_{i\in I} T^o_i\eqo T\backslash B_r$. For all $i\ino I$, let $\sigma_i \ino T$ be such that $T_i\! :=\! \{ \sigma_i \} \cup T^o_i$ is the closure of $T^{o}_i$. We denote by $\widetilde{T}_i$ the isometry class of $(T_i, d, \sigma_i, \mu_i \! :=\! \mu( \, \cdot \, \cap T_i^o))$. Here we suppose that $I$ is finite and that 
$\mu(\{ \sigma_i\})\eqo 0$, for all $i\ino I$ and we note that $\mathtt{Abv}_r(\widetilde{T})\eqo \circledast_{i\in I} \widetilde{T}_i$.  
Then, we set $J\! :=\!  \{ i\ino I: \langle \mu_i \rangle \! >\! h\} $ and $w\! :=\! \big( h(\#J \! -\! 1) + \sum_{i\in I\backslash J} \langle \mu_i \rangle \big)_+$ and we recall that $\Upsilon_{\! w}$ stands for the isometry class of the point-tree $(\{ \rho \}, d, \rho, w\delta_\rho)$. Then we easily check that 
$$ \mathtt{Abv}_r \big(\mathcal E_h (\widetilde{T} )\big)= \circledast_{j\in J} \mathcal E_h (\widetilde{T}_j)  \quad \textrm{and} 
\quad  \mathcal E_h \big(\mathtt{Abv}_r (\widetilde{T}) \big)= \Upsilon_{\! w}  \circledast  \mathtt{Abv}_r \big(\mathcal E_h (\widetilde{T} )\big). $$
Thus $\flambda_{ \mathcal E_h \big(\mathtt{Abv}_r (\widetilde{T}) \big)} (\{ 0\})\eqo w$. We also observe that $\mathtt{Abv}_r \big(\mathcal E_h (\widetilde{T} )\big)= \mathtt{Abv}_0 \big(  \mathcal E_h \big(\mathtt{Abv}_r (\widetilde{T}) \big)\big)$. 
Furthermore, we get 
$$ Z^+_r \big( \mathcal E_h (\widetilde{T} )\big) \eqo \# J \eqo Z^+_0 
 \big(\mathcal E_h ( \mathtt{Abv}_r(\widetilde{T}) )\big)  \quad \textrm{and} \quad \flambda_{\mathcal E_h (\widetilde{T} )} ([0, r])\eqo \flambda_{\widetilde{T} } ([0, r])+w. $$
By Proposition \ref{branchingprop} $(a)$ applied to $Q$, we first get for all integers $p\geqo n \geqo 1$ 
\begin{eqnarray*}
Q\Big[ \un_{\{ Z^+_r = p \, ;\,  Z^+_r \! \circ \, \mathcal E_h= n \}} \!\!\! \! \!\!  & & \!\!  \!\! \!  e^{-y (\flambda \circ \mathcal E_h )([0, r]) } 
G \big( \mathtt{Abv}_r \circ \mathcal E_h \big) \Big]   \\  
\!\!\!\!\!  \!\!\! \!\!  \!\!\! \!\!   \!\!\! \!\!   \!\!\! \!\! & & \!\!\! \!\!  \!\!\! \!\!  \!\!\! \!\!  \!\!\! \!\!  \!\!\!\!\!  =Q \Big[ \un_{\{ Z^+_r = p\}}  e^{-y \flambda ([0, r])} \un_{ \{ Z^+_0 (\mathcal E_h \circ \mathtt{Abv}_r) = n \}} 
e^{-y (\flambda \circ \mathcal E_h \circ \mathtt{Abv}_r ) (\{ 0 \}) } G \big( \mathtt{Abv}_0 \! \circ \!  \mathcal E_h \! \circ \!  \mathtt{Abv}_r  \big) \Big] \\
\!\!\!\!\!  \!\!\! \!\!  \!\!\! \!\!   \!\!\! \!\!   \!\!\! \!\! & & \!\!\!\!\!  \!\!\! \!\!  \!\!\! \!\!   \!\!\! \!\!   \!\!\! \!\! = Q \Big[ \un_{\{ Z^+_r = p\}}  e^{-y \flambda ([0, r])}\Big] Q^{\circledast p} \Big[\un_{ \{ Z^+_0 (\mathcal E_h ) = n \}} 
e^{-y (\flambda \circ \mathcal E_h ) (\{ 0 \}) } G \big( \mathtt{Abv}_0 \! \circ \!   \mathcal E_h \big) \Big]  .
\end{eqnarray*}
To compute the above expectation under $Q^{\circledast p}$, we proceed as follows. For all $i\ino \{1, \ldots,  p\}$, 
let $(T'_i, d_i', \rho_i', \mathrm m_i') \! \equiv\! \widetilde{T}'_i\ino \mathbb{T}^1_c $ such that 
$\mathrm m_i' (\{ \rho_i' \})\eqo 0$. Then set $\widetilde{T}'\eqo \circledast_{1\leq i \leq p} \widetilde{T}'_i$ and we observe that  
\begin{eqnarray*}
& & \!\!\! \!\!\!  \un_{\! \{ Z_{0}^{+} (\mathcal E_h (\widetilde{T}') ) = n\} }  e^{-y \flambda_{\mathcal E_h (\widetilde{T}')} (\{ 0\}) } G (\mathtt{Abv}_0 (\mathcal E_h (\widetilde{T}'))) \\ 
& & =\! \sum_{\substack{J \subset \lgeo 1, p\rgeo \\ \# J =n}} G \big( \! \circledast_{j\in J} \mathcal E_h (\widetilde{T}'_j) \big)e^{-yh(n-1)} \Big( \prod_{j\in J} \un_{\{  \langle \flambda_{\widetilde{T}'_j  } \rangle >h \}}\Big)  \Big( \prod_{i \in \lgeo 1, p\rgeo  \backslash J} \!\!\! e^{-y \langle\flambda_{\widetilde{T}'_i}\rangle }\un_{\{  \langle \flambda_{\widetilde{T}'_i  } \rangle \leq h \}}\Big) \, .
\end{eqnarray*}   
Therefore, 
\begin{equation}
\label{Qpstep}
Q^{\circledast p}\!  \Big[\un_{ \! \{ Z^+_0 \circ \mathcal E_h  = n \}} 
e^{-y \flambda \circ \mathcal E_h  (\{ 0 \}) } G \big( \mathtt{Abv}_0 \circ   \mathcal E_h \big) \Big]= e^{yh} \binom{p}{n} \big( (1 - \alpha)e^{-yh} \big)^{\! n} \!  \falpha_{h} (y)^{p-n}
Q_h^{\circledast n} [G] .
\end{equation}
We then set $g(s,y)\eqo Q[ s^{Z^+_r} e^{-y \flambda ([0, r])}]$ and the previous equalities entail
\begin{eqnarray*}
(1\! -\! \alpha) Q_h \Big[ \un_{\! \{ Z^+_r = n  \}} e^{-y \flambda ([0, r])} G\big( \mathtt{Abv}_r \big) \Big] \!\!\!\! & =&\!\!\!\!  \sum_{p\geq n} Q\Big[ \un_{\! \{ Z^+_r = p\,  ;\,  Z^+_r \! \circ \,  \mathcal E_h= n \}}  e^{-y (\flambda \circ \mathcal E_h) ([0, r]) } 
G \big( \mathtt{Abv}_r \! \circ\!  \mathcal E_h \big) \Big]  \\
\!\!\!\!  & =& \!\!\!\!  \frac{1}{n!} e^{yh}\big( (1\! -\! \alpha)e^{-yh} \big)^n \partial_1^n g \big( \falpha_{h} (y), y\big) Q_h^{\circledast n} [G]. 
\end{eqnarray*}
This shows in particular that $Q_h$ satisfies Proposition \ref{branchingprop} $(a)$: so there are $\xi^*\!$, a proper and (sub)critical offspring distribution, $c^*\ino \bbR_+^*$, 
$\phi^*\! $ and $(\Lambda^*_n)_{n\in \bN}$ such that $Q_h\eqo Q_*\! :=\!  Q_{\xi^*\! , c^*\! , \phi^*\! , \Lambda^*\! }$. 
To complete the proof of $(i)$ we need to characterise $\xi^*\! , c^*\! , \phi^*\!$ and $\Lambda^*\! $ in terms of $\xi, c, \phi$ and $\Lambda$. To this end, for all $n\ino \bN\backslash \{ 0, 1\}$ and all $x,y,z \ino \bbR_+$, we define $L_{n,x,y,z} \colon \mathbb{T}^1_c \! \to \! \bbR_+$ by  
$$ L_{n,x,y,z}\eqo \un_{\{  Z^+_0=1\, ;\, 0<  D <\infty \, ;\,  \mathbf k=n \}}e^{-xD -y\flambda([0, D)) -z\flambda(\{ D\}) }  \; .$$ 
Let $\xi^o, c^o, \phi^o, \Lambda^{\! o}$ be as in Definition \ref{meaGWregedef}. We easily see that 
\begin{equation}
\label{Lectgeneral}
Q_{\xi^o,c^o,\phi^o, \Lambda^o} \big[ L_{n,x,y,z}\big]= \frac{c^o}{c^o+ x + \phi^o(y)} \xi^o(n) \widehat{\Lambda}^o_n (z) \; .
\end{equation}   
We first compute for all integers $n\geqo 2$
$$ f_{n} (x,y,z)= Q_*  \big[ L_{n,x,y,z}\big]= \frac{1}{1\! -\! \alpha} Q \big[ \big(L_{n,x,y,z} \! \circ \!  \mathcal E_h \big)   \un_{\{\langle \flambda \rangle >h  \}}\big] .$$
(\emph{Indeed}, since $n\! \neq \! 0$, $Q$-a.s.~$ \big(L_{n,x,y,z} \! \circ \!  \mathcal E_h \big)   \un_{\{\langle \flambda \rangle >h  \}}\eqo L_{n,x,y,z} \! \circ \!  \mathcal E_h $.)

Let $\mathbf T$ be a $\mathbb{T}^1_c$-valued r.v.~with law $Q$. To simplify, we set $\vartheta (\mathbf T)\eqo \mathtt{Abv}_{{D(\mathbf T)}}(\mathbf T)$, $\mathbf T_{\! h}\eqo  \mathcal E_h (\mathbf T)$ 
and $\mathbf k_h  \eqo 
Z_{{D(\mathbf T )}}^{+} (\mathbf T_{\! h})$, so that the latter is the number of points of 
$\mathbf T_{\! h}$ just above the height of the first branch point of $\mathbf T$. 
Recall that $n\geqo 2$. On the event $ \{ \mathbf k (\mathbf T_{\! h}) \eqo n\}$, a.s.~only two cases occur:   
either $D(\mathbf T_{\! h}) \eqo D(\mathbf T)$ and necessarily $\mathbf k_h \eqo \mathbf k (\mathbf T_{\! h} )\eqo n \leqo \mathbf k (\mathbf T)$, or $D(\mathbf T) \! < \!  D(\mathbf T_h)$ and necessarily $\mathbf k_h \eqo 1$. Let $p\geqo n$.  
In the first case, we observe that 
$$ \un_{\{  \mathbf k (\mathbf T)= p\, ; \,  \mathbf k_h = n \} } L_{n,x,y,z} (\mathbf T_{\! h}) = 
L_{p,x,y,z} (\mathbf T) e^{-z \flambda_{\mathcal E_h (\vartheta (\mathbf T)) } (\{ 0\}) 
 }\un_{\{ Z^+_0  (\mathcal E_h (\vartheta (\mathbf T))) =n\}}.$$
By (\ref{df:GWweight}) in Lemma \ref{measwithmeas}, (\ref{Qpstep}) and (\ref{Lectgeneral}) we get 
\begin{eqnarray*}
\bE \big[  \un_{\{  \mathbf k (\mathbf T)= p\, ; \,  \mathbf k_h = n \} } L_{n,x,y,z} (\mathbf T_{\! h}) \big] & =& Q \big[ L_{p,x,y,z } \big] \!\!  \cdot \! Q^{\circledast p} \Big[ e^{-z (\flambda \circ \mathcal E_h) (\{ 0\} )} \un_{\{ Z^+_0 \!  \circ\,  \mathcal E_h =n  \} }\Big] 
 \\
& =& \tfrac{c}{c+ x + \phi(y)} \xi(p) \widehat{\Lambda}_p (z) \! \cdot \!  e^{zh} \binom{p}{n} ((1\! -\! \alpha)e^{-zh})^n \falpha_{h} (z)^{p-n} .\end{eqnarray*} 
By summing over $p\geqo n$, we get 
\begin{equation}
\label{compuLstep}
\bE \big[  \un_{\{ \mathbf k_h = n \} } L_{n,x,y,z} (\mathbf T_{\! h}) \big]= \frac{e^{zh}((1\! -\! \alpha)e^{-zh})^n }{(c+ x + \phi(y))n! } \partial^n_1 F\big(\falpha_h (z), z \big). 
\end{equation}
 Let $p\geqo 2$. In the second case we get 
\begin{eqnarray*} \un_{\! \{  \mathbf k (\mathbf T)= p\, ; \,  \mathbf k_h = 1 \} } \! \!\! \!\!  & & \!\! \!\! \!   L_{n,x,y,z} (\mathbf T_{\! h}) = \\
\!\! \!\! \!\! & &\!\! \!\! \!\!  L_{p,x,y,y} (\mathbf T) e^{-y \flambda_{\mathcal E_h (\vartheta (\mathbf T)) } (\{ 0\}) } 
L_{n,x,y,z} \big( \mathtt{Abv}_0 (\mathcal E_h (\vartheta (\mathbf T))) \big) \un_{\{ Z^+_0  (\mathcal E_h (\vartheta (\mathbf T))) =1\}}.
\end{eqnarray*} 
Again, by (\ref{df:GWweight}) in Lemma \ref{measwithmeas}, (\ref{Qpstep}) and (\ref{Lectgeneral}), we get 
\begin{eqnarray*}
\bE \big[ \un_{\{  \mathbf k (\mathbf T)= p\, ; \,  \mathbf k_h = 1 \} } L_{n,x,y,z} (\mathbf T_{\! h})\big] \!\!\!\!   & =& \!\!\!\! 
Q \big[ L_{p,x,y,y } \big] \!\!  \cdot \! Q^{\circledast p} \! \Big[ e^{-y (\flambda \circ \mathcal E_h) (\{ 0\} )} L_{n,x,y,z} \big( \mathtt{Abv}_0 \! \circ \! \mathcal E_h \big)\un_{\! \{ Z^+_0 \!  \circ \, \mathcal E_h =1 \}} \! \Big] 
 \\
\!\!\!\!  & =& \!\!\!\! \tfrac{c}{c+ x + \phi(y)} \xi(p) \widehat{\Lambda}_p (y) \! \cdot \!  p(1\! -\! \alpha) \falpha_{h} (y)^{p-1}Q_h \big[ L_{n,x,y,z}\big].
\end{eqnarray*} 
By summing over $p\geqo n$, we get 
\begin{equation}
\label{compuLstep2}
\bE \big[  \un_{\{ \mathbf k_h = 1 \} } L_{n,x,y,z} (\mathbf T_{\! h}) \big]=  \frac{(1\! -\! \alpha) Q_h \big[ L_{n,x,y,z}\big]}{c+ x + \phi(y)} \partial_1 F\big(\falpha_h (y), y \big). 
\end{equation}
By summing (\ref{compuLstep}) and (\ref{compuLstep2}) we get 
\begin{eqnarray*}
 (1\! -\! \alpha) f_n (x,y,z) &=& \bE \big[  \un_{\{ \mathbf k_h = n \} } L_{n,x,y,z} (\mathbf T_{\! h}) \big] + \bE \big[  \un_{\{ \mathbf k_h = 1 \} } L_{n,x,y,z} (\mathbf T_{\! h}) \big] \\
 & =&\frac{e^{zh}((1\! -\! \alpha)e^{-zh})^n }{(c+ x + \phi(y))n! } \partial^n_1 F\big(\falpha_h (z), z \big) + \frac{(1\! -\! \alpha) f_n (x,y,z)}{c+ x + \phi(y)} \partial_1 F\big(\falpha_h (y), y \big). 
\end{eqnarray*} 
To simplify we set $q(y)\eqo  \partial_1 F (\falpha_h (y), y )$ and we rewrite the previous equality as follows. 
$$ (c+x+ \phi (y) \! -\! q(y))  (1\! -\! \alpha) f_n (x,y,z) = \frac{1}{n!} e^{zh} ((1\! -\! \alpha)e^{-zh})^n  \partial^n_1 F\big(\falpha_h (z), z \big) .$$
Thus by (\ref{Lectgeneral}) we obtain for all integers $n\geqo 2$ and all $x,y,z \ino \bbR_+$ 
\begin{equation}
\label{compuparamstep}
\frac{c+x+ \phi (y) \! -\! q(y)}{c^*+x+ \phi^* (y)} c^* (1\! -\! \alpha)  \xi^* (n) \widehat{\Lambda}_n^* (z) = \frac{1}{n!} e^{zh} ((1\! -\! \alpha)e^{-zh})^n  \partial^n_1 F\big(\falpha_h (z), z \big)\; . 
\end{equation}
To simplify we denote by $\varphi_*$ (resp.~$\varphi$) the generating function of $\xi^*$ (resp.~of $\xi$). By definition 
$F (s,0)\eqo c\varphi (s)$, $s\ino [0, 1]$. Thus $q(0) \eqo c\varphi'(\alpha)$. By taking $x\eqo y\eqo z\eqo 0$ in (\ref{compuparamstep}) we get $(1\! -\! \alpha) (1 \! -\! \varphi'(\alpha)) \xi^*(n) \eqo (1\! -\! \alpha)^n \varphi^{{(n)}} (\alpha)/ n!$. Since $\xi^*(1)\eqo 0$, we get for all $s\ino [0, 1]$,  
$$ \varphi_* (s) = s + \frac{\varphi \big( (1\! -\! \alpha) s + \alpha \big) \! -\! \alpha -(1\! -\! \alpha)s }{(1\! -\! \alpha)(1\! -\! \varphi'(\alpha))} \quad \textrm{and thus } \quad \xi^* (0) \eqo \frac{\varphi (\alpha) \! -\! \alpha }{(1\! -\! \alpha)(1\! -\! \varphi'(\alpha)) } .$$  
By taking $z\eqo 0$ in (\ref{compuparamstep}) and by the previous arguments, we get $ \frac{c+x+ \phi (y)  - q(y)}{c^*+x+ \phi^* (y)} c^* \eqo c (1\! -\! \varphi'(\alpha))$ for all $x,y \ino \bbR_+ $. By taking $x\! \neq \! 0 $ and $y\eqo 0$ in this equality, we get $c^* \eqo c(1\! -\! \varphi'(\alpha))$. This entails that $c+ \phi (y) \! -\! q(y)\eqo c^*+ \phi^* (y)$. 
Combined with (\ref{compuparamstep}) again we get 
\begin{eqnarray}
(1\! -\! \alpha) e^{-zh} \big(F_* (s,z) \!\!\!\!\!  & -& \!\!\!\!\! c^* \xi^*(0) \widehat{\Lambda}^*_0 (z)  \big)  = \sum_{n\geq 2} c^* (1\! -\! \alpha) e^{-zh} \xi^* (n) s^n \widehat{\Lambda}^*_n (z) \label{Compustepstep} \\
& =& F \big( (1\! -\! \alpha) e^{-zh}s+ \falpha_h (z)  ,  z\big) \! -\! F(\falpha_h (z) ,z) \! -\!  (1\! -\! \alpha) e^{-zh}sq(z).\nonumber
\end{eqnarray} 

We next need to compute $\xi^*(0) \widehat{\Lambda}^*_0 $. To this end, 
we 
set $\gamma (z)\eqo Q [ e^{-z \langle \flambda \rangle}] $ and $\gamma^* (z) \eqo Q_* [ e^{-z \langle \flambda \rangle}] $. By definition of the function $\falpha_h$ we get $\gamma (z)\eqo  \falpha_h(z)+ (1\! -\! \alpha) e^{-zh} \gamma^* (z)$. We use Lemma \ref{ballGWmass} $(ii$-$a)$ which asserts that $F(\gamma (z),z)\eqo \gamma (z)(c+ \phi(z))$ and $F_*(\gamma^* (z),z)\eqo \gamma^* (z)(c^*+ \phi^*(z))$ for all $z\ino \bbR_+^*$. 
Namely, we take  
 $s\eqo \gamma^*(z)$ in (\ref{Compustepstep}). The member on the left-hand side of (\ref{Compustepstep}) becomes
\begin{eqnarray*}
  (1\! -\! \alpha) e^{-zh} \big( (c^* + \phi^* (z) ) \gamma^* (z)\!\! \!\! & -& \!\! \!\!  
 c^* \xi^*(0) \widehat{\Lambda}^*_0 (z)  \big) \\
\!\! \!\! &= & \!\! \!\! (c+ \phi (z) \! -\! q(z)) (\gamma (z)\! -\! \falpha_h(z)) -(1\! -\! \alpha) e^{-zh} c^* \xi^*(0) \widehat{\Lambda}^*_0 (z).
\end{eqnarray*} 
The member on the right-hand side of (\ref{Compustepstep}) becomes 
\begin{eqnarray*}
F (\gamma (z), z) -F(\falpha_h (z), z)\!\! \!\!  &-&\!\! \!\!  (\gamma (z) \! -\! \falpha_h (z)) q(z) \\
\!\! \!\! & =& \!\! \!\!  (c+ \phi (z) \! -\! q(z))  (\gamma (z)\! -\! \falpha_h(z)) 
\! -\! F(\falpha_h (z), z)\!  +\! \falpha_h (z) (c+ \phi(z)). 
\end{eqnarray*} 
Thus $(1\! -\! \alpha) e^{-zh} c^* \xi^*(0) \widehat{\Lambda}^*_0 (z)\eqo F(\falpha_h (z), z) - \falpha_h (z) (c+ \phi(z)) $, which completes the proof of (\ref{charaQera}) and thus of $(i)$. 

\smallskip

It remains to prove $(ii)$. We next compute the law $P_*$ of $\mathcal E_h $ under $P$ by computing the law of the atom at the root, the law of the degree of the root and the law of the subtrees stemming from the root. To this end, we fix $p \geqo n\geqo1$ and we first get the following (recall that $\mathtt{Abv}_0$ eliminates any mass at the root). 
\begin{eqnarray*}
P \Big[ \un_{\{ Z^+_0 =p \, ; \, Z^+_0 \! \circ \, \mathcal E_h =n \}} \!\!\!\!\!\! \!\!\!\!\!\! \! & &\!\!\!\!\!\! \!\!\!\!\!\! \!e^{-y (\flambda \circ \mathcal E_h) (\{ 0\})} 
G\big(\mathtt{Abv}_0\!  \circ \mathcal \! E_h \big)  \Big] \\
& \overset{\textrm{by (\ref{forestintrin})}}{=}& \varrho (p)
\widehat{\Lambda}_{\varnothing, p} (y) Q^{\circledast p} \Big[ \un_{\{  Z^+_0\!  \circ \, \mathcal E_h =n \}} e^{-y (\flambda \circ \mathcal E_h) (\{ 0\})} 
G\big(\mathtt{Abv}_0 \! \circ \mathcal \! E_h \big)  \Big] \\
& \overset{\textrm{by (\ref{Qpstep})}}{=}&  e^{hy} \varrho (p)
\widehat{\Lambda}_{\varnothing, p} (y) \binom{p}{n} \falpha_h(y)^{p-n} \big((1\! -\! \alpha) e^{-hy} \big)^n Q_*^{\circledast n} \big[ G \big].
\end{eqnarray*} 
By summing over $p\geqo n$, we get 
\begin{eqnarray}
\!\!\!\!\!P_* \Big[ \un_{\{ Z^+_0 =n \}}\!e^{-y \flambda  (\{ 0\})} 
G\big(\mathtt{Abv}_0 \big)  \Big] & \!\!\!=\!\!\!& P \Big[ \un_{\{  Z^+_0\!  \circ\,  \mathcal E_h =n \}} e^{-y (\flambda \circ \mathcal E_h) (\{ 0\})} 
G\big(\mathtt{Abv}_0 \! \circ \! \mathcal E_h \big)  \Big] \nonumber \\
& \!\!\!=\!\!\!& e^{hy} \frac{1}{n!} \big((1\! -\! \alpha) e^{-hy} \big)^n \partial_1^n F^\varnothing_{{\! \varrho, \Lambda^\varnothing}} (\falpha_h (y), y)  Q_*^{\circledast n} \big[ G \big]. \label{groundn}
 \end{eqnarray} 
We next observe that 
$$ P_* \Big[ \un_{\{ Z^+_0 =0 \}}\!e^{-y \flambda  (\{ 0\})} 
G\big(\mathtt{Abv}_0 \big)  \Big] =  P \Big[ \un_{\{  Z^+_0 \! \circ \, \mathcal E_h =0 \}} e^{-y (\langle \flambda \rangle -h)_+}
 \Big] G(\Upsilon_0). $$
To compute the right hand side of this equality, we observe for all $p\in \bN$ that 
\begin{eqnarray*}
P \Big[ \un_{\{ Z^+_0 =p \, ; \, Z^+_0\!  \circ\,  \mathcal E_h =0 \}}e^{-y \langle \flambda  \rangle } 
 \Big] \!\!\! & \overset{\textrm{by (\ref{forestintrin})}}{=}& \!\!\! \varrho (p)
\widehat{\Lambda}_{\varnothing, p} (y) Q^{\circledast p} \Big[ \un_{\{   Z^+_0 \! \circ \, \mathcal E_h =0 \}} 
e^{-y \langle \flambda  \rangle } \Big] \\
\!\!\! &=& \!\!\! \varrho (p)
\widehat{\Lambda}_{\varnothing, p} (y) \Big(Q\big[ e^{-y \langle \flambda  \rangle} \un_{\{ \langle \flambda  \rangle \leq h\}}\big] \Big)^p \! =\!  \varrho (p)
\widehat{\Lambda}_{\varnothing, p} (y) \falpha_h (y)^p. 
\end{eqnarray*} 
By summing over $p\ino \bN$, we get $ P \big[ \un_{\{  Z^+_0 \! \circ \, \mathcal E_h =0 \}} e^{-y \langle \flambda\rangle }
 \big]\eqo  F^\varnothing_{\! \varrho, \Lambda_\varnothing} (\falpha_h (y),y)$. If one takes $y\eqo 0$, then 
 $P \big(Z^+_0 \! \circ\!  \mathcal E_h \eqo 0 \big)\eqo  \varphi_\varrho (\alpha)$ and we obtain that 
 $$   P_* \Big[ \un_{\{ Z^+_0 =0 \}} e^{-y \flambda  (\{ 0\})} 
G\big(\mathtt{Abv}_0 \big)  \Big] =  \varphi_\varrho(\alpha) \bE \big[ e^{-y (M-h)_+} \big] G(\Upsilon_0), $$
where $M$ is distributed as in $(ii)$. Combined with (\ref{groundn}) this completes the proof of (\ref{charaPera}).  \cqfd 

\subsection{Constructions of L\'evy forests} \label{applLevy}
Lévy trees were introduced in Le Gall and Le Jan \cite{LGLJ1} as the continuum random trees corresponding to 
the genealogy of \emph{Continuous-State Branching Processes} (CSBP for short). They extend Aldous's 
Continuum Random Tree \cite{Al1}, which corresponds to the Brownian case (see Le Gall \cite{LG2} and 
Aldous \cite{Al2}). They turn out to be the limits of rescaled Galton--Watson trees (see Aldous \cite{Al2} in the 
Brownian case and see 
Le Gall and Le Jan \cite{LGLJ1} and Duquesne and Le Gall  \cite{DuLG} in more general cases). 
The construction of L{\'e}vy trees given in Le Gall and Le Jan \cite{LGLJ1} and in Duquesne and Le Gall  \cite{DuLG} uses coding by 
a height process. The goal of this section is to provide another construction of Lévy trees as the projective 
limit of a family of measured GW-trees that is consistent under mass erasure: this construction encompasses cases 
where the (sub)critical CSBP is never absorbed at $0$ (i.e.~it does not satisfy the Grey condition: see below). It even includes 
bounded variation cases, which is new.  

 \smallskip

 \noi
 \textbf{Continuous-state branching processes and spectrally positive Lévy processes.}  Before proving these results we recall the basic results on CSBPs that we need in the proofs. 
Let $x\ino \bbR^*_+$. As already mentioned, we want to define the genealogy of a CSBP $(Z^{x}_s)_{s\in \bbR_+}$ with initial population $Z^{x}_0\eqo x$ supposing the CSBP to be critical or subcritical: namely $\bP$-a.s.~$\lim_{s\to \infty} Z^{x}_s\eqo 0$.  (Sub)critical CSBPs are the continuous analogues of (sub)critical GW-Markov chains: namely they are $\bbR_+$-valued Feller--Dynkin Markov processes whose transition probabilities are characterised by their \emph{branching mechanism} $\psi\colon \bbR_+ \! \to \! \bbR_+$ that is the Laplace exponent of a L{\'e}vy process without negative jumps. In the (sub)critical case $\psi$ is necessarily of the following L{\'e}vy--Khintchine form: 
\begin{equation}\label{brmech}
\psi(y)=\mathbf{a}y+\tfrac{1}{2}\mathbf{b}y^2+\int_{\bbR^*_+}(e^{-y z}-1+yz)\pi(dz)\, , 
\end{equation}
where $\mathbf{a}, \mathbf{b}\ino \bbR_+$, and where the 
Borel measure $\pi$ on $\bbR^*_+$ satisfies $\int_{\bbR_{+}^*} \! (z \! \wedge \! z^2)\, \pi(dz)\! <\! \infty$.  

More precisely for all 
$s,r,y\ino \bbR_+$, $\bE [ \exp( -yZ^{x}_{r+s}) | Z^x_r]\eqo \exp (-Z^x_r u_s(y))$ $\bP$-almost surely, where $u_s(y)$ is given by $\int_{u_s(y)}^y dz/ \psi(z)\eqo s$. This entails for all $y\ino \bbR_+$ that $\lim_{s\to \infty} u_s(y)\eqo 0$ and that $\lim_{s\to \infty} Z^{x}_{s}\eqo 0$ in probability (martingale arguments shows that this limits holds a.s.). In many of the cases, 
extinction occurs in finite time. More precisely, we set $\mathtt{Ext} \! :=\! \inf \{ s\ino \bbR_+ \! : \! Z^{x}_{s} \eqo 0 \}$, with the convention $\inf \emptyset \eqo \infty$. Since $0$ is an absorbing state, if $\mathtt{Ext} \leko \infty$, then $Z^{x}_{s}\eqo 0 $ for all $s\geko \mathtt{Ext}$. Namely, $\mathtt{Ext}$ is the extinction time of $Z^{x}_{\cdot}$. We can check  that 
\begin{equation}
\label{Grey} 
\bP (\mathtt{Ext} \leko \infty) \geko 0 \Longleftrightarrow \bP (\mathtt{Ext} \leko \infty) \eqo 1 \Longleftrightarrow \int_1^{\infty} \frac{du}{\psi (u)}  <\infty .
\end{equation}
The latter integral condition is called the \emph{Grey condition} (cf.~Bingham \cite{Bi2} for more details). Let us mention that the Grey condition necessarily implies that the CSBP($\psi, x$) $(Z^{x}_s)_{s\in \bbR_+}$ 
has infinite-variation sample paths, but the converse does not hold. 
Under the Grey condition, a simple argument derived from the integral equation satisfied by $u_s(y)$ implies the following 
\begin{equation}
\label{defandmeaningv}
\bP (\mathtt{Ext} \leko s) \eqo e^{-xv(s)} \quad \textrm{where $v\colon \bbR_+^* \!\! \to \!  \bbR_+^*$ is given by} \quad s= \int_{v(s)}^\infty \frac{du}{\psi (u)}.
\end{equation}

More generally, for all $s,y,z\ino \bbR_+$, 
\begin{equation}
\label{CSBPlvlmss}
 \bE \Big[ \! \exp \Big( -\! zZ^{x}_s  \! -\!  y\! \int_0^s \! \!\!  Z^x_r\,   dr \Big)   \Big] \eqo e^{-xw_s(z,y)} \quad \textrm{where} \quad  \int_{w_s(z,y)}^z \frac{du}{\psi (u)\! -\! y} \eqo  s. 
\end{equation}
When $s \! \rightarrow \! \infty$, this implies 
\begin{equation}
\label{totalmass}
\forall z\ino \bbR_+, \quad \bE \big[e^{ - y \int_0^\infty \! Z^x_r \, dr } \big]=e^{-x \psi^{-1} (y) } \, , 
\end{equation}
where $\psi^{-1}$ is the inverse of $\psi$ which is a one-to-one function from $\bbR_+$ onto itself. Here $\int_0^\infty \! Z^x_r \, dr$ is $\bP$-a.s.~finite and 
represents the total 
size of the population (and ultimately the mass of the underlying L{\'e}vy forest). 
It is distributed as a subordinator with Laplace exponent $\psi^{-1}$ at time $x$. 

More precisely, let $(X_s)_{s\in \bbR_+}\!$ be a c{\`a}dl{\`a}g L{\'e}vy process without negative jumps with initial value $X_0\eqo 0$ and whose Laplace exponent is $\psi$, i.e.~$\psi(y)\eqo \log \bE [e^{-yX_1}]$, $y\ino \bbR_+$. (Sub)criticality for $\psi$ means for $X$ that a.s.~$\liminf_{s\to \infty} X_s \eqo  - \infty$. So for all $x\ino \bbR_+$, it makes sense to define $S_x \eqo \inf \{ s \ino \bbR_+: X_s \! < \! -x\}$ and standard results imply that $(S_x)_{x\in \bbR_+}$ is a  c{\`a}dl{\`a}g subordinator with Laplace exponent $\psi^{-1}\! $, i.e.~$\psi^{-1} (y)\eqo -\! \log \bE [e^{-yS_1}]$, $y\ino \bbR_+$. 

Let us recall the Lamperti transform \cite{La2}: for all $s\ino \bbR_+$ we set $A_s\eqo \inf\{ r \ino \bbR_+\! :\!  \int_0^r \! 
Z^{x}_u \, du \! >\! s \}$ with the convention that $\inf \emptyset\eqo \infty$. We also agree on the convention $Z^{x}_\infty \eqo 0$. Then 
\begin{equation}  
\label{Lamperti}
 \big( Z^{x}_{A_s}\big)_{\! s\in \bbR_+} \overset{\textrm{(law)}}{=} \big( x + X_{s\wedge S_{x}} \big)_{\! s\in \bbR_+} \; .
\end{equation}
For downward skip-free $\mathbb{Z}$-valued L\'evy processes, this transform yields $\bN$-valued continuous-time Galton--Watson processes -- this can be read off from the infinitesimal generator. Then the result (\ref{Lamperti}) for CSBPs can be obtained by limit theorems and by the continuity of Lamperti's time-change proved in Helland \cite{He}.    

Let us denote by $ \mathbf c\ino \bbR_+$ the drift coefficient of the subordinator $(S_x)_{x\in \bbR_+}$ and by $\nu(dz)$ its L{\'e}vy measure. Namely for all $y\ino \bbR_+$, 
\begin{equation}
\label{LKformpsi-1}
\psi^{-1} (y) = \mathbf{c} y + \int_{\bbR^*_+} \!\!\!  \big(1 \! -\! e^{-yz} \big) \nu (dz) , 
\end{equation}
where $\nu$ satisfies $\int_{\bbR^*_+} \! (1\! \wedge \! z) \,  \nu (dz) \! < \! \infty$. We shall need the following properties of $\psi^{-1}$. We recall from Definition \ref{ABCtypesdef} Sato's terminology on branching mechanisms (Types A, B and C).  
\begin{lemma}
\label{psi-1prop} Let $\psi$ be a (sub)critical branching mechanism of the form (\ref{brmech}). We keep the notation from above. Then the following holds true. 
\begin{compactenum}

\smallskip

\item[$(i)$] $ \mathbf c\! >\! 0$ iff $\psi$ is of Type A or B. In that case $1/ \mathbf c\eqo \mathbf a+ \int_{\bbR^*_+} \! z \, \pi(dz) $. 

\smallskip

\item[$(ii)$] If $\psi$ is of Type B or C, then $\nu$ is diffuse and $\nu (\bbR^*_+)\eqo \infty$.  

\smallskip

\item[$(iii)$] $\lim_{h\to 0^+} h\nu ((h, \infty))\eqo 0$. 
\end{compactenum}
\end{lemma}
\noi
\textbf{Proof.} We first prove $(i)$. By elementary arguments, we get the following for all $y\ino \bbR_+^*$. 
$$y^{-1} \psi^{-1} (y)= \mathbf c + \! \int_{0}^\infty\!\!\!  e^{-yu} \nu ((u, \infty)) \, du \quad \textrm{and} \quad  y^{-1}\psi (y) = \mathbf a + \mathbf b y +\!  \int_{\bbR_+^*} \!\!  \pi(dz) \int_0^z \! (1\! -\! e^{-yu}) \, du . $$
This implies that 
$\lim_{y\to \infty} \psi^{-1} (y)/y \eqo  \mathbf c \ino \bbR_+$ 
and that $\lim_{y\to \infty} \psi(y)/y$ is finite iff $\psi$ is of Type A or B and in that case, we see that $1/ \mathbf c \eqo \lim_{y\to \infty} \psi (y)/y\eqo \mathbf a+ \int_{\bbR^*_+} \! z \, \pi(dz)$, by monotone convergence. 

Let us prove $(ii)$. To this end, we assume that $\psi$ is of Type B or C and for all $s\ino \bbR_+$, 
we set 
$I_s\eqo \inf_{r\in [0, s]} X_r$. Then the value $0$ is instantaneous and regular for the 
reflected process $X\! -\! I$ and $-I$ is a local time at $0$ for this process (see Bertoin \cite{Be} Chapter VII). 
We denote by $\mathbf N_\psi$ the corresponding excursion measure of  $X\! -\! I$ above $0$, by 
$(Y_s)_{s\in \bbR_+}$ the canonical c{\`a}dl{\`a}g process and by $\zeta$ the duration of the excursion, namely, 
$\mathbf N_\psi$-a.e.~$Y_0\eqo 0$, $Y_s\ino \bbR^*_+$ for all $s\ino (0, \zeta)$ and $Y_s \eqo 0$ for all 
$s\ino [\zeta, \infty)$. Since $\psi$ 
is a (sub)critical branching mechanism, 
a.s.~$\liminf_{s\to \infty} X_s \eqo-  \infty$ and we get $\mathbf N_\psi$-a.e.~$\zeta \! < \! \infty$. 
Since $0$ is instantaneous and regular, $\mathbf N_\psi$ has an infinite mass and standard arguments 
show that the L{\'e}vy measure 
$\nu(dz)$ of $S$ is $\mathbf N_\psi (\zeta \ino dz)$. This proves that $S$, as a L\'evy process without negative 
jumps, is of Type B. Then Theorem 27.4 in Sato \cite{Sat99} asserts that for all $x\ino \bbR^*_+$, the law of $S_x$ is 
diffuse. We now apply the Markov property under $\mathbf N_\psi$ to get for all $s\ino \bbR^*_+$ that 
$\nu (\{ 2s\}) \eqo \mathbf N_\psi  (\zeta \eqo  2s)\eqo \int_{\bbR^*_+}\!  \mathbf N_\psi( Y_s \ino dx \, ; \zeta \! >\! s) \, 
\bP (S_x\eqo  s)\eqo 0$, which proves $(ii)$. 

 We then prove $(iii)$. For all $\theta\ino (1, \infty)$ and all $h\ino \bbR^*_+$ such that $h\theta \! < \! 1$, we get 
 $$ h\nu ((h, \infty)) =  h\nu ((h, h\theta]) + h\nu ((h\theta, 1])+ h\nu ((1, \infty)) \! \leq \!\! \int_{(h,h\theta]} 
 \!\!\!\!  \!\!\!\! \!\! x \, \nu(dx) + \frac{1}{\theta} \!   \int_{(h\theta, 1]} \!\!\!\! \!\!\!\!  \!\! x \, \nu(dx)+ h\nu ((1, \infty)). $$
Thus we get $\limsup_{h\to 0^+} h\nu ((h, \infty)) \leqo \frac{1}{\theta}    \int_{(0, 1]} x \, \nu(dx)$, which entails the 
desired result since $\theta$ can be taken arbitrarily large.  \cqfd

\smallskip

\noi
\textbf{Unmeasured compact Lévy forests.} Before constructing measured Lévy forests whose branching mechanism is of Type B or C, let us first briefly recall from 
Duquesne and Winkel \cite{DuWi2}  a construction of purely metric compact Lévy forests whose (sub)critical branching mechanism $\psi$ is of the form (\ref{brmech}).
We fix $x_0\ino \bbR_+^*$ and for all $\lambda \ino \bbR_+^*$, we claim that there are $c^\lambda \ino \bbR_+^*$, $\xi^\lambda$, a proper (sub)critical offspring distribution and $\varrho^\lambda$ a law on $\bbN$ such that 
\begin{equation}
\label{lambdaparam}
c^\lambda\eqo \psi'(\lambda), \quad \varphi_{\xi^\lambda} (s) \eqo s+ \frac{\psi(\lambda (1\! -\! s)) }{\lambda \psi'(\lambda)}, \; s\ino [0, 1], \quad  \textrm{and} \quad \varrho^\lambda  \! :=\! \mathtt{Poisson} \, (x_0 \lambda). 
\end{equation}
\noi
\emph{Indeed,} we defined $\xi^\lambda$ by setting 
$\psi'(\lambda) \xi^\lambda  (n)$ $\eqo$ $\mathbf b \lambda \un_{\{ n=2\}} $  $\! + \! \frac{{1}}{{n!}}\int_{\bbR_+^*} \! \pi(dx) $ $ x^n\lambda^{n-1}e^{-x\lambda } $, 
for all integers $n\geqo 2$, and also $\xi^\lambda (1)\eqo 0$ and $\psi'(\lambda) \xi^\lambda (0)\eqo \psi(\lambda) /\lambda$. 
A basic summability argument combined with the Taylor expansion of $\psi$ at $\lambda$, entails the middle equation for $\varphi_{\xi^\lambda}$ in (\ref{lambdaparam}). 
By letting $s\uparrow 1^-$, we can see that $\sum_{n\in \bbN} \xi^\lambda (n)\eqo 1$, and $\xi^\lambda$ 
can therefore be viewed as an offspring 
distribution, which is proper by definition. By taking the derivative at $s\eqo 1$ we get 
$\sum_{n\in \bbN} n\xi^\lambda (n)\eqo (\psi'(\lambda) \! -\! \mathbf a)/ \psi'(\lambda) \leqo 1$, 
which shows that $\xi^\lambda$ is (sub)critical. \cq 

\smallskip

In this theorem we recall from \cite{DuWi2} the definition of the purely metric compact Lévy forests. 
\begin{theorem}
\label{metriccpctLevytrees} Let $\psi$ be a branching mechanism of the form (\ref{brmech}). Let 
$x_0\ino \bbR_+^*$. For all $\lambda \ino \bbR_+^*$, let 
$\xi^\lambda\! , c^\lambda\! , \varrho^\lambda$ be as in (\ref{lambdaparam}) and let 
$\mathbf F^{\lambda}$ be a $\bbT^0_{\! c}$-valued GW($\xi^\lambda\! , c^\lambda\! , \varrho^\lambda$)-forest.  
Then the laws of the $(\mathbf F^{\lambda})_{\lambda \in \bbR_+^*}$ are tight in $(\bbT^0_{\! c} , \dGH)$ iff $\psi$ satisfies the Grey condition (\ref{Grey}). Moreover, in this case $\mathbf F^\lambda \! \! \to \! \mathbf F$ in law in 
 $(\bbT^{0}_{\! c} , \dGH)$ as $\lambda \! \to \! \infty\! :$ 
 $\mathbf F$ is an (unmeasured) $(x_0, \psi)$-Lévy forest, whose law is the unique probability distribution on 
 $(\bbT^0_{\! c} , \dGH)$ such that for all $b\ino \bbR_+^*$, its $b$-length erasure $R_b (\mathbf F)$ is distribued as a GW($\xi^{v(b)}\! , c^{v(b)}\! , \varrho^{v(b)}$)-forest, where $v(\cdot)$ is given by (\ref{defandmeaningv}).  
\end{theorem}
\noi
\textbf{Proof.} See \cite[Thm 3.18]{DuWi2}. \cqfd

\smallskip

\noi
\textbf{Standard measured Lévy forests.} We fix $x_0\ino \bbR^*_+$ and a (sub)critical 
branching mechanism $\psi$ of the form (\ref{brmech}) that is of Type B or C, and we introduce families of GW-forests whose laws are consistent under mass erasure. 
To this end, we first recall notation $\psi^{-1}$, $ \mathbf c$ and $\nu(dz)$ from (\ref{LKformpsi-1}). For all $h,y\ino \bbR^*_+$ we set 
\begin{align} 
\nu(h)\eqo \nu ((h, \infty)), \quad \psi^{{-1}}_h (y) = &   \mathbf c y+\!\! \int_{\bbR^*_+} \!\!\!\big(1 \!  -\!   e^{-yz}\un_{(0, h)} (z) \big) \nu (dz) \nonumber  \\ \;\textrm{and} \qquad & \, \overline{\! \psi}^{{-1}}_h\! (y)\eqo  \mathbf c y +\!\! \int_{(0, h)}  \!\!\! \big(1 \! -\! e^{-yz} \big) \nu (dz) = \psi^{{-1}}_h (y) -\nu(h). \label{psihnuhdef}
\end{align}
We denote by $(S^h_x)_{x\in \bbR_+}$ a subordinator with Laplace exponent $\, \overline{\! \psi}^{{-1}}_h$.
\begin{lemma}
\label{existlaws} Let $x_0\ino \bbR^*_+$ and let $\psi$ a (sub)critical 
branching mechanism of the form (\ref{brmech}) that is of Type B or C. We keep the above notation. 
Then, there are parameters $\xi_h$, $c_h$, $\phi_h$, $(\Lambda^{\! h}_n)_{n\in \bN}$, $\varrho_h$ and 
$(\Lambda^{\! h}_{\varnothing , n})_{n\in \bN}$ as in Definition \ref{meaGWregedef}, with $\xi_h$ proper, such that for all $s\ino [0, 1]$ and all $y\in \bbR_+$ the following transform identities hold
\begin{equation}
 \label{chphih}
 c_h + \phi_h (y) = \psi'( \psi^{{-1}}_h (y))
\end{equation}
\begin{equation}
\label{xihLamh}
\nu(h) e^{-hy} F_{\! \xi_h, c_h,  \Lambda^{\! h}} (s,y) = \psi \big(  \psi^{{-1}}_h (y)\! -\! \nu(h) e^{-hy}s \big) \! -\! y + \nu(h) e^{-hy} s \psi'( \psi^{{-1}}_h (y)),
\end{equation}
\begin{equation}
\label{varrhoh}
 e^{-hy} F^\varnothing_{\!\! \varrho_h, \Lambda^{\! h}_\varnothing} (s,y)= e^{-x_0 (  \psi^{{-1}}_h (y)- \nu(h) e^{-hy}s)}\!\!  -\! 
e^{-x_0 \psi^{{-1}}_h (y)}\! + e^{-hy-x_0\nu(h)}\bE \big[ e^{-y (S^h_{x_0} -h)_+}\!  \big],
\end{equation}
where $F_{{\! \xi_h, c_h, \phi_h, \Lambda^{\! h}} } (s,y)\eqo c_h\sum_{n\in \bN} \xi_h (n) s^n\widehat{\Lambda}^h_n(y) $ and 
$F^\varnothing_{{\!\! \varrho_h,  \Lambda^{\! h}_\varnothing}} (s,y)\eqo \sum_{n\in \bN} \varrho_h (n)s^n \widehat{\Lambda}^h_{\varnothing, n} (y)$, using Laplace transform notation introduced in (\ref{notaLapla}). 
\end{lemma}
\noi
\textbf{Proof.} See Appendix \ref{pfsecexistlaws}. \cqfd 

\smallskip 

To simplify notation we write 
\begin{equation}
\label{simplification}
Q_{\psi}^h = Q_{\xi_h, c_h, \phi_h, \Lambda^h},  \quad F_h \eqo F_{\xi_h, c_h, \Lambda^h}  , \quad P_{x_0, \psi}^h = P_{\xi_h, c_h, \phi_h, \Lambda^h, \varrho_h, \Lambda^h_\varnothing}\quad \textrm{and} \quad   F^\varnothing_h \eqo F^\varnothing_{\!\! \varrho_h, \Lambda^h_\varnothing}. 
\end{equation}
The following lemma shows that the laws $(P_{{x_0, \psi}}^{{h}})_{h\in \bbR_+^*} $ are consistent under mass erasure. 
\begin{lemma}
\label{consistency} Let $x_0\ino \bbR^*_+$ and let $\psi$ a (sub)critical 
branching mechanism of the form (\ref{brmech}) that is of Type B or C. We keep the above notation. 
Let $h_0, h \ino \bbR_+^*$. We set $h_1\! :=\! h_0+h$. Then $P_{{\! x_0, \psi}}^{{h_1}} $ is the law of 
$\mathcal E_h $ under $P_{{\! x_0, \psi}}^{{h_0}} $ . 
\end{lemma}
\noi
\textbf{Proof.} 
By Proposition \ref{prop:eraGW} $(i)$ there are $\xi^*, c^*, \phi^*$ $(\Lambda_n^{\! *})_{n\in \bN}$ as in Definition \ref{meaGWregedef} such that $\mathcal E_h$ under 
$Q_{\psi}^{{h_0}} ( \, \cdot \, | \, \langle \flambda \rangle \! >\! h)$ has law $Q_{*}\! :=\! Q_{\xi^*\! , c^*\! , \phi^*\! , \Lambda^{\! *}\! }$.  
We first prove that $Q_*\eqo Q^{{h_1}}_{\psi}$. 
To this end, we first compute $\gamma_{h_0}(y)\eqo Q^{{h_0}}_{\psi} [ e^{ -y\langle \flambda \rangle} ]$ from  
Proposition \ref{ballGWmass} $(ii$-$a)$, which asserts that $\gamma_{h_0}(y)(c_{h_0}+ \phi_{h_0} (y))\eqo F_{h_0} (\gamma_{h_0} (y) ,y)$. 
By (\ref{chphih}) and (\ref{xihLamh}), this implies that 
\begin{equation}
\label{gammahcomput}
\nu (h_0) e^{-h_0y} \gamma_{h_0} (y) \eqo \psi^{-1}_{h_0} (y) \! -\! \psi^{-1} (y)\eqo \int_{(h_0, \infty)} \!\! e^{-yz} \nu(dz) \; .
\end{equation}
Consequently, $\nu(h_0)Q^{{h_0}}_{\psi} (\langle \flambda \rangle + h_0 \in \cdot  )\eqo \nu (\, \cdot\,  \cap \, (h_0, \infty))$.   
We recall that $\falpha_h(y)\eqo Q^{{h_0}}_{\psi} [ e^{-y  \langle \flambda \rangle} 
\un_{\{  \langle \flambda \rangle \leq h\}} ]$ and that 
$\alpha\eqo \falpha_h(0)\eqo Q^{{h_0}}_{\psi }(\langle \flambda \rangle \leqo h)$. Therefore, we get  
\begin{equation}
\label{alphhycomput}
\nu(h_0) e^{-h_0y}\falpha_h(y)\eqo \psi^{-1}_{h_0} (y)\! -\! \psi^{-1}_{h_1} (y) \quad \textrm{and} \quad  \nu(h_0) (1\! -\! \alpha) \eqo \nu (h_1)\, . 
\end{equation}
By (\ref{charaQera}) in Proposition \ref{prop:eraGW} $(i)$ and by a somewhat long but elementary computation that is left to the reader, 
 we see that $\phi_*\eqo \phi_{h_1}$ and that $F_{\xi^*\! , c^*\! , \Lambda^{\! *}}\eqo F_{h_1}$, which proves that $Q_*\eqo Q_{\psi}^{{h_1}}$.

By Proposition \ref{prop:eraGW} $(ii)$, there are $\varrho^*$ and $(\Lambda^{*}_{{\varnothing, n}})_{n\in \bN}$ such that 
$\mathcal E_h $ under $P^{{h_0}}_{{\! x_0, \psi}}$ has law $P_{*}\! :=\! P_{\xi^*\! , c^*\! , \phi^*\! , \Lambda^{\! *}\! , \varrho^*\! ,
\Lambda^{\! *}_{\varnothing}}$. We set $F_{\! *}^{\varnothing} \eqo F_{\! \varrho^*\! , \Lambda^{\! *}_{\varnothing, \cdot}}$ and we show that 
$F_{\! *}^{\varnothing}\eqo F_{{h_1}}^{\varnothing}$, which entails $P_{*}\! :=\! P^{{h_1}}_{{\! x_0, \psi}}$, since $\phi_*\eqo \phi_{h_1}$ and 
$Q_*\eqo Q_{\psi}^{{h_1}}$. To this end, we  multiply (\ref{charaPera}) by $e^{-h_0y}$ to get 
\begin{align}
e^{-h_1y}F^{{\varnothing }}_{\! *} (s,y)& =\underbrace{e^{-yh_0}F^{{\varnothing }}_{{h_0}} \big( (1\! -\! \alpha)e^{-hy} s+\falpha_h(y) \, , \, y\big)}_{=:T_1 (s,y)} \, - \, \underbrace{e^{-yh_0}F^{{\varnothing }}_{{h_0}} \big( \falpha_h(y)  ,  y \big)}_{=T_1(0,y)} \nonumber \\
& \qquad \qquad +\,  \underbrace{e^{-h_1y}\varphi_{\! \varrho_{h_0}}(\alpha) \bE \big[ e^{-y (M_{h_0} -h)_+}\big]}_{=:T_2} \, , \label{Fstar1}
\end{align}
where $\varphi_{\! \varrho_{h_0}}(\alpha) \bE \big[ e^{-yM_{h_0} }\big] \eqo F^{{\varnothing }}_{{h_0}} ( \falpha_h(y)  ,  y )$. 
Recall that $\varphi_{\varrho_{h_0}}(s)\eqo \exp (-x_0 \nu(h_0) (1\! -\! s))$, since $\varrho_{h_0}$ is a Poisson r.v.~with parameter $x_0 \nu(h_0)$. We then take $h\eqo h_0$ in (\ref{varrhoh}) and we use (\ref{alphhycomput}) to get 
\begin{equation}
\label{DeltaT}
\varphi_{\varrho_{h_0}} (\alpha)\eqo e^{-x_0 \nu(h_1) }\quad \textrm{and} \quad T_1(s,y)\! -\! T_1(0,y)\eqo 
e^{ -x_0 \big( \psi_{h_1}^{-1} (y)  - \nu (h_1)e^{-h_1y} s \big)} \! -\! e^{ -x_0  \psi_{h_1}^{-1} (y) }. 
\end{equation}
We next compute the law of $M_{h_0}$. For any $h'\ino \bbR_+^*$, we recall that $S^{{h'}}$ is a subordinator with 
Laplace exponent $\, \overline{\! \psi}^{{ -1}}_{{h'}}\eqo \psi^{{-1}}_{{h'}} \! -\! \nu(h')$. Then, we take $h\eqo h_0$ in 
(\ref{varrhoh}) and we use (\ref{alphhycomput}) and the first equality in (\ref{DeltaT}) to get 
$$ \bE \big[ e^{-yM_{h_0}}\big]=  \bE \big[ e^{-y(S^{h_1}_{x_0}-h_0)}\big] -e^{-x_0 \nu ((h_0, h_1]) } \bE \big[ e^{-y(S^{h_0}_{x_0}-h_0)}
 \! -\! e^{-y(S^{h_0}_{x_0}-h_0)_+} \big] .$$
Basic results on Poisson point processes allow to couple the subordinateors $S^{h'}$, $h'\ino \bbR_+^*$ in the following way. Let $S$ be a 
subordinator with Laplace exponent $\psi^{-1}$; for all $s, h'\ino \bbR_+^*$, we set $S^{h'}_s\eqo S_s  - 
\sum_{r\in [0, s]}\!  \Delta S_r \un_{\{ \Delta S_r > h'\}}$, where $\Delta S_r\eqo S_r \! -\! S_{r-}$ is the jump of $S$ at time $r$. 
Then we set 
$D(x_0)\eqo S^{h_1}_{x_0} \! -\! S^{h_0}_{x_0}$, which is independent from $S^{h_0}$ and 
$\bP (D(x_0)\eqo 0)\eqo e^{-x_0\nu((h_0, h_1]) }$. 
Moreover, if $D(x_0)\eqo 0$, then $S^{h_0}_{x_0}\eqo S^{h_1}_{x_0}$ and  if $D(x_0)\! \neq \! 0$, 
then $S^{h_1}_{x_0} \geko h_0$. Thus we get 
\begin{eqnarray*}
\bE \big[ e^{-yM_{h_0}}\big] \!\!\! & =& \!\!\! \bE \big[e^{-y(S^{h_1}_{x_0}-h_0)} - \big(e^{-y(S^{h_0}_{x_0}-h_0)} \! -\! 
e^{-y(S^{h_0}_{x_0}-h_0)_+}\big) \un_{\{ D(x_0)=0\}} \big]   \\
 \!\!\! &= &  \!\!\! \bE \big[e^{-y(S^{h_1}_{x_0}-h_0)} - \big(e^{-y(S^{h_1}_{x_0}-h_0)} \! -\! e^{-y(S^{h_1}_{x_0}-h_0)_+}\big) 
\un_{\{ D(x_0)=0\}} \big]  \\
\!\!\! &= &  \!\!\! \bE \big[e^{-y(S^{h_1}_{x_0}-h_0)} \un_{\{ D(x_0)\neq 0\}}+ e^{-y(S^{h_1}_{x_0}-h_0)_+}
\un_{\{ D(x_0)=0\}} \big] = \bE \big[e^{-y(S^{h_1}_{x_0}-h_0)_+ } \big]
\end{eqnarray*}
Consequently, $M_{h_0}$ has the same as $(S^{h_1}_{x_0}-h_0)_+$ and $(M_{h_0}\! -\! h)_+$ has the same law as $((S^{h_1}_{x_0}-h_0)_+\! -h)_+\eqo 
(S^{h_1}_{x_0}-h_1)_+$. Combined with (\ref{Fstar1}) and (\ref{DeltaT}), this implies that 
$F^{{\varnothing}}_{\! *}\eqo F^{{\varnothing}}_{\! h_1}$, which completes the proof as mentioned above. \cqfd 

\smallskip

In the following theorem, we define standard $(x_0, \psi)$-Lévy forests (see the Remark \ref{remconstructLevy} $(\mathbf b)$  for an explanation of the use of the term \emph{standard} and also for other important comments). 
\begin{theorem}
\label{thm:growthprocconvbis} Let $x_0\ino \bbR^*_+$ and let $\psi$ a (sub)critical 
branching mechanism of the form (\ref{brmech}) that is assumed to be of Type B or C. We keep the above notation. 
Then the following holds true. 
\begin{compactenum}

\smallskip

\item[$(i)\! $] There is a $\bT$-valued r.v.~$\mathbf m \! \equiv \! (\mathbf F, d, \rho, \mathrm m)$ such that $\mathcal E_h \mathbf m $ has law $P_{\! x_0, \psi}^h$, for all $h\ino \bbR_+^*$. This uniquely characterizes the law of $\mathbf m$, which is called the \emph{standard measured ($x_0, \psi$)-Lévy forest}: $x_0$ is its \emph{initial population} and $\psi$ its \emph{branching mechanism}.

\smallskip

\item[$(ii)\! $] $\big( \flambda_{\mathbf m} ([0, s]) \big)_{\! s\in \bbR_+} \!\! \!\!\! \!\overset{\textrm{(law)}}{=}\! \! \big( \mathbf c Z^{x_0}_s\!  +\!  \int_0^s \! Z^{x_0}_r \, dr  \big)_{\! s\in \bbR_+}\!$, where 
$Z^{x_0}_{\cdot} \!\! \overset{\textrm{(law)}}{=}\! \!\textrm{CSBP}(x_0, \psi)$ and $\mathbf c$ is as in (\ref{LKformpsi-1}).

\smallskip

\item[$(iii)$] $\bP(\mathbf m \ino \bbT_{\! c} )\geko 0 \! \Longleftrightarrow\! \bP(\mathbf m \ino \bbT_{\! c} )\eqo 1\!  
\! \Longleftrightarrow\! $ (\ref{Grey}),  and in this case, $\Phi_0(\mathbf m)$ is the unmeasured compact $(x_0, \psi)$-Lévy forest defined in Theorem \ref{metriccpctLevytrees}.
\end{compactenum}
\end{theorem}
\noi
\textbf{Proof.} By Lemmas \ref{consistency} and \ref{consistlaw} 
there exists a unique probability distribution $P_{\! x_0, \psi}$,  on $(\bT^*\! , \dera)$ such that 
$P_{\! x_0, \psi} \! \circ \! \mathcal E_{h}^{{-1}}\eqo P^{h}_{{\! x_0, \psi}}$ for all $h\ino \bbR^*_+$. 
Namely, there is a $\bT^*$-valued r.v.~$\mathbf m\equiv (\mathbf F, d, \rho, \mathrm m)$ such that for all 
$h\ino \bbR^*_+$, $\mathcal E_h \mathbf m $ has law $P_{{\! x_0, \psi}}^{h}$. To prove $(i)$, i.e.~$\bP (\mathbf m\ino \bbT)\eqo 1$, and $(ii)$, we claim that it is sufficient to prove 
 \begin{equation}
\label{branchLevybis}
\Big( \langle \mathbf m \rangle \, , \, \big( \flambda_{\mathbf m} ([0, r]) \big)_{\! r\in \bbR_+}\Big)  \overset{\textrm{(law)}}{=} \Big(  \int_0^\infty  \!\! Z^{x_0}_r \, dr \, \, , \, \Big( \mathbf c Z^{x_0}_r + \! \int_0^r \!\! Z^{x_0}_u \, du  \Big)_{\! \! r\in \bbR_+}\, \Big)
\end{equation}
\emph{Indeed,} (\ref{branchLevybis}) implies $(ii)$, trivially. To see that \eqref{branchLevybis} also implies $\bP (\mathbf m \ino \bbT)\eqo 1$, we set $M_r\eqo \mathbf c Z^{x_0}_r+ \int_0^{r} \! Z^{x_0}_u du$, for all $r\ino \bbR_+$. We recall from (\ref{totalmass}) that $ \int_0^{r} \! Z^{x_0}_u du\leko \infty $ a.s.~and we observe that 
$M_\infty \! :=\! \lim_{r\uparrow \infty} \uparrow M_r \eqo \int_0^\infty \! Z^{x_0}_u du $, since $\lim_{r\to \infty} Z^{x_0}_r \eqo  0$ 
almost surely. 
By (\ref{branchLevybis}), $( \flambda_{\mathbf m} ((r, \infty]))_{r\in \bbR_+}$ has the same law as $(M_\infty\! -\! M_r)_{r\in \bbR_+}$. Thus $\mathbf{P}$-a.s.~$\lim_{r\to \infty}$ $\flambda_{\mathbf m} ((r, \infty]) $ $\eqo 0$, 
and $\mathbf m \ino \bT$ a.s. \cq

It remains to show (\ref{branchLevybis}), which is the main part of this proof.  
We recall that (\ref{gammahcomput}) asserts that $\nu(h)e^{-hy}\gamma_h(y)\eqo \psi^{-1}_h(y)\! -\! \psi^{-1} (y) $. This implies $\nu(h) (1\! -\! e^{-hy})$ $ + e^{-hy} \nu(h) (1\! -\! \gamma_h(y)) \eqo \psi^{-1} (y) \! -\! \int_{(0, h)}(1\! -\! e^{-yz}) \nu(dz)$. 
By Lemma \ref{psi-1prop} $(iii)$ we get 
\begin{equation}
\label{gammalimi}
\forall y\ino \bbR_+, \qquad \lim_{h\to 0^+} \nu(h) (1\! -\! \gamma_h(y))  = \psi^{-1} (y) \; .
\end{equation}
We next compute for all $h\ino \bbR^*_+$ and all $r,y,z\ino \bbR_+$ the function  
\begin{equation}
\label{ghdef}
g_r^h (e^{-z}\! , y)\eqo Q_\psi^h \big[e^{-zZ^{+}_r - y \flambda ([0, r]) } \big].
\end{equation} 
By (\ref{discrGrey}) in Lemma \ref{ballGWmass} $(iib)$ and an elementary computation, we get  
\begin{equation}
\label{gcomputbis}
\psi^{-1}_h (y)\! -\! \nu(h) e^{-hy} g^h_r (e^{-z}\! , y)\eqo  w_r \big( \psi^{-1}_h (y)\! -\! \nu(h) e^{-hy -z} \!  , y  \big) \;  \textrm{where} \, \int_{w_r(q,y)}^{q} \! \frac{du}{\psi(u)-y}\eqo r 
\end{equation}
for all $q\ino \bbR^*_+$ distinct from $\psi^{-1}(y)$, by (\ref{CSBPlvlmss}). We next observe that 
$$ \psi^{-1}_h(y)\! -\! \nu(h) e^{-hy}e^{-z/\nu(h)}\eqo \mathbf c y + \! \int_{(0, h)}\!\!\! \!\!\!\!   (1\! -\! e^{-yz'})\nu(dz') + \nu(h)  (1\! -\! e^{-hy})+ e^{-hy} \nu(h)  (1\! -\! e^{-z/\nu(h)}) ,  $$
which implies that $\lim_{h\to 0^+}  \psi^{-1}_h(y)\! -\! \nu(h) e^{-hy}e^{-z/\nu(h)}\eqo \mathbf c y + z$ since $\lim_{h\to 0^+} h\nu(h)\eqo 0$ by Lemma \ref{psi-1prop} $(iii)$. 
This combined with (\ref{gcomputbis}) implies 
\begin{equation}
\label{gcomputter}
\lim_{h\to 0^+} \nu(h) \big(1 \! -\! g^{h}_r \big(e^{-z/\nu(h)} ,y \big) \big) \eqo w_r \big( \mathbf c y +z,y\big) -\mathbf c y \; .
\end{equation}
We next fix $y\ino \bbR_+$. We derive from Lemma \ref{ballGWmass} that $z\ino \bbR_+ \! \mapsto \! J_h(z,y) \! :=\!  
\nu(h) (1 \! -\! g^{h}_r \big(e^{-z/\nu(h)} ,y ) )$ is monotonic and continuous. Similarly, we derive from (\ref{CSBPlvlmss}) that 
$z\ino \bbR_+ \! \mapsto\!  J(z,y)\! :=\!  w_r \big( \mathbf c y +z,y\big) \! -\! \mathbf c y $ is monotonic too. 
By Dini's theorem $\lim_{h\to 0^+}J_h(\cdot, y) \eqo  J(\cdot, y)$ uniformly on every compact subset of $\bbR_+$. 
For all $h\ino \bbR^*_+$, let $q_h \ino [0, 1]$ be such that $q\! :=\! \lim_{h\to 0^+}\nu(h) (1\! -\! q_h)$ exists. 
Then the previous arguments imply that  
\begin{equation}
\label{gcomputquart}
\lim_{h\to 0^+} \nu(h) \big(1 \! -\! g^{h}_r \big( q(h) ,y \big) \big) \eqo w_r \big( \mathbf c y +q,y\big) -\mathbf c y \; .
\end{equation}

We next fix $n\ino \bN^*$, $\mathbf r\eqo (r_1, \ldots, r_n)$ and $\mathbf y\eqo (y_1,  \ldots, y_n)$ 
in $\bbR_+^n$ and $z\ino \bbR_+$. We assume that such that 
$r_1\leqo \cdots \leqo r_n$ and we set $\overline{\mathbf y}\eqo y_1+ \cdots+ y_n$. 
We derive from the Markov property of $(Z^{x_0}_r)_{r\in \bbR_+}$ that there are function 
$\mathtt W_n(\mathbf r, \mathbf y,z)$ such that 
\begin{equation}
\label{limiadditit}
 \bE \Big[\exp \Big( \! \! - z \! \int_0^\infty\!\! Z^{x_0}_r dr -\!\!\!  \sum_{1\leq j\leq n}y_j \Big(\mathbf c Z^{x_0}_{r_j} +\!\! 
 \int_0^{r_j} \! \!Z^{x_0}_r dr \Big) \Big) \Big] = \exp \big( \! -x_0 \mathtt W_n (\mathbf r, \mathbf y,z) \big) , 
\end{equation} 
where $\mathtt W_n (\mathbf r, \mathbf y,z) \eqo w_{r_1} \big( \mathbf c y_1+ \mathtt W_{n-1}
 (\mathbf r', \mathbf y', z) , z + \overline{\mathbf y}\big)$ and $\mathtt W_1 (r_1, y_1, z) \eqo
 w_{r_1} (\mathbf c y_1+ \psi^{-1} (z) , z + y_1)$. Here 
$w_r(\cdot, \cdot)$ is as in (\ref{CSBPlvlmss}), $ \mathbf r'\eqo (r_2-r_1, \ldots, r_n-r_1)$ and $\mathbf y'\eqo (y_2, \ldots, y_n)$. 
Similarly, (\ref{discrGrey}) in Lemma \ref{ballGWmass} $(ii$-$b)$ combined with the regenerative property in Proposition \ref{branchingprop} implies that 
$$ f_n (\mathbf r, \mathbf y,z)\eqo Q_\psi^{h} \Big[ e^{-z \langle \flambda\rangle}
 \!\! \! \prod_{1\leq j\leq n} \!\! e^{-y_j \flambda ([0, r_j])}\Big] $$
satisfies $f_{n} (\mathbf r, \mathbf y,z)\eqo g^h_{r_1} \big( f_{n-1} (\mathbf r', \mathbf y', z) ,z + \overline{\mathbf{y}} \big)$ and 
$f_{1} (r_1, y_1, z) \eqo g^h_{r_1} \big(\gamma_h (z) , z+y_1 \big)$, where $g^h_r$ is given in (\ref{ghdef}). Then we recursively prove that for all $\mathbf r, \mathbf y, z$  
\begin{equation}
\label{fnlimit}
\nu (h) \big(1\! -\! f_n (\mathbf r, \mathbf y, z) \big) \underset{h\to 0^+}{-\!\!\!-\!\!\!-\!\!\! \longrightarrow} \mathtt W_n (\mathbf r, \mathbf y, z) - \mathbf c \big( z+ \overline{\mathbf y}\big) .
\end{equation}
\emph{Indeed,} (\ref{gammalimi}) and (\ref{gcomputquart}) imply (\ref{fnlimit}) for $n\eqo 1$. Suppose that (\ref{fnlimit}) holds for $n-1$. Then we apply (\ref{gcomputquart}) to $q(h)\eqo  f_{n-1} (\mathbf r', \mathbf y', z)$, $r\eqo r_1$ and $y\eqo z + \overline{\mathbf y}$, and we get $\lim_{h\to 0^+} \nu (h)(1\! -\! f_n (\mathbf r, \mathbf y, z) )\eqo 
w_{r_1} (\mathbf c (z + \overline{\mathbf y}) + \mathtt W_{n-1} (\mathbf r', \mathbf y', z) \! -\!  
\mathbf c ( z+ \overline{\mathbf y^\prime} )) \! -\!  \mathbf c ( z+ \overline{\mathbf y})$, which is equal to $\mathtt W_{n} (\mathbf r, \mathbf y, z)\! -\!  \mathbf c ( z+ \overline{\mathbf y})$ and which completes the proof of (\ref{fnlimit}). 

By (\ref{forestintrin}) we then get 
\begin{equation}
\label{forestadditive}
P^h_{x_0, \psi} \Big[  e^{-z \langle \flambda\rangle}
 \!\! \! \prod_{1\leq j\leq n} \!\! e^{-y_j \flambda ([0, r_j])}\Big]= F^\varnothing_h \big( f_n (\mathbf r, \mathbf y, z),  z+ \overline{\mathbf y} \big)
\end{equation}
where we recall $F^\varnothing_h$ from (\ref{simplification}) and (\ref{varrhoh}). We easily check that 
$ \lim_{h\to 0^+} \exp (-x_0 \psi^{-1}_h (z+ \overline{\mathbf y}))\eqo 0$ and that $  \lim_{h\to 0^+} e^{-x_0\nu (h)} \bE \big[ 
\exp (-(z+ \overline{\mathbf y})  (S^{h}_{x_0} \! -\! h)_+ )\big]\eqo 0$ since $\nu(h)\! \rightarrow\! \infty$ by Lemma \ref{psi-1prop} $(ii)$. We then observe that 
\begin{eqnarray}
\label{theleftmemberin}
 \psi^{-1}_h (z+ \overline{\mathbf y}) \! -\! \nu (h) e^{-h (z+ \overline{\mathbf y}) }f_n (\mathbf r, \mathbf y, z) \!\!\! &  = &\!\!\! \mathbf c (z+ \overline{\mathbf y}) + \! \int_{(0, h]} \!\! \big(1 \! -\! e^{-(z+ \overline{\mathbf y})z' } \big)\nu (dz'). \\
\!\!\!\! & & \!\!\!\! +\nu(h) \big(1\! -\!  e^{-h (z+ \overline{\mathbf y}) } \big)+ e^{-h (z+ \overline{\mathbf y}) }\nu(h) \big( 1\! -\!  f_n (\mathbf r, \mathbf y, z)\big). \nonumber
\end{eqnarray}
Since $\lim_{h\to 0^+}h\nu(h)\eqo 0$ (by Lemma \ref{psi-1prop} $(iii)$) and by (\ref{fnlimit}), the right member in (\ref{theleftmemberin}) tends to $\mathtt W_n (\mathbf r, \mathbf y, z)$ and we get that 
$$ \lim_{h\to 0^+}P^h_{\! x_0, \psi} \Big[  e^{-z \langle \flambda\rangle}
 \!\! \! \prod_{1\leq j\leq n} \!\! e^{-y_j \flambda ([0, r_j])}\Big]= e^{-x_0 \mathtt W_n (\mathbf r, \mathbf y, z)} \; .$$
By (\ref{limiadditit}), this shows that the finite-dimensional marginals of $(\langle \flambda \rangle , (\flambda([0, r]))_{r\in \bbR_+})$ under $P^{h}_{{\! x_0, \psi}}$, i.e.~the finite-dimensional marginals of $(\langle \flambda_{\mathcal E_h \mathbf m} \rangle , (\flambda_{\mathcal E_h \mathbf m }([0, r]))_{r\in \bbR_+})$ under $\bP$, 
converge in law to the finite-dimensional marginal laws of $(M_\infty, (M_r)_{r\in \bbR_+ })$ where we recall the notation $M_r\eqo \mathbf c Z^{x_0}_r+ \int_0^{r} \! Z^{x_0}_u du$ for all $r \ino \bbR_+$. 
Now we recall from Lemma \ref{cvlawmeasR+} that there is a deterministic set $C\subseteq\mathbb{R}_+$ such that 
$\mathbb{R}_+\! \setminus \! C$ is countable and $\mathbf{P}$-a.s., $\flambda_\mathbf{m}(\{r\})\eqo 0$ for all $r\ino C$. We then apply Proposition \ref{weakcvera} to $(\mathbf{F},d^*\! ,\rho,\mathrm{m})$ and the balls $B_\mathbf{F}(\rho,r)\eqo \{\sigma\ino \mathbf{F} \! :\! 
d^*(\rho,\sigma)\leqo 1\! -\! e^{-r}\}$, it follows that $\mathbf{P}$-a.s.~for all $r\ino C$, $\lim_{h\to 0^+} \flambda_{\mathcal E_h \mathbf m} ([0, r])
=\flambda_{\mathbf m} ([0, r])$, which entails (\ref{branchLevybis}), by the right-continuity of $r\mapsto\flambda_\mathbf{m}([0,r])$. As already mentioned  (\ref{branchLevybis}) implies $(i)$ and $(iii)$.

 Let us prove $(iii)$. Let us first assume that $\psi$ satisfies the Grey condition (\ref{Grey}). We observe that $\Phi_0 (\era_h \mathbf m) $ has law $\texttt{P}_{\! \xi^{\nu(h)}\! , c^{\nu(h)}\! , \varrho^{\nu(h)}}$. Since $\nu(h)\! \to \! \infty$ as $h\! \to \! 0$, because $\psi$ is of Type B or C (see Lemma \ref{psi-1prop}), Theorem \ref{metriccpctLevytrees} applies and yields that $\Phi_0 (\era_h \mathbf m)$ converges in law to the unmeasured compact ($x_0, \psi$)-Lévy forest as $h\! \to \! 0$. This implies in particular that the laws of the 
 $\Phi_0 (\era_h \mathbf m)$, $h\ino \bbR_+^*$, are tight on $(\bbT^0_{\! c}, \dGH)$, and Lemma \ref{Tcrandomcriterion} implies that $\bP (\mathbf m\ino \bbT_{\! c} )\eqo 1$ and that a.s.~$\lim_{h\to 0} \dGH(\Phi_0 (\era_h \mathbf m), \Phi_0 (\mathbf m) )\eqo 0$. Thus $\Phi_0 (\mathbf m)$ is distributed as an  unmeasured compact ($x_0, \psi$)-Lévy forest. 
 
 Let us now assume that $\psi$ does not satisfy the Grey condition (\ref{Grey}). Since $\Phi_0 (\era_h \mathbf m) $ has law $\texttt{P}_{\! \xi^{\nu(h)}\! , c^{\nu(h)}\! , \varrho^{\nu(h)}}$, (\ref{heightGWexpli}) and a simple computation imply that $\mathtt{Ht} (\era_h \mathbf m) \! \uparrow \! \infty$, in probability as $h\! \downarrow \! 0$. We then easily check that 
 a.s.~$\mathtt{Ht} (\mathbf m) \geqo \mathtt{Ht} (\era_h \mathbf m )$, which a.s.~implies $\mathtt{Ht} (\mathbf m)\eqo \infty$ Thus $\bP (\mathbf m \! \notin \! \bbT_{\! c})\eqo 1$, which completes the proof of $(iii)$.  \cqfd 
\begin{remark}
\label{remconstructLevy} Let $x_0\ino \bbR^*_+$, let $\psi$ be a (sub)critical 
branching mechanism of the form (\ref{brmech}) that is of Type B or C and let $(\mathbf F, d, \rho, \mathrm m)\! \equiv \! \mathbf m \ino \bbT$ be a standard $(x_0, \psi)$-Lévy forest as in Theorem \ref{thm:growthprocconvbis}.
 
\smallskip

\noi
$(\textbf{a})$ When $\psi$ is of type C, i.e.~when the corresponding Lévy process has infinite variation sample paths, we claim that 
$\mathbf m$ has the same law as the $\bbR$-tree encoded by the $\psi$-height process $H$ introduced in Le Gall and Le Jan \cite{LGLJ1} and Duquesne and Le Gall \cite{DuLG}, where it is proved that in compact cases, $H$ is continuous, otherwise it is not continuous but lower semicontinuous. More details will be given in a forthcoming work. 
The definition of Lévy forests when $\psi$ is of type B seems new: an alternative, and somewhat more direct construction is discussed below in $(d)$.

\smallskip

\noi
$(\textbf{b})$ Like measured GW trees, standard $(x_0, \psi)$-Lévy forests satisfy a branching property. In the case of unmeasured compact $\bbR$-trees, i.e.~on $(\bbT^0_{\! c}, \fdelta_{\mathtt{GH}} )$, Weill \cite{Weill} has proved that the branching property 
characterizes the laws of Lévy trees. However in the case of measured trees, i.e.~on $(\bbT, \fdelta_{\mathtt{GP}} )$ or 
$(\bbT^1_{\! c}, \fdelta_{\mathtt{GHP}} )$, this no longer holds. 
Indeed, even in compact cases, a Lévy forest with initial population $x_0$ and branching mechanism $\psi$ can be equipped 
with finite Borel measures in multiple ways, with certain atoms in branch points, while preserving the branching property.  
Even though it seems possible to fully describe all distributions of $\bbT$-valued r.v.s that satisfy the branching property, we leave this question for future work. We however think that \emph{among Lévy forests with initial population $x_0$ and branching mechanism $\psi$, the law introduced in Theorem \ref{thm:growthprocconvbis} is the only one whose total mass has the same distribution as the total population of a CSBP$(x_0, \psi)$. 
This is why we use the term \emph{standard} to refer to them}. This includes branching mechanisms of Type B whose measure is purely atomic at branch points (see $(c)$ and $(d)$ below).

\smallskip

\noi
$(\textbf{c})$ If $\psi$ is of type C, then $\mathbf c\eqo 0$ in (\ref{LKformpsi-1}) by Lemma \ref{psi-1prop} $(i)$ and $(\flambda_{\mathbf m} ([0, s]) )_{s\in \bbR_+}$ has the same law as $( \int_0^s \! Z^{x_0}_r \, dr )_{s\in \bbR_+}$ by Theorem \ref{thm:growthprocconvbis} $(ii)$. Therefore, $\bP$-a.s.~$\mathbf m$ is diffuse. 
On the other hand, if $\psi$ is of Type B, $\mathbf c \geko 0$ by Lemma \ref{psi-1prop} $(i)$, and Theorem \ref{thm:growthprocconvbis} $(ii)$ implies that $\mathbf m$ has atoms, indeed, $(Z^{{x_0}}_{s})_{s\in \bbR_+}$ has jumps (this is clarified by the direct construction provided in $(d)$). 

\smallskip

\noi
$(\textbf{d})$ Let us assume that $\psi$ is of Type B. Then standard $\psi$-Lévy forests can be obtained by the following construction. For all 
$u\ino \bbU$, let $\mathcal N_u \! :=\! \sum_{j\in \bbN^*} \delta_{(\ell_{u\ast (j)} , w_{u\ast (j)}, y_{u\ast (j)} )} $, be i.i.d.~Poisson 
point measures on $(\bbR_+^*)^3$ with intensity $e^{-t/\mathbf c} dt \, \pi (dw) \,  dy$. We also set 
$\ell_\varnothing \eqo 0$, $w_{\varnothing}\eqo  x_0$, and $ y_\varnothing \eqo 0$ (here, $y_\varnothing$ 
plays no role). We then set 
$\tau_{x_0}\! :=\! \{ \varnothing \} \cup \{ u\ino \bbU:  y_v \leqo w_{\overleftarrow{v}}, \, 
\forall v\!  \in \, \rgeo \varnothing , u\rgeo \} $, where $\rgeo\varnothing,u\rgeo$ is the ancestral line of $u$ in $\bbU\setminus\{\varnothing\}$. This random subset $\tau_{x_0}$ is a subtree of $\bbU$ in the sense that 
$\overleftarrow{u}\ino \tau_{x_0}$ 
for every $u\ino \tau_{x_0}\backslash \{ \varnothing\} $. We then 
set $\mathbf m \! \equiv \! (\mathbf F, d, \rho, \mathrm m)\! :=\! \mathtt{Tree} (\tau_{x_0}, (\ell_u, \mathbf c w_u, 0)_{u\in \tau_{x_0}} )$ 
(with an obvious extension of the construction provided in Definition \ref{ellonemodel} to infinite trees). Note that $\mathrm m$ 
is a sum of atoms at branch points. Then we claim, without providing a complete and formal proof, that the equivalence class $\mathbf m$ of 
$(\mathbf F, d, \rho, \mathrm m)$ is a standard $(x_0, \psi)$-Lévy forest. 

To see this, we fix $\varepsilon \ino (0, x_0)$ and set $\tau_{x_0, \varepsilon}\!:=\!
 \{ u\ino \tau_{x_0}: w_v\geko \varepsilon, \forall v\ino \lgeo \varnothing, u\rgeo\}$  and 
$\mathbf m_{[\epp]} \! \equiv \!(\mathbf F_{[\epp]} , d, \rho, \mathrm m_{[\epp]})\! :=\!  
\mathtt{Tree} (\tau_{x_0,\varepsilon}, (\ell_u, \mathbf c w_u, 0)_{u\in \tau_{x_0 ,\varepsilon}} )$, which is a subtree of 
$\mathbf m$ in the sense that $\mathbf F_{[\epp]}\!\subseteq\! \mathbf F$ and 
$\mathrm m\! -\! \mathrm m_{[\epp]}$ is a nonnegative measure. We also observe that
 a.s.~$\mathrm m_{[\epp]}\! \to \! \mathrm m$ weakly on $\mathbf F$ (and here, even in the total variation distance). 
 This implies $\mathbf m_{[\epp]} \! \to \! \mathbf m$ a.s.~in $(\bbT, \dGP)$ and by Theorem \ref{eraGcont} $(i)$, 
 $\era_h \mathbf m_{[\epp]} \! \to \! \era_h \mathbf m$ a.s.~in $(\bbT, \dGP)$ for all $h \ino \bbR^*_+$. 
 On the other hand, basic computations imply that $\mathbf m_{[\epp]}$ is a GW forest whose parameters are denoted by 
 $\xi_{[\epp]}$, $c_{[\epp]}$, $\phi_{[\epp]}$, $(\Lambda^{{\! [\epp]}}_n)_{n\in \bN}$, $\varrho_{[\epp]}$ and 
$( \Lambda^{{\! [\epp]}}_{{\varnothing , n}})_{n\in \bN}$ and are given by the following: $c_{[\epp ]} \eqo 1/\mathbf c$, 
$\phi_{[\epp]}\! \equiv \! 0$ and for all $s\ino [0, 1]$ and $y\ino \bbR_+$ 
$$ \mathbf c \pi(\epp) F_{\! \xi_{[\epp]}, c_{[\epp]} , \Lambda^{\! [\epp]}} (s,y) \! :=\! \int_{(\epp, \infty)} \!\!\!\!\! \!\! 
e^{-\mathbf c z (\pi(\epp) (1-s)+y)} \pi(dz)  \quad \textrm{and} \quad   
F_{\!\!  \varrho_{[\epp]},  \Lambda^{\! [\epp]}_{\varnothing}} (s,y)= e^{-\mathbf cx_0 (\pi(\epp) (1-s)+y)} .
$$
where $\pi(\epp)\!:=\! \pi((\epp, \infty))$. Let us set $\psi_{[\epp]} (y) \! :=\! \frac{1}{\mathbf c}y \! -\! \int_{(\epp, \infty)} (1\!-\!e^{-yz})\pi(dz)$, $y\ino \bbR_+$. 
By Proposition \ref{prop:eraGW}, $\era_h \mathbf m_{[\epp]}$ is a GW forest whose parameter transforms are explicitly obtained by taking $\psi_{[\epp]}$ instead of $\psi$ in (\ref{chphih}), (\ref{xihLamh}) and (\ref{varrhoh}) (the computation, which is simple but somewhat long, is left to the reader). As $\epp \!\to \! 0$, we easily see that parameter transforms of  $\era_h \mathbf m_{[\epp]}$ converge to $(\phi_h, F_h , F^{\varnothing}_h)$, which implies, by 
Proposition \ref{proppoubelle}, that $\era_h \mathbf m_{[\epp]}$ converges in law to $P^{h}_{{\! x_0, \psi}}$ in $(\bbT^1_{\! c} , \dGHP)$, for all $h\ino \bbR^*_+$. Thus $\era_h \mathbf m$ has law $P^{h}_{{\! x_0, \psi}}$ and $\mathbf m$ 
is a standard $(x_0, \psi)$-Lévy forest as introduced in Theorem \ref{thm:growthprocconvbis}.

Let us briefly outline some basic facts about $\mathbf m$. 
We see that the length measure of $\mathbf F$ (i.e.~its one-dimensional Hausdorff measure) attributes an infinite value to any open ball. Moreover, $\mathtt{Lf} (\mathbf F)\eqo \emptyset$, $\mathbf F \backslash (\{ \rho\} \cup \mathtt{Br} (\mathbf F)))$ is a countably infinite union of branches, which are isometric to non-empty finite open intervals,
and $\mathtt m (\mathbf F \backslash (\{ \rho\} \cup \mathtt{Br} (\mathbf F)))\eqo 0$: namely $\mathrm m$ 
is the sum of atoms $w_u$ at the branch points corresponding to $u\ino \tau_{x_0}$. Now let $\mathrm m'$ be any finite measure on $\mathbf F$ that enjoys the branching property: this implies that the  $\mathrm m'$-masses  of the countably infinitely many branches of $\mathbf F$ are i.i.d.~and since $\mathrm m'$ is finite, they must be equal to $0$. Thus $\mathrm m'$ is also a sum of atoms at branch points: namely, $(\mathbf F,d,\rho, \mathrm m') \eqo  
\mathtt{Tree} (\tau_{x_0}, (\ell_u, W_u,  0)_{u\in \tau_{x_0}} )$. The branching property also implies that conditionally given 
$(\tau_{x_0} , (\ell_u, w_u,  0)_{u\in \tau_{x_0}}) $, the $W_u$ are independent and the conditional law of $W_u$ only depends on $w_u$. We claim, without providing a proof, 
that $\mathrm m$, which corresponds to $\mathbf c w_u\eqo W_u$, is the only case where the total mass is distributed as the total population of a CSBP ($ x_0, \psi$).
 
Moreover for all $s\ino \bbR_+$, we set $Z^{x_0}_s\eqo \int_{[0, s]} \mathbf c e^{- (s-r)/\mathbf c} \flambda_{\mathbf m} (dr)$, $s\ino \bbR_+$. We claim  that $Z^{x_0}$ is a CSBP$(x_0, \psi)$, which results from the superposition of the following populations: at every branch point/root $\sigma \ino \mathtt{Br} (\mathbf F)\cup\{\rho\}$, whose height we denote by $r\! :=\! d(\rho, \sigma)$, a population is created with initial size $\mathtt m (\{ \sigma \}) \geko 0$ that decreases deterministically at rate $1/\mathbf c$, i.e.~at time $s\geko r$, its size is $\mathtt m (\{ \sigma \}) e^{- (s - r)/\mathbf c}$. These deterministic evolutions do not represent further branching, but are continuum pure death processes. In aggregate, they can be viewed as spreading the mass $\mathtt m(\{\sigma\})$ as dust across time $[r,\infty)$ according to an exponential density, but the lack of further genealogical structure prevents a genealogical representation of this dust as spread-out mass on a measured separable $\bbR$-tree. Therefore, this dust is responsible for the lack of regularity of  bounded variation cases.  \cq  
\end{remark}

\subsection{An invariance principle}
\label{sec:IPnew}
Before stating our main limit theorem, Theorem \ref{thmIP}, which is proven in Section \ref{pfsecthmIP}, we define the parameters and assumptions of this invariance principle.

Let $x_0\ino \bbR^*_+$ and let $\psi$, a (sub)critical 
branching mechanism of the form (\ref{brmech}), whose parameters are $\mathbf a$, $\mathbf b$ and $\pi$ such that $\int_{(0, \infty)} (z\! \wedge \! z^2) \pi(dz) \leko \infty$. We first make the following assumption. 
\begin{compactenum}
\smallskip
\item[\bf A1:] \emph{$\psi$ is of Type B or C.}
\smallskip
\end{compactenum}
Let $(Z^{x_0}_s)_{s\in \bbR_+}$ be a CSBP($\psi, x_0$), whose law is characterized by (\ref{CSBPlvlmss}). 
We recall the notation 
$\mathtt{Ext} \eqo \inf \{ s\ino \bbR_+\! : \! Z^{x_0}_s \eqo 0\}$ (with the convention $\inf \emptyset \eqo \infty$) 
for the extinction time of $Z^{x_0}$, which is a.s.~finite iff $\psi$ satisfies the Grey condition (\ref{Grey}), and whose law is given by (\ref{defandmeaningv}) in that case. 
We denote by $\mathbf m \! \equiv\! (\mathbf F, d,\rho, \mathrm m)$ a $\bbT$-valued r.v.~that is 
distributed as a standard measured ($x_0, \psi$)-Lévy forest as defined in Theorem \ref{thm:growthprocconvbis}, 
which also asserts that 
$\bP (\mathbf m\ino \bbT_{\! c}) \eqo 1$ iff $\psi$ satisfies the Grey condition (\ref{Grey}). 

For all $p\in\mathbb{N}$, we fix $\xi_p,c_p,\phi_p,\Lambda^{\! (p)},\varrho_p$ and $\Lambda^{\! (p)}_\varnothing$ as in Definition \ref{meaGWregedef}. We assume the following. 
\begin{compactenum}
\smallskip
\item[\bf A2:] \emph{$\xi_p$ is proper and (sub)critical, i.e.~$\xi_p(1)\eqo 0$ and $\sum_{k\in\mathbb{N}}k\xi_p(k)\leqo 1$}.
\smallskip
\end{compactenum}
To simplify, we set 
\begin{equation}
\label{shorthandp}
Q_p:=Q_{\xi_p,c_p,\phi_p,\Lambda^{\!(p)}},\;\,  
P_{\! p}:=P_{\! \xi_p,c_p,\phi_p,\Lambda^{\! (p)} \! ,\, \varrho_p,\Lambda_\varnothing^{\! (p)}}, \;\,  F_{\! p} 
:=F_{\! \xi_p,c_p,\Lambda^{\! (p)}}\quad\text{and}\quad F^\varnothing_{\! p}:=F^\varnothing_{\!\!  \varrho_p,
\Lambda_\varnothing^{\! (p)}}.
\end{equation}
Let $\mathbf{m}_p\! \equiv\! (\mathbf{F}_{\! p},d_p,\rho_p,\mathrm{m}_p)$ be a GW forest with law 
$P_{\! p}$. Here we furthermore suppose: 
\begin{compactenum}
   \smallskip
\item[\bf A3:] \emph{$\mathrm{m}_p$ is a deterministic function of $(\mathbf{F}_{\! p},d_p,\rho_p)$ with no mass at $\rho_p$, and $\mathbf m_p$ is minimal. Namely,}
  \begin{compactenum}
    \smallskip
    \item[$(a)\! $] \emph{$\phi_p(y)\eqo \kappa_py$, $y\ino\mathbb{R}_+$, for some $\kappa_p\ino\mathbb{R}_+$;}
    \item[$(b)\! $] \emph{$\Lambda_{n}^{{\! (p)}}(dx) \eqo \delta_{w_p(n-1)}(dx) $, $n\ino\mathbb{N}$, for some weight function $w_p\colon \{-1\}\cup\mathbb{N}\! \rightarrow \! \mathbb{R}_+$;}
    \item[$(c)\! $] \emph{$\Lambda_{{\varnothing,n}}^{{\! (p)}} (dx) \eqo \delta_0(dx) $, $n\ino \mathbb{N}$, i.e.~there is no weight at the root: $\mathrm{m}_p(\{\rho_p\})\eqo 0$;}
\smallskip
    \item[$(d)\! $] \emph{Either $\kappa_p \geko 0$ or $w_p (-1) \geko 0$.}  
  \end{compactenum}
 \smallskip
\end{compactenum} 
Indeed, by Remark \ref{minimalGW}, $\mathbf{A3}$ $(d)$ implies that a.s.~$\mathbf m_p \ino \bbT^1_{\! c, \mathtt{min}}$ and thus $\bP(\Phi^{-1}_1(\mathbf m_p) \ino \bbT_{\! c} )\eqo 1$. This leads to the following. 

\smallskip
\noi
\textbf{Convention.} \emph{When it is convenient, and to simplify the notation, we shall view $\mathbf m_p$ either as a 
$\bbT^1_{\! c, \mathtt{min}}$-valued r.v.~(as in Definition \ref{meaGWregedef}) or as a $\bbT_{\! c}$-valued r.v. (to avoid notation $\Phi^{-1}_1(\mathbf m_p)$).} \cq  

\smallskip

We next make assumptions on the convergence of the branching process 
$Z_{{ s}}^{{(p)}}\! :=\! Z_{s}^{+}(\mathbf{m}_{ p})$, $s\ino\mathbb{R}_+$. Here, $Z_s^{{(p)}}$ is the number of connected
components of $\{\sigma\ino\mathbf{F}_{\! p}   \colon \! d_p(\rho_p,\sigma)\!>\! s\}$ and 
$s\!\mapsto\! Z_s^{{(p)}}$ is c\`adl\`ag on $\mathbb{R}_+$. We then use the notation 
$\mathtt{Ht}_{p} \eqo \inf \{ s\ino \bbR_+ \! :\! Z_s^{{(p)}}\eqo 0\}$ for the a.s.~finite total height of $(\mathbf F_{\! p}, d_p, \rho_p)$, whose law is characterized by (\ref{heightGWexpli}) in Remark \ref{heightGW}.  
We keep the previous notation and in addition to $\mathbf{A1}$,  $\mathbf{A2}$ and  $\mathbf{A3}$, we furthermore assume the following. 
\begin{compactenum}
\smallskip
\item[\bf A4:]  \emph{There exists a sequence $(a_p)_{p\in\mathbb{N}}$ such that}
$$
 a_p\longrightarrow\infty,\quad 
\text{\emph{and for all $s\ino \bbR_+$,}} \quad 
\tfrac{1}{a_p}Z_s^{{(p)}}\!\!  \xrightarrow[p\to \infty]{\textrm{\emph{(law)}}} \, Z^{x_0}_s.
$$
\end{compactenum}
Let us mention that under $\mathbf{A1}$--$\mathbf{A3}$, Assumption $\mathbf{A4}$ is equivalent to an analytical assumption
$\mathbf{A4}'''$, which we discuss in Section \ref{pfsecthmIP}.
Since $Z_0^{{(p)}}$ has law $\varrho_p$, $\mathbf{A4}$ also entails the convergence in law 
$\varrho_p\big(\frac{dy}{a_p}\big)\! \to \! \delta_{x_0}(dy)$ on $\bbR_+$. 
In general, $\mathbf{A1}$--$\mathbf{A4}$ do not necessarily imply a convergence in law of the total heights $\mathtt{Ht}_{ p}$. We shall need the following lemma. 
\begin{lemma}
\label{Greydis} We keep the previous notation. We assume $\mathbf{A1}$--$\mathbf{A4}$. For all $z\ino [0, a_p]$, we set $\psi_p (z) \! :=\! b_p \big( \varphi_{\xi_p} \big(1\! -\! \frac{1}{a_p}z \big) \! -\! 1 + \frac{1}{a_p}z\big)$, where we set $b_p\!:=a_pc_p$. Then $b_p\rightarrow\infty$, and the laws of $(\mathtt{Ht}_{ p})_{p\in \bbN}$ are tight on $\bbR_+$ iff the following condition is satisfied: \vspace{-0.1cm}
\begin{equation}
\label{Greydiscr}
\lim_{z\to \infty} \limsup_{p\to \infty} \int_{z}^{a_p}\!\!  \frac{du}{\psi_p(u)}= 0 \; .
\end{equation} 
In that case, the Grey condition (\ref{Grey}) is satisfied, $\mathtt{Ext}$ is a.s.~finite and $\mathtt{Ht}_{ p} \! \to \! \mathtt{Ext}$ in law on $\bbR_+$. 
\end{lemma}
\noi
\textbf{Proof.} See Duquesne and Rebei \cite{DuRe26prepub}. See also the proof of \cite[Prop.~7.3]{BrDuWa21}.\cqfd  

\smallskip

We next specify our assumptions on the weights $w_p (\cdot)$ and on $\kappa_p$. To this end, we set $L^{{(p)}}_{0}\eqo 0$ and for all $k\ino \bbN^*$, $L^{{(p)}}_{k}\eqo \sum_{1\leq l\leq k} w_p (Y^{{(p)}}_{l})$, where 
$(Y^{{(p)}}_{k})_{k\in \bbN^*}$ is a sequence of independent r.v.s whose law is given by $\bP (Y^{{(p)}}_{k}\! \eqo j)\eqo \xi_p (j+1)$, for all $j \ino \{ -1\} \cup \bbN$. We then assume the following. 
\begin{compactenum}
\smallskip
\item[\bf A5:] $a_p \kappa_p + L^{{(p)}}_{{\lfloor b_p\rfloor}} \!\!  \rightarrow \! 1$ in law in 
$\mathbb{R}_+$, where we set $b_p\!:=\!a_pc_p$. 
\smallskip
\end{compactenum}
Under $\mathbf{A1}$--$\mathbf{A4}$, $\mathbf{A5}$ is equivalent to analytical assumptions $\mathbf{A5}'$, which we discuss in Section \ref{pfsecthmIP}. Assumption $\mathbf{A5}$ is the minimal assumption under which GW forests whose measures deterministically depend on their metric only, converge to \emph{standard} Lévy forests: see Remark \ref{IPcomments} $(c)$.  
The invariance principle is stated in the following theorem. 
\begin{theorem}
\label{thmIP} We keep the previous notation, we assume $\mathbf{A1}$--$\mathbf{A5}$ and we recall that $\mathbf m_p$ is viewed as a $\bbT_{\! c}$-valued r.v.~when it is convenient. Then the following holds true. 
\begin{compactenum}

\smallskip

\item[$(i)$] $\mathbf m_p \! \to \! \mathbf m$ in law in $(\bbT^*\! , \dera)$. 

\smallskip

\item[$(ii)$] $\mathbf m_p \! \to \! \mathbf m$ in law in $(\bbT, \dGP)$ iff $\psi$ is of Type C. In particular, if $\psi$ is of Type B, $(\mathbf m_p)_{p\in \bbN}$ cannot converge in law in $(\bbT, \dGP)$.

\smallskip

\item[$(iii)$] The laws of $(\mathbf m_p)_{p\in \bbN}$ are tight on $(\bbT^1_{\! c}, \dGHP)$ iff (\ref{Greydiscr}) holds true. In this case (\ref{Grey}) holds true, $\bP (\mathbf m \ino \bbT_{\! c})\eqo 1$ and $\mathbf m_p \! \to \! \Phi_1 (\mathbf m)$ in law in $(\bbT^1_{\! c}, \dGHP)$. \smallskip
\end{compactenum}
Moreover, the convergence in law $\mathbf m_p \! \to \! \mathbf m$ holds jointly with the convergence in law in $\bbR_+$ of the total masses $\langle \mathbf m_p\rangle$ to $\langle \mathbf m \rangle$, and $\langle \mathbf m \rangle$ is furthermore distributed as the total population $\int_0^\infty \! Z^{x_0}_s ds$. If $\psi$ is of Type C in $(ii)$ (resp.~if (\ref{Greydiscr}) is satisfied as in $(iii)$), the convergence $\mathbf m_p \! \to \! \mathbf m$ in law in $(\bbT, \dGP)$ (resp.~in law in $(\bbT^{1}_{\! c}, \dGHP)$) holds jointly with $(\flambda_{\mathbf m_p} ([0, s]))_{s\in \bbR_+}  \!\!\!\!  \longrightarrow \! \!(\flambda_{\mathbf m} ([0, s]))_{s\in \bbR_+}\!\! \! \overset{{\textrm{(law)}}}{=} \!\! (\int_0^s \! Z^{x_0}_r dr)_{s\in \bbR_+} $ in law in $\mathbf C (\bbR_+, \bbR)$. 
\end{theorem}
\noi
\textbf{Proof.} See Section \ref{pfsecthmIP}. \cqfd 
\begin{remark}
\label{IPcomments}
\noi
$\textbf{(a)}$ Even in compact cases (including those involving only metric convergence), this result is new in its generality. It includes more specific results (regarding offspring distributions and branching mechanisms)  
where trees are equipped with their length measure, i.e.~their one-dimensional Hausdorff measure, (in this case, $w_p(n) \eqo 0$, for all $n \geqo -1$ and $a_p\kappa_p\!  \to \! 1$), 
the suitably normalized empirical measure on the vertices (in this case, $\kappa_p \eqo 0$ and $w_p(n) \eqo 1/b_p$, for all 
$n \geqo -1$), 
or the suitably normalized 
empirical measure on the leaves (in this case, $\kappa_p \eqo 0$, $w_p(-1)\eqo 1/ (b_p \xi_p(0))$ and 
$w_p(n) \eqo 0$, for all $n \geqo 0$), which can be proved directly using Gromov--Hausdorff--Prokhorov or coding function methodology. 

A similar result is certainly also valid for discrete measured Galton--Watson trees, i.e.~those with edges of fixed length: we leave this question for future work. Such a result would extend the purely metric convergence of \cite[Thm 3.25]{DuWi2}, and 
in the compact case, the convergence of height processes by Duquesne and Le Gall \cite[Thm 2.3.1]{DuLG}, which corresponds to a $\dGHP$-convergence of the trees equipped with their normalized length measure (or empirical measure on the vertices).  

\smallskip

\noi
$\textbf{(b)}$ 
Assumption \textbf{A5} means that the measures equipping the GW forests are relatively evenly spread: for Type C branching mechanisms, this is also the case for the limiting standard Lévy forest whose measure is diffuse. For Type B branching mechanisms this is no longer the case, the limiting L\'evy forest being purely atomic. Reflecting the lack of regularity of Type B standard Lévy forests, the convergence holds only in the weakest sense, i.e.~in law in $(\bbT^{*}\! , \dera)$: GW forests equipped with deterministic measures  
cannot converge in law in $(\bbT, \dGP)$ to a standard Lévy forest. 
However, the construction provided in Remark \ref{remconstructLevy} $(d)$ exhibits an example where GW forests converge  in law in $(\bbT, \dGP)$ to any Type B standard L\'evy forest, but where their measures show  extra randomness with respect to the metric.  

\smallskip

\noi
$\textbf{(c)}$ In Remark \ref{A5minimal} in Section \ref{pfsecthmIP}, we explain that under 
$\mathbf{A1}$–$\mathbf{A4}$, if total masses 
$\langle \mathbf m_p \rangle $ converge in distribution to the total population $\int_0^\infty \! Z^{x_0}_s ds$ of the limiting CSBP, then $\mathbf{A5}$ is necessarily satisfied. Indeed, we think that the assumptions of the theorem are the weakest under which an invariance principle holds for measured trees whose measures are deterministic functions of the metric and which converge to standard measured $(x_0,\psi) $-Lévy forests, i.e.~whose total mass is distributed like the total population of a CSBP$(x_0, \psi)$. 

To extend Theorem \ref{thmIP} to GW forests whose measures show extra randomness with respect to their  metric, Assumption $\textbf{A5}$ could be modified into a joint convergence \vspace{-0.1cm} of $(\frac{1}{{a_p}}V^{{(p)}}_{{\lfloor b_p\rfloor}} , a_p\kappa_p+L^{{(p)}}_{{\lfloor b_p\rfloor}} )$ where 
$V^{{(p)}}_{{k}}\! \eqo  \sum_{1\leq l\leq k} Y^{{(p)}}_{l}$, $k\ino \mathbb{N}^*$. Here, the limit has to be $(X_1, L_1)$, where $(X_s)_{ s\in \mathbb{R}_+}$ is a spectrally positive Lévy process with Laplace exponent $\psi$, $(L_s)_{ s\in \mathbb{R}_+}$ is a subordinator correlated with the jumps of $X$, and the limit of the $\mathbf m_p$ may no longer be a standard Lévy forest.  \cq
\end{remark}

\subsection{Proof of Theorem \ref{thmIP}}
\label{pfsecthmIP}

\noi
\textbf{Proof of Theorem  \ref{thmIP} $(i)$.} To prove Theorem \ref{thmIP} $(i)$, we first compute the parameter transforms of $\era_h \mathbf m_p$ using Proposition \ref{prop:eraGW}, then we apply Proposition \ref{proppoubelle} to prove that these parameter transforms converge to those of $\era_h \mathbf m$, which are given in (\ref{chphih}), (\ref{xihLamh}) and (\ref{varrhoh}) in Lemma \ref{existlaws}.  

To this end, we first use the Lamperti transform to view the branching processes $(Z^{{(p)}}_s)_{s\in \bbR_+}$ 
as time-changed compound Poisson processes stopped at $0$. We recall that $(Y^{{(p)}}_k)_{k\in \bbN^*}$ are 
independent r.v.s whose laws are given by $\bP (Y^{{(p)}}_k\eqo j)\eqo \xi_p(j+1)$, for all $j\ino \{-1\} \cup \bbN$. 
Lamperti's coupling allows to view the initial part of the sequence $(Y^{{(p)}}_k)_{k\in\bbN^*}$ as the successive jumps of 
$s\ino \bbR_+ \! \mapsto \! Z^{{(p)}}_s$. More precisely, let 
$(\mathtt{N} (s))_{s\in \bbR_+}$ be a linear homogeneous Poisson process with unit intensity. We assume that $Z^{{(p)}}_0$, 
$(Y^{{(p)}}_k)_{k\in \bbN^*}$ and $(\mathtt{N}(s))_{s\in \bbR_+}$ are independent, and for all $p, k\ino \bbN$ and all $s\ino \bbR_+$, we set  
\begin{equation}
\label{VXRpdef} V^{(p)}_0\! = 0, \quad V^{(p)}_k \!\! = \! \!  \sum_{1\leq l\leq k} \!\! Y^{(p)}_l \!  , \quad X^{(p)}_s \! =   Z^{{(p)}}_0\! + V^{(p)}_{\mathtt{N}(c_p s)}\, , 
\end{equation}
\begin{equation}
\label{CApdeff}
R^{(p)}_0 \!\! =  \inf \big\{ s\ino \bbR_+: X^{(p)}_{a_ps}\!\! = 0 \, \big\}  \quad \textrm{and} \quad A^{(p)}_s\!\! =  \inf \! \Big\{ r\ino \bbR_+ \! :\!   \int_0^r \!\!  \tfrac{1}{a_p}Z^{(p)}_u du \geko s  \Big\}, 
\end{equation}
with the convention $\inf \emptyset \eqo \infty$. We observe that 
\begin{equation}
\label{Xppsip}
\textrm{\emph{$\big( \tfrac{1}{a_p} X^{(p)}_s\big)_{\! s\in \bbR_+}\! $ is a Lévy process with Laplace exponent $\psi_p$ as in Lemma \ref{Greydis}}.}
\end{equation}
As already mentioned, we can assume w.l.o.g.~that $(Z^{{(p)}}_s)_{s\in \bbR_+}$ and $(X^{{(p)}}_s)_{s\in \bbR_+}$ are coupled via the Lamperti transform. Namely 
\begin{equation}
\label{Lampertip}
\forall s\ino \bbR_+ , \quad \tfrac{1}{a_p} Z^{{(p)}}_{{\! \!A^{{(p)}}_s}}= \tfrac{1}{a_p}X^{{(p)}} \!\!\!\! \!\!_{{\!\! a_p (s\wedge R^{{(p)}}_0})} \quad \textrm{and} \quad \tfrac{1}{a_p}X^{{(p)}}\!\!\!\! \!\! \!_{a_p C^{(p)}_s}\eqo  
\tfrac{1}{a_p}Z^{{(p)}}_{ s}  
\end{equation}
where $C^{{(p)}}_s\eqo  R^{{(p)}}_0 \! \wedge 
\inf \{ r\ino \bbR_+ \! : \! \int_0^r du/ ( \frac{1}{{a_p}} X^{{(p)}}_{a_p u}) \geko s \}$. Here, we choose to express the Lamperti transforms in terms of the normalized processes $( \frac{1}{{a_p}} X^{{(p)}}_{a_p s})_{s\in \bbR_+}$ and 
$( \frac{1}{{a_p}} Z^{{(p)}}_{s})_{s\in \bbR_+}$. 

We note that $1+Y_{k}^{{(p)}}\! \eqo \bn(\sigma_k,\mathbf{F}_{\! p}) \! -\! 1$, $1\leqo k\leqo K_p$, where for all 
$\sigma\ino \mathbf F_{\! p}$, we recall that $\bn(\sigma,\mathbf{F}_{\! p})$ stands for the number of 
connected components of $\mathbf F_{\! p}\setminus\{\sigma\}$ and where $( \sigma_k)_{ 1\leq k\leq K_p}$ is the indexation 
of the subset $\mathtt{Lf}(\mathbf{F}_{\! p})\cup\mathtt{Br}(\mathbf{F}_{\! p})\backslash\{\rho_p\}$ for which 
$d_p (\rho_p, \sigma_k) \leko d_p (\rho_p, \sigma_{k+1})$, $1\leqo k\leqo K_p-1$, which exists a.s.

Next, we denote by $(X_s)_{s\in \bbR_+}$ a spectrally positive Lévy process with Laplace exponent $\psi$ and initial value $X_0\eqo x_0$ and we recall the following result. 
\begin{lemma}
\label{recallHelland} We keep the previous notation and we assume $\mathbf{A1}$--$\mathbf{A3}$. Then the sequence $(a_p)_{p\in\mathbb{N}}$ in $\mathbf{A4}$ equivalently satisfies the display in $\mathbf{A4}$ or $\mathbf{A4}'$ or $\mathbf{A4}''$ or $\mathbf{A4}'''$ where: 
\begin{compactenum}

\smallskip

\item[$\mathbf{A4}'$] $\;\big(\frac{1}{a_p}X_{{a_p s}}^{{(p)}}\big)_{\! s\in\mathbb{R}_+}\!\!   \xrightarrow[p\to \infty]{\textrm{(law)}} (X_s)_{s\in\mathbb{R}+}$ in $\mathbf{D}(\mathbb{R}_+,\mathbb{R})$.

\item[$\mathbf{A4}''\!$]  $\;\,\frac{1}{a_p} Z^{{(p)}}_{0} \!\!   \xrightarrow[p\to \infty]{\textrm{(law)}} x_0$ and $\frac{1}{a_p} V^{{(p)}}_{{\lfloor b_p\rfloor}}\!  \xrightarrow[p\to \infty]{\textrm{(law)}}  X_1-x_0$ in $\bbR$ (recall the notation $b_p\eqo a_pc_p$). \smallskip

\item[$\mathbf{A4}'''\!\!$]
\begin{compactenum}
\item[$\!\!\!(a)$] $ \varrho_p (dy /a_p) \! \to \! \delta_{x_0} (dy) $ in law in $\bbR_+$;\smallskip 
\item[$\!\!\!(b)$] there is a bounded continuous $f_0\colon  \bbR \! \to \! \bbR$ such that $f_0 (x)\eqo x$ around $0$ and as $p\rightarrow\infty$ 
\end{compactenum}

\vspace{-5mm}

$$\hspace{-8mm}b_p \sum_{k\in \bbN}\!  f_0 \big( \tfrac{k-1}{a_p} \big)\xi_p (k) \! \longrightarrow \! -\mathbf{a}+\!\int_{\bbR_+^*} \!\!\!\! ( f_0(z) \!  -\! z ) \pi (dz) \ \  \textrm{and} \ \ b_p  \sum_{k\in \bbN}\! f^2_0  \big( \tfrac{k-1}{a_p} \big) \xi_p (k)  \!   \longrightarrow \!    \mathbf b  + \! 
\int_{\bbR_+^*} \!\!\!\!  f^2_0(z)\pi  (dz) ;$$

\vspace{-2mm}

\begin{compactenum}
\item[$(c)$] $ b_p \sum_{k\in \bbN} g(\frac{k-1}{a_p} )  \xi_p (k) \!   \xrightarrow[p\to \infty]{\;} \! \int_{\bbR_+^*}  \! g(x) \pi (dx) $, for all bounded continuous $g \colon \bbR\! \rightarrow \! \bbR$ vanishing on a neighbourhood of $0$.\smallskip
\end{compactenum}
\end{compactenum}
Moreover, under $\mathbf{A1}$--$\mathbf{A3}$ and $\mathbf{A4} \! \Leftrightarrow \!\mathbf{A4}' \! \Leftrightarrow \!\mathbf{A4}'' \Leftrightarrow \!\mathbf{A4}''' $, the following joint convergence holds in law in $\mathbf D (\bbR_+, \bbR)^2\! \times \! \mathbf C (\bbR_+, \bbR)^2$
\begin{equation}
\label{jointcv} 
\Big( \big( \tfrac{1}{a_p} X^{(p)}_{a_ps} \big)_{\! s\in \bbR_+},  \big( \tfrac{1}{a_p} Z^{(p)}_{s} \big)_{\! s\in \bbR_+},  
\big( C^{(p)}_{s} \big)_{\! s\in \bbR_+},  \big( A^{(p)}_{s} \big)_{\! s\in \bbR_+} \Big) \xrightarrow[p\to \infty]{\; } 
(X, Z^{x_0}, C, A) 
\end{equation}
where for all $s\ino\bbR_+$, $A_s \eqo \inf \{ r\ino \bbR_+\! :\! \int_0^r\! Z^{x_0}_u du \geko s\}$, 
$C_s \eqo R_0\!  \wedge \! \inf \{ r\ino \bbR_+\! :\! \int_0^r\!  du/\! X_u \geko s\}$, $R_0 \eqo \inf \{ s\ino \bbR_+\! :\! X_s\eqo 0\}$, where $X$ and $Z^{x_0}$ are coupled via Lamperti transforms: 
$X_{C_s}\eqo Z^{x_0}_s$ and $Z^{x_0}_{\! A_s}\eqo X_{s\wedge R_0}$, with the conventions that $\inf \emptyset \eqo \infty $ and a.s.~$Z^{x_0}_\infty \eqo 0$; note that $C_s \eqo \int_0^s \! Z^{x_0}_r dr$. 
\end{lemma}
\noi
\textbf{Proof.} We assume $\mathbf{A1}$--$\mathbf{A3}$. See Helland  \cite[Thm 6.1]{He} for $\mathbf{A4} \! \Leftrightarrow \!\mathbf{A4}' \! \Leftrightarrow \!\mathbf{A4}'' $ and (\ref{jointcv}). See Jacod and Shiryaev \cite[Thm VII.2.35]{JaSh02} for $ \!\mathbf{A4}'\! \Leftrightarrow \!\mathbf{A4}''' $.    \cqfd 

\smallskip

Standard arguments then show that assumptions $\mathbf{A1}$--$\mathbf{A3}$, $\mathbf{A4}''$ and $\mathbf{A5}$  yield
\begin{equation}\label{eqn:IP1}
\textstyle\Big(\frac{1}{a_p}V_{\lfloor b_ps\rfloor}^{(p)},a_p\kappa_ps+L_{\lfloor b_p s\rfloor}^{(p)}\Big)_{\! s\in\mathbb{R}_+} \xrightarrow[p\to \infty]{\textrm{(law)}}  \big( X_s\! -\! x_0 \, ,\, s \big)_{\! s\in\mathbb{R}_+} \quad \textrm{in  $\mathbf{D}(\mathbb{R}_+,\mathbb{R}^2)$.}
\end{equation}
\emph{Indeed}, $\mathbf{A5}$ implies that 
$(a_p\kappa_ps+L_{{\lfloor b_ps\rfloor}}^{{(p)}})_{s\in\mathbb{R}_+} \!\! \rightarrow\! \mathtt{Id}_{\bbR_+}$
in law in $\mathbf{D}(\mathbb{R}_+,\mathbb{R})$, by standard arguments on scaling limits of random walks, where $\mathtt{Id}_{\bbR_+}$ stands for the identity function on $\bbR_+$; 
the convergence in law is a convergence in probability since the limiting process $\mathtt{Id}_{\bbR_+}$ is deterministic.  We then use Slutsky's theorem to show a joint convergence \eqref{eqn:IP1} in law in $\mathbf{D}(\mathbb{R}_+,\mathbb{R})^2$. Since $\mathtt{Id}_{\bbR_+}$ is continuous, standard results in Skorokhod topology (see e.g.~Jacod and Shiryaev \cite[Prop.~VI.2.2 (b)]{JaSh02}) imply 
\eqref{eqn:IP1} in $\mathbf{D}(\mathbb{R}_+,\mathbb{R}^2)$. \cq 

By Jacod and Shiryaev \cite[Thm VII.2.35]{JaSh02}, under $\mathbf{A1}$--$\mathbf{A4}$, assumption
$\mathbf{A5}$ is equivalent to 
\begin{compactenum} 
\smallskip

\item[$\mathbf{A5'}$] 
\begin{compactenum}
\item[$(a)$] 
There is a bounded continuous $f_0\colon  \bbR \! \to \! \bbR$ such that $f_0 (x)\eqo x$ around $0$ and 

\vspace{-2mm}

$$ \hspace{-0mm} a_p\kappa_p+b_p \sum_{k\in \bbN}  f_0 \big( w_p (k-1) \big)\xi_p (k)  \xrightarrow[p\to\infty] \! 1 \quad \textrm{\emph{and}} \quad  b_p  \sum_{k\in \bbN} f^2_0  \big( w_p(k-1) \big) \xi_p (k)  \!   \xrightarrow[p\to \infty]{\;} 0;$$
\item[$(b)$]  $b_p \sum_{k\in \bbN} g(w_p(k\!-\!1))  \xi_p (k) \! \to  \! 0$ for all bounded continuous $g \colon \bbR\! \rightarrow \! \bbR$ vanishing around~$0$.
\end{compactenum}
\end{compactenum}

\smallskip

We next derive analytical consequences of the previous assumptions in several lemmas. Let us recall the short notation $F_{\! p}$ from (\ref{shorthandp}) and the general notation from (\ref{masterfunc}). 
\begin{lemma}\label{lmconvF} We keep the previous notation and assumptions. 
Let $(z,y), (z_p, y_p) \ino\mathbb{R}_+^2$ be such that $(z_p,y_p)\! \rightarrow \! (z,y)$. Then the following holds true.
\begin{compactenum}
\smallskip
\item[$(i)$] $a_p\big(F_p \big( 1\! -\! \frac{1}{a_p} z_p\,  ,\, y_p \big)-\kappa_py_p\big)-b_p\big(1\! -\! \frac{1}{a_p} z_p\big) \! \xrightarrow[p\to\infty]\,  \psi(z)\! -\! y$.
\smallskip
\item[$(ii)$] If $z>0$, then $c_p\! -\! \partial_1F_p \big( 1\! -\! \frac{1}{a_p} z_p\,  ,\, y_p \big)\! \xrightarrow[p\to\infty]\, \psi^\prime(z)$.
\smallskip
\end{compactenum}
(This makes sense since $z_p \leqo a_p$ for all $p$ sufficiently large).
\end{lemma}
\noi{\bf Proof.} By \eqref{eqn:IP1}, $\big(\frac{z_p}{a_p}V_{{\lfloor b_p\rfloor}}^{{(p)}},y_p\big (a_p\kappa_p\!+\!L_{{\lfloor b_p\rfloor}}^{{(p)}}\big)\big)\!\! \overset{\textrm{(law)}}{\longrightarrow}(zX_1,y)$. By standard results on Laplace transforms (see Grimvall \cite{Gr} and easy arguments to handle the second coordinate) we get
$$ \bE \Big[\! \exp\!\Big(\!-\!\tfrac{z_p}{a_p}V_{\lfloor b_p\rfloor}^{(p)}\!-\!y_p\big(a_p\kappa_p\!+\!L_{\lfloor b_p\rfloor}^{(p)}\big)\!\Big)\Big]\xrightarrow[p\to \infty]{\; } e^{\psi(z) - y}$$
and thus
$$b_p \Big( \bE \Big[\!\exp\Big(\!-\!\tfrac{z_p}{a_p}V_{1}^{(p)}\!-\!y_pL_{1}^{(p)}\Big)\Big] \!-1\Big)\!-a_p\kappa_py_p \xrightarrow[p\to \infty]{\; } \psi(z) \!-\! y .$$
Now observe that
$
F_{\! p}(e^{-z_p/a_p},y_p) \eqo c_pe^{-z_p/a_p}\mathbf{E}\big[\exp \! \big(\! -\! \frac{z_p}{a_p}V_{1}^{{(p)}} \! -\! y_p L_{1}^{{(p)}}\big)\big].
$
Therefore, we get 
\begin{equation}\label{eqn:IP2}
b_p \big( \tfrac{1}{c_p}F_{\! p} \big(e^{-z_p/a_p},y_p \big)\! -\! e^{-z_p/a_p}\big)-a_p\kappa_py_p \xrightarrow[p\to \infty]{\;} \psi(z)\! -\! y.
\end{equation}
We apply \eqref{eqn:IP2} to $z_p^\prime \! =\! - a_p\log \big(1\! -\! \frac{1}{a_p} z_p \big)$ to get  $(i)$ since $z_p/a_p\rightarrow 0$. 

To prove $(ii)$, we set $G_p(z)\! :=\! a_p\big(F_p\big( 1-\frac{1}{a_p}z, y_p\big)\!-\!\kappa_py_p\big)\! -\! b_p\big( 1\! -\! \tfrac{1}{a_p} z\big)$, $z\ino [0, a_p]$, so that 
$$
 G_p(z_p)\eqo a_p\Big(\!F_p\Big(\!1-\tfrac{1}{a_p}z_p , y_p\! \Big)\!-\!\kappa_py_p\Big)\!-\! b_p\big(1\!-\!\tfrac{1}{a_p}z_p\big)\xrightarrow[p\to \infty]{\;}\psi(z)\!-\!y\!=:\!G_\infty(z).
$$
Then we note that 
$$
G_p^\prime(z)= c_p- \partial_1F_p \Big( 1-\tfrac{1}{a_p}z ,\, y_p  \Big)
                    = c_p - c_p\mathbf{E}\Big[ \big( Y_{1}^{{(p)}}\!\! +\!1 \big) 
                      \big(1\!-\!\tfrac{1}{a_p}z\big)^{(Y_{1}^{{(p)}})^+}\!\exp\big( \!-\!y_p w_p\big( Y_{1}^{{(p)}} \big) \big)\! \Big] .
$$
$G_p^\prime$ is therefore non-decreasing on $[0,a_p]$: if $z\geko0$, then for all $\varepsilon\in(0,z)$ and $p$ sufficiently large
\begin{equation}\label{eqn:IP3}
\tfrac{1}{\epp} \big( G_p(z)\! -\! G_p(z-\varepsilon)\big)\leq G_p^\prime(z)\leq \tfrac{1}{\epp} 
\big(G_p(z+\varepsilon)\! -\! G_p(z)\big).
\end{equation}
This easily  implies that $\lim_{p\rightarrow\infty} G_p^\prime(z)=G_\infty^\prime(z)=\psi^\prime(z)$. \it (Indeed\rm, \eqref{eqn:IP3} entails that
$$
\tfrac{1}{\epp} \big( G_\infty(z)\! -\! G_\infty(z-\varepsilon)\big)\leq\liminf_{{p\rightarrow\infty}}G_p^\prime(z)\leq\limsup_{{p\rightarrow\infty}}G_p^\prime(z)\leq \tfrac{1}{\epp} 
\big(G_\infty(z+\varepsilon)\! -\! G_\infty(z)\big), $$
where both sides tend to $G_\infty^\prime(z)$ as $\varepsilon\rightarrow 0$.) 
Since $G_p^\prime$ is continuous non-decreasing and $G_\infty^\prime$ is continuous, we get for all $z_p\! \rightarrow\!  z\ino\bbR^*_+$ that $\lim_{p\rightarrow\infty}G_p^\prime(z_p)\eqo\psi^\prime(z)$, as required for $(ii)$, which completes the proof of the lemma. \cqfd

\medskip

We next recall that $\nu(dz)$ is the L\'evy measure of the subordinator with Laplace exponent $\psi^{-1}$. As $X$ is of Type B or C, $\nu$ is diffuse by Lemma \ref{psi-1prop} $(ii)$ and for all $h\ino\bbR^*_+$, we recall from (\ref{psihnuhdef}) the notation $\nu(h)$, $\psi_{h}^{{-1}}$ and $\, \overline{\! \psi}^{{-1}}_{h}$. We also recall from (\ref{shorthandp}) the short hand $Q_p$. 
\begin{lemma} \label{lmIP2} We keep the previous notation and assumptions. Then for all $y\in\mathbb{R}_+$ 
and $h\ino\bbR^*_+$, we get $\lim_{p\to \infty} a_p \big( 1\! -\! Q_p \big[ e^{-y \langle\flambda\rangle } 
\un_{\{ \langle\flambda\rangle\leq h  \}}\big]\big)\eqo \psi_{h}^{{-1}}(y)$. 
\end{lemma}
\noi 
{\bf Proof.} Since $\mathbf m_p$ has no mass at its root, we get $\mathbf{E}[e^{-y\langle\flambda_{\mathbf{m}_p}\rangle}|Z_{0}^{{(p)}}]\eqo Q_p [e^{-y\langle\flambda\rangle}]^{Z_{0}^{{(p)}}}\! $, by construction of GW forests. 
Therefore, the total mass $\langle \mathbf{m}_p\rangle\eqo \langle\flambda_{\mathbf{m}_p}\rangle$ has the same law as $\Sigma_p (Z_{0}^{{(p)}} )$, 
where $Z_{0}^{{(p)}}$ and $(\Sigma_p (k))_{k\in\mathbb{N}}$ are 
independent and where $(\Sigma_p (k))_{k\in\mathbb{N}}$ is a random walk starting at $0$ and 
whose jump law is the law of $\langle\flambda\rangle$ under $Q_p$. 
On the other hand, Lemma \ref{ballGWmass} $(i)$ yields 
\begin{equation}
\label{totmassrepres}
\langle\flambda_{\mathbf{m}_p}\rangle \eqo a_p\kappa_p R_{0}^{{(p)}}+L^{{(p)}}_{
{\mathtt{N} (b_pR_0^{(p)})} } \; ,
\end{equation}
where we recall definition of $R^{{(p)}}_0$ from (\ref{CApdeff}).  
Since $\mathbf{A4}'$ holds true by Lemma \ref{recallHelland}, standard arguments (see e.g.~Jacod and Shiryaev \cite[Chapter VI, Section 2, Prop.~2.11]{JaSh02}) imply 
\begin{equation}
\label{jointcvR0}
R_{0}^{{(p)}}\! \xrightarrow[p\to \infty]{\textrm{(law)}} \! R_0\! :=\! \inf\{s\ino\mathbb{R}_+\! :\! X_s\eqo 0\}, \; \textrm{jointly with (\ref{jointcv}).}
\end{equation}
We also see, by the law of large numbers and a Dini argument that $(\frac{1}{{b_p}} \mathtt{N} (b_ps))_{s\in \bbR_+}\! \to \! \mathtt{Id}_{\bbR_+}$ in probability 
(and even $\bP$-a.s.) in $\mathbf D(\bbR_+, \bbR)$, where $\mathtt{Id}_{\bbR_+}$ stands for the identity function on $\bbR_+$. 
This implies that $\mathtt{N} (b_p R_{0}^{{(p)}})/b_p \! \to \!  R_0$, in law. 
We then deduce from \eqref{eqn:IP1} that 
$(a_p\kappa_ps+L^{{(p)}}_{{\lfloor b_p s\rfloor}})_{s\in \bbR_+}$ actually converges to $\mathtt{Id}_{\bbR_+}$ in probabilty  in $\mathbf D(\bbR_+, \bbR)$. 
Before completing the proof, we recall the following result from Whitt \cite{Wh80}.
\begin{lemma}
\label{composko} Let $f_n\! \rightarrow \! f$ and $\ell_n \! \rightarrow \! \ell $ in $\mathbf D (\bbR_+, \bbR)$. We assume that 
$\ell_n (0) \! = \! 0$, that $\ell$ is strictly increasing and continuous. Then $f_n \! \circ\!  \ell_n \! \to \! f\!\circ \! \ell$ in $\mathbf D (\bbR_+, \bbR)$.
\end{lemma}
\noi
\textbf{Proof:} See Whitt \cite[Thm.~3.1]{Wh80}. \cqfd 

\smallskip
\noi
To continue the proof of Lemma \ref{lmIP2}, we deduce from Lemma \ref{composko} that 
 \begin{equation}
\label{LandPoisson}
\big( a_p\kappa_ps+L^{\!(p)}_{\mathtt{N} (b_ps)} \big)_{\! s\in \bbR_+}\!  \xrightarrow[p\to \infty]{\; } \mathtt{Id}_{\bbR_+}
\end{equation}
where  $\mathtt{Id}_{\bbR_+}$ stands for the identity function on $\bbR_+$ (of course here, the monotonicity and the continuity of the processes allow for a direct proof; we however find convenient to state here Lemma \ref{composko}, which is used crucially in the middle the proof of Theorem \ref{thmIP} $(ii)$).     
We next deduce from \eqref{eqn:IP1} that 
$(a_p\kappa_ps+L^{{(p)}}_{{\lfloor b_p s\rfloor}})_{s\in \bbR_+}$ actually converges to $\mathtt{Id}_{\bbR_+}$ in probability for the norm of 
the uniform convergence on every compact intervals of $\bbR_+$: combined with the previous arguments and 
$\mathbf{A5}$, this yields $\langle\flambda_{\mathrm{m}_p}\rangle \! \to \!  R_0$ in law.  
We now recall that $R_0$ has the same law as a subordinator $(\Sigma_s)_{s\in\mathbb{R}_+}$ with Laplace exponent $\psi^{-1}$ at time $x_0$. Standard arguments on scaling limits of random walks imply that 
\begin{equation}
\label{Skocvsubo}
\big( \Sigma_{p} (\lfloor a_p x\rfloor) \big)_{\! x\in\mathbb{R}_+}\xrightarrow[p\to \infty]{\textrm{(law)}} \! \big(\Sigma_x\big)_{\! x\in\mathbb{R}_+} \quad \textrm{on $\mathbf D (\bbR_+, \bbR)$.}
\end{equation} 
This first entails $Q_p[e^{-y\langle \flambda \rangle} ]^{\lfloor a_p\rfloor} \eqo \bE [e^{-y\Sigma_p(1)}]^{\lfloor a_p\rfloor}\eqo \bE [e^{-y\Sigma_p (\lfloor a_p\rfloor)}]\! \to \! \bE [e^{-y\Sigma_1}]\eqo \exp(-\psi^{-1} (y))$, which implies for all $y\ino \bbR_+$, that 
\begin{equation}
\label{uncondpsi-1}
a_p \big( 1\! -\! Q_p \big[ e^{-y\langle \flambda \rangle} \big]\big)\xrightarrow[p\to\infty]{\; } \psi^{-1} (y).
\end{equation}
For all $x\ino \bbR_+$ we set $\Delta \Sigma_x \! :=\! \Sigma_x \! -\! \Sigma_{x-}$ and 
$\mathbf x (h)\! :=\!  \inf \{ x\ino \bbR_+\!  :\!  \Delta \Sigma_x\geko h\}$. Since 
$\{ (x,\Delta \Sigma_x)\,  ; \, x\ino \bbR_+\! :\! \Delta \Sigma_x \geko 0 \}$ is a Poisson point process with intensity 
$dx\nu(dz)$, $\mathbf x (h)$ is an exponentially distributed r.v.~with parameter $\nu(h)$ and 
$(\Delta \Sigma)(h)\! :=\! \Delta \Sigma_{\mathbf x(h)}$ has law $\nu (\, \cdot \, \cap (h, \infty))/ \nu(h)$. 
We also set $\mathbf x_p(h)\! :=\! \frac{1}{a_p}\inf \{k\ino \bbN^*\!  :\!  \Sigma_p (k) \! -\! \Sigma_p(k\! -\! 1) \geko h\}$ and 
$(\Delta \Sigma_p) (h) \! :=\!  \Sigma_p (a_p \mathbf x_p(h)) \! -\! \Sigma_p(a_p \mathbf x_p(h)\! -\!1)$. 
We easily check that $a_p \mathbf x_p(h)$ is a geometric r.v.~with parameter $Q_p (\langle \flambda \rangle \geko h)$ 
and that 
$\Delta \Sigma_p (h)$ has law $Q_p (\langle \flambda \rangle  \ino \cdot | \langle \flambda \rangle \geko h)$.   
Since $\nu$ is diffuse by Lemma \ref{psi-1prop} $(ii)$, we get $\bP$-a.s.~$h\! \notin \{ \Delta \Sigma_x ; x\ino \bbR_+\}$ and 
by standard arguments on the Skorokhod topology 
(see e.g.~Jacod and Shiryaev \cite[Chapter VI, Section 2, Prop.~2.7]{JaSh02}), (\ref{Skocvsubo}) implies that 
$(\mathbf x_p(h), (\Delta\Sigma_p )(h)) \! \to \! (\mathbf x (h), (\Delta \Sigma) (h)) $ in law on $\bbR_+^2$. In particular, 
we get 
$$ a_p Q_p \big(\langle \flambda \rangle \geko h \big)\xrightarrow[p\to\infty]{\; } \nu(h) \quad \textrm{and} \quad 
Q_p \big[ e^{-y\langle \flambda \rangle }\, \big| \, \langle \flambda \rangle \geko h\big] \xrightarrow[p\to\infty]{\; } \, 
\tfrac{1}{\nu(h)}\! \int_{(h,\infty)} \!\!\!  e^{-yz} \nu(dz) , $$
which implies the desired result by (\ref{uncondpsi-1}). \cqfd 
  
\begin{remark}
\label{A5minimal} Let us assume only $\mathbf{A1}$--$\mathbf{A4}$. A more general assumption than $\mathbf{A5}$ would be the convergence $(\frac{1}{{a_p}}V_{{\lfloor b_ps\rfloor}}^{{(p)}},a_p\kappa_ps+L_{{\lfloor b_p s\rfloor}}^{{(p)}} )_{\! s\in\mathbb{R}_+} \! \! \longrightarrow \!   \big( X_s\! -\! x_0 \, , L_s \big)_{\! s\in\mathbb{R}_+} \!$ in law in $\mathbf D (\bbR_+, \bbR^2)$ where $L$ is a subordinator whose jumps are a deterministic function of those of $X$: namely, for all $s\ino \bbR_+$, 
$L_s\eqo \alpha s+ \sum_{r\in [0, s]} f(\Delta X_r)\un_{\{\Delta X_r >0\}}$, where $\Delta X_r\eqo X_r \! -\! X_{r-}$ and where $f\colon \bbR_+^*\! \to \! \bbR^*_+$ is such that $\int_{\bbR_+^*} \! \min (1, f(z)) \pi(dz) \leko \infty $.  
Let us assume that, jointly with the preceding convergences, the total masses $\langle \flambda_{\mathbf m_p} \rangle$ converge to the total population of the CSBP($x_0, \psi$), i.e.~to $R_0$. Then (\ref{totmassrepres}) and standard results entail that $L_{R_0}\eqo R_0$ almost surely, which in turns implies that $L_{s}\eqo s$, $s\ino \bbR_+$, by an explicit 
computation of the Laplace exponent $k(y_1,y_2) \! :=\!\log \bE [e^{-y_1 (X_1-x_0)-y_2 L_1} ] \eqo \psi(y_1) -\alpha y_2 -\int_{\bbR^*_+} e^{-y_1z} (1\! -\! e^{-y_2 f(z)})\pi(dz) $ and stopping arguments at time $R_0$ for the positive martingale $s\! \mapsto e^{-y_1 (X_s-x_0)-y_2L_s-sk(y_1,y_2)}$. \cq  
\end{remark}
 
We now compute the parameter transforms of 
$\era_h \mathbf m_p$ using Proposition \ref{prop:eraGW}, which asserts that  $\era_h\mathbf m_p$ is distributed as a GW-forest associated with the parameters $\xi_{p,h}$, $c_{p,h}$, $\phi_{p,h}$, $(\Lambda^{\! p,h}_{n})_{n\in \bbN}$, $\varrho_{p,h}$, $(\Lambda_{\varnothing, n}^{\! p,h})_{n\in \bbN}$ that are characterised as follows. First, to simplify, we set 
$$ F_{\! p, h}\! :=\! F_{\! \xi_{p,h}, c_{p,h}, \Lambda^{\! p,h}_{\cdot}} \quad \textrm{and} \quad F^{\varnothing}_{\! p,h} \!:=\! F^\varnothing_{\!\! \varrho_{p,h}, \Lambda_{\varnothing ,\cdot}^{\! p,h}}. $$
We also define 
the function $\falpha_{p,h} (\cdot)$, the real number $\alpha_{p,h}$ and the $\bbR_+$-valued r.v.~$M_{p,h}$ by 
$$ \falpha_{p,h}(y)\! := \!  Q_p[ e^{-y \langle\flambda\rangle}\un_{\{\langle\flambda\rangle\leq h\}}], \quad  \alpha_{p,h}\! :=\! \falpha_{p,h}(0) \quad   \textrm{and} \quad  \bE \big[e^{-yM_{p,h} }\big]\! :=\! \frac{F^\varnothing_{\! p} ( \falpha_{p,h}(y) , y )}{\varphi_{\! \varrho_p} (\alpha_{p,h}) } .$$
Then, for all $s\ino [0, 1]$ and for all $y\ino \bbR_+$, 
\begin{eqnarray}
\label{charaQera2IP}
c_{p,h}+ \phi_{p,h}(y) \!\!\!\! &\eqo &\textstyle  \!\!\!\! c_p+ \phi_p(y) -\partial_1 F_{\! p} ( \falpha_{p,h}(y),y)   \\
(1\! -\! \alpha_{p,h}) e^{-hy}   F_{\! p,h} (s,y) \!\!\!\! & \eqo&  \!\!\!\!  \textstyle 
F_{\! p}  \big( (1\! -\! \alpha_{p,h}) e^{-hy}s +  \falpha_{p,h}(y)\, , \, y \big)
 - \falpha_{p,h}(y)(c_p+ \phi_p(y)) \nonumber \\
 & & \textstyle \quad  - (1\! -\! \alpha_{p,h}) e^{-hy}s \,  \partial_1 F_{\! p} \big( \falpha_{p,h}(y)\, ,\, y \big) \label{charaQera3IP} \\
 \text{and}\quad e^{-hy}   F^\varnothing_{\! p,h} (s,y) \!\!\!\! & \eqo&  \!\!\!\!  \textstyle F^\varnothing_{\! p}  \big( (1\! -\! \alpha_{p,h}) e^{-hy}s +  \falpha_{p,h}(y)\, , \,y \big)- F^\varnothing_{\! p}  \big(   \falpha_h(y)\, ,\,   y \big)\nonumber
  \\
 & & \quad \textstyle  +e^{-hy} \varphi_{\! \varrho_p} (\alpha_{p,h}) \bE \big[ e^{-y(M_{p,h}-h)_+}\big]  .\label{charaPeraIP}
\end{eqnarray}  

We set $z_p\! :=\!  z''_p -z_p'se^{-hy}$ where $z_p'\!\!  :=\! \eqo a_p(1\! -\! \alpha_{p,h})$ and 
$z_p''\!\! :=\! a_p(1\! -\! \falpha_{p,h}(y))$. By Lemma \ref{lmIP2}, 
we see that for all $(s,y)\in[0,1]\times\mathbb{R}_+$ $z_p'\! \to \! \nu(h)$, $z_p''\! \to \! \psi_{h}^{{-1}}(y)$ and thus 
$z_p\! \to \! \psi_{h}^{{-1}}(y)\! -\! \nu(h) se^{-hy} $. Then (\ref{charaQera2IP}) can be rewritten using $z_p''$: the previous limits and 
Lemma \ref{lmconvF} then yield 
\begin{equation}
\label{primelimit}
\textstyle  c_{p,h}+\phi_{p,h}(y)=\kappa_p y+c_p-\partial_1 F_{\! p} \big(1\! -\! \frac{1}{a_p} z_p'',\,  y\big)\; \xrightarrow[p\to \infty]{\; }\; \psi^\prime(\psi_h^{-1}(y)), 
\end{equation}
since $\kappa_p \! \to \! 0$ by $\mathbf{A4}$ and $\mathbf{A5}$. 
By multiplying it by $a_p$, we rewrite \eqref{charaQera3IP} in terms of $z_p$, $z_p'$ and $z_p''$: 
\begin{align*}
   z_p^\prime e^{-hy}F_{\! p,h}(s,y)=& \textstyle a_p\big(F_{\! p} \big( 1\! -\! \frac{1}{a_p}z_p \, ,\, y\big)-\kappa_py\big)-b_p\big(1\! -\! \frac{1}{a_p}z_p\big) \\
	&\textstyle\quad+z_p'se^{-hy}(c_p\! -\! \partial_1 F_{\! p}
	 \big(1\! -\! \frac{1}{a_p} z_p^{\prime\prime},\, \, y\big) \big) 
	 + \kappa_p z_p'' y .
\end{align*}
By the previous limits and Lemma \ref{lmconvF}, and as $\kappa_pz_p''\rightarrow 0$, we thus obtain for all $(s,y)\ino[0,1]\times\mathbb{R}_+$, 
\begin{equation}
\nu(h)e^{-hy}\lim_{p\rightarrow\infty}F_{\! p,h}(s,y)=
								  \, \psi \big( \psi_h^{-1}(y)\! -\! \nu(h)se^{-hy} \big)-y+\nu(h)se^{-hy}\psi^\prime \big( \psi_h^{-1}(y) \big).\label{limit2}
\end{equation}

  Since $Z_0^{{(p)}}$ has law $\varrho_p$, $\mathbf{A4}$ also entails the convergence in law 
$\varrho_p (dy/a_p)\!\!  \to \! \delta_{x_0}(dy)$ on $\bbR_+$. 
Thus for any nonnegative sequence $(\texttt z_p)_{p\in \bbN}$ such that $\texttt z_p \! \rightarrow\!  z\ino\mathbb{R}_+$, we get 
$\varphi_{\! \varrho_p}(e^{-\texttt z_p/a_p})\! \rightarrow \! e^{-zx_0}$. We apply this to $\texttt z_p\! :=\! -a_p\log(\falpha_{p,h}(y))\rightarrow\psi_h^{-1}(y)$ to get 
$\varphi_{\! \varrho_p}(\falpha_{p,h}(y)) \! \rightarrow \! e^{-x_0\psi_h^{-1}(y)}$, and we recall that $\psi^{{-1}}_{h}(0)\eqo \nu(h)$ and $\alpha_{p,h}\eqo \falpha_{p,h}(0)$. 

We next observe that $F^\varnothing_{\! p} (s,y)\eqo \varphi_{\! \varrho_p}(s)$, for all $s\ino[0,1]$ and all $y\ino\mathbb{R}_+$,  because $\mathbf m_p$ has no mass at its root. For any nonnegative sequence $(\mathtt z_p)_{p\in \bbN}$ such that $\mathtt z_p \! \rightarrow\!  z\ino\mathbb{R}_+$, we then observe that $F_{\! p}^\varnothing(e^{-\mathtt z_p/a_p},y)\! \rightarrow \! e^{-zx_0}$, $(z,y)\ino\mathbb{R}^2_+$. Therefore,\vspace{-0.1cm}
$$
 \bE \big[e^{-yM_{p,h} }\big] \xrightarrow[p\to \infty]{\; }
e^{-x_0\psi_h^{-1}(y)+x_0\nu(h)}=e^{-x_0 \, \overline{\! \psi}^{-1}_h(y)},
$$
where we recall that $\, \overline{\! \psi}_{h}^{{-1}}(y)\eqo \int_{(0,h)}(1 - e^{-zy})\nu(dz) $, is the Laplace exponent of the subordinator $(S^{h}_s)_{s\in \bbR_+}$. 
Then we pass to the limit in \eqref{charaPeraIP} and for all $(s,y)\in[0,1]\times\mathbb{R}_+$ we get 
\begin{equation}
\label{limit3}
e^{-hy} \!\!  \lim_{p\rightarrow\infty}\! F_{\! p,h}^\varnothing(s,y)\! =\! e^{-x_0(\psi_h^{-1}(y)-\nu(h)se^{-hy})}\! - e^{-x_0\psi_h^{-1}(y)}\! +e^{-hy}e^{-x_0\nu(h)}\mathbf{E} \big[ e^{-y(S^h_{x_0}-h)_+} \big].
\end{equation}

Let $\mathbf m$ be a standard measured ($x_0, \psi$)-Lévy forest as introduced in 
Theorem \ref{thm:growthprocconvbis}. By definition $\era_h \mathbf m$ is a 
GW forest whose parameter transforms are given in (\ref{chphih}), (\ref{xihLamh}) and (\ref{varrhoh}). 
We observe that (\ref{primelimit}) converges to (\ref{chphih}), (\ref{limit2}) to (\ref{xihLamh}) and (\ref{limit3}) to (\ref{varrhoh}). Therefore, the parameters $\xi_{p,h}$, $c_{p,h}$, $\phi_{p,h}$, $(\Lambda^{{\! p,h}}_{n})_{n\in \bbN}$, $\varrho_{p,h}$, $(\Lambda_{\varnothing, n}^{{\! p,h}})_{n\in \bbN}$ satisfy Assumptions $(a)$, $(b)$ and $(c)$ of Proposition \ref{proppoubelle}, which yields that $\era_h \mathbf m_p \! \to \! \era_h \mathbf m$ in law in $(\bbT^{1}_{\! c}, \dGHP)$ and thus also in 
$(\bbT, \dGP)$. Since this holds true for all $h\ino \bbR_+^*$, this shows that $\mathbf m_p \! \to \! \mathbf m$ in law in $(\bbT^*\!, \dera)$, which is Theorem \ref{thmIP} $(i)$. \cqfd  

\medskip

\noi
\textbf{Proof of Theorem \ref{thmIP} $(ii)$.} We recall that we assume that the branching mechanism is either of Type B or C. We recall from (\ref{LKformpsi-1}) that $\mathbf c$ is the drift coefficient of $\psi^{-1}$. By Lemma \ref{psi-1prop}, 
$\mathbf c\geko 0$ iff $\psi$ is of Type B and in this case, $1/\mathbf c \eqo  \mathbf{a}+\int_{\bbR^*_+}\! z\pi(dz) $. 
We deduce from Lemma \ref{ballGWmass} that
\begin{equation}
\label{represdiscr} \Big( \tfrac{1}{a_p}Z_s^{(p)}, \flambda_{\mathbf{m}_p}([0,s]) \Big)_{\! \!s\in\mathbb{R}_+} \!\!\!= 
\Big( \tfrac{1}{a_p}X^{(p)}\!\!\!\!\!\!\!   \!_{a_p C_s^{(p)}} \, ,\, a_p\kappa_p\,  C_s^{(p)}\!\! +L^{\!(p)} \!\!\! \!\!\!_{\mathtt{N} (b_p C_s^{(p)}) } \Big)_{\! \! s\in\mathbb{R}_+},
\end{equation}
where we recall $C^{(p)}_\cdot$ from (\ref{Lampertip}). 
We also recall from (\ref{jointcv}) in Lemma \ref{recallHelland} and from (\ref{jointcvR0}) that 
$( (\tfrac{1}{{a_p}}X^{{(p)}}_{a_ps} )_{s\in \bbR+} , (C_s^{{(p)}} )_{s\in \bbR_+}, R^{{(p)}}_0)$ converges in law jointly in  
$\mathbf D (\bbR_+, \bbR) \! \times \! \mathbf C (\bbR_+, \bbR)\! \times \! \bbR_+$ to 
$(X_\cdot  , C_\cdot, R_0)$. Here, $Z^{x_0}$, $X$ are coupled by the Lamperti transform 
$Z^{x_0}_s\eqo X_{C_s}$ and where $C_s\eqo \int_0^s \! Z^{x_0}_r dr$, for all $s\ino \bbR_+$, as recalled 
in Lemma \ref{recallHelland}. 

To take the limit in (\ref{represdiscr}), we next want to use the Skorokhod continuity of composition of functions as recalled in 
Lemma \ref{composko} and we need to extend $ C_\cdot^{{(p)}}$ and $C_\cdot$ in order to have increasing processes on $\bbR_+$. 
To this end, we observe that $C^{{(p)}} \! (R^{{(p)}}_0)\eqo R^{{(p)}}_{0}$ and 
$C_{\! R_0}\eqo R_{0}$. For all $s\ino \bbR_+$, we set $\overline{C}_s^{{(p)}}\eqo C_s^{{(p)}}$ if $s\leko R^{{(p)}}_0$ and $\overline{C}_s^{{(p)}}\eqo s$ if $s\geqo R^{{(p)}}_0$ and $\overline{C}_s\eqo C_s$ if $s\leko R_0$ and 
$\overline{C}_s\eqo s$ if $s\geqo R_0$. It is easy to derive from the previous joint convergence that 
$( (\tfrac{1}{{a_p}}X^{{(p)}}_{a_ps} )_{s\in \bbR+} , (\overline{C}_s^{{(p)}} )_{s\in \bbR_+}, R^{{(p)}}_0)$ converges in law jointly in  
$\mathbf D (\bbR_+, \bbR) \! \times \! \mathbf C (\bbR_+, \bbR)\! \times \! \bbR_+$ to 
 $(X_\cdot  , \overline{C}_\cdot, R_0)$.
By the convergence in probability (\ref{LandPoisson}) and Slutsky's theorem, we thus get the joint convergence 
$$ \Big( \big( \tfrac{1}{a_p} X^{(p)}_{a_ps} \big)_{\! s\in \bbR_+},   
\big( \overline{C}^{(p)}_{s} \big)_{\! s\in \bbR_+},  \big(a_p\kappa_ps+ L^{\!(p)}_{\mathtt{N} (b_ps)} \big)_{s\in \bbR_+}, 
R^{(p)}_0 \Big) \xrightarrow[p\to \infty]{\; } 
(X_\cdot, \overline{C}_\cdot , \mathtt{Id}_{\bbR_+}, R_0), $$ 
in law in $\mathbf D(\bbR_+, \bbR)^3\! \times \! \bbR_+$. 
This convergence and Lemma \ref{composko}, which applies since $\overline{C}_\cdot$ is strictly increasing and continuous, then entail the following convergence in law in $\mathbf D(\bbR_+, \bbR)^2\! \times \! \bbR_+$:
\begin{equation}
\label{extentedcv}
\Big( \big(\tfrac{1}{a_p}X^{(p)}\!\!\!\!\!\!\!   \!_{a_p \overline{C}_s^{(p)}} \big)_{\! s\in \bbR_+} \, ,\, \big( a_p\kappa_p\,  \overline{C}_s^{(p)}\!\! +L^{\!(p)} \!\!\! \!\!\!_{\mathtt{N} (b_p \overline{C}_s^{(p)}) } \big)_{\! s \in \bbR_+} , R^{(p)}_0 \Big)  \xrightarrow[p\to \infty]{\; } \big( X_{\overline{C}_\cdot}, \overline{C}_\cdot, R_0 \big).
\end{equation}
We next use the following elementary fact on the Skorokhod topology: if $f_p \! \to \! f$ in $\mathbf D (\bbR_+, \bbR)$, $r_p \! \to r$ in $\bbR_+$ and $f$ is continuous at $r$, then the stopped paths $f_p (\cdot \wedge r_p)$ converge to $f(\cdot \wedge r)$ in $\mathbf D (\bbR_+, \bbR)$. Then (\ref{extentedcv}) combined with  (\ref{represdiscr}) shows that 
$$ \Big( \big(\tfrac{1}{a_p}Z_s^{(p)} \big)_{\! s\in \bbR_+} , \big( \flambda_{\mathbf{m}_p}([0,s])\big)_{\! s\in \bbR_+} \Big)\! \xrightarrow[p\to \infty]{\; }\!  \big( X_{C_\cdot}, C_\cdot \big)\! = \! \Big( Z^{x_0}_\cdot\, , \Big(\!  \int_0^s \!\! Z^{x_0}_r dr\!  \Big)_{\! \! s\in \bbR_+} \Big). $$
We now recall from Theorem \ref{thm:growthprocconvbis} that $(\flambda_{\mathbf{m}}([0,s])_{s\in\mathbb{R}_+}\eqo(\mathbf c Z^{x_0}_s+\int_0^s\! Z^{x_0}_rdr)_{s\in\mathbb{R}_+}$. Therefore, if $\psi$ is of Type C, then $\mathbf c\eqo 0$ and the previous convergence implies that $\flambda_{\mathbf m_p}\! \to \! \flambda_{\mathbf m} $ in law in $\cM_f(\bbR_+)$. Since $\mathbf m_p \! \to \! \mathbf m$ in law in $(\bbT^*\! , \dera)$, Theorem \ref{condilawGP} applies and shows that $\mathbf m_p\! \to \! \mathbf m$ in law in $(\bbT, \dGP)$. 

If $\psi$ is of Type $B$, then $\mathbf c\geko 0$, 
$\flambda_{\mathbf m_p}\! \not \to \! \flambda_{\mathbf m} $ and Theorem \ref{condilawGP} asserts that $(\mathbf m_p)_{p\in \bbN}$ does not converge in law to $\mathbf m$ in $(\bbT, \dGP)$. Moreover, 
let us prove that $(\mathbf m_p)_{p\in \bbN}$ cannot converge in law in $(\bbT, \dGP)$ by arguing by contradiction:  
if $\mathbf m_p\! \to \! \mathbf m'$ in law in $(\bbT, \dGP)$, for some $\bbT$-valued r.v.~$\mathbf m'$, then the $\dGP$-continuity of mass erasure proved in Theorem \ref{eraGcont} $(i)$ implies that $\era_h \mathbf m_p \! \to \! \era_h \mathbf m'$ in law in $(\bbT, \dGP)$, for all $h\ino \bbR_+^*$, which entails 
$\mathbf m_p \! \to \! \mathbf m'$ in law in $(\bbT^*\! , \dera)$ by Proposition \ref{traducvlawera} $(i)$. Thus $\mathbf m'$ would have the same law as $\mathbf m$ by Theorem \ref{thmIP} $(i)$, which yields a contradiction. This completes the proof of Theorem \ref{thmIP} $(ii)$. \cqfd

\medskip

\noi
\textbf{Proof of Theorem \ref{thmIP} $(iii)$}. Since total height is $\dGH$-continuous, the tightness of the laws of 
$(\Phi_0 (\mathbf m_p))_{p\in \bbN}$ implies the tightness of the laws of $(\mathtt{Ht}_{ p})_{p\in \bbN}$ on $\bbR_+$,  
and Lemma \ref{Greydis} entails (\ref{Greydiscr}).

Conversely, let us assume (\ref{Greydiscr}). It implies the Grey condition (\ref{Grey}) by Lemma \ref{Greydis}. 
By Theorem \ref{thm:growthprocconvbis} $(iii)$, $\bP (\mathbf m\ino \bbT_{\! c}) \eqo 1$ and $\Phi_0 (\mathbf m)$ is the unmeasured compact ($x_0, \psi$)-Lévy forest defined in Theorem \ref{metriccpctLevytrees}.   

Next, from the definition of the total height of $\mathbf F_{\! p}$ and from 
(\ref{heightGWexpli}) in Remark \ref{heightGW} we get for all 
$b\ino \bbR_+^*$ and all $p\ino \bbN$,  
$\bP (\mathtt{Ht}_{ p} \leqo b)$  $\eqo$ $\bE [ \exp( -Z^{{(p)}}_{0} \!  \log (\, 1\! -\! \frac{1}{{a_p}}v_p (b))\,  ) ]$, $v_p (b)$ $\! :=\!$ $ a_p Q_{\! p} (\mathtt{Ht} \geko b)$. 
By Lemma \ref{Greydis}, under (\ref{Greydiscr}), the Grey condition (\ref{Grey}) holds and $\bP(\mathtt{Ht}_{ p} \leqo b)\!\to \! \bP (\mathtt{Ext} \leqo b)\eqo e^{-x_0 v(b)}$, where $v(b)$ is given by (\ref{defandmeaningv}), which 
implies that $v_p (b)\!  \to \! v(b)$, for all $b\ino \bbR_+^*$. 

  We recall $b$-length erasure $R_b \colon \bbT^0_{\! c} \! \to \! \bbT^0_{\! c}$ from (\ref{deflengthera}). 
By Proposition \ref{GWlengthera}, $R_b (\mathbf F_{\! p})$ is an unmeasured 
GW forest whose parameters are denoted by $c_{p,b}$, $\xi_{p,b}$, $\varrho_{p,b}$ and characterized by 
\begin{equation}
\label{lengtherapara}
\!\!c_{p, b} \eqo \psi'_p (v_p(b)), \; \varphi_{\xi_{p,b}}\!(r) \eqo r\!+\! \frac{\psi_p \big(v_p(b) (1\! -\! r) \big)}{v_p(b) 
\psi_p'(v_p(b))} \; \textrm{and} \; \varphi_{\! \varrho_{p,b}}\! (r)\eqo \bE \big[ e^{Z^{(p)}_0 \log (1-\frac{1}{{a_p}} v_p(b) (1-r))}\big] ,
\end{equation}
for all $r\ino [0, 1]$, where we recall that $\psi_p (z)\eqo b_p \big( \varphi_{\xi_p} \big(1\! -\! 
\frac{1}{{a_p}}z \big) \! -\! 1 + \frac{1}{a_p}z\big)$ from Lemma \ref{Greydis} (indeed, (\ref{lengtherapara}) 
follows from (\ref{lengtheraparam}) 
with $\alpha 
\eqo Q_{ p} (\mathtt{Ht} \leqo b)\eqo 1\! -\! \frac{1}{{a_p}} v_p(b)$ and a straightforward computation). 
Under $\mathbf{A1}$--$\mathbf{A4}$ and by Lemma \ref{lmconvF}, 
$\psi_p (z) \eqo \frac{1}{{c_p}} F_{\! p} ( 1\! -\! \frac{1}{{a_p}} z_p\,  ,\, 0)\!  -\! 
b_p(1\! -\! \frac{1}{{a_p}} z) \! \! \to\!   \psi(z)$, and $\psi_p'(z) \eqo c_p \! -\! \partial_1 F_{\! p} ( 1\! -\! 
\frac{1}{{a_p}} z\,  ,\, 0) \!\!  \to \! \psi'(z)$, uniformly on every compact interval of $\bbR_+$. Therefore 
$(c_{p, b}, \varphi_{\xi_{p,b}}, \varphi_{\! \varrho_{p, b} }) \! \to \! (c^{v(b)}\! , \varphi_{\xi^{v(b)}} ,
\varphi_{\! \varrho^{v(b)}})$, where we recall the notation $(c^\lambda, \xi^\lambda, \varrho^\lambda)$, 
$\lambda\ino \bbR_+^*$ from 
(\ref{lambdaparam}).  
Proposition \ref{proppoubelle} then asserts that the laws of $R_b(\mathbf F_{\! p})$ converge to 
$\texttt{P}_{ \! \xi^{v(b)}\! , c^{v(b)} \! ,\varrho^{v(b)}}$ in law in $(\bbT^{0}_{\! c}, \dGH)$. By Theorem 
\ref{metriccpctLevytrees} the law of the $b$-length erasure $R_b (\Phi_0 (\mathbf m))$ of the unmeasured compact ($x_0, \psi$)-Lévy forest is 
equal to $\texttt{P}_{ \! \xi^{v(b)}\! , c^{v(b)} \! ,\varrho^{v(b)}}$. Since by definition, $\Phi_0 (\mathbf m_p)\eqo \mathbf F_{\! p}$, 
we have proved for all $b\ino \bbR_+^*$ that $R_b (\Phi_0 (\mathbf m_p))\! \to \! R_b (\Phi_0 (\mathbf m))$ in law in 
$(\bbT^{0}_{\! c}, \dGH)$. Lemma \ref{lengtheraprop} $(v)$ implies that $\Phi_0 (\mathbf m_p)\! \to \! \Phi_0 (\mathbf m)$ 
in law in $(\bbT^{0}_{\! c}, \dGH)$. 
Since we have already proved Theorem \ref{thmIP} $(ii)$, which asserts that $\mathbf m_p\! \to \! \mathbf m$ in law in 
$(\bbT, \dGP)$, Theorem \ref{GPGH+GHP} $(ii)$ applies and shows that $\Phi_1 (\mathbf m_p)\! \to \! \Phi_1 (\mathbf m)$ 
in law in $(\bbT^{1}_{\! c}, \dGHP)$, which completes the proof of Theorem \ref{thmIP} $(iii)$. \cqfd 

\smallskip

\noi
\textbf{End of the proof of Theorem \ref{thmIP}}. The joint convergence $(\mathbf m_p, \langle \mathbf m_p \rangle) \! \to \! (\mathbf m, \langle \mathbf m \rangle)$ in law in $\bbT^*\! \times \! \cM_f(\bbR_+)$ is a straightforward consequence of Theorem \ref{thmIP} $(i)$ and the $\dera$-continuity of total mass stated in 
Theorem \ref{cveracomplet} $(iv)$. 

Let us assume that $\psi$ is of Type C. Then the convergence $(\mathbf m_p, \flambda_{\mathbf m_p} ) \! \to \! (\mathbf m, \flambda_{\mathbf m})$ in law in $\bbT \! \times \! \cM_f(\bbR_+)$ (and in law in $\bbT^1_{\! c} \times \! \cM_f(\bbR_+)$ under the assumptions Theorem \ref{thmIP} $(iii)$)
is a straightforward consequence of Theorem \ref{thmIP} $(ii)$ and the $\dGP$-continuity (resp.~$\dGHP$-continuity) of 
the function $\fmu \! \mapsto \! \flambda_{\fmu}$ stated in Theorem \ref{eraGcont} $(v)$. By Skorokhod's representation theorem, up to a slight abuse of notation, we may assume that $(\mathbf m_p, \flambda_{\mathbf m_p} ) \! \to \! (\mathbf m, \flambda_{\mathbf m})$ holds a.s.~By Theorem \ref{thm:growthprocconvbis}, in Type C cases $(\flambda_{\mathbf m} ([0, s]))_{s\in \bbR_+}$ has the same law as $(\int_0^s\! Z^{x_0}_r dr)_{s\in \bbR_+}$: $\flambda_{\mathbf m}$ is then diffuse, which a.s.~implies for all $s\ino \bbR_+$, that 
$\flambda_{\mathbf m_p} ([0, s])\! \to \! \flambda_{\mathbf m} ([0, s])$. By Dini's theorem, $(\flambda_{\mathbf m_p} ([0, s]))_{s\in \bbR_+}\! \to \! \flambda_{\mathbf m} ([0, s]))_{s\in \bbR_+}$ a.s.~in $\mathbf C (\bbR_+, \bbR)$. This completes the proof of Theorem \ref{thmIP}. \cqfd

\appendix 
\renewcommand{\eqnsection}{
\renewcommand{\theequation}{\Alph{section}.\arabic{equation}}
    \makeatletter
    \csname  @addtoreset\endcsname{equation}{section}
    \makeatother}
\eqnsection

\section{Proofs of 
Lemmas 
\ref{edgelength} and \ref{fchangemetric} 
}
\label{ProofRtreeapp}

\subsection{Proof of Lemma \ref{edgelength}}
\label{Proofedgelength} 
Let us assume that $\ttt$ is of finite type. Then there are $N\ino \bbN^*$ and $\gamma_1, \ldots, \gamma_N \ino \ttt$ such that $\ttt\eqo \bigcup_{1\leq i\leq N} \lgeo \rho, \gamma_i\rgeo$, which implies that $\ttt$ is compact (and thus closed) and  $\sup_{\sigma \in \ttt} d(\rho, \sigma) \leko \infty$.  
Suppose that $\sigma_1, \ldots, \sigma_n\ino \ttt$ are such that $\sigma_i  \wedge \sigma_j \! \notin \! \{\sigma_i ,\sigma_j \}$ for all distinct $i,j\ino \{ 1, \ldots, n\}$. For all $i\ino \{ 1, \ldots, n\}$, there is $k_i\ino \{ 1, \ldots, N\}$ such that $\sigma_i \ino \lgeo \rho, \gamma_{k_i} \rgeo$. If $k_j \eqo k_i$ and $j\! \neq \! i$, then $\sigma_j \ino \lgeo \rho, \gamma_{k_i} \rgeo$ which contradicts $\sigma_i \wedge \sigma_j \! \notin \!  \{ \sigma_i, \sigma_j \}$. 
Therefore, $i\!\mapsto \! k_i$ is injective, $n\leqo N$, and $\ttt$ satisfies (\ref{bbound}). 

Conversely, we suppose that $\ttt$ satisfies the three conditions. We denote by $N\ino \bbN^*$ the supremum 
in (\ref{bbound}). We can find $\sigma_1, \ldots, \sigma_N \ino \ttt$ such that 
$\sigma_i \! \wedge \! \sigma_j \! \notin \! \{\sigma_i , \sigma_j \}$, for all distinct $i,j\ino \{ 1, \ldots, N\}$. 
Let us suppose that there are $x, x' \ino \ttt \cap \theta_{\sigma_N} T$ which are distinct from $\sigma_N$ and 
such that $x\wedge x'\! \notin \! \{ x,x'\}$: we see that $\{ \sigma_1, \ldots, \sigma_{N-1}, x, x'\}$ satisfies the conditions in (\ref{bbound}) which contradicts the maximality of $N$. 
Therefore $x \wedge x' \ino  \{ x, x'\}$. Since $\sup_{\sigma \in \ttt} d(\rho, \sigma) \leko \infty$, there exists a sequence $x_p \ino  \ttt \cap \theta_{\sigma_N} T$, $p\ino \bbN$, such that $d(\rho, x_p) \leqo d(\rho, x_{p+1}) \! \to \! \ell \! :=\! \sup_{\sigma \in \ttt\cap \theta_{\sigma_{\! N}}\! T } d(\rho, \sigma)\leko \infty$. The previous argument implies that $\lgeo \rho, \sigma_N \rgeo \!\subseteq\! \lgeo \rho, x_p \rgeo \!\subseteq\! \lgeo \rho, x_{p+1} \rgeo$. Thus $d(\rho, x_{p+1} )\eqo  d(\rho, x_{p} )+d(x_p, x_{p+1} )$ which implies that $\ell \eqo d(\rho, x_0)+ \sum_{p\in \bbN} d(x_p, x_{p+1})$. Consequently, $(x_p)_{p\in \bbN}$ is a Cauchy sequence. Since $\ttt$ is closed, there is $\gamma_N\ino \ttt$ 
such that $\lim_{p\to \infty} d(x_p, \gamma_N)\eqo 0$. We necessarily get $\gamma_N \ino  \ttt\cap \theta_{\sigma_{\! N}}\! T $ and if there is $x\ino  \ttt\cap \theta_{\sigma_{\! N}}\! T $, then the previous argument implies that $x\wedge \gamma_N \ino \{ x, \gamma_N\}$. Thus $x\ino \lgeo \rho, \gamma_N\rgeo$ since $d(\rho, \gamma_N)\eqo \ell$. This proves that $\ttt\cap \theta_{\sigma_N} T\eqo \lgeo \sigma_N, \gamma_N \rgeo$. Similarly, for all $i\ino \{ 1, \ldots, N\}$, 
there is $\gamma_i\ino \ttt$ such that $\ttt\cap \theta_{\sigma_i} T\eqo \lgeo \sigma_i, \gamma_i \rgeo$. We note that $\gamma_i$ is necessarily a leaf of $\ttt$ (i.e.~$\ttt\backslash \{ \gamma_i\}$ is connected). 
Then, for all distinct $i,j\ino \{ 1, \ldots N\}$, we get $\sigma_i \wedge  \sigma_j \eqo \gamma_i  \wedge  \gamma_j $ and thus $ \gamma_i  \wedge  \gamma_j \! \notin \! \{\gamma_i , \gamma_j \}$. We now suppose that there is $\sigma\ino \ttt$ such that $\sigma \! \notin \! \ttt'\! :=\! \bigcup_{1\leq i\leq N} \lgeo \rho, \gamma_i\rgeo$. Then, for all $i\ino \{ 1, \ldots, N\}$, we see that $\sigma \wedge \gamma_i \! \neq \! \sigma$ because 
$\sigma \wedge \gamma_i\ino \lgeo \rho, \gamma_i\rgeo \!\subseteq\! \ttt'$, and we also get $\sigma \wedge \gamma_i \! \neq \! \gamma_i $ because $\gamma_i$ is a leaf of $\ttt'$. Thus $\{ \gamma_1, \ldots, \gamma_N, \sigma\}$ satifies the condition in (\ref{bbound}), which contradicts the maximality of $N$. This implies that $\ttt\eqo \ttt'$ and $\ttt$ is therefore of finite type. \cqfd

\subsection{Proof of Lemma \ref{fchangemetric}}
\label{Prooffchangemetric}
We let the reader check that $d_g$ is a distance satisfying the four-point condition (\ref{4ptscondi}). To simplify, we set 
$|\sigma|\eqo d(\rho, \sigma)$ for all $\sigma\ino T$. 
Let us prove that $(T^o\! , d_g)$ is connected. To this end, we fix $\sigma_n, \sigma\ino T^o$, $n\ino \bbN$, such that 
$\lim_{n\to \infty} d(\sigma_n, \sigma)\eqo 0$. We first see that $|\sigma_n | \! \to \! |\sigma| \leko a$. 
Since $(T^o\! , d)$ is an $\bbR$-tree, we also 
get $|\sigma_{n} \! \wedge\!  \sigma|\eqo \frac12 (|\sigma_{n} |+ |\sigma| \! -\! 2 d(\sigma_n, \sigma)) \! \to \! |\sigma|$, which implies $d_g (\sigma_n, \sigma)\! \to \! 0$, since $g$ is continuous. This shows that the identity function on $T^o$ is $(d,d_g)$-continuous. This implies that $(T^o\! , d_g)$ is connected because $(T^o\! , d)$ is connected. Theorem \ref{4ptsth} then yields that $(T^o\! , d_g)$ is an $\bbR$-tree and by the very definition of $\bbR$-trees, geodesic arcs in $(T^o\! , d)$ and $(T^o\! , d_g)$ have to be the same, hence branch points, too. 

Let $\sigma,\sigma_n \ino T^o$, $n\ino \bbN$, be such that $\lim_{n \to \infty} d_g (\sigma_n, \sigma)\eqo 0$. We first get $d_g(\rho, \sigma_n) \! \to \! d_g (\rho , \sigma) \leko b$ and  
since $(T^o\! , d_g)$ is an $\bbR$-tree, 
$d_g(\rho, \sigma_n \! \wedge \! \sigma) \eqo \frac12 (d_g(\rho, \sigma_n )+d_g(\rho, \sigma)\! -\! 2 d_g(\sigma_n,\sigma))\! \to d_g(\rho, \sigma)$. 
 Since $g$ is a homeomorphism, this implies that $\lim_{n\to \infty} |\sigma_n|\eqo \lim_{n\to \infty} 
 | \sigma_n \! \wedge \! \sigma| \! \to |\sigma|$ and thus $d(\sigma, \sigma_n)\eqo   |\sigma_n |+|\sigma |\! -\! 2|\sigma_n \! \wedge \!  \sigma|\! \to \! 0$. This shows that the identity function on $T^o$ is $(d_g,d)$-continuous, which completes the proof of the lemma. \cqfd

\section{Proof of Lemma \ref{measwithmeas}}
\label{pfsecmeaswithmeas}

   Lemma \ref{measwithmeas} generalises \cite[Lemmas 2.3--2.5]{DuWi2} to include finite mass measures. 
We first establish some auxiliary results.
\begin{lemma}\label{lememb} For any $\widetilde{T},\widetilde{T}'\ino\bbT_c^1$ 
with $\dGHP(\widetilde{T},\widetilde{T}')\!<\!\varepsilon$, there is a Polish 
  pointed metric space $(X,\delta,\rho)$ with subsets $T,T'\!\subseteq\! X$ and measures 
  $\mu,\mu^\prime\ino\mathcal{M}_f(X)$ with 
  $\mathtt{Span}(\supp\mu)\!\subseteq\! T$ and $\mathtt{Span}(\supp\mu')\!\subseteq\! T'$ such that 
  $(T,\delta,\rho,\mu)$ and $(T',\delta,\rho,\mu')$ are 
  representatives of $\widetilde{T}$ and $\widetilde{T}'$ with $\delta_{\mathtt{Haus}}(T,T')\vee
 \delta_{\mathtt{Pro}}(\mu,\mu')\!<\!2\varepsilon$.
\end{lemma} 
 \noi
 \textbf{Proof.} By definition, there are $(\hat{X},\hat{\delta})$ with subsets $\hat{T},\hat{T}'\subseteq \hat{X}$ and measures 
  $\hat{\mu},\hat{\mu}^\prime\ino\mathcal{M}_f(\hat{X})$, as well as $\hat{\rho}\ino\hat{T}$ and $\hat{\rho}'\ino\hat{T}'$ such that 
  $(\hat{T},\hat{\delta},\hat{\rho},\hat{\mu})$ and $(\hat{T}'\! ,\hat{\delta},\hat{\rho}'\! ,\hat{\mu}')$ are representatives of $\widetilde{T}$ and $\widetilde{T}'$ 
  with $\hat{\delta}(\hat{\rho},\hat{\rho}')\vee\hat{\delta}_{\mathtt{Haus}}(\hat{T},\hat{T}')\vee\delta_{\mathtt{Pro}}(\hat{\mu},\hat{\mu}')\leko \varepsilon$, without loss of generality 
 $\hat{T}$ and $\hat{T}'$ are disjoint with union $\hat{X}$. We then set $X\eqo \hat{X}\! \setminus\! \{\hat{\rho}'\}$, $\rho\eqo \hat{\rho}$, $T\eqo \hat{T}$, 
  $T'\! \eqo (\hat{T}'\!  \setminus\! \{\hat{\rho}'\}) \cup\{\rho\}$, $\mu\eqo \hat{\mu}$, $\mu^\prime\! \eqo \hat{\mu}'|_X+\hat{\mu}(\{\hat{\rho}'\})\delta_{\! \rho}$ and consider the 
 unique metric (triangle inequality easily checked) $\delta$ which extends $\hat{\delta}$ on $T^2$ and $(T'\! \setminus\! \{\rho\})^2$ and which satisfies 
 $$ \delta(x,y)=\big(\hat{\delta}(x,y)+\hat{\delta}(\hat{\rho},\hat{\rho}')\big)\wedge\big(\hat{\delta}(x,\hat{\rho})+\hat{\delta}(\hat{\rho}',y)\big) \quad \textrm{for all $(x,y)\ino T \! \times \! (T'\! \setminus \! \{\rho\})$.} $$
  Then $\delta_{\mathtt{Haus}}(T,T')\le\hat{\delta}(\hat{\rho},\hat{\rho}')+\hat{\delta}_{\mathtt{Haus}}(\hat{T},\hat{T}')$ and 
  $\delta_{\mathtt{Pro}}(\mu,\mu')\le\hat{\delta}(\hat{\rho},\hat{\rho}')+\hat{\delta}_{\mathtt{Pro}}(\hat{\mu},\hat{\mu}')$. It is therefore easily checked that
  this setup has all the properties required.\cqfd
\smallskip

\begin{lemma}\label{aecont} Let $\widetilde{T}_n$, $n\ino\bN$, be a $\bbT_c^1$-valued sequence that 
$\dGHP$-converges to $\widetilde{T}$, and let 
  $(r_n)_{n\in\bN}$ be an $\bbR_+$-valued sequence that converges to $r\ino\bbR_+$. 
  If furthermore $\flambda_{\widetilde{T}}(\{r\})\eqo 0$, then 
  $(\mathtt{Abv}_{r_n}(\widetilde{T}_n))_{n\in\bN}$ $\dGHP$-converges to $\mathtt{Abv}_r(\widetilde{T})$ and  
  $(\mathtt{Blw}_{r_n}(\widetilde{T}_n))_{n\in\bN}$ $\dGHP$-converges to $\mathtt{Blw}_r(\widetilde{T})$.
\end{lemma}

\noi
\textbf{Proof.} First note that $\flambda_{\widetilde{T}_n}\! \rightarrow\! \flambda_{\widetilde{T}}$ weakly in 
$\cM_f(\bR_+)$. Now let 
$\varepsilon\geko 0$. As $\flambda_{\widetilde{T}}(\{r\})\eqo 0$, there is $\Delta\geko 0$ such that 
$\flambda_{\widetilde{T}}([r\! -\! \Delta,r+\Delta])\leko \varepsilon/8$. Let $\gamma\eqo(\Delta/4)\wedge(\varepsilon/8)$. 
Then by the Portmanteau theorem, there is $N\ino\bN$ such that for all $n\geqo N$, we have 
$\flambda_{\widetilde{T}_n}([r\! -\! \Delta,r+\Delta])\leko \varepsilon/4$, $|r\! -\! r_n|\leko \gamma$ and 
$\dGHP(\widetilde{T}_n,\widetilde{T}) \leko \gamma$. Fix $n\geqo N$. By Lemma \ref{lememb}, 
we can find $(X,\delta,\rho)$, $T$, $T_n$, $\mu$, $\mu_n$ such that $(T,\delta,\rho,\mu)$ and $(T_n,\delta,\rho,\mu_n)$ 
are representatives of 
$\widetilde{T}$ and $\widetilde{T}_n$ with $\delta_{\mathtt{Haus}}(T,T_n)\vee\delta_{\mathtt{Pro}}(\mu,\mu_n)\!<\!2\gamma$. 

Now it is easily checked that we have 
$\delta_{\mathtt{Haus}}(\mathtt{Blw}_r(T),\mathtt{Blw}_{r_n}(T_n))\!<\!|r\!-\!r_n|\!+\!2\delta_{\mathtt{Haus}}(T,T_n)$ and
$\delta_{\mathtt{Pro}}(\mu(\cdot\cap\mathtt{Blw}_r(T)),\mu_n(\cdot\cap\mathtt{Blw}_{r_n}(T_n)))
\leqo\delta_{\mathtt{Pro}}(\mu,\mu_n)+\flambda_{\widetilde{T}}([r\!-\!\Delta,r\!+\!\Delta])\vee 
\flambda_{\widetilde{T}_n}([r\!-\!\Delta,r\!+\!\Delta])\leqo\delta_{\mathtt{Pro}}(\mu,\mu_n)\!+\!\varepsilon/4$. Hence 
$\dGHP\big(\mathtt{Blw}_r(\widetilde{T}),\mathtt{Blw}_{r_n}(\widetilde{T}_n)\big)\le 5\gamma+\varepsilon/4<\varepsilon$.

  Also, consider the quotient metric space $X^\bullet\eqo X/\!\!\sim$, where the equivalence relation $\sim$ is given by $x\!\sim\! y$ iff $x,y\ino B_X(\rho,r)$, 
with the quotient metric $\delta^\bullet$ and root $\rho^\bullet\eqo[\rho]_\sim$. Then the projections/push-forwards 
of $(T,\mu)$ and 
$(T_n,\mu_n)$ under the natural projection map $X\!\rightarrow\! X^\bullet$ are such that 
$\big( T^\bullet\! ,\delta^\bullet\! ,\rho^\bullet\! ,\mu^\bullet(\cdot\cap(X^\bullet \! \setminus\! \{\rho^\bullet\}))\big)$ and 
$\big(T_n^\bullet,\delta^\bullet \! ,\rho^\bullet \! ,\mu_n^\bullet(\cdot\cap(X^\bullet\! \setminus\! \{\rho^\bullet\}))\big)$ 
are representatives of 
$\mathtt{Abv}_r(\widetilde{T})$ and $\mathtt{Abv}_r(\widetilde{T}_n)$, with 
$\delta_{\mathtt{Haus}}^\bullet(T^\bullet,T_n^\bullet)\leqo\delta_{\mathtt{Haus}}(T,T_n)$ and
$$\delta_{\mathtt{Pro}}^\bullet \big( \mu^\bullet|_{X^\bullet\setminus\{\rho^\bullet\}},\mu_n^\bullet|_{X^\bullet\setminus\{\rho^\bullet\}} \big)
\leq \delta_{\mathtt{Pro}}(\mu,\mu_n)+\flambda_{\widetilde{T}}([r\!-\!\Delta,r\!+\!\Delta])\vee
\flambda_{\widetilde{T}_n}\! ([r\!-\!\Delta,r\!+\!\Delta]) \; .$$
Hence $\dGHP\big(\mathtt{Abv}_r(\widetilde{T}),\mathtt{Abv}_{r}(\widetilde{T}_n)\big)\leqo 2\gamma+\varepsilon/4$. 

  Let $(T_n^-\! ,d^-\! ,\rho^-\! ,\mu_n^-)$ be a representative of $\mathtt{Abv}_{r\wedge r_n}(\widetilde{T}_n)$. 
Let us define the equivalence relation with $x\!\sim\! y$ iff
$x,y\ino B_{T_n^-\! , d^-_n}(\rho^-\! ,|r_n\!-\!r|)$. Let $T_n^+\eqo T_n^-/\!\!\sim$ with the quotient metric denoted by $d^+$, root  
$\rho^+\! =\! [\rho^-]_\sim$ and push-forward $\mu_n^+$ of $\mu_n^-$ under the natural projection map 
$p\colon T_n^-\!\rightarrow T_n^+$. Then 
$(T_n^+,d^+,\rho^+,\mu_n^+(\cdot\cap T_n^+\setminus\{\rho^+\})$ is a representative of 
$\mathtt{Abv}_{r\vee r_n}(\widetilde{T}_n)$. Consider 
$X^\pm=T_n^-\sqcup(T_n^+\setminus\{\rho^+\})$ equipped with the metric $d^\pm$ which extends $d^-$ on $(T^-_n)^2\! $, 
$d^+$ on $ (T^+_n\!\setminus\!\{\rho^+\})^2\! $ and which satisfies for all $(x,y)\ino T^-_n\!\times\! (T^+_n\!\setminus\!\{\rho^+\})$, 
$$d^\pm(x,y) \eqo d^\pm(y,x) \eqo \big(d^-_n(x,p^{-1}(y))\!+\!|r\!-\!r_n|\big)\!\wedge\!\big(d^-_n(x,\rho^-)\!
+\!d^+_n(\rho^+,y)\big) . $$
  Then $d^\pm_{\mathtt{Haus}}(T_n^-,T_n^+)\leqo|r-r_n|<\gamma$ and $d^\pm_{\mathtt{Pro}}(\mu_n^-,\mu_n^+)\leqo|r-r_n|+\flambda_{\widetilde{T}_n}([r-\Delta,r+\Delta])<\gamma+\varepsilon/4$. Hence
  $\dGHP\big(\mathtt{Abv}_r(\widetilde{T}_n),\mathtt{Abv}_{r_n}(\widetilde{T}_n)\big)
                     \!<\!2\gamma\!+\!\varepsilon/4$,
  and $\dGHP\big(\mathtt{Abv}_r(\widetilde{T}),\mathtt{Abv}_{r_n}(\widetilde{T}_n)\big)\!<\!4\gamma\!+\!\varepsilon/2\leqo\varepsilon$.\cqfd

\medskip

\noi
\textbf{Proof of Lemma \ref{measwithmeas} $(i)$.}  We first prove the measurablility of $\flambda\colon\mathbb{T}_c^1\to\cM_f(\bbR_+)$. Let $\widetilde{T},\widetilde{T}^\prime\ino\bbT_c^1$ with $\dGHP(\widetilde{T},\widetilde{T}^\prime)\!<\!\varepsilon$. Consider the embedding of Lemma \ref{lememb}. Then clearly $d_{\mathtt{Pro}}(\flambda_{\widetilde{T}},\flambda_{\widetilde{T}^\prime})\leqo \delta_{\mathtt{Pro}}(\mu,\mu^\prime)<2\varepsilon$. Hence $\flambda$ is Lipschitz continuous.

We next prove the measurability of  $\mathtt{Blw},\mathtt{Abv}\colon\bbR_+\!\times\!\bbT_c^1\!\rightarrow\!\bbT_c^1$. 
We adapt the argument for the measurability of $f$ in the proof of 
 Lemma \ref{cpctsubsets}. We first $r\ino\bbR_+$, $\widetilde{T}\ino\bbT_c^1$ and $h\ino\bbR^*_+$. Then we observe that 
 $s\ino(0,1) \! \mapsto\!  F(sh,\widetilde{T})\! :=\! \mathtt{Blw}_{r+sh}(\widetilde{T})$ is $\dGHP$-c\`adl\`ag and 
 hence Borel-measurable. Consider a uniform r.v.~$U\colon \Omega\! \rightarrow \! (0,1)$. Then Lemma \ref{aecont} yields that $\bP$-a.s.~$\widetilde{T}\! \mapsto\!  F(Uh,\widetilde{T})$ is continuous and 
$\lim_{h\downarrow 0} F(Uh,\widetilde{T})\eqo \mathtt{Blw}_r(\widetilde{T})\! =:\! f(\widetilde{T})$ by right-continuity. Then the argument based on convergence of laws and Kuratowski's theorem discussed in Lemma \ref{cpctsubsets} establishes the measurability of $\widetilde{T} \! \mapsto \! f(\widetilde{T})\eqo \mathtt{Blw}_r(\widetilde{T})$ for each $r\ino\bbR_+$. Since 
  $\widetilde{T}\! \mapsto \! (\mathtt{Blw}_r(\widetilde{T}))_{r\in \bbR_+}$ takes values in the Skorokhod space $\mathbf{D}(\bbR_+,\bbT_c^1)$, whose Borel 
  sigma-field is generated by the evaluation maps, this entails the measurability of $\widetilde{T}\mapsto(\mathtt{Blw}_r(\widetilde{T}))_{r\in \bbR_+}$, and
  a standard argument yields the claimed measurability of $(\widetilde{T},r)\mapsto\mathtt{Blw}_r(\widetilde{T})$. The argument for $\mathtt{Abv}$ is similar
  and left to the reader.

We finally derive the measurability of the functions $D\colon \bbT_c^1\! \to \! [0,\infty]$, $\mathbf{k} \colon \bbT_c^1\! \rightarrow\! \bN\cup\{\infty\}$ and $Z_r^+\colon \bbT_c^1\! \rightarrow\! \bN\cup\{\infty\}$ from the corresponding results on  $\bbT_c^0$ in \cite[Lemmas 2.4--2.5]{DuWi2} by composition with the 1-Lipschitz function 
$\Phi^1_0\colon \bbT_c^1\! \rightarrow \! \bbT_c^0$ which assigns to the isometry class  of $(T,d,\rho,\mu)$ in $\bbT_c^1$, the isometry class of $(T,d,\rho)$ in $\bbT_c^0$. \cqfd

\medskip

\noi
\textbf{Proof of Lemma \ref{measwithmeas} $(ii)$.}  Consider $\mathbf{T}\!:=\!\widetilde{\mathtt{T}}\mathtt{ree}(\tau,\ell(\cdot),w_\cdot,\mu_\cdot))$ with law $Q_{\xi,c,\phi,\Lambda}$ in the setting of Definition 
\ref{meaGWregedef}. Then $\mathbf{k}_\varnothing(\tau)$ has 
law $\xi$. Let $n\ino \bbN\backslash \{ 1\}$ be such that $\xi (n) \geko 0$. Then under $\bP(\,\cdot\,|\,\mathbf{k}_\varnothing(\tau)\eqo n)$, the subtrees $(\theta_{ (i)}\tau)_{1\leq i\leq n}$ are independent and distributed as $\tau$
under $\bP$. Furthermore, the conditional independence and the conditional laws of the r.v.s $(\ell(u),w_u,\mu_u)$, $u\ino\tau$, show that under $\bP(\,\cdot\,|\,\mathbf{k}_\varnothing(\tau)\eqo n)$, 
the vector $(\ell(\varnothing),w_\varnothing,\mu_\varnothing)$ and the marked 
trees $(\theta_{(i)}\tau,(\ell((i)*u),w_{(i)*u},\mu_{(i)*\mu})_{u\in\theta_{(i)}\tau})$, $1\leqo i\leqo n$, are independent and the latter are distributed as 
$(\tau,\ell(\cdot),w_\cdot,\mu_\cdot)$ under $\bP$, while the former have the law given in \eqref{lawatu} with $\mathbf{k}_\varnothing(\tau)\eqo n$. 

We next observe that $\bP$-a.s.~$\mathtt{Ht}(\mathbf{T})\geqo D(\mathbf{T})\geqo\ell(\varnothing)\!>\!0$. Hence, 
$\mathbf{k}(\mathbf{T})\eqo\mathbf{k}_\varnothing\!\not=\! 1$, 
$D(\mathbf{T})\eqo\ell(\varnothing)$, 
$Z_r^+(\mathbf{T})\eqo 1$ for $0\leqo r\leqo D(\mathbf{T})$, 
$\flambda_\mathbf{T}|_{[0,D(\mathbf{T}))}\eqo\mu_\varnothing$ and 
$\flambda_\mathbf{T}(\{D(\mathbf{T})\})\eqo w_\varnothing$.
Furthermore, these all are necessarily independent of 
$$\mathtt{Abv}_{D(\mathbf{T})}(\mathbf{T})\eqo\circledast_{1\leq i\leq \, \mathbf{k}(\mathbf{T})}\widetilde{\mathtt{T}}\mathtt{ree}\Big( \theta_{(i)}\tau, \big( \ell((i)\!*\!u),w_{(i)*u},\mu_{(i)*u} \big)_{u\in\theta_{(i)}\tau} \Big)\; , $$ 
which has law $Q_{\xi,c,\phi,\Lambda}^{\circledast n}$ under $\bP(\,\cdot\,|\,\mathbf{k}(\mathbf{T})\eqo n)$. Since this holds true for all $n\ino\bN\!\setminus\!\{1\}$, this shows that 
$Q_{\xi,c,\phi,\Lambda}$ satisfies \eqref{df:GWweight}.

Conversely, we suppose that $Q$ and $Q'$ are two laws on $\bbT_c^1$ that satisfy \eqref{df:GWweight}. 
and we recall from the end of the proof of Lemma \ref{measwithmeas} $(i)$ the function 
$\Phi_0^1\colon  \bbT_c^1\! \rightarrow \! \bbT_c^0$. 
Then the push-forwards $Q_*$ and $Q_*^\prime$ of $Q$ and $Q'$ under $g$ satisfy the corresponding equation for $\bbT^0_{\! c}$-valued (unmeasured) Galton--Watson trees and we obtain from \cite[Lemma 2.15]{DuWi2} that $Q_*=Q_*^\prime$. The idea of proof there was to 
use a measurable map to associate with members of $\bbT_c^0$ equipped a discrete branching structure 
a marked ordered tree in $\bbT_{\mathtt{dis}}^0:=\bigsqcup_{\mathbf{t}\in\bbT_{\mathtt{dis}}}\{\mathbf{t}\}\times[0,\infty]^\mathbf{t}$ and then to use a shuffling kernel to associate with $Q_*$ and $Q_*^\prime$ two laws on $\bbT_{\mathtt{dis}}^0$. These laws on $\bbT_{\mathtt{dis}}^0$ were then identified as the laws of Galton--Watson trees with i.i.d.\ marks, seen to be identical entailing also that $Q_*=Q_*^\prime$. The same idea applies here, with members of $\bbT_c^1$ 
of finite type and marked ordered trees in $\bbT_{\mathtt{dis}}^1$, and the steps of the proof are the same. We omit the details.\cqfd
\section{Proof of Lemma \ref{existlaws}}
\label{pfsecexistlaws}
We set $(\xi_h, c_h, \varrho_h)\eqo (\xi^{\nu(h)} \! , c^{\nu(h)}\! , \varrho^{\nu(h)})$, the 
$(\xi^\lambda \! , c^\lambda \! , \varrho^{\lambda})$ being as in (\ref{lambdaparam}) and $\nu(h)$ as in (\ref{psihnuhdef}). 
Since $\psi'$ is the Laplace exponent of a subordinator, so is  
$f_h \! :=\! \psi' (\, \cdot +\nu (h))\! -\! \psi'(\nu(h))$. We recall that $\, \overline{\! \psi}^{{-1}}_{h}$ is the Laplace exponent of a subordinator: so is $\phi_h \! :=\! f_h\!  \circ \, \overline{\! \psi}^{{-1}}_{h}$ by subordination and 
$\phi_h$ satisfies  (\ref{chphih}). 
For all $n\geq 2$, we next denote by $\mathbf s_n$ a r.v.~whose law is 
$$\tfrac{1}{\xi_h(n)n!} \Big( \tfrac{\nu(h)}{\psi'(\nu(h)) }\mathbf b\, \un_{\{ n=2\}} \delta_0 (dx)  
+\un_{\bbR_+^*} (x) \tfrac{\nu(h)^{n-1}x^n}{\psi'(\nu(h))} e^{-x\nu(h)} \pi(dx) \Big).$$ 
In particular $\bP (\mathbf s_n\eqo 0)\eqo 0$ if $n\geqo 3$ and $\bP (\mathbf s_2\eqo 0)\eqo (\mathbf b \nu(h))/ 
(2\xi_h(2) \psi'(\nu(h)))$. We assume that $\mathbf s_n$ is independent of a subordinator $(S^h_s)_{s\in \bbR_+}$ of Laplace exponent $\overline{\!\psi}_h^{-1}$, and 
we define $\Lambda^{\! h}_n$ as the law of $(n\! -\! 1) h +S^h_{\mathbf s_n}$. Then we easily check that 
\begin{equation}
\label{Lambdahenndef}
\forall n\geqo 2, \quad \forall y\ino \bbR_+ , \quad \nu(h) e^{-hy} c_h \xi_h (n) \widehat{\Lambda}_n^{h} (y) 
\eqo \tfrac{1}{n!} \big(\! -\! \nu(h) e^{-hy} \big)^n \psi^{(n)} \big( \psi^{-1}_h(y)\big), 
\end{equation} 
which is consistent with (\ref{xihLamh}). 

We next set $\mathcal L (y) \! :=\! e^{hy}(\psi(\psi_h^{-1}(y))\! -\! y)/\psi(\nu(h))$, for all 
$y\ino \bbR_+$ and we prove that it is the Laplace transform of a law on $\bbR_+$, which is not obvious, 
and we will employ a limiting procedure. 
For all $\lambda\ino \bbR_+^*$, we denote by $Q^\lambda$ the law of a 
GW tree whose parameters are $\xi^\lambda$, $c^\lambda$, $\phi^\lambda$ and $(\Lambda_{\lambda,n} )_{n\in \bbN}$ 
where $\phi^\lambda \! \equiv\!  0$ (no mass on edges), 
$\Lambda_{\lambda,n}\eqo \delta_0$ for all $n\ino \bbN^*$ (no weight on branch points) and 
$\Lambda_{\lambda,0} \eqo \delta_{1/\psi(\lambda)}$: namely the measure which equips the tree is 
$1/\psi (\lambda)$ times the empirical measure of the leaves of the GW tree. Then under $Q^\lambda$ the total 
mass of the tree 
$\langle \flambda \rangle$ is equal to 
$L/\psi(\lambda)$, where $L$ is the number of leaves. 
We then observe that 
$F_{\xi^\lambda\! , c^\lambda \! , \Lambda_{\lambda , \cdot}} (s,y)\eqo 
c^\lambda (\varphi_{\xi^\lambda} (s) \! -\! \xi^\lambda (0) (1\! -\! e^{-y/\psi(\lambda)}))$, $s\ino [0, 1]$, $y\ino \bbR_+$. 
By Proposition \ref{prop:eraGW} $(i)$, $\mathcal E_h$ under 
$Q^\lambda ( \, \cdot  \, | \,   \langle \flambda \rangle \! >\!  h)$ has law 
$Q_{\xi^\lambda_h , c^\lambda_h , \phi^\lambda_h,  \Lambda^h_{\lambda, \cdot}\!\!}$, where 
$\xi^\lambda_h , c^\lambda_h , \phi^\lambda_h$ and 
$(\Lambda^h_{\lambda, n})_{n\in \bbN}$ are characterized by (\ref{charaQera}). More precisely, 
let us set $\falpha^{\lambda}_{h} (y)\! :=\! 
Q^{\lambda} [\exp (-y\langle \flambda \rangle) \un_{\{ \langle \flambda \rangle \leq h\}} ]$ and 
$\alpha^\lambda_h\! :=\!  \falpha^\lambda_h (0)$. By (\ref{charaQera}) we get $(c^\lambda_h, \xi^\lambda_h)\eqo (c^{\alpha^\lambda_h}\! , \xi^{\alpha^\lambda_h})$, where the 
$(c^\lambda\! , \xi^\lambda)$ are as in (\ref{lambdaparam}) and $(1\! -\! \alpha^\lambda_h) e^{-hy} c^\lambda_h \, \xi^\lambda_h (0) \, \widehat{\Lambda}^h_{\lambda, 0} (y) 
\eqo c^\lambda \big( \varphi_{\xi^\lambda} (\falpha^{\lambda}_h (y)) \! -\! \xi^\lambda (0) \big( 1\! -\! e^{-y/\psi(\lambda)} \big) \big)\! -\! 
c^\lambda \falpha^{\lambda}_h (y), $
 for all $s\ino [0, 1]$ and all $y\ino \bbR_+$, 
 which implies 
\begin{equation}
\label{approxLambdazero}\widehat{\Lambda}^h_{\lambda, 0} (y) = e^{hy}\! \cdot \! 
\frac{\psi \big( \lambda (1\! -\! \falpha^\lambda_h (y) )\big)  - 
\psi (\lambda) \big(1\! -\! e^{-y/\psi (\lambda) } \big)}{\psi \big( \lambda (1\! -\! \alpha^\lambda_h) \big)} . 
\end{equation}
We next compute the limit of $\lambda (1\! -\! \falpha^\lambda_h (y) )$ as $\lambda \! \to \! \infty$. 
To this end, we set 
$g_\lambda (s)\eqo Q^\lambda [s^{L}]$, $s\ino [0, 1]$. 
We easily see that $g_\lambda$ satisfies the following functional equation: 
$g_\lambda(s)\eqo \varphi_{\xi^\lambda} (g_\lambda (s)) \! -\! \xi^\lambda(0) (1\! -\! s)$, which implies that 
$1\! -\! g_\lambda (s) \eqo \tfrac{1}{\lambda}\psi^{-1} (\psi (\lambda) (1\! -\! s))$. Namely 
\begin{equation}
\label{leavescomput}
\lambda \big( 1\! -\! Q^\lambda  \big[ e^{-y \langle \flambda \rangle} \big] \big) = \lambda \big( 1\! -\! 
g_\lambda \big( e^{-y/\psi (\lambda)}\big) \big)\eqo \psi^{-1} \big( \psi (\lambda) (1\! -\! e^{-y/\psi (\lambda)} \big) \big) . 
\end{equation}
Let $(L_{\lambda, p})_{p\in \bbN^*}$, be independent r.v.s with the same law as $\langle \flambda \rangle$ under 
$Q^\lambda$ and let $(N_s)_{s\in \bbR_+}$ be a homogeneous Poisson process on $\bbR_+$ with unit intensity, 
which is assumed to be independent of the $(L_{\lambda, p})_{p\in \bbN^*}$. We set 
$W_{\lambda, s} \! :=\!  \sum_{1\leq p\leq N_{\lambda s}} L_{\lambda, p}$ for all $s\ino \bbR_+^*$, with the convention $W_{\lambda , s}\eqo 0$ if $N_{\lambda s }\eqo 0$. Then, $W_{\lambda, \cdot}$ is a compount Poisson process with nonnegative 
jumps 
whose Laplace exponent $-\frac{1}{s}\log \bE \big[ \exp (-yW_{\lambda, s})\big]$ is the member on the left hand side of (\ref{leavescomput}), which tends to $\psi^{-1} (y)$ as $\lambda\! \to \! \infty$. 
Standard arguments on scaling limits 
of Lévy processes imply that $(W_{\lambda, s})_{s\in \bbR_+}$ converges in law in the Skorokhod space $\mathbf D(\bbR_+, \bbR)$ 
to the subordinator $(S_s)_{s\in \bbR_+}$ with Laplace exponent $\psi^{-1}$. Since the Lévy measure $\nu$ of $S$ 
is diffuse (by Lemma \ref{psi-1prop} $(ii)$ because $\psi$ is of type B or C ), the parameter of the first jump of 
$W_{\lambda, \cdot}$ larger that $h$ and the law of this jump converge respectively to $\nu(h)\eqo \nu ((h,\infty))$ and 
$\nu (\, \cdot \, \cap (h, \infty))/\nu(h)$. Namely 
$\lim_{\lambda \to \infty} \lambda (1\! -\! \alpha^\lambda_h)\eqo  \lim_{\lambda \to \infty} 
\lambda Q^{\lambda} (\langle \flambda \rangle \geko h )\eqo \nu(h)$ and $\lim_{\lambda \to \infty} 
Q^\lambda [ \exp (-y \langle \flambda \rangle)| \langle \flambda \rangle \geko h ]\eqo\int_{(h, \infty)} e^{-yx}\nu(dx)/ \nu(h)$. 
Therefore we get $\lim_{\lambda \to \infty}$ $ \lambda (1\! -\! \falpha^\lambda_h(y)) \eqo \psi^{-1} (y) \! -\! \int_{(h, \infty)} \! e^{-yx}\nu(dx)\eqo 
\psi^{-1}_h(y)$ and by (\ref{approxLambdazero}), $\lim_{\lambda \to \infty} \widehat{\Lambda}^h_{\lambda, 0} (y)\eqo \mathcal L(y)$, 
for all $y\ino \bbR_+$. Since $\lim_{y\to 0^+}\mathcal L (y)\eqo \mathcal L(0)\eqo 1$, this implies that the laws 
$\Lambda_{\lambda, 0}^{\! h} $ converge in distribution on $\bbR_+$ to a limiting law that is taken as the definition of $\Lambda^{\! h}_0$ whose Laplace transform is $\mathcal L$. 
Combined with (\ref{Lambdahenndef}), this easily entails (\ref{xihLamh}).

  Finally, we recall that $\varrho_h\eqo \mathtt{Poisson} (x_0 \nu(h))$. For all $n\ino \bbN$, we then take $ \Lambda^{\!h}_{\varnothing , n}$ as the law of $( (n\! -\! 1) h + S^{h}_{x_0} )_+$. 
We easily check (\ref{varrhoh}). This completes the proof of the lemma. \cqfd

\section*{Acknowledgements.} The authors are grateful to Mie Gl\"uckstad for useful feedback on earlier drafts.

For the purpose of Open Access, the authors have applied a CC BY public copyright licence to any Author Accepted Manuscript (AAM) version arising from this submission.

{\small
\bibliographystyle{acm}
\bibliography{biblio}
}
\end{document}